\documentclass{article}
\pdfoutput=1
\usepackage{fullpage}
\usepackage[breaklinks,hidelinks]{hyperref}
\hypersetup{
  pdftitle={Coinduction in Grounded Arithmetic},
  pdfauthor={Bryan Ford}}

\usepackage{amsmath}
\usepackage{amssymb}
\usepackage{amsthm}
\usepackage{booktabs}
\usepackage{array}
\usepackage{longtable}
\usepackage{tikz}

\usepackage[T1]{fontenc}	%

\usepackage{times}
\usepackage{xspace}
\usepackage{cite}		%
\usepackage{bussproofs}
\usepackage{stmaryrd}
\usepackage[noabbrev,capitalize]{cleveref}
\usepackage{thmtools}
\usepackage{xstring}
\usepackage{graphicx}
\usepackage{subcaption}
\usepackage{nicematrix}		%
\usepackage{longtable}
\usepackage{pifont}

\usepackage{isabelle}
\usepackage{isabellesym}
\usepackage{pdfpages}
\usepackage{ragged2e}
\newcommand{\SNIP}[2]{\expandafter\newcommand\csname snippet--#1\endcsname{#2}}
\IfFileExists{snips.tex}{\input{snips}}{}

\newcommand{\GetSnip}[1]{%
    \ifcsname snippet--#1\endcsname%
        \csname snippet--#1\endcsname%
    \else%
        \PackageWarning{snips}{Snippet ``#1'' is undefined.}%
        \emph{Warning: Snippet ``#1'' is undefined.}%
    \fi%
}

\newcommand{\Snippet}[1]{{%
  \newcount\i
  \i=0
  \loop
    \GetSnip{#1-\the\i}%
    \advance \i 1
  \ifcsname snippet--#1-\the\i\endcsname
  \repeat
}}

\newcommand{\SnippetPart}[3]{{%
  \newcount\i
  \i=#1
  \loop
    \ifnum \i=#2
      \renewcommand{\isanewline}{}%
    \fi
    \GetSnip{#3-\the\i}%
    \advance \i 1
    \ifnum \i>#2 {}
    \else \repeat
}}

\hypersetup{breaklinks=true}

\newcommand{\com}[1]{}

\long\def\note#1{}	%

\newcommand{\rga}{RGA\xspace}		%

\com{	%

}%

\makeatletter
\newlength{\doublefracgap}
\DeclareRobustCommand{\doublefrac}[2]{%
  \mathinner{\mathpalette\doublefrac@{{#1}{#2}}}%
}
\newcommand{\doublefrac@}[2]{\doublefrac@@#1#2}
\newcommand{\doublefrac@@}[3]{%
  \ooalign{%
    \raisebox{\doublefracgap}{$\m@th#1\frac{#2}{\phantom{#3}}$}\cr
    \raisebox{-\doublefracgap}{$\m@th#1\frac{\phantom{#2}}{#3}$}\cr
  }%
}
\newcommand{\ddoublefrac}[2]{{\displaystyle\doublefrac{#1}{#2}}}

\makeatother

\newcommand{\irl}[1]{\ensuremath{\mathit{#1}}}		%

\newcommand{\infrule}[3][]{\cfrac{#2}{#3}\IfStrEq{#1}{}{}{\ \irl{#1}}}
\newcommand{\infeqv}[3][]{\ddoublefrac{#2}{#3}\IfStrEq{#1}{}{}{\ \irl{#1}}}
\newcommand{\infceqv}[4][]{\cfrac{#2\qquad}{}\ddoublefrac{#3}{#4}
				\IfStrEq{#1}{}{}{\ \irl{#1}}}

\newtheorem{thm}{Theorem}[section]
\newtheorem{cor}[thm]{Corollary}

\newtheorem{defi}[thm]{Definition}
\newtheorem{rem}[thm]{Remark}

\crefname{thm}{Theorem}{Theorems}
\crefname{cor}{Corollary}{Corollaries}
\crefname{lem}{Lemma}{Lemmas}
\crefname{defi}{Definition}{Definitions}
\crefname{rem}{Remark}{Remarks}

\providecommand{\oga}{OGA\xspace}
\providecommand{\rga}{RGA\xspace}
\providecommand{\iga}{IGA\xspace}
\newcommand{\oiga}{OIGA\xspace}
\newcommand{\inga}{INGA\xspace}
\newcommand{\oinga}{OINGA\xspace}
\newcommand{\dinga}{DINGA\xspace}
\newcommand{\cA}{\mathcal A}             %
\newcommand{\fix}{\mathrm{fix}}          %
\newcommand{\Acc}{\mathrm{Acc}}
\newcommand{\Path}{\mathrm{Path}}

\newcommand{\Dec}{\mathrm{Dec}}
\newcommand{\Tot}{\mathrm{Tot}}
\newcommand{\Term}{\mathrm{Term}}
\newcommand{\Div}{\mathrm{Div}}
\newcommand{\WinI}{\mathrm{Win}_{\mathrm I}}
\newcommand{\WinII}{\mathrm{Win}_{\mathrm{II}}}
\newcommand{\lfp}{\mathrm{lfp}}
\newcommand{\gfp}{\mathrm{gfp}}
\newcommand{\Mstar}{M^{*}}               %
\newcommand{\dual}[1]{{#1}^{\partial}}  %
\newcommand{\isn}[1]{\textsf{#1}}        %
\newcommand{\artifacttag}{\texttt{coind-}\allowbreak\texttt{2026-11-11}\xspace}	%
\newcommand{\ent}{\vdash}

\title{Coinduction in Grounded Arithmetic}
\author{Bryan Ford\\EPFL}
\date{}

\begin{document}
\maketitle

\begin{abstract}
A grounded arithmetic counts a sentence true only when a
terminating process establishes it; a sentence no such process
settles has no value.  Such a logic can state that an event happens,
but not that it never happens under unbounded choice.  This paper
extends such a logic so that it can.  The base is Reflective Grounded
Arithmetic
(\rga) with Inductive Grounded Arithmetic (\iga)'s least-fixed-point
reading and ground induction.  It adds a
separate greatest-fixed-point former $\nu$ with Park's coinduction rule, proves
that one former cannot carry both ground induction and coinduction,
and adds a dichotomy rule: the classical disjunction between a
least-fixed-point atom and the greatest-fixed-point atom of its dual,
under a derivable premise that the body's leaves are decided.  The
systems are sound for an explicit model reading every $\nu$-body at
its greatest fixed point, hence consistent; they agree with the least
model where it has a value and prove sentences it leaves without
one.  The dichotomy trades the underivable hypothesis that an order
is well-founded for the derivable one that its steps are decided, at
the price of a change of model.  With ground induction the
dichotomy derives the omniscience rule of the sibling system
\oga, so the system has no semantics whose certificates for universals
are its own proofs; without the dichotomy it has one, though unsound
for the least model.  The systems prove, as single sentences, that
every decided tree is well-founded or has an infinite branch given by a
term, total in the new model; that random assignment terminates while its
divergent twin runs forever; that a fair loop terminates and the alternating-bit
protocol delivers every message; that open games are determined; that
the termination gap in Nakata and Uustalu's trace semantics closes;
and that Cinderella wins the bucket game at capacity three.
Isabelle/HOL checks every result.
\end{abstract}

\section{Introduction}
\label{sec:intro}

Five buckets of capacity three stand in a circle.  Every round the
Stepmother pours one unit of water, divided among the buckets as she
likes; then Cinderella empties two adjacent buckets.  The Stepmother
wins if a bucket ever overflows.  Does Cinderella have a strategy that
keeps the water in the buckets forever?  At capacity three she does,
and a single round-robin invariant certifies it; smaller capacities
are not attempted here.  The game is a standard benchmark for
infinite-state safety games~\cite{bodlaender12cinderella,beyene14constraint}.
Look at what the two players claim.  The Stepmother claims that an
overflow happens: a terminating event, witnessed by a finite play.
Cinderella claims that it never happens, whatever the Stepmother
pours, and a unit can be divided in infinitely many ways, so the plays
branch without bound at every round.  A logic that counts as true only
what a terminating process establishes can state the Stepmother's
claim.  It can also state that one given play never overflows, a
universal over the play's rounds.  It cannot state Cinderella's claim,
a ``never'' over unboundedly branching plays, except as the negation
of a claim that nobody has shown false.

This paper is about such logics, and about the step that lets them
say what Cinderella says.  A \emph{grounded} arithmetic is one in
which a sentence is true only when a terminating computation, or a
well-founded tree of them, makes it so; a sentence that no such
process settles is not false but has no value.  Picture a logician
leafing through a book of computations and coloring some of them
true.  The grounded tradition colors only terminating ones.  What if
the logician also colored certain computations that never terminate?
Not all of them: the Liar never terminates, and coloring the Liar
true is the inconsistency that grounded logic exists to avoid.  The
question is which non-terminating computations deserve the color.
The answer this paper gives is the coinductive one.  A computation
deserves the color when it is an instance of a recursive definition
declared for the greatest fixed point and lies in that fixed point,
so that the definition's own unfolding never refutes it.  The paper's
thesis is that a grounded arithmetic supports that reading exactly,
with its rule, its model and its limits stated as theorems.

The setting is Reflective Grounded Arithmetic (\rga)~\cite{rga-paper},
a \emph{paracomplete} arithmetic, one in which a sentence may be
neither true nor false, where every term is a computation and
self-reference is native: a fixed-point former turns any body into a
term that unfolds to itself.  Its extension \iga, Inductive Grounded
Arithmetic~\cite{iga-paper}, adds the elimination rule for every
positive self-referential term, called ground induction, and with it
the grounded meaning of a positive recursive definition is its least
fixed point: accessibility, well-foundedness, the winning positions
of a game, and \rga's own truth predicate are single terms with
induction principles.  The least fixed point is the inductive reading,
written $\mu$ below.  Cinderella's claim needs the other one.

\paragraph{What the paper adds.}
Two rules, each with its semantics, a theorem about what they cost,
and a body of examples:
\begin{itemize}
\item a greatest-fixed-point former $\nu$ with Park's coinduction rule,
  a five-line theorem that one former cannot carry both ground
  induction and coinduction, and an explicit model $\Mstar$ for which
  the systems with the rule are sound and consistent, and semantically
  conservative over the least model (\S\ref{sec:nu});
\item a dichotomy rule, a characterization of exactly when it is valid
  in $\Mstar$, and its soundness there (\S\ref{sec:dich});
\item a theorem that, with ground induction, the dichotomy derives the
  omniscience rule of \iga's sibling \oga, and that no system with that
  rule, a recursively enumerable set of theorems and no false claim of
  nontermination is sound for any semantics whose certificates for
  universals are its own proofs
  (\S\ref{sec:dich});
\item a theorem that the system with $\nu$ but without the dichotomy
  keeps such a semantics: it is sound and complete for one, although
  it is not sound for the least model, so the loss comes with ATI and
  not with coinduction (\S\ref{sec:dich}, Appendix~\ref{app:inga});
\item examples: streams, runs and branches; random assignment and fair
  termination through it, up to the alternating-bit protocol; Nakata
  and Uustalu's gap between termination and non-divergence; open
  determinacy and the bucket game; the decidedness of Kleene's
  ordinal notations (\S\S\ref{sec:streams}--\ref{sec:showpieces});
\item a map, produced by build-time checks, of which rules each proof
  applies (\S\ref{sec:rungmap}).
\end{itemize}

\paragraph{The coinduction rule.}
The former $\nu$ is separate in syntax from $\mu$, with Park's
coinduction rule as its introduction.  To show that a point lies in
$\nu\Phi$, exhibit a certificate: a predicate that holds at the point
and that, wherever it holds, makes the body true with the certificate
in place of the recursive variable.  The rule's closure premise is a
hypothetical derivation at an undecided hypothesis, which a grounded
system admits as a premise sequent and would not admit as an
implication inside a term.  Elimination is unfolding.  Ground
induction cannot be extended to $\nu$-atoms: at a $\nu$-term it is
inconsistent in five lines, which is why the two formers are separate.
The resulting systems, \inga (\iga with $\nu$) and \oinga (the union
\oiga of \iga with its sibling \oga~\cite{oga-paper}, with $\nu$), are
sound for an explicit model $\Mstar$ in which every body of the $\nu$
class is read at its greatest fixed point and everything else as in
the least model.  The systems are consistent.  On every sentence the
least model values, they agree with the least model.  And they prove
sentences the least model leaves without a value, among them the
all-zeros atom of the five-line argument, which is the coloring made
precise.

\paragraph{The dichotomy.}
Classically $\neg\mu\Phi$ is $\nu(\neg\dual\Phi)$, the greatest fixed
point of the \emph{co-body} $\neg\dual\Phi$ (the De Morgan dual of
the body).  A grounded system has no rule that turns a failed
verification of the $\mu$-atom into a proof of the $\nu$-atom, so
the two live apart.  The dichotomy rule supplies the classical
disjunction $\mu\Phi\cdot x \vee \nu(\neg\dual\Phi)\cdot x$ at one
body and one point, under a premise that is a sentence of the system
rather than a semantic condition: every non-recursive leaf of the
body is decided.  With
its premise the rule is sound for $\Mstar$, and a theorem says exactly
when the rule is valid there: when the $\nu$-atom has a value at every
instance.  The \iga paper priced the rule of bivalence of
accessibility, and the price was a hypothesis on the order, that it is
well-founded.  The dichotomy's premise is a hypothesis on the order
too, that its steps are decided, but one that is derivable whenever
the steps are computed by decided terms, as in every example here;
the change of model pays for the move (\S\ref{sec:dich}).  A second cost comes with
ATI, the omniscience rule of \oga.  \rga and \iga are
\emph{proof-grounded}: each is sound and complete for a semantics in
which a universal is true only through a certificate that is one of
the system's own proofs (\S\ref{sec:dich} defines the term and its
counterpart, \emph{omega-grounded}).  With ground induction the
dichotomy derives ATI, and in a system with ATI whose theorems are
recursively enumerable and include no false claim of nontermination,
no semantics of that kind is sound.  The system with both new rules,
\dinga (\inga with the dichotomy), is sound for $\Mstar$
(\cref{thm:price}) and, by \cref{thm:nocert}, for no semantics of
that kind, and neither is \oinga.  \inga, by contrast,
is proof-grounded (\cref{thm:inga-pg}): sound and complete for a
semantics of that kind, and the first system of this paper and its
companion papers to be so while not sound for the least model; that
semantics alone, with no $\omega$-model (no classical model of the kind
the least model and $\Mstar$ are), also proves \inga consistent.  Read
together, the two theorems say that Park's rule keeps
proof-groundedness and that the dichotomy, through ATI, destroys it.
So \dinga's verdicts on $\nu$-atoms are grounded only in
an extended sense: they rest on $\Mstar$, which reads a $\nu$-body at
its greatest fixed point and a universal through all its instances,
and not on a terminating or well-founded process.
\Cref{fig:systems} shows the six systems.

\begin{figure}[t]
\centering
\begin{tikzpicture}[x=1cm,y=1cm,
    sys/.style={draw, rounded corners, inner sep=3pt, align=center, font=\small},
    lab/.style={font=\scriptsize, fill=white, inner sep=1pt}]
  \fill[gray!12] (-5.6,3.75) rectangle (5.6,7.6);
  \node[font=\scriptsize, anchor=north east, align=right] at (5.5,7.5)
    {sound for $\Mstar$;\\unsound for the least model by design};
  \node[font=\scriptsize, anchor=south east, align=right] at (5.5,-0.4)
    {sound for the least model\\and for $\Mstar$};
  \node[sys] (rga)   at (0,0)     {\rga\\[-1pt]{\scriptsize proof-grounded}};
  \node[sys] (iga)   at (-2.8,1.4) {\iga\\[-1pt]{\scriptsize proof-grounded}};
  \node[sys] (oga)   at (2.8,1.4)  {\oga\\[-1pt]{\scriptsize omega-grounded}};
  \node[sys] (oiga)  at (0,2.8)   {\oiga\\[-1pt]{\scriptsize omega-grounded}};
  \node[sys] (inga)  at (-3.6,4.7) {\inga\\[-1pt]{\scriptsize proof-grounded}};
  \node[sys] (oinga) at (1.8,5.1)  {\oinga\\[-1pt]{\scriptsize no certificate-grounded semantics}};
  \node[sys] (dinga) at (-1.6,6.9) {\dinga\\[-1pt]{\scriptsize no certificate-grounded semantics}};
  \draw (rga) -- node[lab] {$\mu$} (iga);
  \draw (rga) -- node[lab] {ATI} (oga);
  \draw (iga) -- node[lab] {ATI} (oiga);
  \draw (oga) -- node[lab] {$\mu$} (oiga);
  \draw (iga) -- node[lab] {$\nu$} (inga);
  \draw (oiga) -- node[lab] {$\nu$} (oinga);
  \draw (inga) -- node[lab] {ATI} (oinga);
  \draw (inga) -- node[lab] {Dich} (dinga);
  \draw[dashed] (oinga) -- node[lab, align=center] {Dich: same\\theorems} (dinga);
\end{tikzpicture}
\caption{The six systems over \rga, ordered by inclusion of their
rules; an edge is labeled with the rule it adds ($\mu$, ground
induction; ATI, \oga's rule; $\nu$, coinduction; Dich, the
dichotomy).  \dinga has no ATI rule but derives it
(\cref{thm:ati-dich}), so \oinga with the dichotomy, the dashed edge,
has \dinga's theorems (\cref{thm:price}).  The shaded systems are
sound for $\Mstar$ and unsound for the least model by design; the
others are sound for both.  Under each name: whether the system has a
semantics whose certificates for universals are its own proofs
(\S\ref{sec:dich}); \oga's label is an instance of \cref{thm:nocert}
with soundness for the least model, and \inga's is \cref{thm:inga-pg}.}
\label{fig:systems}
\end{figure}

\paragraph{The examples.}
The examples are the reason for the rules.  Streams are terms the
system proves total, their prefix sets the derived view, and the round trip between the
two is proved.  Runs of a transition system are streams, and ``some
run is infinite'' is a $\nu$-atom.  Every decided tree is well-founded
at its root or has an infinite branch, with no bound on the
branching; the branch is a term, total in $\Mstar$, and K\H{o}nig's
lemma is the bounded corollary.  Random assignment terminates by an
arithmetic argument and has no natural-number variant, and its
divergent twin runs forever.  Fairness enters through Apt and
Olderog's transformation of a fair program into one with random
assignment: a fair loop terminates, a twin of the loop separates
strong from weak fairness as accessibility against the path atom, and
the alternating-bit protocol over lossy channels delivers every
message.  For a search, Nakata and Uustalu's gap between termination
and non-divergence closes in two independent ways.  The winning
regions of an open game are determined, and in the Baire instance
player II's strategy is a term, total in $\Mstar$.  Cinderella wins
at capacity three by one coinductive step whose certificate carries
her round-robin strategy, while the Stepmother's loss is a
refutation.  The decidedness of $\mathcal O$, Kleene's set of
notations for computable ordinals, and the split of a tree into its
scattered part and its perfect kernel, the shape of the
Cantor--Bendixson theorem, close the list.

\paragraph{Coinduction and reactive systems.}
A reader from verification will want to know what coinduction buys
for reactive systems.  On a given trace, safety and liveness are both
arithmetic, and the proof of a liveness claim is a ranking argument.
The new rules enter with unbounded choice: when the trace is not
given, the controller's forced eventuality is a $\mu$-atom with a
universal inside its body, its negation is the $\nu$-atom of the
dual, and the dichotomy is determinacy of the reachability game.
\S\ref{sec:liveness} says this with a table, after the examples have
supplied the cases.

\paragraph{Where each rule is spent.}
Isabelle/HOL checks every theorem in the paper, and the theories of
the examples end with a build-time check that records which rules a
proof applies and which it does not.  The map that results is a
finding in its own right.  In the examples, asserting a termination
claim, a $\mu$-atom, is arithmetic: folding plus number induction;
the examples reach ranks below $\omega^\omega$, and the \iga paper's
ladder of ordinal notations exists because at higher ranks this
fails.  Using a region, by unfolding or searching inside it, is
arithmetic too.  A $\nu$-region costs one coinductive step when it is
established.  The proofs of the dichotomies themselves (accessibility
or a path, termination or divergence, determinacy, the bucket game,
the scattered part or the kernel) and of the decided $\nu$-regions
apply neither the coinduction rule nor ATI as a rule; the decidedness
of $\mathcal O$ applies ATI, a rule the calculus derives.
\S\ref{sec:rungmap} gives the map and says exactly what the checks
certify.

\paragraph{Reading the metatheory.}
Five kinds of claim appear.  A \emph{grounded} claim is a derivation
in one of the systems.  A \emph{classical} claim is a HOL theorem
about a system or about the least model.  A claim \emph{in $\Mstar$}
is a classical claim about the specific model that reads $\nu$-bodies
at their greatest fixed points; the paper never presents it as truth
without that qualifier, because $\Mstar$ is not the least model.  A
\emph{prediction} is a claim the paper expects but does not prove.
And a \emph{reading} interprets a proved statement, a rank or a count
read off a proof, and is not formalized.  Every theorem header in the
main text names its kind; the appendices, which restate the theorems
with every hypothesis, name it where a clause's kind is not clear from
its wording.

\paragraph{AI disclosure.}
The machine-checked proofs and the first drafts of this paper's text
were produced with substantial assistance from Claude models from
Anthropic, working under the direction and review of the author (the
artifact paragraph of Section~\ref{sec:concl}).  The definitions that a
reader must trust are the rows marked ``def.''\ in
Appendix~\ref{sec:table}, which lists every result with its Isabelle
name; the author reviewed these definitions and theorems.

\paragraph{Road map.}
\S\ref{sec:bg} recalls \rga and \iga.  \S\ref{sec:regions} reads a
body at both fixed points and finds the four regions a point can lie
in.  \S\ref{sec:nu} adds $\nu$: the former, the rule, the model, the
systems.  \S\ref{sec:dich} adds the dichotomy, prices it, and shows that \inga,
without it, keeps a certificate-grounded semantics.
\S\ref{sec:streams} builds streams, runs and branches.
\S\ref{sec:showpieces} works the examples, fairness among them.
\S\ref{sec:liveness} says where coinduction adds truth.
\S\ref{sec:rungmap} is the map.  \S\ref{sec:related} places the work
and \S\ref{sec:concl} lists what is open and where the artifact is.
Theorem headers carry at most two Isabelle names and the kind of
claim; Appendix~\ref{sec:table} carries every name.

\section{Background: \rga and \iga}
\label{sec:bg}

This section recalls what the rest of the paper rests on.  The
systems are in the \rga paper~\cite{rga-paper} and the \iga
paper~\cite{iga-paper}; nothing here is new.

\paragraph{Terms as computations.}
\rga has one syntactic category: terms.  A term is built from
numerals, variables, successor, a recursor, a single binder
$\lambda$, application under call-by-name, equality, negation,
disjunction, and the universal combinator $\cA$.  There are no
separate formulas.  A term computes by head reduction; a closed term
may reach a numeral, which is then its value, or run forever.  The
logical constants are read on values: $0$ and $1$ are the truth
values, so $a = b$ has value $1$ when both sides reach the same
numeral, $0$ when they reach different ones, and no value when
either side diverges.  Negation flips the value; a disjunction is
true when either disjunct is and false when both are.  A universal
$\cA\cdot u$ is true when every instance $u\cdot\bar n$ is and false
when some instance is; this reading is infinitary, a well-founded
tree of evaluations rather than one evaluation, and it is what
``grounded'' means for a universal.  A term with no value is neither
true nor false.  Conjunction is defined, $a \wedge b$ as $\neg(\neg a
\vee \neg b)$, and so is the existential, $\exists$ as $\neg\cA\neg$.
The constant $\bot$ names a divergent term, a term with no value,
never a falsehood.

Three senses of ``has a value'' recur, and the paper keeps them apart.
A closed term \emph{computes} a numeral when head reduction reaches
it.  It has a value \emph{in the least model} below when the
infinitary reading of universals gives it one: every computed value
is one, but a term that tests a universal can have a value that no
reduction reaches.  And from \S\ref{sec:nu} on, a term may have a
value only \emph{in $\Mstar$}, the model of the coinductive systems.
A term is \emph{total} in a model when it has a numeral value there
at every numeral argument, and the paper always says which model is
meant; ``total'' never means ``computes''.  Likewise a term is
\emph{decided} in a model when it has a truth value there, which is a
semantic fact and not an algorithm.

\paragraph{Judgments and rules.}
A judgment is a sequent $\Gamma \ent c$ with a finite set of
hypothesis terms on the left; a context is \emph{well-formed} when
its hypotheses are.  Three forms recur: $\Gamma \ent c$, that $c$ is
true; $\Gamma \ent c\,\mathbf N$, that $c$ reaches a numeral; and
$\Gamma \ent c\,\mathbf B$, that $c$ is \emph{decided}, reaching $0$
or $1$.  The proof system has forty-two rules.  Universal
introduction is schematic: to prove $\cA\cdot u$ one proves $u\cdot v$
for a fresh variable $v$ under the hypothesis that $v$ is a number.
Conversion is a rule in both directions: a term and its head reduct
are interderivable.  Disjunction elimination is available into any
judgment.  Primitive recursive functions enter through a verified
compiler from \emph{index trees}, explicit syntax trees naming a
function by its constructors, and \rga proves every compiled term
total.  Two defined forms recur: $\Dec\,x$ says that $x$ is decided at
every number, $\cA\cdot\lambda n.\,(x\cdot n)\,\mathbf B$, and $\Tot\,f$
that $f$ is total, $\cA\cdot\lambda n.\,(f\cdot n)\,\mathbf N$.  A
binary relation $R$ is decided when $R(y,x)$ is decided at every pair
of numbers; the paper writes this $\Dec\,R$ as well.  Throughout,
$R(y,x)$ reads ``$y$ is one step from $x$'': $y$ is an $R$-predecessor
of $x$, and a chain of steps from $x$ descends.

\paragraph{Self-reference.}
The fixed-point former is a term: $\fix\,f$ is $(\lambda x.\,f(x\,x))
(\lambda x.\,f(x\,x))$ and head-reduces to $f(\fix f)$.  The Liar is
$\fix(\lambda x.\,\neg x)$.  It has no value in the least model, and,
classically, \rga derives neither the Liar nor its negation.

\paragraph{The least model.}
The \emph{least model} reads a universal as true when every instance
is and false when some instance is, and reads everything else by
evaluation.  It is the least fixed point of Kripke's
jump~\cite{kripke75outline} for \rga's language, the operator that
extends a partial assignment of values by one round of evaluation,
and the name comes from that; \rga is sound for it~\cite{rga-paper}.
The \iga and G\"odel papers~\cite{iga-paper,goedel-paper} call the
same model the omega model, which is where the omega in
omega-grounded (\S\ref{sec:dich}) comes from.
Every value the least model assigns persists in the model of \S\ref{sec:nu}.

\paragraph{What \rga lacks.}
The deduction theorem holds exactly at decided hypotheses: from $a,
\Gamma \ent b$ one may conclude $\Gamma \ent \neg a \vee b$ when $a$
is decided, and only then~\cite{ford-conrga}; that paper develops
the reasoning \rga carries out about its own derivations.  Conditional proof at
an arbitrary hypothesis is inconsistent, by Curry's paradox, since
the fixed-point former builds the Curry term~\cite{iga-paper}.  The
consequence that shapes this paper: a premise may be a hypothetical
derivation at an undecided hypothesis, but no term expresses that
conditional.  Both the coinduction rule and the dichotomy rule below
are built around that asymmetry.

\paragraph{Bodies and the inductive reading.}
A \emph{body} is a term $\Phi$ with two distinguished free variables,
the recursive predicate variable $a$ and the point $x$; any further
free variables are \emph{parameters}.  The mechanization uses
de~Bruijn indices; the paper writes named variables.  The
\emph{$\mu$-atom} of a body is the closed term $\mu\Phi \equiv
\fix(\lambda a.\,\lambda x.\,\Phi)$, and $\mu\Phi\cdot x$ unfolds in
three head steps to $\Phi$ with $\mu\Phi$ for $a$ and $x$ for the
point.  (The \iga paper calls it the ground fixed point; the paper
writes $\mu$ because the greatest fixed point is about to appear.)  A
body is \emph{positive} when $a$ occurs only under an even number of
negations and every recursive atom $a\cdot u$ has as its argument a
numeral expression in the bound variables: the point, a
quantifier-bound variable, zero, or a successor of one of these.  A
parameter may not appear in the argument of a recursive atom; a
computed argument $t$ is written through an equality guard, $\cA\cdot
\lambda y.\,\neg(y = t) \vee a\cdot y$ or its existential form.  The
restriction is forced by the semantics: after substitution a
parameter is an arbitrary term, possibly without a value, and no
reading of the body can consult the recursive predicate at it.

\iga adds one rule, \emph{ground induction}.  From $\Gamma \ent
\mu\Phi\cdot x$, a closure argument at a fresh point $i$ (from $i\,
\mathbf N$ and the body with a motive $c$ in place of $a$, derive the
motive at $i$), and $\Gamma \ent x\,\mathbf N$, conclude the motive at
$x$.  It is Park's fixed-point induction~\cite{park69fixpoint} in
grounded form, with parameters allowed under the grammar above.
Semantically, the truth set of $\mu\Phi$ in the least model is the
least fixed point of the operator that reads the body with the
recursive atoms interpreted as membership in a set~\cite{iga-paper}:
a positive self-referential term is an inductive definition.  \iga is
sound for the least model and consistent.  Three facts about the
$\mu$-atom recur below.  A $\mu$-atom is introduced by folding: one
backward unfolding, by the conversion rule, plus whatever number
induction the body needs, which is arithmetic.  It is consumed by
ground induction.  And it is false only through a well-founded
refutation; at a point where no finite unfolding verifies it and none
refutes it, as at the start of a search that never succeeds, it has no
value, not the value $0$.

\paragraph{The sibling \oga and the lattice of rules.}
\oga~\cite{oga-paper} is \rga with one rule, \emph{ATI}, the \oga
paper's name for omega-grounded universal decidedness introduction
(\S\ref{sec:dich} defines omega-grounded):
from a schematic proof that every instance of a universal's matrix is
decided, conclude that the universal is decided.  Together with the
other rules it acts as Bishop's limited principle of
omniscience~\cite{bishop67foundations}: a decided search either
succeeds at some number or fails everywhere, and the system derives
that disjunction.  The system need not derive either disjunct: by a
classical theorem (\isn{oga\_godel\_fict}), \oga derives the
decidedness of a universal that it neither derives nor refutes.
\oiga is the union of \oga's and \iga's rules.
The development states each rule once, as an assumption of an
abstract calculus (an Isabelle locale), and builds the systems by
combining assumptions; a theorem proved in such a calculus holds in
every system that has its rules.  The paper calls the family of
calculi so obtained the \emph{lattice of rules} and names, for each
theorem, the rules its proof uses.  ``Arithmetic'' below always means
``using \rga's rules alone'': number induction, conversion, folding
and unfolding, and no fixed-point rule.

\paragraph{Notation.}
$\bar n$ is the numeral for $n$.  $\Acc_R$ is the accessibility term
of the \iga paper, the $\mu$-atom of the body $\cA\cdot\lambda y.\,
\neg R(y,x) \vee a\cdot y$: every chain of steps from $x$ is finite.
$c[t]$ is the certificate or motive $c$ with $t$ for the point.
Isabelle names appear as \isn{name} in theorem headers and in
Appendix~\ref{sec:table}.

\section{A Body Read at Both Fixed Points}
\label{sec:regions}

Before any new syntax, the least model already has a coinductive
reading of every body, and comparing it with the inductive one
motivates the design.  Every claim in this section is classical: a
HOL theorem about operators on sets of numbers.

\paragraph{Two operators.}
Fix a positive body $\Phi$ with $\mu\Phi$ closed.  Read the body at a
point $n$ with its recursive atoms interpreted through a set $S$: an atom $a\cdot u$ is true when the instance of $u$ is a numeral $\bar m$ with $m \in S$, and has no value otherwise.  The \emph{truth operator} $F_\Phi$ sends $S$ to the
set of $n$ at which the body is true under this reading, and the \iga
paper proves that $\mu\Phi$'s truth set is $\lfp\,F_\Phi$.  Now read
the body with an atom \emph{false} when the instance of $u$ is a numeral $\bar m$ with $m \in X$, and valueless otherwise; the \emph{falsity operator} $G_\Phi$ sends
$X$ to the set of $n$ at which the body is false under that reading.
Under these two three-valued readings both operators are monotone,
for every body, positive or not, because each reading consults the
set in one direction only.

\begin{thm}[\isn{gfix\_omv\_false\_eq\_lfp}, \isn{gfp\_lfp\_disjoint}; classical]
\label{thm:disjoint}
The set of $n$ at which $\mu\Phi\cdot\bar n$ is false in the least
model is $\lfp\,G_\Phi$.  For every body, $\gfp\,F_\Phi$ and
$\lfp\,G_\Phi$ are disjoint.
\end{thm}

So a point $n$ lies in one of four regions.  In $\lfp\,F_\Phi$ the
$\mu$-atom is true, by a well-founded verification.  In
$\lfp\,G_\Phi$ the atom is false, by a well-founded refutation.  In
$\gfp\,F_\Phi \setminus \lfp\,F_\Phi$ the atom has no value in the
least model, yet reading it as true at every point of that region is
consistent with the body: that is the third region, where the
coinductive reading adds truth.  Every other point is in the fourth
region, where nothing can be added.  A point of $\gfp\,F_\Phi$ is
never false in the least model, by the two halves of
\cref{thm:disjoint}.  \Cref{fig:regions} draws the four regions.

\begin{figure}[t]
\centering
\begin{tikzpicture}[scale=0.9, every node/.style={font=\small}]
  \draw[rounded corners] (0,0) rectangle (11,4);
  \node[anchor=north west] at (0.1,3.95) {$\mathbb N$};
  \draw (3.2,2) ellipse (2.9 and 1.6);
  \node at (3.2,3.25) {$\gfp\,F_\Phi$};
  \draw[fill=gray!20] (2.3,1.7) ellipse (1.6 and 0.95);
  \node at (2.3,1.95) {$\lfp\,F_\Phi$};
  \node[font=\scriptsize] at (2.3,1.45) {(1) true};
  \node[font=\scriptsize, align=center] at (4.85,1.75) {(3) no value;\\coinductive\\truth};
  \draw[fill=gray!20] (8.4,2.1) ellipse (1.9 and 1.15);
  \node at (8.4,2.45) {$\lfp\,G_\Phi$};
  \node[font=\scriptsize] at (8.4,1.95) {(2) false};
  \node[font=\scriptsize, align=center] at (6.6,0.45) {(4) no value; nothing to add};
\end{tikzpicture}
\caption{The four regions of a body $\Phi$ in the least model
(\cref{thm:disjoint}): (1) verified, (2) refuted, (3) in the greatest
fixed point of the truth operator but not the least, where the
coinductive reading adds truth, and (4) the rest.}
\label{fig:regions}
\end{figure}

\paragraph{Why $\neg\mu$ is not $\nu$ of the dual.}
Take the search body $q(x) \vee a\cdot(x+1)$ for a decided test $q$
that is false at every number.  Its $\mu$-atom at $0$, ``the search
from $0$ terminates'', has no value in the least model: no finite
unfolding makes it true, and nothing refutes it, since refuting it
would need a finite reason.  So $\neg\mu\Phi\cdot 0$ is not true.  The
greatest fixed point of the body $\neg q(x) \wedge a\cdot(x+1)$, the
co-body of the search body, contains $0$: ``the search from $0$ runs
forever''.  Classically the
two are the same statement.  In a paracomplete logic the first has no
value and the second is a new truth.  That is the gap the dichotomy
rule of \S\ref{sec:dich} will close, and only for bodies whose leaves
are decided.

\paragraph{When the two readings are complementary.}
Classically, $\gfp\,F_\Phi$ is the complement of $\lfp\,G_\Phi$.  In
the paracomplete setting that is a property of the body, and it is
exact.

\begin{thm}[\isn{gfp\_eq\_compl\_iff}; classical]
\label{thm:compl}
$\gfp\,F_\Phi$ is the complement of $\lfp\,G_\Phi$ if and only if, for
every $n$, either the body is true at $n$ with its atoms read as true
on the complement of $\lfp\,G_\Phi$, or $n \in \lfp\,G_\Phi$.
\end{thm}

A sufficient condition is that the body is \emph{definite}: at every
set $S$ and point $n$, the body is true with the atoms read in $S$ or
false with the atoms read in the complement of $S$.  A syntactic
grammar guarantees definiteness: positivity with every non-recursive
leaf a term that has a value.  Definiteness is not necessary: the body
$a\cdot 0 \vee L$ with $L$ the Liar is complementary without being
definite.

\paragraph{Three examples.}
The \emph{truth-teller} body $a\cdot x$ has $\lfp\,F = \emptyset$ and
$\gfp\,F = \mathbb N$: nothing is verified, nothing is refuted, and
every point lies in the third region.  The coinductive reading
colors every instance true.  Kripke's construction of
truth~\cite{kripke75outline} leaves the truth-teller ungrounded in its
least fixed point, and other fixed points make it true or false at
will; the reading here makes a choice, and \S\ref{sec:nu} says which.
A body with the Liar as a conjunct, $L \wedge a\cdot x$, has nothing
refuted and nothing coinductively true: every point lies in the
fourth region.  This body shows why ``color the never-refuted true''
and ``read at the greatest fixed point'' differ: the first would make
$L \wedge a\cdot x$, and so the Liar, true at every point, while the
second adds nothing.  The two coincide exactly under the condition of
\cref{thm:compl}.  Finally, the accessibility body: $\gfp\,F_{\Acc_R}$
is everything, for every closed relation, and its falsity set is
empty.  So the coinductive reading of accessibility has no content,
and the content lies in the $\nu$-atom of the co-body, ``there is
an infinite descent from $x$'', which \S\ref{sec:dich} calls
$\Path_R$.  The same holds of games: the winning-position body of the
\iga paper read at its greatest fixed point is the survival game, the
positions from which a player is never caught (the player
\emph{survives} there), and for the body that moves up by one, every
position survives, none is won, and every position has no value in
the least model.

\section{The Greatest Fixed Point: Former, Rule, Model}
\label{sec:nu}

\subsection{Separate syntax}

The $\nu$-atom of a body is a term, as the $\mu$-atom is, but a
different one.  $\nu\Phi$ is the same self-application as $\mu\Phi$
with an identity redex inserted under the abstraction, $\nu\Phi = s\cdot
s$ with $s = \lambda z.\,I\cdot((\lambda a.\lambda x.\,\Phi)\cdot(z\cdot
z))$ and $I = \lambda y.\,y$, so that $\nu\Phi$ and $\mu\Psi$ are
distinct terms for all bodies $\Phi$ and $\Psi$, and $\nu\Phi\cdot x$ head-reduces in four steps to the
unfolding of the body at $x$ with $\nu\Phi$ for the recursive
variable.  So $\nu$-elimination is unfolding, by the conversion rule,
exactly as for $\mu$.  A reader from programming languages will ask
why a tag that reduction can see through should change meaning.  The
model of \S\ref{sec:nu-model} does not reduce a $\nu$-atom.  It
recognizes the atom syntactically, as a head of its own, and reads
every application whose head evaluates to it from an oracle.  The
oracle is built so that the atom's value equals the value of its
unfolding, which is why the conversion rule stays sound.  And no
sequence of head steps turns a $\nu$-atom of the class into a
$\mu$-atom, provided the body is shaped as a formula, its outermost
constructor a negation, disjunction, equation, universal or recursive
atom; a body that is itself a $\mu$-atom applied to the point reaches
one in four steps.

The separation is forced.  Suppose one former carried both
principles: ground induction, and the coinduction rule below.

\begin{thm}[\isn{nu\_gfix\_clash}; classical]
\label{thm:clash}
Any proof system whose steps include \rga's rules, the coinduction
rule for $\nu$-terms and ground induction at $\nu$-terms, and whose
step relation is closed under extending the proof list, derives
$0 = 1$.
\end{thm}

The derivation is five lines, and it is the first worked use of the
rule.  Its witness is the all-zeros atom rather than the shorter
truth-teller because the same atom is the valueless theorem of
\cref{thm:unsound}.  Let $Z$ be the $\nu$-atom of the body $(\lambda y.\,0)\cdot x
= 0 \wedge a\cdot(x+1)$: ``the constant-zero function is zero at $x$,
and again from $x+1$ on''.  Take the certificate $c \equiv (x = x)$.
It holds at $0$.  At a fresh point $i$ with $c[i]$ as hypothesis, the
body with $\lambda c$, the certificate abstracted over its point, for
$a$ reads $(\lambda y.\,0)\cdot i = 0 \wedge
(i + 1 = i + 1)$, which conversion and reflexivity prove.  So the
coinduction rule gives $Z\cdot 0$.  Now apply ground induction at $Z$
with the point-free motive $0 = 1$: the closure premise asks, at a
fresh $i$ with the body for $Z$ in which the motive replaces the
recursive atom, that is, with $(\lambda y.\,0)\cdot i = 0 \wedge (\lambda z.\,0 = 1)\cdot(i+1)$ as hypothesis, for $0 = 1$, which the hypothesis' second conjunct gives by one conversion step.  The rule concludes $0 = 1$ at $0$.  The theorem uses
no $\omega$-model, and it is the strongest argument for the design: the
least and the greatest fixed point act on distinct terms, each read
at its own extremal point in one model, and no stratification of
the rules is needed.  Isabelle, Coq and Agda keep inductive and
coinductive definitions apart, in Coq and Agda for reasons of
normalization and guardedness as well; the theorem makes one reason
for the separation exact for this calculus.

\subsection{The coinduction rule}

Let $\Phi$ be a positive body in which the recursive variable occurs,
$c$ a well-formed certificate (a term in which the point is the one
distinguished free variable), and $i$ a variable fresh for the hypotheses, the certificate,
and the body's parameters.  The rule reads:
\[
  \frac{
    \Gamma \ent c[x]
    \qquad
    i\,\mathbf N,\ c[i],\ \Gamma \ent \Phi[\lambda c/a][i/x]
    \qquad
    \Gamma \ent x\,\mathbf N
  }{
    \Gamma \ent (\nu\Phi)\cdot x
  }
\]
In words: if the certificate holds at $x$, and wherever the
certificate holds at a fresh point the body holds there with the
certificate in place of the recursive predicate, then $x$ lies in the
greatest fixed point.  The certificate set $W = \{y \mid c[y]\}$
contains $x$ and is a post-fixed point, $W \subseteq \Phi(W)$, so this
is the Knaster--Tarski principle in the form of Park's coinduction
rule~\cite{park69fixpoint,park81concurrency}, the exact dual of
ground induction, where the motive is a pre-fixed point.  In a proof
assistant it is Isabelle's
\texttt{weak\_coinduct}~\cite{paulson97mechanizing}, the plain
Knaster--Tarski form, and not proof by guarded corecursion:
the coinductive hypothesis is never available as a hypothesis, so no
guardedness condition arises.

Two features of the rule matter.  The closure premise is a
hypothetical derivation at the hypothesis $c[i]$, which is in general
undecided; the system admits it because it is a premise sequent,
never an implication inside a term, which Curry forbids
(\S\ref{sec:bg}).  And the side conditions are those of ground
induction with one addition: the body must mention the recursive
variable.  The addition makes the body a formula, a negation,
disjunction, equation, universal or recursive atom, which is the
shape the model's reading of a $\nu$-head requires.

Two derived principles follow from the rule and the replacement of
one predicate by a pointwise weaker one inside a positive body
(Appendix~\ref{sec:table}): monotonicity, that a $\nu$-atom implies
the $\nu$-atom of a pointwise weaker body, and coinduction up to
union, where the certificate need only be closed under the body up to
the $\nu$-atom itself.  The second is the form of Isabelle's
\texttt{coinduct} and of the coinduction rules its coinductive package
generates, which users know as strong coinduction; it is the simplest of the up-to
techniques~\cite{pous12enhancements}.  The accumulation of the
hypothesis in parameterized coinduction~\cite{hur13power}, the
companion~\cite{pous16companion}, and coinduction up to equivalence or
up to context are not developed.  They are undeveloped rather than
known to be blocked; the expected obstacle, a prediction, is that a
context applied to a $\nu$-atom must be a term the body reads, so up
to context would need codes read by an evaluator.

\paragraph{A second worked use.}
Let $R'$ step from $0$ to every $m+1$, from $m+2$ to $m+1$, and from
$1$ back to $0$ (the divergent twin of \S\ref{sec:showpieces}).  Write
$\Path_R$ for the $\nu$-atom of the body $\exists y.\,R(y,x) \wedge
a\cdot y$, ``some chain of steps from $x$ is infinite'', which
\S\ref{sec:dich} obtains as the co-body of accessibility.  To show
$\Path_{R'}\cdot 0$,
take the certificate ``every number'', $c \equiv (0 = 0)$.  It holds
at $0$.  At a fresh point $i$, the body of $\Path$ with $\lambda
c$ for $a$ reads ``some $R'$-predecessor $y$ of $i$ has $c[y]$'', and
the predecessor is given by cases on $i$: $1$ if $i = 0$, $0$ if $i =
1$, $i - 1$ otherwise; each case is decided arithmetic.  One
coinduction step concludes.

\subsection{The model}
\label{sec:nu-model}

Soundness needs a model in which a body of the right class is read at
its greatest fixed point and everything else as before.

\paragraph{The $\nu$ class.}
A body is in the $\nu$ \emph{class} when it mentions the recursive
variable and is positive in a wider sense, in which a recursive atom
may take any argument not mentioning the recursive variable.  The
rule's narrower grammar is what its soundness proof needs, and the
wider class is invariant under value congruence, the replacement of a
subterm by a well-formed term of the same value, which one of \rga's
rules requires (Appendix~\ref{app:nu}).  Writing $\nu\Phi$ for a body
of the class is the only act of declaration there is; the model reads
every such head at its greatest fixed point.

\paragraph{The oracle model.}
An \emph{oracle} $\Theta$ assigns to each body of the class a set of
argument terms at which its $\nu$-atom is true and a set at which it
is false.  The oracle model $M_\Theta$ evaluates as the least model
does except at an application whose head evaluates to a $\nu$-atom of the class: there it reads $\Theta$ at the argument term.  It is monotone
in the oracle, and deterministic when the oracle's two sets are
disjoint.  Each body of the class has a truth operator and a falsity
operator relative to $M_\Theta$, defined as in \S\ref{sec:regions}
with $M_\Theta$ for the least model and argument terms for numerals.  The oracle is
\emph{coherent} when its two sets are disjoint, at every body of the
class its truth set is a fixed point of the truth operator and its
falsity set a fixed point of the falsity operator, and both sets are
empty at every other body; in the model of a coherent oracle a
$\nu$-atom has exactly the value of its unfolding.  The \emph{oracle
map} sends $\Theta$ to the oracle that reads each body at the
greatest fixed point of its truth operator and at the least fixed
point of its falsity operator; the map is monotone in the information
order (both sets grow), and $\Theta^*$ is its least fixed point, a
coherent oracle.  Which non-minimal fixed point of Kripke's
construction this is, the polarity of the declaration decides, that
is, whether a term is written with $\nu$'s tag or $\mu$'s: the
greatest for truth, the least for falsity, body by body.

\begin{defi}
$\Mstar$ is the oracle model at $\Theta^*$.
\end{defi}

Every value of the least model persists in the model of every
coherent oracle, $\Mstar$ among them: a coherent oracle only adds
verdicts.  The soundness proof of \rga's forty-two rules is carried
out once against an interface of twenty-four facts about a model, and
the least model satisfies the interface.  So does the model of every
coherent oracle whose sets are closed under substitution and under
value congruence, of the argument and of the body itself
(Appendix~\ref{app:nu}): the facts, among them value congruence and the preservation of
values by head steps in both directions, are proved once for every
such oracle.  $\Theta^*$ is one.  The closed-off oracle, which reads
each body at the least fixed point of its truth operator and the
greatest of its falsity operator, the exact dual of $\Theta^*$, is
another, and the \iga paper's closed-off model is built on
it~\cite{iga-paper}.

\begin{rem}
\label{rem:kind}
$\Mstar$ is a non-minimal fixed point of the kind Kripke's
construction allows, and not an intrinsic one, a fixed point whose
verdicts no other fixed point contradicts~\cite{kripke75outline}.  The
$\nu$ truth-teller, the $\nu$-atom of the body $a\cdot x$, is true at
every numeral in $\Mstar$, whose truth set for that body is the
greatest fixed point of its truth operator.  The model of the
closed-off oracle also satisfies the interface, and there the truth
set is the least fixed point, empty, and the falsity set the
greatest, so the same atom is false at every numeral; the all-zeros
atom of \cref{thm:clash} behaves the same way.  And the value depends
on the declaration: the $\mu$ truth-teller, the same self-application
without the identity redex, has the same infinite unfolding and no
value in $\Mstar$, since its unfolding never reaches a value and no
oracle reads a $\mu$-head.
So the question of the introduction, which non-terminating
computations deserve the color, is answered for declared bodies: a
computation is colored when it is written as a $\nu$-atom of the
class and lies in the greatest fixed point of its body.  These facts
are classical, by an argument from the two oracles' fixed-point
equations, and are not separate theorems of the development.
\end{rem}

\begin{thm}[\isn{rga\_nu\_sound}; classical]
\rga's rules are sound for $\Mstar$.
\end{thm}

For \inga's two further rules: ground induction is sound for every model of the interface that satisfies one more fact, the induction property of its $\mu$-atoms,
and $\Mstar$ satisfies it by an induction on the model's evaluation.
The coinduction rule's case is the content of the construction: the
set of well-formed, valued argument terms at which the certificate is
true is a post-fixed point of the body's truth operator relative to
$\Mstar$, hence below its greatest fixed point, hence inside
$\Theta^*$.

\subsection{The systems}

\inga is \iga with the coinduction rule; \oinga is \oiga with the
coinduction rule.  The names are compositional: an N for $\nu$ is
inserted before the GA of the system extended.  \S\ref{sec:dich} adds a third
system, \dinga, which is \inga with the dichotomy rule.  Two notions are the paper's
own.  A system is \emph{semantically conservative} over the least
model when a sentence that has a value in the least model, and is
derivable from hypotheses true in the least model, is true in the
least model; no proof-theoretic conservativity over a fragment is
claimed.  A system is \emph{unsound for the least model by design}
when it derives a sentence with no value in the least model; by
conservativity it never derives a false one.

\begin{thm}[\isn{inga\_nu\_sound}, \isn{oinga\_nu\_sound}; classical]
\label{thm:inga}
\inga and \oinga are sound for $\Mstar$, consistent, and
semantically conservative over the least model.
\end{thm}

\begin{thm}[\isn{inga\_Z}, \isn{nuZ\_true}; classical]
\label{thm:unsound}
\inga derives the all-zeros atom $Z\cdot 0$ of \cref{thm:clash}, which
is true in $\Mstar$ at every numeral and has no value in the least
model at any.
\end{thm}

The least model leaves valueless every instance of a body of the
class that lies in its greatest fixed point but not its least; the
systems derive some of them, and by \cref{thm:inga} nothing false
about the arithmetic world enters with them, only verdicts about
never-refuted infinite behavior.  Both theorems hold for \oinga, and
\oiga itself is sound for $\Mstar$, so every theorem of \iga and
\oga is true there.

\paragraph{What totality means after $\nu$.}
In the three senses of \S\ref{sec:bg}: in \rga and \iga a derivable
$\Tot\,f$ says that $f$ is total in the least model, and the \rga
paper shows what that buys.  After $\nu$ and the dichotomy, the
judgments $\mathbf N$ and $\mathbf B$ are claims in $\Mstar$.  A term
that searches through a $\nu$-atom may be total in $\Mstar$ and have
no value in the least model, because the atom has none there.
Sentences valued in the least model keep their meaning, by
conservativity.  For every term the examples construct, the paper
says in which model it is total.

\section{The Dichotomy}
\label{sec:dich}

\subsection{The dual, the co-body and the decidedness premise}

The \emph{dual} of a body negates every recursive atom literally and
nothing else: $\dual{\Phi}$ is $\Phi$ with each $a\cdot u$ replaced
by $\neg(a\cdot u)$.  It exists for every positive body, is unique,
commutes with substitution for parameters, and its negation
$\neg\dual\Phi$ is again a positive body, which the paper calls the
\emph{co-body} of $\Phi$.  Classically, $\neg\mu\Phi$ is the
$\nu$-atom of the co-body, $\nu(\neg\dual\Phi)$.  For the
accessibility body, the co-body is $\neg\cA\cdot\lambda y.\,\neg
R(y,x) \vee \neg(a\cdot y)$: ``some $R$-predecessor of $x$ is in the
set''.  Its $\nu$-atom is the term $\Path_R$: there is an infinite
chain of steps from $x$.

A grounded system has no rule that turns a failed verification of the
$\mu$-atom into a proof of the $\nu$-atom (\S\ref{sec:regions}).  The
rule that joins them needs a premise.  A \emph{decidedness premise}
$D$ of a body is a sentence saying that every non-recursive leaf of
the body is decided at all numeral values of the point and the bound
variables.  For the accessibility body it is $\Dec\,R$, with a
trivial conjunct $0=0$ for the recursive atom, and for the game body
of \S\ref{sec:showpieces} it is the decidedness of the three
parameters.  The premise is a term, derivable like any other.  The
dual is computed from the body by a syntactic recursion, and a
premise by a second one, which at a compound part without the
recursive variable may stop at that part as one leaf or descend into
it, so a body can have several premises and the rule accepts each of
them.  Because the premise is computed from the body, the rule needs
no semantic side condition, and a proof checker can check every
application of it.

\subsection{The rule and its characterization}

For a positive body $\Phi$ mentioning its recursive variable, with
dual $\Psi$ and decidedness premise $D$:
\[
  \frac{\Gamma \ent D \qquad \Gamma \ent x\,\mathbf N}
       {\Gamma \ent (\mu\Phi)\cdot x \vee (\nu(\neg\Psi))\cdot x}
\]
The rule asserts the classical disjunction at one body and one
point.  Which side holds is not decided by the rule, and the rule
never refutes the $\mu$-atom: the falsity of a $\mu$-atom needs a
well-founded refutation, which the rule does not supply.  For
accessibility no system in the paper refutes the atom at a closed
relation from no hypotheses at all (the end of \S\ref{sec:excl}).
Two theorems about $\Mstar$ say what the rule is.

\begin{thm}[\isn{swap\_identity}; in $\Mstar$]
\label{thm:swap}
For every positive body $\Phi$ of the rule's class whose co-body is
in the $\nu$ class, and every $n$: $\mu\Phi\cdot\bar n$ is true in
$\Mstar$ if and only if $\nu(\neg\dual\Phi)\cdot\bar n$ is false in
$\Mstar$.
\end{thm}

The falsity half of the co-body's $\nu$-atom is exactly the
$\mu$-term, with no decidedness needed.  So in $\Mstar$ the rule's
conclusion is excluded middle for one class of $\nu$-atoms, with the
negative half renamed to the $\mu$-term: a semantic excluded middle,
gated by a derivable premise and not valid in the least model, and
not unrestricted excluded middle (\S\ref{sec:price}).  Call the rule
\emph{valid} at a body when its conclusion is true in $\Mstar$ at
every well-formed term $x$, in every well-formed context in which $x$
is a number; the premise $D$ plays no part in validity.

\begin{thm}[\isn{thmC}, \isn{dich\_step\_mjt}; in $\Mstar$]
\label{thm:C}
The rule is valid at a body if and only if the $\nu$-atom of every
substitution instance of the co-body has a value in $\Mstar$ at
every numeral.  The rule, with its premise, is sound for $\Mstar$.
\end{thm}

Soundness holds because, wherever the derivable premise $D$ is true
in $\Mstar$, it secures the bivalence of the $\nu$-atom that validity
demands by the characterization above.  Two controls show that the
premise does the work.  At the body $\bot \vee a\cdot x$ the rule is
invalid at every numeral, and the premise, which asks $\bot$ to be
decided at every point, is true under no substitution, so it blocks
every instance.  And the premise can fail in $\Mstar$ for a reason
the system cannot see.  At a body whose non-recursive leaf is itself
a $\mu$-atom, the leaf can have no value at some instance; in the
trace example of \S\ref{sec:showpieces} the leaf is a reachability
atom, valueless at a target state that an infinite run never
reaches.  The premise then has no value there, and by soundness
neither it nor its negation is derivable.  A decidedness atom is never
false in $\Mstar$, so a premise of the dichotomy is either true or
valueless, never false.

\subsection{The system and its price}
\label{sec:price}

\dinga (D for dichotomy) is \inga with the rule.  Its derivations are
\iga steps, coinduction steps and dichotomy steps, nothing else.
\oga's rule ATI is not among its rules, and it need not be.

\begin{thm}[\isn{g\_ati\_of\_dich}; grounded]
\label{thm:ati-dich}
In every calculus with ground induction and the dichotomy, ATI is a
derived rule: from $\Gamma \ent \cA\cdot\lambda v.\,(u\cdot v)\,\mathbf
B$, $\Gamma \ent (\cA\cdot u)\,\mathbf B$.
\end{thm}

The proof is the dichotomy at the search body $\neg(u\cdot x) \vee
a\cdot(x+1)$, whose computed premise is ATI's premise with a conjunct
$0 = 0$.  The $\mu$-atom at $0$ says that a search for an instance at
which $u$ is false terminates, and ground induction then refutes the
universal; the $\nu$-atom of the co-body, of $u\cdot x \wedge
a\cdot(x+1)$, says that the search runs forever, and number induction
with unfolding then proves every instance.  So adding ATI to \dinga
adds no theorem: \oinga with the dichotomy, where ATI is primitive,
has \dinga's derivability relation (\cref{thm:price}), and results
first proved for that system, the examples of \S\ref{sec:showpieces}
among them, are results about \dinga.
ATI is \oga's omniscience about decided matrices, so, with ground
induction, this omniscience is the dichotomy at one body, and the
dichotomy is the stronger rule.  Over ground induction alone the
inclusion is strict: \oiga, which is \iga with ATI, derives less than
\iga with the dichotomy.  That is a classical fact, by combining two
theorems: \oiga is sound for the least model (\isn{oiga\_omega\_sound}),
and the dichotomy fails there at a closed, decided, ill-founded order
whose premise \rga derives (\isn{dich\_not\_omega\_valid\_decRel}), so
\oiga does not derive that instance.  Whether \oinga derives every
instance of the dichotomy is open: $\Mstar$ validates the rule
wherever its premise holds, so it cannot separate the two systems,
and a separation would need a different kind of model.  Whether the
dichotomy derives ATI without ground induction is also open.  Read
informally, and not as a theorem, the dichotomy is the general
coinductive closing-off of a positive definition with decided leaves,
the assertion that its least fixed point or the greatest fixed point
of its co-body holds, and ATI is its arithmetical special case.
The \iga paper's \emph{class rung}~\cite{iga-paper} closes off a
similar kind of body from the other side.  Here $\mu$-bivalence is
the decidedness of the $\mu$-atom at every number.  The class rung is
\iga with $\mu$-bivalence at every positive body without parameters
that uses its recursive variable and whose other parts have a value
in the least model at every assignment of numbers.  By classical
theorems the class rung is consistent and sound for the closed-off
model, which is built on the closed-off oracle of \S\ref{sec:nu}
(\isn{mubivD\_consistent}, \isn{imb\_closed\_off\_sound}), while
$\mu$-bivalence at every positive body derives $0=1$
(\isn{g\_mubiv\_contra}).  So the two rules close off similar bodies
from opposite sides.  The dichotomy concludes that the $\mu$-atom or
the $\nu$-atom of the co-body holds, which in $\Mstar$ is excluded
middle for the $\nu$-atom (\cref{thm:swap}); the dichotomy never
refutes the $\mu$-atom, and its warrant is $\Mstar$.  The class rung
concludes the decidedness of the $\mu$-atom itself, and its warrant is
the closed-off model, built on the closed-off oracle, the exact dual
of $\Theta^*$.  This comparison is a reading.  The class rung's
condition on the leaves is semantic, a value in the least model, while
the dichotomy's premise is a sentence the system derives.  The class
rung also admits only bodies without parameters, while the
dichotomy's bodies may have parameters.  So the two classes of bodies
are similar, not the same.  No theorem of the development relates the
two rules, and a system with both has not been studied.

The \emph{price} of a rule, in the sense of the \iga paper, is what
has to hold for the rule to be sound: there, bivalence of
accessibility for a fixed order was sound for the least model exactly
when the order is well-founded, a hypothesis on the order.

\begin{thm}[\isn{dinga\_dich\_price}, \isn{dinga\_eq\_oingad}; classical]
\label{thm:price}
\dinga is sound for $\Mstar$ and consistent.  In it, for every
relation $R$ and well-formed context, from $\Dec\,R$ and $x\,\mathbf
N$: $\Acc_R\cdot x \vee \Path_R\cdot x$, and $\Path_R\cdot x$ is
decided.  It is semantically conservative over the least model, and
it derives a sentence with no value in the least model.  Its
derivability relation is that of \oinga with the dichotomy.
\end{thm}

Read the dichotomy for accessibility in the least model, with the
$\nu$-atom replaced by the $\mu$-atom of the same co-body, which is
what a $\nu$-term amounts to there, since the least model simply
unfolds it: the dichotomy is then true at every point exactly when
the order is well-founded, the same price as bivalence of
accessibility.  Read in $\Mstar$, it is true at every decided
relation, well-founded or not.  At the relation that steps from every
number to its successor, bivalence of accessibility is unsound for
the least model while the dichotomy is true in $\Mstar$ everywhere.
Bivalence of accessibility asked $\lfp\,F_{\Acc} = \mathbb N$; the
dichotomy asks $\lfp\,F_{\Acc} \cup \gfp\,F_{\Path} = \mathbb N$ with
both operators read in $\Mstar$, which holds at every decided relation
because the body is then definite, the sufficient condition that
follows \cref{thm:compl}.  The dichotomy's premise, $\Dec\,R$, is
still a hypothesis on the order, and \cref{thm:C} shows that it does
work.  What changes is its kind: well-foundedness is a $\Pi^1_1$
property, not derivable in general, while decidedness of the steps is
derivable, in \rga, at every relation whose steps are computed by
decided terms (every example here), and is otherwise a sentence like
any other.  The change of model pays for that move.
\Cref{thm:price} states the cost, a reading of its clauses rather
than a further theorem: soundness for $\Mstar$, nothing changed on
the sentences the least model values, and new verdicts only where it
has none.  This is the price of the dichotomy.

The dichotomy is semantic excluded middle for one class of
$\nu$-atoms (\cref{thm:swap}), not unrestricted excluded middle, in
three ways.  It is not valid for the least model: at every closed
ill-founded order whose premise \rga derives, the conclusion, with the
$\nu$-atom itself, is valueless at some point, so the rule is not
sound for the semantics \rga and \iga are sound for.  It needs a
derivable premise, which can fail by being valueless.  And it decides
the $\nu$-atom but not the $\mu$-atom, which for accessibility is
never false in $\Mstar$, so that no system in the paper decides it
at a point where it is not verified, where classical logic decides
both.

\paragraph{The price in proofs: no certificate-grounded semantics.}
A semantics of a system is \emph{certificate-grounded} when a
universal is true in it only through a certificate, a derivation in
the system of the universal's open instance (a different use of the
word from the post-fixed point of \S\ref{sec:nu}), and false only
through a false instance, a disjunction is true only through a true
disjunct, a negation only when its argument is false, and nothing is
both true and false.  The conditions are necessary conditions, ``only
through'', and that is what makes the theorem below strong: it
excludes every semantics that meets them, whatever else the semantics
does.  A system is \emph{proof-grounded} when it is sound and complete
for such a semantics: \rga is~\cite{rga-paper}, and so is
\iga~\cite{iga-paper}.  It is \emph{omega-grounded} when it is sound
for the least model but for no certificate-grounded semantics, as
\oga and \oiga are.  None of these terms is one of the kinds of claim
of \S\ref{sec:intro}; each names a classical property of a system.
\inga, \oinga and \dinga are not sound for the least model, so none
is omega-grounded; \oinga and \dinga are not proof-grounded
(\cref{cor:nocert}), while \inga is (\cref{thm:inga-pg} below).
The \emph{halting family}
has one universal for each $x$, ``machine $x$ does not halt on input
$x$''.

\begin{thm}[\isn{ati\_cert\_unsound}, \isn{dich\_cert\_unsound}; classical]
\label{thm:nocert}
Let a calculus have ATI, derive a universal of the halting family only
where it holds, and derive a recursively enumerable set of them.  For
every certificate-grounded semantics over the calculus's own
derivations under which \rga's closed theorems are true, the calculus
derives the decidedness of some universal of the family, which the
semantics does not make true.  By \cref{thm:ati-dich} the same holds
for every calculus with ground induction and the dichotomy.
\end{thm}

The reason is a count of complexity.  Suppose the semantics made
ATI's conclusion true at every universal of the family, so that each
of them had a value.  A universal of the family is true only through
a certificate, a derivation, so the true ones are among the derivable
ones, and the calculus derives one only where it holds, so the
derivable ones are true non-halting claims.  A universal of the
family is false only through a false instance, a halting computation.
Every $x$ would then either have a derivable universal, a recursively
enumerable condition that holds exactly when machine $x$ does not halt
on $x$, or halt, also recursively enumerable; so the halting set would
be decidable, which it is not.  With ground induction the dichotomy
derives ATI, and so no certificate-grounded semantics can follow the
dichotomy at the halting family either.  For \dinga the hypotheses are
classical theorems.  A primitive-recursive proof checker for \oinga
with the dichotomy, carried over by the equality of
\cref{thm:price}, makes \dinga's derivability recursively enumerable
(\dinga's own step relation has a primitive-recursive decider), and
semantic conservativity makes it derive a universal of the family
only where it holds.

\begin{cor}[\isn{dinga\_not\_cert\_sound}, \isn{oinga\_not\_cert\_sound}; classical]
\label{cor:nocert}
\dinga and \oinga are sound for $\Mstar$ (\cref{thm:price,thm:inga})
and for no certificate-grounded semantics: their warrant is a model
that reads a universal through its instances, as the least model
does, and no certificate.  \oiga is omega-grounded.
\end{cor}

\noindent For \oinga the hypotheses of \cref{thm:nocert} are again
classical theorems: its derivability is recursively enumerable, and
its conservativity makes it derive a universal of the family only
where it holds; the corollary for \oinga goes through ATI's theorem
directly, not through the dichotomy.  The corollary for \oiga
(\isn{oiga\_not\_cert\_sound}) rests on the same theorem, with
soundness for the least model in place of conservativity.

\paragraph{Park's rule keeps a certificate-grounded semantics.}
\inga has the coinduction rule but neither ATI nor the dichotomy, and
it has a certificate-grounded semantics.  The reading that serves is
derivability itself: a term is true when \inga derives it and false
when \inga derives its negation.  With truth read so, soundness and
completeness are immediate, and the content of the next theorem is
that this reading meets the conditions above.

\begin{thm}[\isn{inga\_pg}, \isn{inga\_pg\_not\_omega}; classical]
\label{thm:inga-pg}
Call a term true when \inga derives it from no hypotheses, and false
when \inga derives its negation from no hypotheses.  This reading is
a certificate-grounded semantics over \inga's own derivations under
which \rga's closed theorems are true, the kind of semantics that
\cref{thm:nocert} quantifies over: a universal is true only when \inga
derives its open instance, and false only when some numeral instance
is false; a disjunction is true only when a disjunct is true; a
negation is true only when its argument is false; and no term is both
true and false.  So \inga is proof-grounded.  \inga is not sound for
the least model.
\end{thm}

Of the six conditions, four hold outright at \inga
(\isn{inga\_cert\_allT}, \isn{inga\_cert\_excl}): a derivable universal
has a derivable open instance by universal elimination, a derivable
negation is a refutation by definition, \rga's theorems embed, and
nothing is both derivable and refutable, by \inga's consistency.  The
arguments use only \rga's rules and consistency.  The two conditions
that remain, on disjunctions and on false universals, are the
substance.

\begin{cor}[\isn{inga\_disj\_prop}, \isn{inga\_witness\_prop}; classical]
\label{cor:inga-dpwp}
If \inga derives $a \vee b$ from no hypotheses, it derives $a$ or
$b$.  If \inga derives $\neg(\cA\cdot u)$ from no hypotheses, it
derives $\neg(u\cdot\bar n)$ for some $n$.  Hence a sentence whose
decidedness \inga derives is derived or refuted by \inga.
\end{cor}

The corollary's last clause is its first clause at the disjunction
$c \vee \neg c$.  Conversely, the two properties are all that the
conditions ask beyond the four that hold outright, so the corollary's
first two clauses are equivalent to the theorem's claim that the
reading is a certificate-grounded semantics (Appendix~\ref{app:inga}).

The proof does not argue from derivability.  It builds a semantics in
which a derivable term is true for a structural reason, as the \rga
and \iga papers do for their systems, and extends it to $\nu$-atoms in
the manner of $\Mstar$.  The evaluation is \rga's certificate-grounded
one: a universal is true when the system derives its open instance
and every numeral instance is true, and false when the system refutes
some instance and that instance is false.  To it the construction
adds the $\nu$-head of $\Mstar$ and an oracle that the evaluation
consults on both sides: a $\nu$-atom is true when its argument lies in
the oracle's truth set and false when the argument lies in the
falsity set, and a $\nu$-term is never unfolded as a head.  The oracle
map reads each body at the greatest fixed point of its truth operator
and the least fixed point of its falsity operator, as for $\Mstar$,
but it filters both sides: a point enters the truth side only if
\inga derives the atom there, and the falsity side only if \inga
derives the atom's negation.  The \emph{certified oracle} is the
least fixed point of that map, and the \emph{certified evaluation} is
the evaluation at the certified oracle.  Because of the filter, every
well-formed term true in the certified evaluation is derivable and
every well-formed false term is refutable.  Soundness then goes rule
by rule, as for $\Mstar$.  For Park's rule it is a coinduction against
the filtered greatest fixed point.  The region is the set of
well-formed terms with a value at which the rule's certificate is
true, and a point of the region passes the filter because
completeness turns the truths at that point into derivations, to
which one application of the rule adds the atom.  Once \inga is sound
for the certified evaluation, a well-formed term is derivable exactly
when it is true there and refutable exactly when it is false there,
the oracle consists exactly of the derivable and the refutable
$\nu$-atoms, and a derivable disjunction has a true disjunct, hence a
derivable one.  The same construction applies to \rga with the
coinduction rule alone.  Appendix~\ref{app:inga} gives the
construction in detail.

Both sides of a $\nu$-atom are read from the oracle, not one side
from the oracle and the other by unfolding.  The reason is
determinism, the fact that no term gets two values.  An evaluation
that read only truth from the oracle and still unfolded $\nu$-atoms to
make them false would meet, in a proof of determinism by induction on
the evaluation, a $\nu$-atom true by the oracle, with no
subderivation, and false through its unfolding; the two verdicts are
related only through the fixed point, not by the induction.  With both
sides in the oracle, the clash is a clash inside the oracle, which is
consistent.  This is a reading of why the construction has the shape
it has, not a theorem that the alternative fails.

Read together, \cref{thm:inga-pg,thm:nocert} locate the loss of a
certificate-grounded semantics.  Coinduction is not where the loss
comes from: \inga is unsound for the least model by design, as \oinga
and \dinga are, and it is nonetheless proof-grounded.  The loss comes
with ATI, primitive in \oinga and derived from the dichotomy with
ground induction in \dinga.  The count of complexity behind
\cref{thm:nocert} needs the system to derive the decidedness of every
universal of the halting family.  By \cref{cor:inga-dpwp}, \inga
derives the decidedness of a sentence only where it derives the
sentence or the sentence's negation, so the count has nothing to work
on.  In this reading of the two theorems together, Park's rule keeps
proof-groundedness, and the dichotomy, through ATI, destroys it; no
theorem of the development states that ATI is the cause.

\paragraph{Proof-grounded without the least model.}
In the \rga and \iga papers a proof-grounded system is also sound for
the least model, because the certificate-grounded evaluation there
refines the least model: a term true in the evaluation is true in the
least model.  The development states that kind of semantics as a
locale, \isn{gsem\_sound}, whose instances are therefore sound for the
least model.  \inga is an instance for no step relation that contains
\inga's, since \inga derives the all-zeros atom, which has no value
in the least model (\cref{thm:unsound}).  So \inga is the first
system, among those of this paper and its companion papers, that is
proof-grounded and not sound for the least model.  The definition of this section asks for the necessary
conditions and for soundness and completeness, and nothing more; it
does not ask that the least model confirm every verdict.  That
consequence belongs to the locale, not to proof-groundedness.
\inga's verdicts on $\nu$-atoms rest on its derivations, and the
certified evaluation and $\Mstar$ confirm them; the least model does
not.

\paragraph{Consistency from the certified evaluation.}
\Cref{thm:inga-pg} takes nothing from $\Mstar$, not even \inga's
consistency.  The certified oracle is consistent, no point lying on
both of its sides, for every calculus and with no hypothesis: at a
consistent oracle the truth and the falsity readings of one body never
give it opposite values, so the truth side of the oracle map's image,
a greatest fixed point, and its falsity side, a least fixed point, are
disjoint; the map therefore preserves consistency, so do suprema of
chains, and induction along the iteration that reaches the least
fixed point makes the certified oracle consistent.  The filters play
no part in that argument.  Soundness for the certified evaluation then
gives consistency, because $0 = 1$ is false there.

\begin{cor}[\isn{inga\_consistent\_mf}, \isn{rganu\_consistent\_mf}; classical]
\label{cor:inga-con-mf}
\inga does not derive $0 = 1$ from no hypotheses, and neither does
\rga with the coinduction rule.  The proof is soundness for the
certified evaluation, and no fact about $\Mstar$ lies under it (a
build-time dependency check, \S\ref{sec:concl}).
\end{cor}

\noindent \inga's consistency thus has two proofs: the proof through
$\Mstar$ of \cref{thm:inga}, which also gives soundness and
conservativity, and this one, from \inga's certificate-grounded
semantics alone.  Exclusion, the one condition
of \cref{thm:inga-pg} that needs consistency, takes it from
\cref{cor:inga-con-mf}.  So the proof that \inga is proof-grounded,
and \cref{cor:inga-dpwp} with it, is free of $\omega$-models: it rests on an
evaluation built from \inga's own derivations, and on no classical model; only
the theorem's last sentence, unsoundness for the least model, speaks
of one (Appendix~\ref{app:inga}).

\subsection{Exclusion and decided regions}
\label{sec:excl}

Two further principles close the circle, and they need the rule less
than one would guess.

\begin{thm}[\isn{g\_mu\_not\_nu}, \isn{g\_dec\_nu}; grounded]
\label{thm:excl}
In every calculus with ground induction, in a well-formed context
with $x$ a number: from $\mu\Phi\cdot x$, $\neg\nu(\neg\dual\Phi)\cdot
x$.  In every calculus with ground induction and the dichotomy: from
the decidedness premise, $\nu(\neg\dual\Phi)\cdot x$ is decided.
\end{thm}

Exclusion uses ground induction and no other rule beyond \rga, at the
motive ``not in the $\nu$-atom'': in the closure, one conversion at
each recursive atom turns the body at that motive into the co-body at
the $\nu$-atom, and the $\nu$-atom is then folded under negation; it
needs no decidedness.  For accessibility it says $\Acc_R\cdot x \ent
\neg\Path_R\cdot x$, and with the dichotomy, $\Path_R$ is decided at
every decided relation.  The $\mu$-side stays undecided.  No system in
the paper proves $\neg\Acc_R\cdot x$ from no hypotheses at a closed
relation $R$, and none can: accessibility is never false in $\Mstar$,
at any closed $R$ and any point, decided or not, as it is never false
in the least model, and non-derivability follows by soundness
(\cref{thm:app-accO}).

\paragraph{The lattice.}
\iga adds ground induction ($\mu$) to \rga, \oga adds ATI, \oiga
both; \inga and \oinga add $\nu$; \dinga adds the dichotomy to \inga
(\cref{fig:systems}).  The development's abstract calculi are all
sixteen combinations of these four rules, one locale each, ordered by
inclusion; by \cref{thm:ati-dich} every calculus with ground
induction and the dichotomy has ATI as a derived rule.  Each theorem
is stated in the least of them that its proof uses, so the absence of
a primitive rule from a proof is carried by the calculus, and the
theorem holds in every system that has the calculus' rules.  For ATI
in a calculus with ground induction and the dichotomy, the calculus
carries only the absence of ATI's rule, since it derives ATI's
conclusions; \S\ref{sec:rungmap} says what certifies the rest, and,
theorem by theorem, which rules the proof uses within its calculus.
\begin{center}
\begin{tabular}{@{}llll@{}}
\toprule
system & rules beyond \rga & sound for & proof-grounded \\
\midrule
\iga & $\mu$ & least model, $\Mstar$ & yes~\cite{iga-paper} \\
\oga & ATI & least model, $\Mstar$ & no (\cref{thm:nocert}) \\
\oiga & $\mu$, ATI & least model, $\Mstar$ & no (\cref{cor:nocert}) \\
\inga & $\mu$, $\nu$ & $\Mstar$, the certified evaluation & yes (\cref{thm:inga-pg}) \\
\oinga & $\mu$, ATI, $\nu$ & $\Mstar$ & no (\cref{cor:nocert}) \\
\dinga & $\mu$, $\nu$, dichotomy (ATI derived) & $\Mstar$ & no (\cref{cor:nocert}) \\
\bottomrule
\end{tabular}
\end{center}

\section{Streams, Runs and Branches}
\label{sec:streams}

\subsection{Streams as total terms}

An infinite stream is a term $s$ with $\Tot\,s$: every $s\cdot\bar n$
is a number.  The paper takes this representation as primary and the
alternative, a stream as the set of its finite prefixes, as the
derived view; the two are interdefinable.  Finite sequences are coded
as numbers with a decided parent function, the empty sequence coded
as $1$.  The prefix function $\mathrm{pfx}_s$ sends $n$ to the code of
the first $n$ elements of $s$, and the \emph{prefix set}
$\mathrm{pset}_s$ is the predicate ``some $n$ has $\sigma =
\mathrm{pfx}_s\cdot n$''.  A \emph{branch} is a decided set of sequence codes that is prefix-closed, has a member of every length, and has at most one member of each length.

\begin{thm}[\isn{g\_pset\_Dec}, \isn{g\_pset\_Branch}; grounded]
\label{thm:branch}
For a total $s$, its prefix set is decided, by ATI, and is a branch.
\end{thm}

\begin{thm}[\isn{g\_RT\_s}, \isn{g\_RT\_p}; grounded]
\label{thm:roundtrip}
For a total $s$, the stream recovered from the prefix set agrees with
$s$ at every $n$, and the prefix set of the stream recovered from a
branch has exactly the branch's members.
\end{thm}

The round trip uses ATI to decide the prefix set; the stream is then
recovered by a search over that decided set, a search which is a
derived rule of \rga, and no fixed-point rule is used.  ATI is spent
only because the prefix set is defined with an unbounded existential,
``some $n$''.  The witness is in fact bounded, since
$\mathrm{pfx}_s\cdot n$ has length $n$, and the elementwise view, that
$\sigma$ is a prefix of $s$ exactly when every position $i$ below its
length has $\sigma_i = s\cdot i$, is arithmetic without ATI; with the
prefix set defined by that bounded form, the round trip would need no
ATI (a prediction; the development keeps the unbounded definition).

Streams carry no coinductive content here, and a reader from
programming languages should read them that way.  A stream is a total
function on the numbers, not an element of a coinductive type.
Productivity is totality, the judgment $\Tot\,s$ derived in the
system, proved by number induction; no $\omega$-model enters.  A
guarded corecursive definition such as $s = 1 :: s$ is an ordinary
recursive function on indices, written with the fixed-point former.
Bisimilarity of two streams is pointwise equality, a $\Pi^0_1$
sentence; \rga has no extensionality for terms, and none is wanted.
Every $\nu$ in the paper is a predicate on numerically coded states,
such as ``some chain of steps from this state is infinite'', and never
a type of infinite objects.

\paragraph{The prefix tree.}
A prefix set is a tree, and a tree has the child relation: $y$ is a child of $x$ when $y$ is a member and the parent of $y$ is $x$.  The
$\nu$-atom $\Path$ of that relation at a node says there is an
infinite branch through the node.

\begin{thm}[\isn{g\_pset\_path}, \isn{g\_path\_pset\_hol\_iff}; grounded]
\label{thm:pset-path}
For a total $s$, every member of its prefix set has an infinite
branch through it, by one coinduction step with the prefix set as
certificate.  Conversely, at any code that has a length, the branch
predicate and membership in the prefix set are interderivable, with
$\nu$ and no ATI; with ATI added to $\nu$ the equivalence is one
sentence and the branch predicate is decided, without the dichotomy.
\end{thm}

So at a single stream's prefix tree, the coinductive predicate
``there is an infinite branch'' is arithmetic.  The dichotomy earns
its keep only on trees with many branches, which is where the next
subsection goes.

\subsection{Runs}

A \emph{run} of a relation $R$ is a total term $s$ with $R(s\cdot(n+1),
s\cdot n)$ at every $n$: an infinite chain of steps.  ``Every chain of
steps from $x$ is finite'' is $\Acc_R\cdot x$, and ``some chain of
steps from $x$ is infinite'' is $\Path_R\cdot x$.  Three theorems
connect runs with the two atoms, and their proofs spend three
different rules.

\begin{thm}[\isn{g\_run\_acc}, \isn{g\_run\_path}; grounded]
\label{thm:runs}
If $s$ is a run of $R$ and $x$ is accessible, then $s$ never visits
$x$ (ground induction).  If $s$ is a run of $R$, every state of $s$
lies on an infinite chain (one coinduction step, with the image of
$s$ as certificate).
\end{thm}

\begin{thm}[\isn{g\_run\_of\_path}; grounded]
\label{thm:runofpath}
For decided $R$, from $\Path_R\cdot x$ there is a run from $x$, a
term total in $\Mstar$, every state of which lies on an infinite
chain.  The proof uses the dichotomy and ground induction, and applies
neither $\nu$ nor ATI.
\end{thm}

The run is a generic construction, a \emph{replay}: given a
decided relation $C$, a decided predicate $Q$ with a progress
property (every $Q$-state has a $C$-successor in $Q$) and a
$Q$-state $a$, the term $\mathrm{path}_{C,Q}(a)$, defined by
primitive recursion with a least-witness search at each step, is
total, starts at $a$, is a $C$-run and stays in $Q$.  That is
arithmetic.  What the dichotomy supplies is the decidedness of $Q =
\Path_R$, which the search needs, and the progress property is one
unfolding.  So non-termination is stated positively, as an object:
the run.  The run is total in $\Mstar$; its search reads a
$\nu$-atom, and the paper does not claim that it reaches numerals by
reduction.

\subsection{Well-founded, or an infinite branch}

K\H{o}nig's lemma says that a finitely branching infinite tree has an
infinite path, and the \oga development proves it for decided,
prefix-closed trees with a bound on the branching, by a path term
that at each node picks the first child with an infinite subtree; the
pigeonhole argument that some child has one is where the bound is
spent.  Without a bound the lemma is false as stated.  What holds
without a bound is a dichotomy.

\begin{thm}[\isn{g\_tree\_acc\_or\_branch}, \isn{g\_tree\_branch}; grounded]
\label{thm:koenig}
Let $T$ be a decided, prefix-closed set of sequence codes.  Then
either the root is accessible for the child relation of $T$, that
is, $T$ is well-founded, or there is a branch given by a term, total
in $\Mstar$, whose values are nodes of $T$, each the parent of the next, the
$n$-th of length $n$.  The proof uses the dichotomy and ground
induction, and applies neither $\nu$ nor ATI.  The branch itself needs neither prefix closure nor sequence codes: for any decided $T$ whose root is a member with $\Path$ at the root, the replay is total, starts at the root and stays in $T$ on nodes with $\Path$.
\end{thm}

The branch is the replay at the child relation and $Q = \Path$, which
the dichotomy decides.  K\H{o}nig's lemma for bounded branching is the
same replay at the parent relation of the full tree, with progress
from the pigeonhole, and the \oga theorem is literally a corollary of
the replay.  The two differ in where the path is total.  The bounded
path term is total in the least model.  The unbounded one is total in
$\Mstar$, by soundness; its search reads a $\nu$-atom, and whether it
has a value in the least model is not established.  Neither is
computable in general, and the classical landmark says why: the
Kleene tree is a decided, binary, infinite tree no computable
function is a path through, and the \oga development has it as a
theorem pair, the tree decided and infinite, and no path recursive.
The bounded path term through it is total in the least model, as an
infinitary object, and is not a computable function; the theorem
above would give the same path from the dichotomy, and the formal
instance, which needs a recoding of that tree into sequence codes,
is left for later work.

\section{Examples}
\label{sec:showpieces}

Each example gives the program or game, its bodies, the theorems,
and the rules the proofs use.  The examples sit where
\S\ref{sec:liveness} will say the new rules matter: unbounded
nondeterminism, an adversary, or a tree with many branches.

\subsection{Random assignment}

The program $x := ?;\ \mathbf{while}\ x > 0\ \mathbf{do}\ x := x - 1$
is the standard example of countable
nondeterminism~\cite{apt86countable}: it always terminates, but no
bound on its running time exists, and its termination has rank
$\omega$.  Code the state before the assignment as $0$ and the state
with $x = m$ as $m+1$.  The step relation $R$ then sends $0$ to $m+1$
for every $m$ ($\omega$-branching) and $m+2$ to $m+1$, with no
successor at $1$.  The divergent twin $R'$ adds a step from $1$ to
$0$.

\begin{thm}[\isn{g\_acc\_R1\_all}, \isn{g\_R1\_no\_variant}; grounded]
\label{thm:rand}
Every state is accessible for $R$: every execution terminates.  The
proof uses no fixed-point rule.  No total term is a strictly
decreasing natural-number ranking along $R$: the variant hypothesis
entails $0=1$ arithmetically, and the proof spends ATI to decide the
hypothesis, which turns that entailment into its negation.
\end{thm}

The first claim is the finding of \S\ref{sec:rungmap} at its
simplest.  Accessibility at $0$ is introduced by one backward
unfolding whose successors are all loop states, and each loop state
is accessible by number induction; folding and conversion are not
fixed-point rules, and no ground induction is needed to assert
termination, only to use it.  The second claim says that no
natural-number measure certifies termination: a variant would have to
assign to $0$ a number above every $m+1$.  That the rank is exactly
$\omega$ is a reading of the two claims, not a formalized ordinal
statement.

\begin{thm}[\isn{g\_R1p\_div}, \isn{g\_R1p\_run}; grounded]
\label{thm:rand-div}
$\Path_{R'}\cdot 0$ holds, by one coinduction step with the
certificate ``every number'' (\S\ref{sec:nu}).  A run from $0$ along
$R'$, a term total in $\Mstar$, exists, every state of which lies on
an infinite chain; the proof uses $\nu$ for the path and the
dichotomy with ground induction for the replay.  For both relations,
every state is accessible or lies on an infinite chain, by the
dichotomy, and the chain predicate is decided, by the dichotomy with
ground induction.
\end{thm}

The replay is the generic construction of \cref{thm:runofpath}, and
at $R'$ it is not the only witness: the run $0, 1, 0, 1, \ldots$ is
primitive recursive and total by \rga's rules alone, and with
\cref{thm:runs} it gives $\Path_{R'}\cdot 0$ as well (an observation;
the development states no theorem for this particular run).  A nested
loop with a random reset at each level terminates at rank $\omega^2$,
again without a fixed-point rule.  And neither \oinga nor \dinga
proves $\Acc_{R'}\cdot 0$: from the divergence, the exclusion theorem
and the system's consistency, which is itself proved through
$\Mstar$.

\subsection{Fairness through random assignment}

Fair termination is accessibility over an $\omega$-branching relation
of the same kind.  A variant of Apt and Olderog's transformation of a
fair program into one with random assignment~\cite{apt83fairness},
with natural-number counters and a strict guard, turns the
two-component loop that either increments a counter or falsifies its
own guard into a step relation whose infinite chains project onto
strongly fair infinite runs of the loop, and every strongly fair
infinite run of the loop lifts to such a chain.  A component may run
only while every other enabled component's counter is positive;
running it resets its own counter to any number, decrements every
other enabled component's counter by one and keeps the counter of
every disabled component.

\begin{thm}[\isn{g\_acc\_ao\_all}, \isn{aomod\_fair\_term}; grounded and classical]
\label{thm:fair}
Grounded: every state of the transformed loop is accessible, by one
introduction and one number induction on the postponed component's
counter, with no fixed-point rule; from accessibility, no state lies
on an infinite chain (ground induction); every state is accessible or
lies on an infinite chain (the dichotomy); and the chain predicate is
decided (the dichotomy with ground induction).  Classically, through
\rga's soundness for the least model and the correspondence of
infinite runs, the loop has no infinite strongly fair run from its
initial state.
\end{thm}

The postponed component's counter is the measure, and the reset
belongs to the component that just ran, so the reset at every
iteration does not nest: the rank at the initial state is $\omega$
(a reading).  The correspondence is a classical theorem about runs,
a projection of chains onto runs and a lifting of runs to chains, and
it holds for infinite runs only, since the strict guard can block
where the loop does not.

A twin of the loop, the \emph{fairness twin}, separates strong from
weak fairness, and the separation is exactly accessibility against the
path atom.  Let the first component also toggle a bit $c$, and let the
second component, which exits the loop by falsifying the guard, be
enabled only when $c$ holds.  The second component is then enabled at
every other step, so a strongly fair run schedules it, while no run
keeps it continuously enabled, so weak fairness asks nothing of it.
The two transformations differ in one conjunct: the counter of a
disabled component is kept, for strong fairness, or reassigned at
random, for weak.

\begin{thm}[\isn{g\_acc\_aoTs\_all}, \isn{g\_aoTw\_div}; grounded and classical]
\label{thm:fair-twin}
Grounded: under the strong transformation every state of the fairness
twin is accessible, with no fixed-point rule; under the weak
transformation the path atom holds at the initial state, by one
coinduction step at a three-part certificate.  Classically: neither
\oinga nor \dinga proves the initial state accessible under the weak
transformation, from the exclusion theorem and the system's
consistency; through \rga's soundness for the least model, the
fairness twin has no strongly fair infinite run from its initial
state; and the run that always schedules the first component is
weakly fair, a fact about the HOL program alone.
\end{thm}

Read with $c$ as ``channel full'', the first component as ``send when
empty, lose when full'' and the second component as ``deliver'', the
fairness twin is already the one-message lossy channel: strong
fairness is what defeats a channel that may drop messages, and weak
fairness is not enough.

The development also builds the transformation as a term former, for
any list of guarded components with a flag saying whether the
scheduler owes each one fairness, with a monotone rather than exact
counter update; whenever every guard and step is decided, the
transformed relation is decided, so the dichotomy holds at every such
program (Appendix~\ref{app:fair}).  Generic accessibility is not on
offer, by Harel's theorem~\cite{harel86effective}: fair termination is
$\Pi^1_1$-complete.

\paragraph{Eventualities and leads-to.}
For a step relation $T$, read forward, and a target $G$, the
\emph{eventually atom} $\mathrm{Ev}_G$ is the $\mu$-atom of the body
$G\cdot x \vee \forall y.\,\neg T(x,y) \vee a\cdot y$: every run from
$x$ reaches $G$.  At a state with no successor the universal holds
vacuously, so the atom is true there; read as ``every run reaches
$G$'', it presumes a relation in which every state has a successor,
and on one with deadlocks it says that every run reaches $G$ or
stops.  ``$P$ leads to $Q$'' is the sentence $\forall x.\,\neg
P\cdot x \vee \mathrm{Ev}_Q\cdot x$.  The classical leads-to rules
are derived rules.

\begin{thm}[\isn{g\_lattice}, \isn{g\_well}; grounded]
\label{thm:lattice}
In every calculus with ground induction, for closed, well-formed $T$,
$\prec$, $F$ and $G$, Lamport's Lattice
rule~\cite{lamport94temporal}, which descends from the proof lattices
of Owicki and Lamport~\cite{owicki82proving}: if the sentence ``for
every $c$, either $c$ is not accessible along $\prec$ or $F\cdot c$
leads to $G$ or to some $F\cdot d$ with $d \prec c$'' is derivable,
then at every accessible $c$, $F\cdot c$ leads to $G$.  The
Manna--Pnueli WELL rule~\cite{manna92temporal} is the instance at the
family ``$P$ holds and the rank is $c$'' (\cref{thm:app-lattice}).
\end{thm}

Both rules are one ground induction along accessibility, with the
motive carrying accessibility and the leads-to claim together, after
one ground induction for the transitivity of ``eventually''; the check
forbids $\nu$, the dichotomy and ATI under them.  The first disjunct of
the hypothesis is never derivable from no hypotheses at a closed
order, since $\neg\Acc_\prec\cdot c$ is not (\cref{thm:app-accO}); it is there so that the hypothesis can be one sentence over all
$c$ while being needed only at accessible ones.  For the same reason
the conclusion is a rule and not a sentence.  A countdown that
decrements to zero is the worked instance, its hypothesis and the
accessibility of every counter value derived without any fixed-point
rule (Appendix~\ref{app:programs}).

The main instance of the term former is the alternating-bit protocol over
two lossy one-slot channels, six components: send, lose, deliver on
the forward channel and acknowledge, lose, deliver on the return
channel.  The two loss components owe no fairness, so the adversary
drops messages as often as it likes.

\begin{thm}[\isn{aoAb\_every}, \isn{aomod\_ab\_deliver}; grounded and classical]
\label{thm:abp}
Grounded: every message is eventually delivered and acknowledged.  For
every $K$, the eventually atom whose target is ``message index at
least $K$'' holds at every code of the transformed protocol whose
program part satisfies the protocol's invariant, with any counters;
as one sentence, with $K$ a parameter of the atom's target, the atom
holds at the codes with the invariant and index $0$.  Both forms are
derived with no fixed-point rule.  Classically, through \rga's
soundness for the least model, every infinite run of the protocol
from its initial state that is strongly fair to send, acknowledge and
the two deliveries reaches every index.
\end{thm}

Read through the atom, the grounded claim says that every maximal
execution of the transformed protocol from such a code reaches an
index at least $K$ or ends at a code with no successor; the classical claim, about
infinite runs, has no such proviso.  The proof uses seven number
inductions (a count read off the proof script; the check records only
that no fixed-point rule is present).  The Lattice rule is not on the
path: it consumes a leads-to premise, and a concrete protocol proves
its round lemma rather than assuming it.  The rank is $\omega^2$ per
phase (a reading), since while one channel is empty the other
channel's components run, so each phase needs an outer induction on
the delivery counter and an inner one on the refill counter, and a
stale delivery resets the counter at random.  Read informally, strong
fairness, not weak, is what the lossy channel requires: a loss empties
the channel and keeps each delivery only intermittently enabled.  No
theorem of the development shows the protocol failing under weak
fairness; the fairness twin is the formal separation.

The fair-termination theorem also has a grounded form, in which
fairness is a sentence of the calculus rather than a side condition
read in the metatheory.  A run is a stream parameter whose values
carry a program state and the component scheduled at each step.  The
run sentence says that at every step some component is scheduled,
enabled at the current state, and steps to the next; the fairness
sentence says, for every component the scheduler owes fairness, that
it is enabled only finitely often or scheduled infinitely often, one
disjunction per component.

\begin{thm}[\isn{g\_ao\_fair\_rule}, \isn{g\_ao\_run\_ev}; grounded]
\label{thm:fair-grounded}
For every program whose guards and steps are decided: if the
transformed relation is accessible at every code over the initial
state, then from the run sentence and the fairness sentence for a
stream parameter the calculus derives that the stream does not start
at the initial state.  For every program, take an eventually atom
along the transformed relation, read forward; every run of the
transformed relation that starts at a code where that atom holds
meets the atom's target.  In the first claim ground induction is spent only in
the no-run theorem; the second claim is one ground induction on the
eventually atom; neither claim uses ATI, the dichotomy rule or
$\nu$.
\end{thm}

The derivation builds the counters of the transformation inside the
branches of each fairness disjunction, by a search or as a constant,
so it never decides whether a component is enabled infinitely often
or scheduled again (Appendix~\ref{app:fair}).  The conclusion is a
rule about a stream parameter, not a negated sentence.

\subsection{Termination and non-divergence}

Nakata and Uustalu's coinductive trace semantics of a while
language~\cite{nakata15hoare} separates termination from
non-divergence, and in their constructive setting the gap is
Markov's principle.  The grounded form of the gap is different.  Let
$q$ be a decided test and consider the search from $x$: the body
$q(x) \vee a\cdot(x+1)$.  Its $\mu$-atom $\Term_q$ says the search
terminates; the $\nu$-atom of the co-body, $\Div_q$, with body $\neg
q(x) \wedge a\cdot(x+1)$, says it diverges.

\begin{thm}[\isn{tr\_Q1\_iff}, \isn{tr\_Q2\_iff}; grounded]
\label{thm:trace}
In a well-formed context $\Gamma$ with $x$ a number: $\Term_q\cdot x$
and $\exists m.\,q(x+m)$ are interderivable in $\Gamma$, for decided
$q$ (ground induction for one direction, nothing for the other);
$\Div_q\cdot x$ and $\forall m.\,\neg q(x+m)$ are interderivable in
$\Gamma$ ($\nu$ for one direction); and $\Term_q\cdot x \ent
\neg\Div_q\cdot x$ (ground induction).
\end{thm}

Markov's principle proper, from $\neg\neg\exists$ to $\exists$, is
free: the existential is $\neg\cA\neg$ and double-negation
elimination is a rule, since a grounded negation flips a value rather
than asserting refutability.  What separates termination from
non-divergence is excluded middle for the $\Pi^0_1$ sentence $\forall
m.\,\neg q(x+m)$, and for the search body the gap closes in two
independent ways, with a third route that uses both.

\begin{thm}[\isn{tr\_term\_of\_ndiv}, \isn{tr\_dich}; grounded]
\label{thm:trace-gap}
From $\neg\Div_q\cdot x$, $\Term_q\cdot x$: the proof uses ATI, which
here acts as the limited principle of omniscience, and $\nu$; it uses
neither ground induction nor the dichotomy.  Independently, $\Term_q
\cdot x \vee \Div_q\cdot x$ by the dichotomy, with ATI, $\nu$ and
ground induction all absent from the proof.  A third proof of the
disjunction uses ATI and $\nu$ and no dichotomy.
\end{thm}

Whether \inga, with neither ATI nor the dichotomy, proves the first
claim is open.  The disjunction is the grounded counterpart of Leroy
and Grall's coinductive big-step semantics~\cite{leroy09coinductive},
which proves, by classical excluded middle in Coq, that a program
either evaluates or diverges.  For a deterministic language with a
decided step function, take $q(n)$ to say that the state after $n$
steps is final: $\Term_q\cdot 0 \vee \Div_q\cdot 0$ is then their
``evaluates or diverges'', obtained from the dichotomy under the
derivable premise that the steps are decided, and for a
nondeterministic decided step relation the same role falls to
$\Acc \vee \Path$ (a reading; the development states no instance for
a while language).

Nakata and Uustalu's later section concerns the predicate ``the trace
reaches $n$ and then goes up'', and its grounded form is a $\nu$-body
whose non-recursive leaf is a $\mu$-atom, reachability.  Take the
counting program, whose state is a pair $\langle x, y\rangle$:
\[
  \langle x, 0\rangle \to \langle x+1, x+1\rangle, \qquad
  \langle x, y+1\rangle \to \langle x, y\rangle,
\]
so the first component never decreases.  Write $\mathrm{Reach}_t$ for
the $\mu$-atom of $x = t \vee \exists s.\,(x \to s) \wedge a\cdot s$,
``the run from the point reaches $t$''.  The up body reads, at a point
$\langle n, \sigma\rangle$:
\[
  \exists \sigma'.\,\mathrm{Reach}_{\sigma'}\cdot\sigma \wedge
  \mathrm{fst}(\sigma') = n \wedge \exists\sigma''.\,(\sigma' \to
  \sigma'') \wedge a\cdot\langle n+1, \sigma''\rangle .
\]
The counting program's trace satisfies the $\nu$-atom of the up body
at $0$ by one coinduction step.  This is also the body at which the
dichotomy's premise fails in $\Mstar$.  Its leaf
$\mathrm{Reach}_{\sigma'}\cdot\sigma$ at the target $\sigma' =
\langle 0,0\rangle$ and the state $\sigma = \langle 1,0\rangle$ has
no value: it is not true, because every state from which the target
is reachable has first component zero, and the run from $\langle
1,0\rangle$ never returns there; and it is not false, because every
state has a successor and a falsity would need a well-founded
refutation.  So the premise, which asks every leaf to be decided, has
no value at that instance, and neither it nor its negation is
derivable in \dinga from no hypotheses (\cref{thm:app-premise}).
That is the control named in \S\ref{sec:dich}.

\subsection{Open determinacy}

A game between players I and II on positions is given by three
parameters: a decided set $U$ of positions already won by I, a
decided relation $\mathit{Move}$ for I's moves, and a decided
relation $\mathit{Reply}$ for II's replies.  Player I's winning body
is
\[
  U(x) \;\vee\; \exists q.\,\bigl(\mathit{Move}(q,x) \wedge
  \forall r.\,(\neg \mathit{Reply}(r,q) \vee a\cdot r)\bigr),
\]
with the existential as a negated universal; $\WinI$ is its
$\mu$-atom, I's winning region, and $\WinII$ the $\nu$-atom of its
co-body, II's.  The game is open for I: a win is witnessed by a
finite play.

\begin{thm}[\isn{gm\_det}, \isn{gm\_decII}; grounded]
\label{thm:det}
For decided $U$, $\mathit{Move}$, $\mathit{Reply}$ and a position
$x$: $\WinI\cdot x \vee \WinII\cdot x$ (the dichotomy alone);
$\WinI\cdot x \ent \neg\WinII\cdot x$ (ground induction alone, with no
decidedness); and $\WinII\cdot x$ is decided (both).
\end{thm}

The decidedness premise is derived from the three parameters by
\rga's rules, with no ATI.  Classically this is Gale--Stewart
determinacy for open games~\cite{gale53infinite}, and over a weak base
theory open determinacy for games on Baire space, where a move
appends a number to a finite sequence, is equivalent to arithmetical
transfinite recursion~\cite{steel77determinateness}, the subsystem of
second-order arithmetic that iterates arithmetical comprehension
along any well-order.  That theorem is about the existence of
strategies; \cref{thm:det} is about regions, and the paper claims no
calibration.  The Baire instance is one substitution.  Two concrete
games fix the two sides: I wins in one round by two folds, and, when
I's target is the set of positions at which II's reply repeats I's
last move, II wins forever by answering every move $m$ with $m+1$,
with the complement of $U$ as certificate.  The grounded setting
also yields a strategy, a term total in $\Mstar$.

\begin{thm}[\isn{gm\_sigII}; grounded]
\label{thm:sigII}
In the Baire instance with $U$ decided, there is a term $\sigma$,
total at every position, such that from any position in $\WinII$ and
any move of I, the reply $\sigma$ chooses keeps the play in $\WinII$.
The proof uses the dichotomy, ground induction and least-witness
search, and no coinduction step.
\end{thm}

The strategy searches, at each position, for the least reply that
keeps II in her region; the search is decided because $\WinII$ is,
and the region is a hypothesis, only unfolded.  The term is total in
$\Mstar$, where $\WinII$ has a value; the paper does not claim that
it computes.  That is as it must be, since winning strategies in
computable open games can be non-computable: the open player, when he
wins, has a hyperarithmetic strategy (computable from an iteration of
the halting problem along a computable well-order), of arbitrarily
high level~\cite{blass72complexity}; the closed player's may be
non-hyperarithmetic (Kleene's tree).  Player II is the closed player
here.  Player I's strategy, whose extraction would need ordinal ranks
rather than one search, is not built.

\subsection{The bucket game}

The game of the introduction, with buckets of capacity $\kappa$, a
natural number, is an instance of the same body.  The Stepmother is
player I, the open player; she wins when a pour overflows a bucket,
since Cinderella may reply only where no bucket exceeds the capacity,
so the target $U$ is empty and an overflow leaves Cinderella without
a legal reply.  Cinderella is player II.  Amounts are rational, coded
over a common denominator that the Stepmother may refine by any
factor $k \ge 1$ at each pour, so her move relation is
$\omega$-branching; restricting her to rational pours is the one
simplification of the game of the introduction.  The refinement
factor is stored in the successor state, which removes every bounded
existential from the move relations and makes their decidedness a
matter of \rga's rules alone.

\begin{thm}[\isn{cs\_cind\_init}, \isn{cs\_stepmother\_refuted}; grounded]
\label{thm:bucket}
For $\kappa \ge 3$, Cinderella's region contains the initial state:
she has a strategy that keeps the water in the buckets forever.  The
proof uses one coinduction step and otherwise \rga's rules; ground
induction, ATI and the dichotomy are absent from it.  For every
$\kappa$, every state is won by one of the players, by the dichotomy,
and Cinderella's region is decided, by the dichotomy and ground
induction.  For $\kappa \ge 3$, the hypothesis that the Stepmother
wins from the initial state refutes itself: from
$\WinI\cdot\mathit{init}$, $0 = 1$.
\end{thm}

The certificate is the round-robin invariant: with a pointer $r$ to
the pair Cinderella empties next, each bucket holds at most the entry
of a fixed table of bounds indexed by $r$ and the bucket, in units of
the current denominator.  The proof factors through one strategy
term, uniform in the state: the coinduction's closure premise asks,
at a certificate state and any pour, for a reply that re-establishes
the certificate, and the paper's proof supplies the round-robin reply
as a single term that evaluates to the right pair in each of the five
cases, so that the strategy's step theorem and Cinderella's region
come out of one argument.  The strategy term is total by \rga's rules
alone, hence in the least model.

The last clause of \cref{thm:bucket} is the grounded way to say that
a player loses, and it is the only one available: the negation
$\neg\WinI\cdot\mathit{init}$ is not derivable from no hypotheses in
\dinga, nor in any system sound for $\Mstar$, at any capacity
(\cref{thm:app-notrefuted}).  A $\mu$-atom is false only through a
well-founded refutation, and the Stepmother always has a legal move,
the pour of the whole unit into one bucket, so no refutation of her
atom ever bottoms out.

The certificate is for capacity three and above, and that bound is
attained against the strategy it certifies.  Let the Stepmother pour
the whole unit into bucket $0$ every round while Cinderella follows
the round-robin strategy.  Every pour of that play is legal, the
first five replies are legal at capacity two, and after the sixth
pour bucket $0$ holds three units, the denominator being one.  So at
capacity two or less no reply is legal there: the hypothesis that any
well-formed term is a legal reply of Cinderella's at that position
entails $0 = 1$, and the Stepmother has won that play (\cref{thm:app-tight}, in
\rga).  The threshold capacity of the game
is studied in the literature~\cite{bodlaender12cinderella,hurkens11cinderella};
the paper makes no claim about it.

\subsection{Decidedness of $\mathcal O$ and the Cantor--Bendixson shape}

Two results at a body whose $\mu$-atom is $\Pi^1_1$-complete, stated
briefly.  Let $\mathcal O$ be the $\mu$-atom of the \iga development's
body for Kleene's constructive ordinal notations, the codes of
well-founded computable ordinals, and $\mathrm{co}\mathcal O$ the
$\nu$-atom of its co-body.  At the greatest fixed point the body adds
non-notations: a limit code lies there exactly when its sequence never
produces a refutable code, and the constant-zero sequence gives a code
in the greatest fixed point that is not a notation.

\begin{thm}[\isn{g\_dec\_O}, \isn{g\_O\_not\_coO}; grounded]
\label{thm:decO}
For every number $a$: $\mathcal O\cdot a \vee \mathrm{co}\mathcal O
\cdot a$, and $\mathrm{co}\mathcal O\cdot a$ is decided.  From
$\mathcal O\cdot a$, $\neg\mathrm{co}\mathcal O\cdot a$.
\end{thm}

The leaves of the decidedness premise are equations between pairing
codes, decided without ATI, and graph atoms of the partial recursive
functions, each an existential over a step bound of an interpreter of
index trees, which ATI decides; the step inside each graph atom is
total because a decoder of index trees is total, which \rga proves
alone.  The theorem does not make the system prove $\mathcal O\cdot
a$ at any new $a$; it decides the complement, in $\Mstar$.
``Decided'' here is the judgment $\mathbf B$ in $\Mstar$, not an
algorithm.  The decision has a term.  Let $\chi$ be $\lambda
a.\,\neg\,\mathrm{co}\mathcal O\cdot a$, the negated co-set atom with
its point abstracted.

\begin{thm}[\isn{O\_Mstar\_iff}, \isn{dinga\_chiO\_char}; grounded and in $\Mstar$]
\label{thm:chiO}
For every number $a$, $\mathcal O\cdot\bar a$ is true in $\Mstar$
exactly when $a$ is in the least model's $\mathcal O$, Kleene's set of
notations.  \dinga derives $\Tot\,\chi$ in every well-formed context
(\isn{dinga\_Tot\_chiO}), and in $\Mstar$, $\chi\cdot\bar a$ has the
value $1$ when $a \in \mathcal O$ and $0$ otherwise
(\isn{dinga\_chiO\_Mstar}): $\chi$ is the characteristic function of
$\mathcal O$ in $\Mstar$ (\cref{thm:app-chiO}).
\end{thm}

The reading as a characteristic function is $\Mstar$'s, not the least
model's: off $\mathcal O$ the co-set atom, and with it $\chi\cdot\bar
a$, has no value in the least model.  No search for a least witness is
needed.  An \rga formula is a term whose values are $1$ and $0$, so
the decidedness of the co-set atom (\cref{thm:decO}) makes $\chi\cdot
a$ decided at every number $a$, hence a number, which is
$\Tot\,\chi$.  The first clause
is the substantial one.  One direction is the persistence of the least
model's values in $\Mstar$.  For the other, ground induction holds in
$\Mstar$ at a closed set of numbers, and at the least model's
$\mathcal O$ the body's reading relative to $\Mstar$ agrees with its
reading relative to the least model: every arithmetic leaf of the body
at numerals is decided in \oiga, so it has a value in the least model,
which persists to $\Mstar$, where values are unique.  So $\Mstar$
decides membership in $\mathcal O$ by a term, and since $\mathcal O$
is $\Pi^1_1$-complete, a classical fact of the literature not
mechanized here, $\Mstar$'s truth set is not hyperarithmetic.  Beyond
that lower bound the paper says nothing about the strength of $\Mstar$
in reverse-mathematical terms.

The same pattern over three abstract decided parameters, a set of
nodes, an extension relation and an incompatibility relation, gives
the shape of the Cantor--Bendixson theorem, which splits a tree into a
scattered part and a perfect kernel.  A node is \emph{scattered} when
it is outside the set or every incompatible pair of its extensions
has a scattered member (the $\mu$-atom), and in the \emph{kernel}
when it is in the set and splits into two kernel points, forever (the
$\nu$-atom of the co-body).  Every node is scattered or in the kernel
(the dichotomy), the kernel is decided (with ground induction), a
kernel point splits by one unfolding, and induction along the
scattered rank is ground induction.  Countability of the scattered
part and a concrete tree coding are out of scope.

\section{Where Coinduction Adds Truth}
\label{sec:liveness}

A reader from verification expects a slogan: safety is $\nu$,
liveness is $\mu$.  With the examples in hand the division can be
drawn more finely.  The classifications in this section are
classical, and the table marks which of its rows the development
proves as theorem instances.  The step relations of this section are
read forward, $T(x,y)$ for a step from $x$ to $y$, the orientation of
the eventually atom of \S\ref{sec:showpieces}, while accessibility and
the path atom are written child first, $R(y,x)$, as in
\S\ref{sec:bg}.

\paragraph{On a given trace, everything is arithmetic.}
Let $\tau$ be a total term, a stream of events.  A safety property
of $\tau$ is a universal over its finite prefixes, $\Pi^0_1$ when the
prefix predicate is decided.  A liveness property in the
response form, ``whenever an obligation arises at position $n$ it is
met at some $m \ge n$'', is $\Pi^0_2$ over $\tau$; ``infinitely
often'' is $\Pi^0_2$ and ``eventually always'' is $\Sigma^0_2$.  No
fixed point appears, and the claims are decided leaf by leaf.  A
proof of such a claim is a ranking argument, a bound or a
well-founded measure on how long the obligation stays open, and in a
grounded arithmetic that is number induction when the rank is finite
or $\omega$ and ground induction along an order beyond that.  For a
deterministic system with its trace given, safety and liveness both
live in \rga, \oga or \iga, and the $\nu$ layer adds no truth.  This
is the setting of trace-based verification in B\"urki's
thesis~\cite{burki26grounded}, where liveness is such a $\Pi^0_2$
statement proved per program and no coinduction is needed because the
trace is given; it is also the setting of ranking-function synthesis
in the classical literature~\cite{manna92temporal,lamport94temporal},
whose leads-to rules are derived rules here (\cref{thm:lattice}).  The
eventually atom used for them is true vacuously at a state without a
successor, so ``leads to'' means what it classically means under the
assumption that no state deadlocks.

\paragraph{$\nu$ enters with unbounded choice.}
When the trace is not given, because the environment chooses or the
schedule is open, a quantifier over infinite behaviors is
second-order.  Consider forced eventuality: ``from state $q$ the
controller can force the good event, whatever the environment does''.
As a recursive definition it is the attractor, a $\mu$-atom whose
body has an existential over the controller's move and, inside it, a
universal over the environment's reply.  Its negation, ``the
environment can delay the event forever'', is the $\nu$-atom of the
co-body.  This is accessibility against the infinite chain in game
clothing, and the dichotomy of \S\ref{sec:dich} is determinacy of the
reachability game read as a fixed-point identity, exactly as in
\S\ref{sec:showpieces}.  So the controller's forced eventuality is
inductive even under nondeterminism.  Coinduction is for its
negation, and for safety under unbounded choice: ``the controller can
keep the system safe forever'' is a $\nu$-atom of a body with the
same shape, and its certificate is an invariant region $W$ containing
$q$, consisting of safe states, each of which has a move keeping
every environment reply in $W$.  Under finitely and computably
bounded control the region is already arithmetic: $q$ wins exactly
when, for every $n$, the environment cannot force a violation within
$n$ rounds, and each round is a bounded search; under unbounded control the attractor
does not close at $\omega$, and the region is $\Sigma^1_1$ in
general~\cite{harel86effective}.  The coinduction rule gives the
classical certificate, and under a decided step relation the
dichotomy makes the two regions jointly exhaustive, $\WinI \vee
\WinII$ at every position, with $\WinII$ decided and the
reachability side $\WinI$ left undecided, since it is refuted only
through a well-founded refutation: this is the grounded form of Gale--Stewart determinacy for
closed games~\cite{gale53infinite}.  What is new is not the invariant
method, which classical verification applies under any branching, but
that a grounded system, which could state neither the region nor the
negation of a liveness claim, now states and decides both.

\paragraph{The table.}
On a decided transition relation in which every state has a
successor, and a decided state predicate $p$, the properties of
computation-tree logic sort as follows.  The last column gives each
row's status, explained after the table.
\begin{center}
\small
\begin{tabular}{@{}p{0.27\textwidth}p{0.17\textwidth}p{0.27\textwidth}p{0.19\textwidth}@{}}
\toprule
property & classical form & fixed point & status \\
\midrule
EF\,$p$: some finite run reaches $p$ & $\Sigma^0_1$ & arithmetic; the $\mu$/$\nu$ readings add nothing & classical \\
AG\,$p$: every reachable state satisfies $p$ & $\Pi^0_1$ & arithmetic; the $\mu$/$\nu$ readings add nothing & classical \\
AF\,$p$: every run eventually hits $p$ & $\Pi^1_1$ under $\omega$-branching & $\mu$: accessibility of $R$ restricted to $\neg p$ & classical (identification) \\
EG\,$\neg p$: some infinite run stays in $\neg p$ & $\Sigma^1_1$ under $\omega$-branching & $\nu$: the path atom of $R$ restricted to $\neg p$ & classical (identification) (decided: thm.) \\
AF\,$p$ $\vee$ EG\,$\neg p$ & & the dichotomy & thm. \\
AG\,AF\,$p$, with AF\,$p$ as $\neg$EG\,$\neg p$ & & $\nu$ of a body whose leaf is a path atom & classical (identification) (decided: thm.) \\
the controller can force $p$ & $\Pi^1_1$ under $\omega$-branching & $\mu$ of an $\exists\forall$ body & classical (identification) \\
the controller can keep $p$ & $\Sigma^1_1$ under $\omega$-branching & $\nu$ of an $\exists\forall$ body & classical (identification) (decided: thm.) \\
\bottomrule
\end{tabular}
\end{center}
A ``thm.'' marks a theorem instance at the relation or game named:
the dichotomy row is the dichotomy rule at the accessibility body of
$R$ restricted to $\neg p$, sound by \cref{thm:C}, and the
decidedness of the EG\,$\neg p$ and keeping rows is the second clause
of \cref{thm:excl} at that body and at the game body (for games,
\cref{thm:det}); the decidedness of AG\,AF\,$p$ is \cref{thm:agaf}
below.  ``Classical (identification)'' marks a row whose
match of the property with the atom in the third column is a
classical identification, not a theorem of the development; the two
arithmetic rows are a classification only.  Under computably bounded branching
K\H{o}nig's lemma collapses AF and EG to arithmetic.  Unbounded
nondeterminism is exactly where $\mu$, $\nu$ and the dichotomy earn
their keep, and the examples of \S\ref{sec:showpieces} are programs
and games at that boundary.

\paragraph{The alternation boundary.}
The bodies of this paper are closed in their recursive variable: a
nesting is of closed bodies, so AG\,AF\,$p$ is $\nu X.\,(\mathrm{AF}\,p
\wedge \mathrm{AX}\,X)$ with AF\,$p$ a separate closed $\mu$-body.
That covers the properties of computation-tree logic expressible
with closed bodies and the coinduction they need, safety and reachability games, open and closed
determinacy, and the trace predicate of \S\ref{sec:showpieces}.  It
does not cover the dichotomy at every nesting: a nested $\mu$-leaf,
such as AF\,$p$ inside AG\,AF\,$p$, can be valueless, and then the
outer body's decidedness premise fails, as at the trace example's up
body (\cref{thm:app-premise}); a $\nu$-over-$\nu$ nesting whose
inner body has decided leaves has an inner leaf the dichotomy decides
(\cref{thm:excl}).  Writing the inner AF\,$p$ as
$\neg\mathrm{EG}\,\neg p$, a negated path atom, removes the valueless
leaf.  Let $R|_{\neg p}$ be $R$ restricted to the states where $p$
fails, $R|_{\neg p}(y, x)$ when $R(y, x)$ and not $p\cdot x$
(\isn{restrRel}), so that $\Path_{R|_{\neg p}}\cdot x$ is
EG\,$\neg p$ at $x$.  Let $\Phi$ be the body ``$\Path_{R|_{\neg
p}}\cdot x$, or some successor $y$ of $x$ is in the recursive
variable'' (\isn{agafPhi}), whose $\mu$-atom is EF\,EG\,$\neg p$ and
whose dual's $\nu$-atom $\nu(\neg\dual\Phi)$ is AG\,AF\,$p$.

\begin{thm}[\isn{g\_agaf\_premise}, \isn{g\_dec\_agaf\_restr}; grounded]
\label{thm:agaf}
In every calculus with ground induction and the dichotomy, in a
well-formed context, at closed $R$ and $p$: from $\Dec\,R$ and the
decidedness of $\neg p\cdot a$ at every number $a$ (\isn{decSet}),
$R|_{\neg p}$ is decided (\isn{g\_restrRel\_decRel}), the
decidedness premise of $\Phi$ is derivable, and from $x\,\mathbf N$,
$\nu(\neg\dual\Phi)\cdot x$ is decided.
\end{thm}

\noindent The inner leaf is the path atom, which the dichotomy decides at
a decided relation (\cref{thm:excl}), so the outer premise has no
valueless leaf, and one more application of the dichotomy decides the
outer atom (\cref{thm:app-agaf}).  The identification of the atoms
with AF, EG and AG\,AF is classical, as in the table.  A further claim is a prediction: $\Mstar$,
which reads nested closed bodies through the information-order fixed
point of its oracle map, a stratified reading, should agree with the
standard $\mu$-calculus semantics on alternation-free nestings.
Genuine alternation,
$\nu X.\,\mu Y.\,\varphi(X,Y)$ with $X$ free in the inner body, is not
supported: parity conditions, E\,GF\,$p$, and the fair fragments of
the $\mu$-calculus~\cite{emerson86efficient}.  Supporting it would
need bodies with predicate parameters and a priority semantics, and
the paper does not propose it.  Two classical routes around
opposite-polarity dependence fit the design.  Fairness by
transformation to random assignment~\cite{apt83fairness}, so that fair
termination is accessibility over an $\omega$-branching relation, is
carried out in \S\ref{sec:showpieces} for a two-component loop
(\cref{thm:fair}) and, through a term former over lists of guarded
components, for the alternating-bit protocol over lossy channels
(\cref{thm:abp}).  Progress
measures~\cite{klarlund91progress,jurdzinski00small}, which cast
B\"uchi and parity winning regions as $\nu$-certificates over pairs of
a state and a rank, are future work.  In both, the price moves into
the ordinal height of a closed $\mu$-atom, the ladder of the \iga
paper, though on infinite arenas that height need not be bounded by
any fixed notation.

\section{Where Each Rule Is Spent}
\label{sec:rungmap}

\paragraph{The method, and what it certifies.}
Each theory of the examples in the session \isn{OIGA\_Coind}, and the
theories of the dichotomy and its systems, ends with a block of ML
that runs at build time; the two \oga theories that hold the replay
and the bounded K\H{o}nig lemma (\isn{OGA\_Coind\_Kit},
\isn{OGA\_Koenig\_FB}) carry none, and their results enter the map
only through results of \isn{OIGA\_Coind} whose checks follow the
dependencies into them.
For each listed theorem the block collects the transitive
dependencies of the proof, through every proof body, and applies
three tests.  It fails the build if any dependency is an unchecked
assumption (an Isabelle oracle).  It fails if the name of any
dependency contains one of the strings \isn{omv}, \isn{omh},
\isn{Ostar}, \isn{nu\_model}, \isn{mjt} and \isn{mtrue}.  This
test for $\omega$-model-freeness is a naming convention of the development, and
it works because the closure is transitive: a proof that uses an $\omega$-model
reaches, at the end of its dependencies, the theorems that define
the model's relations, whose names contain those strings, so a model
fact under another name, such as the truth of the all-zeros atom in
$\Mstar$, is caught through them.  And the block fails unless named
rules occur under results that are meant to use them (a \emph{needs}
list), and if named rules occur under results that are meant not to
(a \emph{forbids} list); both are substring tests on the names of the
rules' defining facts.  The forbids are live: the same substrings
occur under the needs of other results, and a control that must be
flagged, the system's own soundness theorem, keeps the model test
from being blind.  A further ML control, in the junction session
\isn{Coind\_Dinga} over the two sibling sessions, confirms that every
axiom of the chain beyond those of Isabelle/HOL's \isn{Main} is a
recorded definition, a type definition or a type arity; the sessions
of the coinduction chain proper (\isn{OIGA}, \isn{OIGA\_Coind},
\isn{IGA\_PGSem}, \isn{IGA\_Closing}, \isn{Coind\_Dinga}) add none
of the last two kinds, and the build fails if one appears.

The abstract calculus a theorem is stated in is the other
certificate: a primitive rule the calculus lacks cannot occur in the
proof, and every theorem sits in the least calculus its proof uses.
One case needs both certificates.  In a calculus with ground
induction and the dichotomy, ATI is a derived rule
(\cref{thm:ati-dich}), so the calculus certifies only that ATI is not
a primitive rule there.  The theories of \isn{OIGA\_Coind} and
\isn{OIGA} are built where that derivation is not in scope: the
sublocale edge that carries it lives in \isn{IGA\_Closing}, which is
not their ancestor, so in those theories no fact of ATI's is available
to such a calculus at all.  The check then certifies that the proof
does not apply ATI's rule, and only the re-exports in
\isn{IGA\_Closing} rest on the check alone.  The junction's \dinga
restatements inherit the accounting of their originals; through the
interpretation their dependencies include the derivation of ATI.  An ATI mark in the map
therefore records the rule a proof applies; every row that applies
ATI is one where the dichotomy and ground induction would do, and the
ATI column records the weaker rule, weaker in that \oga and \oiga,
which have it, are sound for the least model, where the dichotomy is
not valid.

A rank, or a count of steps or of inductions, read off a proof script
is called \emph{a reading} in the text and is not formalized.

\paragraph{The map.}
Rules: $\mu$, ground induction; $\nu$, the coinduction rule; Dich,
the dichotomy; ATI, \oga's rule.  A result whose calculus has no rule
beyond \rga's is arithmetic.  ``The search'' is \rga's least-witness
search, a derived rule of \rga, not a fifth rule.

\begin{center}
\begin{minipage}{\textwidth}
\footnotesize
\renewcommand{\thempfootnote}{\alph{mpfootnote}}
\begin{tabular}{@{}p{0.70\textwidth}cccc@{}}
\toprule
result & $\mu$ & $\nu$ & Dich & ATI \\
\midrule
\multicolumn{5}{@{}l}{\emph{Arithmetic: no rule beyond \rga's}}\\
asserting a $\mu$-atom: termination of nested random assignment at every fixed depth $k$, one-round wins, reachability & $\times$ & $\times$ & $\times$ & $\times$ \\
fair termination: the transformed loop, and the fairness twin under strong fairness, every state accessible & $\times$ & $\times$ & $\times$ & $\times$ \\
the alternating-bit protocol: every message eventually delivered & $\times$ & $\times$ & $\times$ & $\times$ \\
the transformation as a term former: decidedness, one helpful direction & $\times$ & $\times$ & $\times$ & $\times$ \\
the decidedness premises of the games and of random assignment & $\times$ & $\times$ & $\times$ & $\times$ \\
Cinderella's strategy total and invariant-preserving; the round-robin play lost below capacity three & $\times$ & $\times$ & $\times$ & $\times$ \\
\midrule
\multicolumn{5}{@{}l}{\emph{Ground induction}}\\
the Lattice and WELL rules (leads-to) & \checkmark & $\times$ & $\times$ & $\times$ \\
the grounded fairness rule: no fair run from an accessible start (with the search) & \checkmark & $\times$ & $\times$ & $\times$ \\
the eventuality rule: a run from a code satisfying an eventually atom meets its target & \checkmark & $\times$ & $\times$ & $\times$ \\
exclusion $\mu \ent \neg\nu$: $\Acc \ent \neg\Path$, $\Term \ent \neg\Div$, $\WinI \ent \neg\WinII$ & \checkmark & $\times$ & $\times$ & $\times$ \\
no run from an accessible state & \checkmark & $\times$ & & $\times$ \\
\midrule
\multicolumn{5}{@{}l}{\emph{Coinduction}}\\
establishing a $\nu$-region: prefix branches, runs give paths, divergence, player II's copying region, Cinderella's region, the weak fairness twin's path atom & $\times$ & \checkmark & $\times$ & $\times$ \\
\midrule
\multicolumn{5}{@{}l}{\emph{The dichotomy}}\\
dichotomy and determinacy: $\Acc \vee \Path$, $\Term \vee \Div$, $\WinI \vee \WinII$, the bucket game, the Cantor--Bendixson split & $\times$ & $\times$ & \checkmark & $\times$ \\
decided $\nu$-regions: $\Dec(\Path)$, $\Dec(\WinII)$, $\Dec(\Div)$ & \checkmark & $\times$ & \checkmark & $\times$ \\
a run from a path; well-founded or an infinite branch & \checkmark & $\times$ & \checkmark & $\times$ \\
player II's strategy (with the search) & $\circ$ & $\times$ & \checkmark & $\times$ \\
\midrule
\multicolumn{5}{@{}l}{\emph{Several rules}}\\
the runs of the divergent twin of random assignment and of the weakly transformed fairness twin & $\circ$ & \checkmark & \checkmark & \\
the Stepmother's loss as a refutation & \checkmark & \checkmark & $\times$ & $\times$ \\
$\Dec(\mathcal O)$ & \checkmark & & \checkmark & \checkmark \\
\midrule
\multicolumn{5}{@{}l}{\emph{ATI}}\\
no natural-number variant & $\times$ & $\times$ & $\times$ & \checkmark \\
termination from non-divergence (ATI as omniscience) & $\times$ & \checkmark & $\times$ & \checkmark \\
termination or divergence, by ATI & $\times$ & $\circ$ & $\times$ & \checkmark \\
prefix sets decided; the round trip (with the search) & $\times$ & $\times$ & $\times$ & \checkmark \\
\bottomrule
\end{tabular}

\smallskip
\checkmark: in the result's calculus, and the check requires it.
$\circ$: in the calculus, and the check neither requires nor forbids
it; the proof may use it.\footnote{A check certifies which rules a
given proof applies.  It does not show that a weaker calculus cannot
prove the statement, so the paper says ``the proof uses'' throughout,
except where a non-derivability theorem exists; and the required
marks are a lower bound on what the proof text spends.}
$\times$: the check forbids it.\footnote{An ATI cross in a row whose
calculus has ground induction and the dichotomy means that the proof
does not apply ATI's rule, not that ATI's conclusions are out of
reach (\cref{thm:ati-dich}); and a check records presence, not
multiplicity, so ``one coinduction step'' is a reading.}
Blank: not in the calculus, and not tested.
\end{minipage}
\end{center}

\paragraph{What the map says.}
Four lines carry the content.  In the examples, asserting a
$\mu$-atom is arithmetic: folding plus number induction, with ground
induction needed only to consume the atom; the termination of nested
random assignment at rank $\omega^k$, for every fixed depth $k$, is an
\rga theorem, with one number induction per level (the rank a
reading), confirmed by a check that forbids every fixed-point rule
under it.  The examples stay below $\omega^\omega$; at higher ranks
accessibility is not arithmetic in general, which is why the \iga
paper builds a ladder of ordinal notations; its proof-grounded rungs
reach Buchholz's notation for the Bachmann--Howard ordinal, of
which Buchholz's system $\mathrm{OT}(1)$ is a second coding, and the
countable part of each of the larger systems $\mathrm{OT}(k)$ for
$k \geq 2$~\cite{iga-paper}.  Using a region, by
unfolding or searching inside it, is arithmetic as well.  A
$\nu$-region costs one coinduction step when it is established, and
nothing afterward.  The proofs of the dichotomy results, the
dichotomies themselves and the decided $\nu$-regions, apply neither
$\nu$ nor ATI; the proof for $\Dec(\mathcal O)$ applies ATI, which
the calculus derives.  With ground induction the dichotomy decides
every $\nu$-region whose body has decided leaves.  The run from a
path, the branch of a tree and player II's strategy apply no
coinduction step: the path or region is a hypothesis, only unfolded.
ATI is applied for the absence of a natural-number variant, for the
round trip of streams, where it decides the prefix set, and in the
trace example as the limited principle of omniscience; nowhere in the
bucket game.

\paragraph{Non-derivability.}
One kind of negative result is available without a new model: from a
proved divergence, the exclusion theorem and the system's
consistency, the accessibility claim is not derivable.  For the
divergent twin of random assignment, neither \oinga nor \dinga proves
$\Acc_{R'}\cdot 0$; under the weak transformation of the fairness
twin, neither proves the initial code accessible, by the same argument
from the twin's path atom.  The statements at \dinga are checked in
the junction theory (\isn{CD\_Junction}).  The consistency used is
the one proved through $\Mstar$.  For the fairness twin, \rga's
non-derivability of the same accessibility is also proved
classically, through the least model: the weakly fair run lifts to an
infinite chain of edges true there.

\section{Related Work}
\label{sec:related}

\paragraph{Two readings of a recursive definition.}
The inductive and coinductive readings of a definition, $\mu$
against $\nu$, go back to Park's fixed-point
induction~\cite{park69fixpoint} and bisimulation~\cite{park81concurrency}
and to Kozen's modal $\mu$-calculus~\cite{kozen83results}; Sangiorgi
traces the history of bisimulation and coinduction across modal
logic, set theory and concurrency~\cite{sangiorgi09origins}, and
Aczel's non-well-founded sets made coinduction the proof principle of
a greatest fixed point and bisimulation its equality~\cite{aczel88nonwellfounded}.

\paragraph{Fixed points of truth.}
Kripke's least fixed point leaves the truth-teller ungrounded while
other fixed points settle it, and Kripke singled out the intrinsic
fixed points, whose verdicts no fixed point contradicts, and the
maximal intrinsic one~\cite{kripke75outline}.  $\Mstar$ is neither
minimal nor intrinsic (\cref{rem:kind}): it commits to the greatest
fixed point for truth at every declared $\nu$-body, a choice made by
the polarity of the declaration, body by body.  Visser's survey
treats the structure of the fixed points of the truth
jump~\cite{visser89semantics}, and Burgess's analysis of their
complexity~\cite{burgess86truth} is the classical reference for what
calibrating $\Mstar$ would involve.  The axiomatic theories of
Kripkean truth take the other route: the Kripke--Feferman theory
KF~\cite{feferman91reflecting}, whose standard models are the fixed
points of the jump, and the partial-logic theory PKF of Halbach and
Horsten~\cite{halbach06axiomatizing}, sound for the fixed points as
partial models, are sound for a whole class of fixed points, where
\inga commits to one non-minimal fixed point and proves what holds
there.  Field's paracomplete theory~\cite{field08saving} builds on
Kripke's least fixed point with a conditional of its own; the grounded
conditional of \S\ref{sec:bg} is instead a premise sequent at an
undecided hypothesis.  Fitting's Kripke--Kleene semantics for logic
programs~\cite{fitting85kripke} and the well-founded
semantics~\cite{vangelder91wellfounded} work with the same split of a
three-valued operator into a truth and a falsity part that
\S\ref{sec:regions} uses.

\paragraph{Least and greatest fixed points in proof theory.}
Baelde's $\mu$MALL~\cite{baelde12least}, linear logic with least and
greatest fixed points, has $\neg\mu$ equal to $\nu$ of the dual by
construction, which is the duality the dichotomy restores at bodies
with decided leaves; Tiu and Momigliano's Linc$^-$~\cite{tiu12cut}
has both Park rules on defined atoms and proves cut elimination.  Both
are proof-search frameworks rather than arithmetics, and both have
the duality that a grounded system lacks.  Over arithmetic, Lubarsky
characterized the sets definable in the $\mu$-calculus over the
integers~\cite{lubarsky93mu}, the definability side of the
alternation boundary of \S\ref{sec:liveness}, and Möllerfeld showed
that arithmetic with arbitrarily nested least and greatest fixed
points, where $\nu$ is the classical dual of $\mu$, has the strength
of $\Pi^1_2$-comprehension~\cite{moellerfeld02generalized}; Curzi and
Das give the computational counterpart for intuitionistic and
circular proofs with fixed points~\cite{curzi25computational}.  The
theories of iterated inductive definitions~\cite{buchholz81iterated}
are the inductive side of the same picture.  The paper claims no
calibration of \dinga or $\Mstar$ against these theories.  Cyclic
proof systems~\cite{brotherston11sequent}, and cyclic arithmetic,
which Das calibrates by quantifier complexity~\cite{das20logical},
obtain induction and coinduction from a global condition on a finite
cyclic proof; whether that condition can be stated and checked in a
grounded setting is an open question the paper does not touch.

\paragraph{Coinduction in proof assistants.}
The coinduction rule here is Park's rule, stated in Isabelle as
\texttt{weak\_coinduct}, which Paulson developed definitionally from
Knaster--Tarski greatest fixed points~\cite{paulson97mechanizing}; the
codatatype package of Isabelle/HOL~\cite{traytel12foundational,blanchette14truly}
and Lochbihler's library of lazy lists and streams, of the kind used
for Nakata--Uustalu-style traces~\cite{lochbihler10coinductive},
supply coinductive types, which the paper does without
(\S\ref{sec:streams}).  Proof by guarded
corecursion~\cite{coquand94infinite}, productivity
checking~\cite{sijtsma89productivity}, parameterized coinduction's
accumulation of the hypothesis~\cite{hur13power}, the
companion~\cite{pous16companion} and up-to
techniques~\cite{pous12enhancements} refine what a certificate may
use; the paper has only coinduction up to union, the form of
Isabelle's \texttt{coinduct} and of its coinductive package's
generated rules, known to users as strong coinduction.  The separation of inductive and coinductive definitions
in Isabelle, Coq and Agda is the folklore of which \cref{thm:clash}
makes one reason exact for this calculus.  Capretta's delay monad and
the treatments of termination against non-divergence built on it are
the type-theoretic analog of the search example, and of
non-termination stated positively as a run~\cite{capretta05general}.
Leroy and Grall's coinductive big-step
semantics~\cite{leroy09coinductive} and the trace semantics of Nakata
and Uustalu~\cite{nakata15hoare} are the sources of the
termination-and-divergence example, whose disjunction is their
``evaluates or diverges'' for the search body; coinductive logic
programming~\cite{simon06coinductive} reads a logic program at its
greatest fixed point in the same spirit.

\paragraph{Verification.}
Recurrent sets, the certificates of non-termination in Gupta and
others' work~\cite{gupta08proving}, are exactly the certificates of
the path atom: a set closed under one step.  Fair termination by
transformation to random assignment is Apt and
Olderog's~\cite{apt83fairness}, and Francez's monograph is the
standard account of fairness and its proof methods~\cite{francez86fairness};
the development's two variants of the transformation are described in
\S\ref{sec:showpieces}, and the alternating-bit protocol over lossy
channels is the standard benchmark for fairness under message loss.
The rank-$\omega$ example is Apt and Plotkin's~\cite{apt86countable},
and the $\Pi^1_1$-completeness of termination under unbounded
nondeterminism Harel's~\cite{harel86effective}.  Progress measures
for B\"uchi and parity conditions are Klarlund's and
Jurdzi\'nski's~\cite{klarlund91progress,jurdzinski00small}.  Lamport's
Lattice rule~\cite{lamport94temporal}, which descends from the proof
lattices of Owicki and Lamport~\cite{owicki82proving}, and the
Manna--Pnueli well-founded rules~\cite{manna92temporal} are the
classical statements that liveness is ranking; \cref{thm:lattice}
derives the Lattice and WELL rules from ground induction.  The bucket
game is from Bodlaender and others~\cite{bodlaender12cinderella} and
Hurkens and others~\cite{hurkens11cinderella}, and Beyene and others
used it as a benchmark for infinite-state safety-game
solving~\cite{beyene14constraint}.  B\"urki's thesis on grounded
verification of recursive trace systems~\cite{burki26grounded} is the
motivating direction: traces produced with finite fuel, safety as a
decided predicate on prefixes, liveness as a $\Pi^0_2$ statement over
the given trace, no coinduction.  \S\ref{sec:liveness} says what that
treatment can state and what needs the new rules.

\paragraph{Determinacy and reverse mathematics.}
Open determinacy is Gale and Stewart's~\cite{gale53infinite}; its
calibration at arithmetical transfinite recursion is
Steel's~\cite{steel77determinateness}, the finer calibrations
Tanaka's~\cite{tanaka90weak}, and Simpson's treatment of the
Cantor--Bendixson theorem at $\Pi^1_1$
comprehension~\cite{simpson09subsystems} is the classical landmark
behind the last example.  Blass's bounds on the complexity of winning
strategies~\cite{blass72complexity} say what a strategy term can and
cannot be.  The Kleene tree is Kleene's~\cite{kleene52introduction}.
The paper claims no calibration of $\Mstar$ in these terms.

\paragraph{The grounded family.}
\rga~\cite{rga-paper}, with its internal reasoning and
consistency~\cite{ford-conrga} and G\"odel's theorems for
it~\cite{goedel-paper}, \oga~\cite{oga-paper} and
\iga~\cite{iga-paper} are the systems extended.  The \iga paper's
derivation of the consistency sentence of a checker for Peano
Arithmetic, in \iga, and the
derivation of \rga's consistency in \oiga that the \iga and G\"odel
papers give, carry over to the systems here whose derivations include
those of \iga or of \oiga (\S\ref{sec:concl}).  The grounded reading
of implication at an undecided hypothesis, and why it is forced, is
the subject of the \iga paper; both rules here are shaped by it.

\section{Limits, Open Questions, and the Artifact}
\label{sec:concl}

Two rules separate \dinga from \iga.  The coinduction rule reads a
body of the $\nu$ class at its greatest fixed point, and the systems
that have it are sound for an explicit model, consistent,
semantically conservative over the least model, and unsound for the
least model by design: they derive never-refuted instances that the
least model leaves without a value.  The dichotomy rule supplies the
classical disjunction between a least-fixed-point atom and the
greatest-fixed-point atom of its co-body, under a derivable premise;
it trades the hypothesis that an order is well-founded for the
derivable one that its steps are decided, and the change of model
pays for the trade (\S\ref{sec:price}).  A second cost, inherited from
ATI (\cref{thm:ati-dich}), is the loss of every certificate-grounded
semantics: \rga and \iga are proof-grounded, while \oinga and \dinga,
like \oiga, have no certificate-grounded semantics; unlike \oiga they
are unsound for the least model by design, and their warrant is
$\Mstar$.  The loss comes with ATI and not with coinduction: \inga,
which has the coinduction rule but neither ATI nor the dichotomy, is
proof-grounded although it is unsound for the least model
(\cref{thm:inga-pg}), and that semantics alone, with no $\omega$-model, gives
its consistency (\cref{cor:inga-con-mf}).  With the two rules, streams, runs, branches, the regions of
a game, and the termination or divergence of a search are single
sentences, and the map of \S\ref{sec:rungmap} says what each costs.

\paragraph{Expressiveness.}
Genuine alternation is not supported, and the paper does not propose
it.  Of the two classical routes around opposite-polarity dependence,
random assignment for fairness is carried out, up to the
alternating-bit protocol's eventual delivery (\cref{thm:abp}), and
progress measures for B\"uchi and parity conditions are the natural
next example.  Streams enter bodies as parameters, not as points.
Stream bisimilarity for systems with infinite branching, and up-to
context techniques, would need codes read by an evaluator (a
prediction, \S\ref{sec:nu}).  Player I's strategy in an open game is
not extracted.

\paragraph{Open questions about the rules.}
Whether \inga proves termination from non-divergence without ATI is
open.  Whether \oinga
derives every instance of the dichotomy is open; $\Mstar$ cannot
separate the two systems.  The calculus of ground induction with the
dichotomy but without $\nu$ is already unsound for the least model,
so it pays the dichotomy's price without $\nu$, and every rule of
\oiga holds in it, ATI as a derived rule; what model and strength
that calculus has is open.  Whether \dinga, or that calculus, derives
\iga's consistency sentence is open as well.  The \iga paper's
\emph{truth term} is truth for \rga's own language as one positive
fixed point, whose body transcribes the least model's evaluation rules
onto codes, and a \emph{truth atom} is the truth term applied at a
point~\cite{iga-paper}.  \iga with bivalence of the truth atom, the
\iga paper's smallest \emph{closing-off} of the truth term, which adds
the decidedness of that one predicate, derives \iga's consistency
sentence, a grounded claim of the \iga paper.  The dichotomy, where it
applies to the truth term's body, concludes only that the truth atom or
the $\nu$-atom of the co-body holds, which is not that bivalence.

\paragraph{Calibration.}
The strength of $\Mstar$ in reverse-mathematical terms, and the
complexity of the $\nu$-truth sets of the controller's region, are
not established beyond a lower bound: $\Mstar$'s truth set is not
hyperarithmetic (\cref{thm:chiO} with the $\Pi^1_1$-completeness of
Kleene's $\mathcal O$, a classical fact of the literature).
The bucket game's threshold capacity is not settled here; the
certificate is for capacity three and above, and that bound is
attained against the strategy it certifies (\cref{thm:app-tight}).
Two consistency sentences of the companion papers carry over by
inclusion of derivations, and read as lower bounds in consistency
strength.  \iga derives, from no hypotheses, the consistency sentence
of a guarded proof checker for Peano Arithmetic
(\isn{gi\_ConPA})~\cite{iga-paper}.  That sentence is a consistency
sentence for Peano Arithmetic, since by a classical theorem, for every
derivation of a formula in Peano Arithmetic, \rga derives that the
checker accepts a code of a proof of that formula
(\isn{pa\_ChkPAr\_complete}).  \oiga derives \rga's consistency
sentence (\isn{oiga\_ConRGA})~\cite{iga-paper,goedel-paper}.
What \iga derives from given hypotheses, \inga, \oinga and \dinga
derive as well, and what \oiga derives, \oinga derives, and so does
\dinga through the equality of \cref{thm:price}
(\isn{gi\_inga\_embed}, \isn{inga\_oinga\_embed},
\isn{inga\_dinga\_embed}, \isn{oiga\_oinga\_embed},
\isn{oiga\_oingad\_embed}).  So \inga, \oinga and \dinga derive the
consistency sentence for Peano Arithmetic, and \oinga and \dinga
derive \rga's.  These are grounded claims, derivations in the
systems, and they are established classically by composing the named
theorems; no single theorem of the development states them at the
coinductive systems.

\paragraph{What is not claimed.}
No calibration of $\Mstar$, of \dinga, or of the calculus of ground
induction with the dichotomy is claimed beyond the lower bounds above
for \inga, \oinga and \dinga, and no threshold capacity of the bucket game.  The run, branch and
strategy terms are total in $\Mstar$; the paper does not claim that
they compute.  No separation of \oinga from \dinga is claimed.  Of the consistency sentences,
only the two above are claimed: G\"odel's second theorem at \inga,
\oinga or \dinga, that the system does not derive its own consistency
sentence, is not claimed, and nothing is claimed about \iga's
consistency sentence in these systems.

\paragraph{The question of the introduction.}
Which non-terminating computations deserve to be colored true?  The
paper's answer, at the level of declared bodies
(\cref{rem:kind}): the instances of a body of the $\nu$ class,
written with $\nu$'s tag, that lie in the greatest fixed point of its
truth operator, relative to the least fixed point of the oracle map.
They are never refuted; nothing the least model values changes; and
reading at the greatest fixed point is not the same as coloring the
never-refuted true, since the two coincide only for bodies that
satisfy \cref{thm:compl}.

\paragraph{The artifact.}
Isabelle/HOL (Isabelle2025-2) checks every theorem in this paper; the
remarks are arguments and are not mechanized.  For the coinduction
chain, \textsf{Coind\_Dinga} with the sessions it loads, no session
enables \texttt{quick\_and\_dirty}, the build-time checks fail
on any oracle, and an axiom check (\isn{CD\_Junction}) confirms that every axiom beyond
those of Isabelle/HOL's \isn{Main} is a recorded definition, a type
definition or a type arity; the sessions of the coinduction chain
proper (\textsf{OIGA}, \textsf{OIGA\_Coind}, \textsf{IGA\_PGSem},
\textsf{IGA\_Closing}, \textsf{Coind\_Dinga}) add none of the last
two kinds, and the build fails if one appears.  The sessions, over the \rga
development (\textsf{RGA}, with the model interface, and
\textsf{RGA\_Kit}), are \textsf{GCalc\_Core} (the base abstract
calculus), \textsf{AGA}, \textsf{OGA} (the replay, K\H{o}nig's lemma
and the Kleene tree: \isn{OGA\_Coind\_Kit}, \isn{OGA\_Koenig\_FB},
\isn{OGA\_Kleene\_Tree}), \textsf{IGA} (the former, the model and
the dichotomy rule: theories \isn{IGA\_Nu}, \isn{IGA\_Nu\_Sound},
\isn{IGA\_Nu\_Rule}, \isn{IGA\_Dich}, \isn{IGA\_Coind} and
companions), \textsf{OIGA} (the abstract calculi with the rules and
the systems: \isn{GNu}, \isn{GDich}, \isn{GShell} and their
companions), \textsf{OIGA\_Coind} (streams, games, fairness and the
other examples, and AG\,AF\,$p$ in the theory \isn{GAgAf}),
\textsf{IGA\_PGSem} (the theorems on
certificate-grounded semantics), \textsf{IGA\_Closing} (ATI from the
dichotomy, the closed-off oracle, \dinga with a primitive-recursive
step decider and its recursively enumerable derivability through the
equality with \oinga with the dichotomy, whose proof checker is
there, and the characteristic term of Kleene's $\mathcal O$ in the
theory \isn{GChiO}), and \textsf{Coind\_Dinga}, the junction over
\textsf{OIGA\_Coind} and \textsf{IGA\_Closing} (the examples restated
at \dinga, the axiom check, and \oinga's recursively enumerable
derivability).  The axiom check runs over the junction theory's
ancestry, and two theories the paper cites are not in it:
\isn{GAgAf} (\cref{thm:agaf}) and \isn{GChiO} (\cref{thm:chiO}),
which no theory of the chain imports.  Neither has a build-time
dependency check of its own; their sessions do not enable
\texttt{quick\_and\_dirty}, and neither theory declares an axiom or
an oracle.  The two consistency sentences of the calibration paragraph come from
two further sessions, of the \iga and G\"odel papers, outside the
chain: \textsf{IGA\_ConPA\_Main} (the consistency sentence for Peano
Arithmetic in \iga) and \textsf{OIGA\_Fold} (\rga's consistency
sentence in \oiga).  The chain's axiom check does not run over them.  Each of the two
headlines has instead its own session's final build-time dependency
check, which fails, among other conditions, on an oracle under the
headline: the theory \isn{ConPA\_Main} for \isn{gi\_ConPA} and the
theory \isn{ST\_Final} for \isn{oiga\_ConRGA}.  No session that the
two further sessions load enables \texttt{quick\_and\_dirty}, and no
theory of those sessions declares an axiom or an oracle.
\inga's certificate-grounded semantics (\cref{thm:inga-pg},
Appendix~\ref{app:inga}) comes from a third further session,
\textsf{INGA\_PGSem}, a child of the chain's session
\textsf{IGA\_PGSem}, in eight theories (\isn{INGA\_PG0},
\isn{INGA\_PG1}, \isn{NG\_Eval}, \isn{NG\_Kit}, \isn{NG\_Closure},
\isn{NG\_Step}, \isn{NG\_GfixStep} and \isn{NG\_Sound}); its results
are classical, HOL theorems about \inga and about abstract calculi,
and they include \inga's consistency from the certified evaluation
alone (\cref{cor:inga-con-mf}).
\textsf{Coind\_Dinga} does not load it, so the chain's axiom check
does not run over it.  The theory \isn{NG\_Sound} ends instead in a
build-time dependency check of its own, which fails on an oracle
under any of sixteen headlines, \isn{inga\_pg} and
\isn{inga\_consistent\_mf} among them, or on any fact about $\Mstar$
under one of them (the model's defining and consistency theorems by
name, and every theorem whose name mentions the model's readings),
with the model's own consistency theorem as a live control.  The
session does not enable \texttt{quick\_and\_dirty}, and no theory of
it declares an axiom or an oracle.
At the tag \artifacttag, the whole development, the chain and the
three further sessions, builds with \texttt{isabelle build -d
hol Coind\_Dinga IGA\_ConPA\_Main OIGA\_Fold INGA\_PGSem}.  A clean build of the
chain alone, fifteen sessions at that tag (\textsf{Coind\_Dinga}, its
parent \textsf{OIGA\_Coind} with that session's ancestors, and
\textsf{IGA\_PGSem} and \textsf{IGA\_Closing}, which
\textsf{Coind\_Dinga} loads through \texttt{sessions}), takes about seven minutes with ten threads on an Apple-silicon
workstation.  Rebuilt on their parent sessions' stored heaps with ten
threads, the three further sessions take about ten seconds
(\textsf{IGA\_ConPA\_Main}), about a minute and a half
(\textsf{OIGA\_Fold}) and about a minute (\textsf{INGA\_PGSem}),
beyond the time of building those parents.
The development was produced by the author directing Claude models
from Anthropic, which mechanized the $\nu$-former and its model, the
coinduction and dichotomy rules and their soundness proofs, the
theorems on certificate-grounded semantics and the examples, with the
author reviewing the definitions and the theorem statements.
Appendix~\ref{sec:table} lists the theorems with their names.  The
development checked is the one at the tag \artifacttag of the
project repository.

\bibliographystyle{plainurl}
\bibliography{logic}

\clearpage
\appendix
\section{The Rules in Full}
\label{sec:rules}

The mechanization uses de~Bruijn indices: in a body, index $0$ is the
point and index $1$ the recursive variable, and under $k$ binders
they shift to $k$ and $k+1$.  $\Phi[\lambda c/a]$ substitutes the
abstraction of the certificate or motive for the recursive variable,
so that a recursive atom $a\cdot u$ becomes $(\lambda c)\cdot u$, one
head step from $c[u]$; $\Phi[t/x]$ substitutes $t$ for the point.

\paragraph{The positivity grammar (\isn{posP}).}
A term in which the recursive variable does not occur is positive
and negative.  $a\cdot u$ is positive when $u$ is a numeral expression
in the bound variables (\isn{pnexp}): a variable bound inside the body
or the point, $0$, or the successor of one.  $\neg b$ is positive when
$b$ is negative and conversely.  $b \vee c$ is positive when both are,
negative when both are.  $\cA\cdot\lambda\,b$ is positive when $b$ is,
and negative when $b$ is, with indices shifted.

\paragraph{The coinduction rule (\isn{nuS}).}
\sloppy Premises: $\Phi$ positive and mentioning the recursive
variable; $\Gamma \ent c[x]$; $i\,\mathbf N,\ c[i],\ \Gamma \ent
\Phi[\lambda c/a][i/x]$; $\Gamma \ent x\,\mathbf N$; $\Phi$ and $c$
well-formed; $i$ not free in $\Gamma$, in $\lambda c$, or in $\lambda
a.\lambda x.\Phi$.  Conclusion: $\Gamma \ent (\nu\Phi)\cdot x$.

\paragraph{The dual (\isn{ndual}).}
A term not mentioning the recursive variable is its own dual.  The
dual of a recursive atom $a\cdot u$, with $u$ a numeral expression in
the bound variables, is $\neg(a\cdot u)$.  The dual commutes with
negation, disjunction and the universal.

\paragraph{The decidedness premise (\isn{nfdec}).} The premises of a body are generated by a relation, not a function.  For a term $t$ not mentioning the recursive variable under $k$ binders, a premise is the universal closure, over the point and the $k$ bound variables, of $t\,\mathbf B$; the clauses below may also descend into such a term, which is why a body can have several premises, and the rule accepts each of them.  For a recursive atom the premise is $0 = 0$.  For $\neg b$ a premise of $b$ is one; for $b \vee c$ the conjunction of one for each; for $\cA\cdot\lambda\,b$ a premise of $b$ one binder deeper.  For the control body $\bot \vee a\cdot x$ every premise is therefore the universal closure of $\bot\,\mathbf B$
conjoined with $0 = 0$, which no substitution makes true.

\paragraph{The dichotomy rule (\isn{dichS}).}
Premises: $\Phi$ positive and mentioning the recursive variable;
$\Psi$ the dual of $\Phi$; $D$ the decidedness premise of $\Phi$;
$\Gamma \ent D$; $\Gamma \ent x\,\mathbf N$; $\Phi$ well-formed.
Conclusion: $\Gamma \ent (\mu\Phi)\cdot x \vee (\nu(\neg\Psi))\cdot x$.

\paragraph{ATI (\isn{pr\_allTI}).}
From $\Gamma \ent \cA\cdot\lambda v.\,(u\cdot v)\,\mathbf B$, which is
proved schematically at a fresh $v$, conclude $\Gamma \ent (\cA\cdot
u)\,\mathbf B$.

\paragraph{The soundness case of the coinduction rule.}
In $\Mstar$, under an arbitrary well-formed substitution for the free
variables, let the certificate set $W$ be the set of well-formed
argument terms $t$ that have a value and at which $c[t]$ is true.
The first premise makes $c[x]$ true and the third gives $x$ a value, so $x \in W$.  At every term $t$ of $W$, the closure premise, with the fresh point instantiated at $t$ itself (its hypothesis $t\,\mathbf N$ holds because $t$ has a value), makes the body true at $t$ with $\lambda c$ for the recursive variable.  The positive transport lemma for the oracle model (\isn{omvO\_pos\_neg\_transport\_P}) turns this into truth with the recursive atoms read in the set of terms at which $\lambda c$ applies truly; because every atom argument is a numeral expression in the point and the bound variables, it is well-formed and valued at the instance, so the reading restricts to $W$.  That gives $W \subseteq F(W)$ for the body's truth operator relative to $\Mstar$.
So $W$ is below the greatest fixed point, which is the oracle's
truth set for the body, and the $\nu$-atom at $x$ is true.  No value
congruence is used; the certificate set is indexed by terms.

\clearpage
\section*{Appendices B to H: the headline proofs in detail}
\addcontentsline{toc}{section}{Appendices B to H}

The seven appendices below give, for each headline theorem, the
statement as it stands in the Isabelle sources with every hypothesis,
the definitions it rests on, the proof as numbered steps that follow
the mechanized proof, with each induction's motive and each
coinduction's region written out, and a list of the lemmas a reader
re-mechanizing the proof in another assistant would need.  They are
written for a competent Master's or early PhD student without logic
expertise.  Isabelle names appear throughout, in \isn{sans serif}, so
that each step can be found in the sources.  Appendix~\ref{sec:rules}
states the rules; Appendix~\ref{sec:table} indexes every theorem.

\paragraph{Common vocabulary.}
Variables are de~Bruijn indices in the sources; the appendices write
named variables and say which index is meant where it matters.  In a
body, index $0$ is the point and index $1$ the recursive variable.  A
term is \emph{well-formed} (\isn{wft}) when it is built correctly from
the constructors, \emph{closed} (\isn{closed}) when no index occurs
free; lifting (\isn{lift}) moves a term under a binder.  $b[u]$
substitutes $u$ for a binder's variable (\isn{bsub}); $\sigma
\downarrow t$ applies a substitution to every free variable at once
(\isn{subf}), and a substitution is \emph{well-formed} (\isn{wsubf})
when its images are.  A sequent $\Gamma \ent c$ is a derivation in
\rga, or in the system or abstract calculus named; $\Gamma$ is
well-formed (\isn{wfc}) when its members are; $c\,\mathbf N$ says
$c$ is a number and $c\,\mathbf B$ that it is decided.  A derivation
is a \emph{proof list}, each entry a sequent obtained by one rule from
entries before it; the list \emph{contains} $H \ent c$ when $H \ent c$
is one of its entries.  $H \ent c$ is \emph{valid} in a model $M$
when, for every well-formed substitution $\sigma$ that makes every
member of $H$ true in $M$ ($\sigma$ \emph{satisfies} $H$), $\sigma
\downarrow c$ is true in $M$.  An \emph{eigenvariable} is the fresh
variable that a rule binding a variable (universal introduction,
ground induction, coinduction) introduces in its premise, fresh for
everything else in the premise, so that the premise holds at every
instance of the variable.  ATI in its \emph{open form}
(\isn{ir\_allTI\_open}) derives the decidedness of a universal from
the decidedness of its body at a fresh number variable.

\emph{The least model} is the pair of mutually inductive relations
of \cite{rga-paper}: $t$ reaches the head form $w$ (\isn{omh}) and
$t$ has the value $r$ (\isn{omv}), with $1$ for true and $0$ for
false, twenty-one rules in all, and a term has at most one value.
\emph{The oracle model} at an oracle $\Theta$ is the same pair of
relations with a $\nu$-head for each body of the $\nu$ class, a guard
that keeps the beta clause from unfolding a $\nu$-term of the class,
and two clauses that read an application whose head is a $\nu$-head at
the oracle's verdict on the argument term (\isn{omhO}, \isn{omvO}).  The \emph{oracle map}
(\isn{Om}) sends $\Theta$ to the oracle that reads each body of the
class at the greatest fixed point of its truth operator and the least
fixed point of its falsity operator, both relative to the model at
$\Theta$; $\Theta^*$ (\isn{Ostar}) is its least fixed point and
$\Mstar$ is the oracle model at $\Theta^*$, so that ``$t$ is true in
$\Mstar$'' is \isn{omvO}~\isn{Ostar}~$t$ $\{\mathit{True}\}$.  Every value
of the least model persists in $\Mstar$ (\isn{omv\_omvO\_Ostar}).

\emph{The oracle readings.}  A body with recursive index $k$ is read
at an environment with its recursive atoms interpreted as membership
in a set: the truth reading over the least model (\isn{osem}) says
when the body reads true, and requires the atom's argument to be a
literal numeral.  Over the oracle model at $\Theta$ the readings come
in a truth form and a falsity form (\isn{osemO}, \isn{fsemO}), and
the reading relative to $\Mstar$ with numeral membership is
\isn{osemM}.  The truth operator of a body sends a set to the set of
points at which the body reads true (numerals for the least-model
reading, argument terms for the oracle readings); the falsity operator
likewise with the falsity reading.

\emph{The dual and the premise.}  The dual of a body (\isn{ndual}) is
a relation between bodies that negates each recursive atom and
nothing else; it is functional (\isn{ndual\_fun}) and exists for every
body in the positive or the negative grammar (\isn{ndual\_exists}).
The decidedness premise (\isn{nfdec}) is likewise a relation between
a body and a term: it holds of a body and any term that conjoins
decidedness claims for the body's non-recursive leaves, so a body may
have several premises, and the rule accepts any of them.

\emph{The abstract calculi.}  The development states each grounded
theorem in a locale that has exactly the rules the proof uses, over
\rga: \isn{gcalc} has none beyond \rga's; the suffixes \isn{gfix},
\isn{ati}, \isn{nu} and \isn{dich} add ground induction, \oga's rule,
the coinduction rule and the dichotomy, in every combination
(\isn{gcalc\_gfix\_nu}, \isn{gcalc\_ag} for ground induction with
ATI, \isn{gcalc\_ag\_nu\_dich} for all four, and so on).  A result
stated in one calculus holds in every system that interprets the
calculus; \inga interprets \isn{gcalc\_gfix\_nu} and \oinga
\isn{gcalc\_ag\_nu}.  Every calculus with ground induction and the
dichotomy interprets \isn{gcalc\_ati} (\cref{thm:ati-dich}).  \dinga
(\isn{dinga}) interprets \isn{gcalc\_ag\_nu\_dich}, and so does
\isn{oinga\_dich}, \oinga with the dichotomy, which has the same
derivability relation (\cref{thm:price}); results stated for
\isn{oinga\_dich}, under \isn{oingad\_} names, are results about
\dinga, and those the paper cites at \dinga are restated there by name
in the junction session \isn{Coind\_Dinga} (Appendix~\ref{app:dich}).  Facts about
$\Mstar$ are HOL theorems in the locale \isn{rga}.

\section{The model and the soundness of the coinduction rule}
\label{app:nu}

Files: \isn{RGA/RGA\_Omega\_Sound.thy} (the model interface and the generic
soundness theorem), \isn{IGA/IGA\_Nu.thy} (the $\nu$-former, the oracle
model, the oracle map, coherence, head-step invariance, the all-zeros
atom), \isn{IGA/IGA\_Nu\_Sound.thy} (the value congruence of $\Mstar$ and
\isn{rga\_nu\_sound}), \isn{IGA/IGA\_Nu\_Rule.thy} (the coinduction rule,
\inga, its soundness and conservativity, the clash), \isn{IGA/IGA\_Coind.thy} (a body read at both fixed points), and
\isn{IGA/IGA\_KleeneO\_Meas.thy} (the least model's falsity reading and its
transport lemma).  All results
are classical: HOL theorems in the locale \isn{rga}, or in \isn{inga}
for the system.

\paragraph{Definitions.}
The \emph{$\nu$-former} (\isn{trmNu}) is $\nu\Phi = s_\Phi\cdot s_\Phi$
with $s_\Phi = \lambda z.\,I\cdot((\lambda a.\lambda x.\Phi)\cdot(z\cdot
z))$ and $I = \lambda y.\,y$; the identity redex makes $\nu\Phi$ a
different term from every $\mu\Psi$ (\isn{trmNu\_neq\_trmGfix}), and
$\nu\Phi\cdot x$ reaches $\Phi[\nu\Phi/a][x/x]$ in four head steps
(\isn{nu\_step1} to \isn{nu\_step4}).  Three positivity grammars share
the clauses for non-occurrence, disjunction, negation (which swaps the
positive and negative halves) and the universal (one index deeper).
They differ at the recursive atom $a\cdot u$: \isn{pos} asks $u$ to be
a numeral expression in the variables other than $a$; the rule's
grammar \isn{posP} asks for a numeral expression in the variables
below $a$ (\isn{pnexp}: the point and the quantifier-bound variables,
never a parameter); and the model's grammar \isn{posA} accepts any $u$
in which $a$ does not occur (\isn{posP\_posA}: \isn{posP} implies
\isn{posA}).  A body is \emph{in the $\nu$ class} (\isn{nuGood}) when it
is \isn{posA}-positive at index $1$ and index $1$ occurs in it
(\isn{nuGood\_of\_posP}: a \isn{posP}-positive body that mentions $a$ is
in the class).  A body of the class is a disjunction, negation,
equation, universal or recursive atom (\isn{nuGood\_fshape}).

An \emph{oracle} $\Theta$ assigns to each body $\Phi$ two sets of terms,
$\Theta(\Phi, 1)$, the arguments at which $\nu\Phi$ is to be true, and
$\Theta(\Phi, 0)$, those at which it is to be false; oracles are ordered
pointwise by inclusion of both sets.  An oracle is \emph{consistent}
(\isn{ocons}) when the two sets are disjoint for every body.  The
\emph{oracle model} (\isn{omhO}, \isn{omvO}) is a pair of simultaneous
inductive relations, the head judgment $\mathrm{hd}_\Theta(t, w)$ and
the value judgment $\mathrm{val}_\Theta(t, r)$, with nine head clauses
and fifteen value clauses.  Its heads are those of the least model
plus one $\nu$-head $\nu_\Phi$ per body.  The clauses are the least
model's with three changes.
\begin{enumerate}
\item The head-level beta clause (\isn{omhO\_appLam}) carries the guard
  that $t\cdot u$ is not $\nu\Phi$ for a body $\Phi$ of the class, so a
  $\nu$-term of the class is never unfolded as a head.
\item A new head clause (\isn{omhO\_nu}): for $\Phi$ in the class,
  $\mathrm{hd}_\Theta(\nu\Phi, \nu_\Phi)$.
\item Two oracle clauses (\isn{omvO\_nuT}, \isn{omvO\_nuF}): if
  $\mathrm{hd}_\Theta(t, \nu_\Phi)$ and $u \in \Theta(\Phi, b)$, then
  $\mathrm{val}_\Theta(t\cdot u, b)$.  The oracle is consulted at the
  argument term $u$, not at a value of $u$, and at every application
  whose function part evaluates to a $\nu$-head.
\end{enumerate}
A $\nu$-term of the class has only the $\nu$-head (\isn{omhO\_trmNu\_inv}),
so its atoms take values only through the oracle
(\isn{omvO\_nu\_atom\_inv}, \isn{omvO\_nu\_atomI}), and a $\nu$-term itself
has no value (\isn{omvO\_trmNu\_none}).

The \emph{oracle readings} of a body are inductive relations like the
least model's readings, with the set indexed by terms.
$\mathrm{rd}^+_\Theta(S, k, \sigma, t, b)$ (\isn{osemO}) reads a part
in which $k$ does not occur by $\mathrm{val}_\Theta$ of its
$\sigma$-instance, reads an atom $v_k\cdot u$, with $k$ not occurring
in $u$, true when the term $\sigma\downarrow u$ lies in $S$ and never
false, and reads disjunction,
negation and the universal (over numerals) clause by clause.
$\mathrm{rd}^-_\Theta(X, k, \sigma, t, b)$ (\isn{fsemO}) is the dual: an
atom reads false when $\sigma\downarrow u \in X$ and never true.  The
environment of a body $\Phi$ at the point $u$ (\isn{nuenv}) sends $0$
to $u$, $1$ to $\nu\Phi$ and $v + 2$ to $v$, so that its instance is
$\Phi[\nu\Phi/a][u/x]$ (\isn{nuenv\_gsub}).  The truth and falsity
operators relative to $\Theta$ (\isn{FO}, \isn{GO}) are $F^\Theta_\Phi(S)
= \{u \mid \mathrm{rd}^+_\Theta(S, 1, \mathrm{env}_\Phi(u), \Phi, 1)\}$
and $G^\Theta_\Phi(X) = \{u \mid \mathrm{rd}^-_\Theta(X, 1,
\mathrm{env}_\Phi(u), \Phi, 0)\}$; the two operators are monotone in the
set and in the oracle, for every body (\isn{FO\_mono2}, \isn{GO\_mono2},
\isn{FO\_mono}, \isn{GO\_mono}).  The \emph{oracle map} (\isn{Om}) is
$\mathrm{Om}(\Theta)(\Phi, 1) = \gfp\,F^\Theta_\Phi$ and
$\mathrm{Om}(\Theta)(\Phi, 0) = \lfp\,G^\Theta_\Phi$ for $\Phi$ in the
class, both empty otherwise; the map is monotone (\isn{Om\_mono}), and
$\Theta^* = \lfp\,\mathrm{Om}$ (\isn{Ostar}).  Hence, for $\Phi$ in the
class, $\Theta^*(\Phi, 1) = \gfp\,F^{\Theta^*}_\Phi$ and $\Theta^*(\Phi,
0) = \lfp\,G^{\Theta^*}_\Phi$ (\isn{Ostar\_T}, \isn{Ostar\_F}).  $\Mstar$
is $\mathrm{val}_{\Theta^*}$ with its head judgment.

The \emph{model interface} (\isn{omega\_model}) is a locale over a
relation $M(t, r)$ with twenty-four assumptions: uniqueness of values;
the inversions and introductions of the least model's clauses for zero,
successor, equation, negation, disjunction and the universal; the
recursor's successor clause (an introduction only); no value for
$\bot$; beta, $M((\lambda d)\cdot u, r)$ iff $M(d[u], r)$; head-step
invariance in both directions (\isn{M\_hstep}, \isn{M\_hstep\_rev}); and
two value congruences, \isn{M\_vc} (if $\sigma$ and $\sigma'$ are well-formed and
agree at every variable up to replacing a term by a well-formed term
with the same value, then $M(\sigma\downarrow p, r)$ implies
$M(\sigma'\downarrow p, r)$) and \isn{M\_vcb} (if $u$ and $u'$ are
well-formed with the same value, $M(b[u], r)$ implies $M(b[u'], r)$).

The \emph{value congruence} of $\Mstar$ (\isn{ovcrdO}, at depth $d$;
\isn{ovcrO} at depth $0$) relates two terms that are equal except that,
at some positions, they hold well-formed terms with the same value in
$\Mstar$ (lifted past $d$ binders).  The congruence image of a
$\nu$-term $\nu\Psi$ is a \emph{mixed self-application} $s_{\Psi_1}\cdot
s_{\Psi_2}$ (\isn{mix}) with $\Psi_1$ and $\Psi_2$ congruent to $\Psi$ at
depth $2$ (\isn{ovcr\_trmNu\_E}); when the halves differ the mixed term
is not a $\nu$-term, and $\mathrm{mix}(\Psi_1, \Psi_2)\cdot u$ has exactly the
value of $\Psi_1[\nu\Psi_2/a][u/x]$ (\isn{mix\_val}).

\subsection{\rga is sound for $\Mstar$}

\begin{thm}[\isn{rga\_nu\_sound}]
\label{thm:app-rga-nu-sound}
If $\Gamma \ent c$ is derivable in \rga, $\Gamma$ is well-formed, and
every member of $\Gamma$ is true in $\Mstar$, then $c$ is true in
$\Mstar$.
\end{thm}

Truth here is the value judgment on $c$ itself, the identity
substitution; the form for every well-formed substitution is
\isn{nu\_model.rga\_sound\_mjt}.  The theorem is one line: the generic
soundness theorem of the interface, instantiated through the
sublocale $\isn{rga} \subseteq \isn{omega\_model}(\mathrm{val}_{\Theta^*})$,
whose justification \isn{nu\_model\_omega\_model} carries the work.

\paragraph{Proof.}
\begin{enumerate}
\item \emph{Generic soundness} (in \isn{omega\_model}, nothing specific to
  $\Mstar$).  \isn{rga\_M\_sound} is \isn{rga\_sound\_mjt} at the identity
  substitution (\isn{wsubf\_var}); \isn{rga\_sound\_mjt} passes to a
  valid proof list containing the judgment.  \isn{prf\_mjt} is an induction on that list,
  with motive:
  for every $H$ and $c$, if the list contains $H \ent c$ and $H$ is
  well-formed, then $H \ent c$ is valid in $M$.  At a list $j :: P'$, a
  judgment is either already in $P'$, or is $j$, obtained by one rule
  from judgments of $P'$; \isn{rga\_step\_mjt} then treats the
  forty-two rules case by case.  Each case fixes a well-formed $\sigma$
  satisfying $H$, applies the premises' validity at $\sigma$ or at an
  updated $\sigma$, and closes with interface facts; for instance the
  conversion rules use the head-step facts, and \isn{pr\_eq\_subs},
  \isn{pr\_nat\_ind} and the universal rules use the value congruence.
\item \emph{Reduction to the value congruences}
  (\isn{nu\_model\_iff\_vc}): $\mathrm{val}_{\Theta^*}$ satisfies the
  interface if and only if it satisfies \isn{M\_vc} and \isn{M\_vcb}.
  Only the ``if'' direction is used.  Of the twenty-two other facts, the
  inversions hold at every oracle by case analysis on the inductive
  definition and the introductions are clauses; uniqueness is step~3 and
  head-step invariance step~4.
\item \emph{Determinism} (\isn{omvO\_Ostar\_unique}).
  \begin{enumerate}
  \item At every consistent oracle, \isn{omhO\_det\_excl} is a
    simultaneous rule induction with motive: a term with head $w$ has no
    other head and no value; a term with value $r$ has no other value
    and no head.  In the oracle cases a true and a false reading at one
    term would put the argument in both of the oracle's sets.  In the
    beta case the guard of clause~1 and \isn{omhO\_trmNu\_inv} separate a
    $\lambda$-head from a $\nu$-head.  The consequence is
    \isn{omvO\_unique}.
  \item $\Theta^*$ is consistent (\isn{ocons\_Ostar}), by fixed-point
    induction along the transfinite iteration of $\mathrm{Om}$
    (\isn{ccpo\_class.fixp\_induct}).  The base, the empty oracle, is
    consistent.  Consistency is admissible (\isn{ocons\_admissible}): a
    term in the union of the truth sets and the union of the falsity
    sets of a chain lies in both sets of the larger of the two members
    involved.  The step (\isn{ocons\_Om}) is \isn{gfpO\_lfpO\_disjoint}:
    at a consistent $\Theta$, $\gfp\,F^\Theta_\Phi$ and
    $\lfp\,G^\Theta_\Phi$ are disjoint.  Write $N = \gfp\,F^\Theta_\Phi$.
    By \isn{lfp\_lowerbound} it suffices that the complement of $N$ is
    closed under $G^\Theta_\Phi$.  A $u$ in $G^\Theta_\Phi(-N) \cap N$
    would satisfy $\mathrm{rd}^+_\Theta(N, \ldots, 1)$, since $N$ is a
    fixed point, and $\mathrm{rd}^-_\Theta(-N, \ldots, 0)$.
    \isn{osemO\_fsemO\_excl}, an induction on the truth reading at
    disjoint sets, refutes the pair: at a non-recursive leaf by
    \isn{omvO\_unique}, at an atom by disjointness, and at connectives by
    the inversions of the falsity reading.
  \end{enumerate}
\item \emph{Coherence and head-step invariance.}
  \begin{enumerate}
  \item Two transports, at any oracle.  \isn{osemO\_term}, an induction
    on the reading, turns $\mathrm{rd}^+_\Theta(S, k, \sigma, t, b)$, for
    $\Phi$ in the class with $S \subseteq \Theta(\Phi, 1)$ and $\sigma k
    = \nu\Phi$, into
    $\mathrm{val}_\Theta(\sigma\downarrow t, b)$; the atom case is
    \isn{omvO\_nu\_atomI}.  \isn{fsemO\_term} is the dual.
    \isn{omvO\_pos\_neg\_transport} goes back, by induction on
    \isn{posA}/\isn{negA}: a true instance of a positive part reads true
    at $\Theta$'s own truth set, a false one reads false at the falsity
    set, and dually for a negative part.  The atom case is
    \isn{omvO\_nu\_atom\_inv}; the model reads a $\nu$-atom only through
    the oracle, so the induction is structural.  Together:
    \isn{omvO\_body\_T\_iff}, for $\Phi$ in the class,
    $\mathrm{val}_\Theta(\Phi[\nu\Phi/a][u/x], 1)$ iff $u \in
    F^\Theta_\Phi(\Theta(\Phi, 1))$, and \isn{omvO\_body\_F\_iff} likewise
    with $G^\Theta_\Phi$.
  \item \emph{Coherence} (\isn{nu\_coherent\_T}, \isn{nu\_coherent\_F}):
    for $\Phi$ in the class, $u \in \Theta^*(\Phi, b)$ iff
    $\Phi[\nu\Phi/a][u/x]$ has value $b$ in $\Mstar$.  By \isn{Ostar\_T}
    and \isn{gfp\_fixpoint}, $F^{\Theta^*}_\Phi(\Theta^*(\Phi, 1)) =
    \Theta^*(\Phi, 1)$; apply (a).  Falsity likewise, with
    \isn{lfp\_fixpoint}.
  \item \isn{nu\_head\_val}: if $\mathrm{hd}_{\Theta^*}(t, \nu_\Phi)$, then
    $t\cdot u$ and $\Phi[\nu\Phi/a][u/x]$ have the same values in
    $\Mstar$.  Forward: by determinism the value of $t\cdot u$ came from
    an oracle clause; apply coherence.  Backward: an instance of a body
    of the class has only truth values (\isn{omvO\_body\_bool}), so
    coherence and an oracle clause apply.
  \item \isn{omO\_hstep\_sim} is an induction on the head step from a
    redex to its reduct, with a two-part motive: every value of the
    reduct is a value of the redex; every head $w$ of the reduct has a
    head $w_0$ of the redex with $w_0 = w$, or one of $w_0, w$ the
    $\nu$-head $\nu_\Phi$ and the other the $\lambda$-head of
    $\Phi[\nu\Phi/a]$ (\isn{heqN}).  At a beta step whose redex is a
    $\nu$-term $\nu\Phi$ of the class, the reduct is $I\cdot((\lambda
    a.\lambda x.\Phi)\cdot\nu\Phi)$, whose only head is the $\lambda$-head
    of the unfolding (\isn{omhO\_Ostar\_nuI\_redex\_eq}), paired with
    $\nu_\Phi$; any other beta step uses the beta clause.  The recursor
    steps are direct.  At a step inside the function part of an
    application, related heads $\nu_\Phi$ and the unfolding's
    $\lambda$-head close by \isn{nu\_head\_val}, and by
    \isn{omhO\_body\_none} (an instance of a body of the class has no
    head).  \isn{omO\_hstep\_rev} is the mirror, from redex to reduct.
  \end{enumerate}
\item \emph{The value congruence of $\Mstar$} (\isn{omvO\_vcr\_ex}): if
  $\mathrm{val}_{\Theta^*}(X, r)$ and $X$ is congruent to $Z$, then
  $\mathrm{val}_{\Theta^*}(Z, r)$.  The proof evaluates at an auxiliary
  oracle, the \emph{congruence interior} $\Theta^\circ$ (\isn{Ocg}):
  $u \in \Theta^\circ(\Psi, b)$ when $u \in \Theta^*(\Psi, b)$ and, for
  all $\Psi_1, \Psi_2$ congruent to $\Psi$ at depth $2$ and every $u'$
  congruent to $u$, the term $\Psi_1[\nu\Psi_2/a][u'/x]$ has value $b$
  in $\Mstar$.  Directly $\Theta^\circ \le \Theta^*$ (\isn{Ocg\_le}).
  \begin{enumerate}
  \item The class is invariant under the congruence at depth $2$
    (\isn{ovcrdO\_nuGood}), by induction on \isn{posA}/\isn{negA}, using
    that congruence preserves whether a low index occurs
    (\isn{ovcrdO\_free\_low}).
  \item \isn{omO\_vcr\_transport} is a simultaneous rule induction on
    $\mathrm{hd}_{\Theta^\circ}$ and $\mathrm{val}_{\Theta^\circ}$, with
    motive: if $X$ has head $w$, every $Z$ congruent to $X$ has, in
    $\Mstar$, a head related to $w$ (\isn{nhrel}: the least model's head
    relation, plus a $\nu$-head related to a $\nu$-head of a congruent
    body or to the $\lambda$-head of a mixed unfolding, and conversely);
    if $X$ has value $r$, every such $Z$ has value $r$ in $\Mstar$.  The
    least model's cases replay \rga's congruence proof with the
    congruence inversions.  Four cases are new.  At $\nu\Psi$
    (\isn{tr\_case\_nu}) the image is a mixed self-application: a
    $\nu$-term of the class with the $\nu$-head when the halves agree
    (by (a)), otherwise a term with the $\lambda$-head of the mixed
    unfolding.  The head-level beta case with a $\nu$-term image is
    \isn{tr\_case\_appLam\_h}; the value-level beta case is
    \isn{tr\_case\_appLam\_v}.  At an oracle clause
    (\isn{tr\_case\_nuTF}), $u \in \Theta^\circ(\Phi, b)$: if the image
    is a swapped leaf, determinism closes; otherwise the image is
    $t'\cdot u'$, $t'$ has head $\nu_{\Psi_1}$ or the $\lambda$-head of
    $\mathrm{mix}(\Psi_1, \Psi_2)$, the definition of $\Theta^\circ$
    supplies the value $b$ of $\Psi_1[\nu\Psi_2/a][u'/x]$, and
    \isn{nu\_head\_val} or the beta clause gives the value of
    $t'\cdot u'$.
  \item \isn{osemO\_vcr} and \isn{fsemO\_vcr}, by induction on the reading
    at $\Theta^\circ$: from $\mathrm{rd}^\pm_{\Theta^\circ}(S, k, \rho, t,
    b)$, a congruent $t'$, related environments (\isn{envrel}) and a set
    map (whenever $x \in S$ and $x$ is congruent to $x'$, then $x' \in
    S'$), conclude $\mathrm{rd}^\pm_{\Theta^*}(S', k, \rho', t', b)$.
    Non-recursive leaves go through (b).
  \item \emph{A coinduction} (\isn{Ocg\_gfp\_cong}).  Region: $R$, the
    set of $u'$ congruent to some $u \in \gfp\,F^{\Theta^\circ}_\Psi$,
    for a $\Psi$ of the class congruent to $\Psi_1$.  Closure obligation:
    $R \subseteq F^{\Theta^*}_{\Psi_1}(R)$, which follows by unfolding the
    greatest fixed point and applying (c) with the environments
    $\mathrm{env}_\Psi(u)$ and $\mathrm{env}_{\Psi_1}(u')$
    (\isn{envrel\_nuenv}, \isn{ovcrO\_trmNu}).  So $R \subseteq
    \gfp\,F^{\Theta^*}_{\Psi_1} = \Theta^*(\Psi_1, 1)$.
  \item \isn{Ocg\_prefix}: $\mathrm{Om}(\Theta^\circ) \le \Theta^\circ$.
    Let $u \in \mathrm{Om}(\Theta^\circ)(\Phi, b)$; monotonicity of
    $\mathrm{Om}$ and $\Theta^\circ \le \Theta^*$ give $u \in
    \Theta^*(\Phi, b)$.  The truth half: (c) with the set map of (d)
    gives $\mathrm{rd}^+_{\Theta^*}(\Theta^*(\Psi_2, 1), 1,
    \mathrm{env}_{\Psi_2}(u'), \Psi_1, 1)$, and \isn{osemO\_term} turns
    that into the truth of $\Psi_1[\nu\Psi_2/a][u'/x]$
    (\isn{nuenv\_gsub2}).  The falsity half is a least-fixed-point
    induction (\isn{lfp\_induct}) with motive: every mixed unfolding of
    congruent bodies at a congruent argument is false at $x$.  Its step
    applies \isn{fsemO\_vcr} with the set map that motive supplies and
    \isn{nu\_coherent\_F}, then \isn{fsemO\_term}.
  \item By \isn{lfp\_lowerbound}, $\Theta^* \le \Theta^\circ$
    (\isn{Ostar\_le\_Ocg}), so $\Theta^\circ = \Theta^*$
    (\isn{Ocg\_eq\_Ostar}).  \isn{omvO\_vcr\_ex} composes monotonicity of
    the model in the oracle, which moves a value at $\Theta^*$ to
    $\Theta^\circ$, with the value half of (b).
  \item \isn{M\_vc} is \isn{omvO\_vc\_var}: the two substitutions are
    congruent pointwise (\isn{ovcrO\_refl}, \isn{ovcrO\_swap}), hence on
    every term (\isn{ovcrO\_subf1}).  \isn{M\_vcb} is \isn{omvO\_vc\_bsub},
    through \isn{ovcrO\_bsub}.
  \end{enumerate}
\item \emph{Close.}  \isn{nu\_model\_omega\_model} combines step~2 with
  the two instances of step~5(g); the sublocale and the theorem follow.
\end{enumerate}

\paragraph{The proof for every coherent oracle.}
Of $\Theta^*$ the proof above uses only that the oracle is consistent,
that its truth and falsity sets are fixed points of the two operators
at every body of the class and empty elsewhere, and two inclusions:
$\Theta^* \le \Theta^\circ$ of step~5, and the analogous inclusion in
the oracle's substitution interior, which makes the oracle's sets
closed under substitution (\isn{Ostar\_le\_Ocg},
\isn{Ostar\_le\_Osub}).  The development states these as the
assumptions of three nested locales over an arbitrary oracle,
\isn{cohO} (consistency and coherence, in \isn{IGA/IGA\_Nu.thy}),
\isn{cohO2} (the substitution inclusion) and \isn{cohO3} (the
congruence inclusion, in \isn{IGA/IGA\_Nu\_Sound.thy}), and proves the
coherence facts, the persistence of the least model's values
(\isn{cohO.omv\_le\_omvO\_Oz}) and the interface
(\isn{cohO3.nu\_model\_omega\_model}) once in them.  $\Theta^*$ is the
instance (\isn{cs}, \isn{cs2}, \isn{cs3}), and every fact above keeps
its name in the locale \isn{rga}.  The closed-off oracle
(\isn{OstarC}, the least fixed point of the map that reads each body
at the least fixed point of its truth operator and the greatest fixed
point of its falsity operator) is a second instance, with its own
proofs of the two inclusions, and its model satisfies the interface
(\isn{closed\_off\_omega\_model}, in
\isn{IGA\_Closing/CM\_Closed.thy}).

\subsection{The coinduction rule and the systems}

The rule (\isn{nuS}) has the side conditions: $\Phi$ \isn{posP}-positive
at index $1$; index $1$ free in $\Phi$; $\Phi$ and $c$ well-formed; $i$
not free in $\Gamma$, in $\lambda c$, or in $\lambda a.\lambda x.\Phi$.
Its premises are, in the proof list, $\Gamma \ent c[x]$, the closure
premise $i\,\mathbf N,\ c[i],\ \Gamma \ent \Phi[\lambda c/a][i/x]$, and
$\Gamma \ent x\,\mathbf N$; its conclusion is $\Gamma \ent
\nu\Phi\cdot x$.  \inga's step relation (\isn{inga\_derives}) is \iga's,
\rga's forty-two rules and ground induction (\isn{pr\_gfix}), or a
coinduction step.  The locale \isn{omega\_model\_gfix} adds to the
interface the ground-induction property \isn{M\_gfix\_ind\_par}: for
$\Phi$ \isn{posP}-positive and a term $P$ with $M(P\cdot\bar n, 1)$
whenever $M(\Phi[P/a][\bar n/x], 1)$, the truth of $\mu\Phi\cdot\bar n$
implies the truth of $P\cdot\bar n$.

\begin{thm}[\isn{inga\_nu\_sound}, \isn{inga\_consistent}]
\label{thm:app-inga-sound}
If $\Gamma \ent c$ is derivable in \inga, $\Gamma$ is well-formed and
every member of $\Gamma$ is true in $\Mstar$, then $c$ is true in
$\Mstar$.  $\emptyset \ent 0 = 1$ is not derivable in \inga.
\end{thm}

\paragraph{Proof.}
\begin{enumerate}
\item \emph{Skeleton.}  \isn{step\_M\_sound} is \isn{prf\_mjt} abstracted
  over the step relation (\isn{step\_prf\_mjt}), followed by the
  identity substitution.  Its inputs are \isn{nu\_model\_omega\_model},
  well-formedness of steps (\isn{inga\_wft}) and step soundness
  (\isn{inga\_step\_mjt}), which splits on the disjunction defining the
  step relation.  \isn{inga\_consistent} is \isn{step\_M\_consistent}
  with the same inputs, since $0 = 1$ is never true.
\item \emph{\iga steps} (\isn{nu\_model\_gfix.gi\_step\_mjt}).  An \iga step
  is an \rga step, handled by the preceding theorem's step~1, or a
  ground-induction step, handled by \isn{gfix\_step\_mjt}, written once
  for every model of \isn{omega\_model\_gfix}: from the premises,
  $\sigma\downarrow x$ has a value $v$; \isn{M\_vc\_arg} moves the
  $\mu$-atom's truth to $\bar v$; the closure premise at the
  substitution sending $i$ to each numeral $\bar n$ discharges the
  hypothesis of \isn{M\_gfix\_ind\_par}; beta and \isn{M\_vcb} move back
  to $\sigma\downarrow x$.  $\Mstar$ has the property
  (\isn{nu\_omega\_model\_gfix}) by \isn{gfix\_omvO\_ind\_par}, through the
  set form \isn{gfix\_omvO\_ind\_set\_par} at $S = \{m \mid P\cdot\bar m$
  true in $\Mstar\}$ and \isn{osemM\_term}.  The set form is the first
  component of the motive of \isn{gfix\_transportO\_par}, one rule
  induction on the derivations of $\Mstar$ with the five-part motive of
  the \iga paper's transport lemma (\isn{gfInvM}), read in the
  numeral-set reading relative to $\Mstar$ (\isn{osemM}).  The five
  parts say: a true $\mu$-atom at a numeral $\bar m$ has $m \in S$
  (\isn{gfA}); a true instance of a positive part reads true at $S$ and
  a false instance of a negative part reads false at $S$ (\isn{gfPM},
  \isn{gfNM}); and the same two for a body under a $\lambda$ applied to
  a numeral (\isn{gfRPM}, \isn{gfRNM}).  The two new cases, the oracle clauses, close because a
  $\nu$-head is none of the heads the motive inspects: it is not the
  head of $\mu\Phi_0$ (\isn{omhO\_Ostar\_det}, \isn{omhO\_trmGfix}), not
  $\cA$ and not a $\lambda$-head.
\item \emph{The coinduction step} (\isn{nu\_step\_mjt}).  Invert the step
  to obtain $\Phi, x, c, i$; by the outer induction the three premises
  are valid in $\Mstar$.  Fix a well-formed $\sigma$ satisfying $H$, and
  write $\Phi_\sigma$ and $c_\sigma$ for the body and the certificate
  with the substitution pushed under two and one binders, $P = \lambda
  c_\sigma$ and $x_\sigma = \sigma\downarrow x$.  $\Phi_\sigma$ is
  \isn{posP}-positive (\isn{posP\_subf2}) and in the class
  (\isn{nuGood\_of\_posP}, \isn{free1\_subf2}).  The certificate set is
  \[
    W = \{u \mid u \text{ well-formed},\ u \text{ has a value in }
    \Mstar,\ c_\sigma[u] \text{ true in } \Mstar\}.
  \]
  \begin{enumerate}
  \item $x_\sigma \in W$: well-formed because $x$ is (the premise
    $\Gamma \ent x\,\mathbf N$ is well-formed) and $\sigma$ is
    (\isn{wsubf\_img\_wft}); valued by the validity of $\Gamma \ent
    x\,\mathbf N$ (\isn{mtrue\_NE}); and $c_\sigma[x_\sigma]$ true by the
    validity of $\Gamma \ent c[x]$ (\isn{bsub\_subf\_commute}).
  \item \emph{Region and closure obligation.}  The coinduction's region
    is $W$ and its obligation is $W \subseteq F^{\Theta^*}_{\Phi_\sigma}(W)$.
    Take $u \in W$ and $\tau = \sigma[i \mapsto u]$, a well-formed
    substitution; the eigenvariable is instantiated at the term $u$
    itself, not at a numeral.  By freshness, $\tau$ leaves $\Phi_\sigma$
    and $c_\sigma$ unchanged (\isn{subf\_ext\_free}).  $\tau$ satisfies
    the closure premise's context: $i\,\mathbf N$ because $u$ has a value
    (\isn{mtrue\_NI}); $c[i]$ because its instance is $c_\sigma[u]$; $H$
    by freshness.  The closure premise's validity at $\tau$ gives the
    truth of $\Phi_\sigma[P/a][u/x]$ (\isn{bsub\_asub},
    \isn{subf\_gsub\_par}), the instance of $\Phi_\sigma$ at the
    environment $\rho$ that sends $0$ to $u$ and $1$ to $P$.
  \item \isn{omvO\_pos\_neg\_transport\_P}, an induction on
    \isn{posA}/\isn{negA} with the recursive index sent to an arbitrary
    term $Q$, turns that truth into $\mathrm{rd}^+_{\Theta^*}(\{v \mid
    P\cdot v \text{ true}\}, 1, \rho, \Phi_\sigma, 1)$.
  \item \isn{osemO\_restrict}, an induction on the reading, shrinks the
    set to $W$.  Its hypotheses: $\Phi_\sigma$ is \isn{posP}-positive or
    negative, and $\rho$ is well-formed and valued below index $1$
    (\isn{envV}; here only the point $u$).  At an atom the argument is a
    numeral expression in the variables below $1$, so its instance is
    well-formed and valued (\isn{pnexp\_envV}), and the truth of
    $P\cdot v$ is the truth of $c_\sigma[v]$ (\isn{omvO\_beta}).  The
    restriction is the one step that needs the narrower grammar
    \isn{posP} rather than \isn{posA}; elsewhere \isn{posP} serves only
    to place $\Phi_\sigma$ in the class.
  \item \isn{osemO\_env\_off} replaces $\rho$ by $\mathrm{env}_{\Phi_\sigma}(u)$:
    the two differ only at index $1$, which the reading never consults.
    So $u \in F^{\Theta^*}_{\Phi_\sigma}(W)$.
  \item \emph{Close.}  \isn{gfp\_upperbound} gives $W \subseteq
    \gfp\,F^{\Theta^*}_{\Phi_\sigma} = \Theta^*(\Phi_\sigma, 1)$, so
    $x_\sigma \in \Theta^*(\Phi_\sigma, 1)$, and \isn{omvO\_nu\_atomI}
    gives the truth of $\nu\Phi_\sigma\cdot x_\sigma$, which is
    $\sigma\downarrow(\nu\Phi\cdot x)$ (\isn{subf\_trmNu\_gen}).
  \end{enumerate}
  No value congruence enters the coinduction step: the certificate set
  is indexed by terms.  Ground induction, by contrast, passes through
  a numeral and needs \isn{M\_vc\_arg} and \isn{M\_vcb}.
\end{enumerate}

\begin{thm}[\isn{inga\_conservative}]
\label{thm:app-inga-cons}
If $\Gamma \ent c$ is derivable in \inga, $\Gamma$ is well-formed, every
member of $\Gamma$ is true in the least model, and $c$ has some value
$r$ in the least model, then $c$ is true in the least model.
\end{thm}

\paragraph{Proof.}
\begin{enumerate}
\item \isn{omv\_le\_omvO\_Ostar}: every value in the least model is a value
  in $\Mstar$.  The inclusion is half of \isn{omv\_omvO\_Ostar}, a simultaneous
  rule induction on the least model, with motive: if $\mathrm{hd}(t, w)$,
  then $\mathrm{hd}_{\Theta^*}(t, w)$, or $w$ is the $\lambda$-head of
  $\Phi[\nu\Phi/a]$ for a $\Phi$ of the class and
  $\mathrm{hd}_{\Theta^*}(t, \nu_\Phi)$; if $\mathrm{val}(t, r)$, then
  $\mathrm{val}_{\Theta^*}(t, r)$.  At the head-level beta case, if $t\cdot
  u$ is a $\nu$-term of the class, the least model's head is the
  unfolding's $\lambda$-head (\isn{omh\_trmNu}, \isn{omh\_det}) and
  $\Mstar$ gives $\nu_\Phi$ (\isn{omhO\_nu}); when the function part
  already sat at a $\nu$-head, the body instance has no head
  (\isn{omhO\_body\_none}).  The value-level beta case at a $\nu$-head
  closes by \isn{nu\_head\_val}.  The other cases use the matching
  clauses.
\item \Cref{thm:app-inga-sound}, with the hypotheses lifted by step~1,
  makes $c$ true in $\Mstar$.
\item Step~1 also gives $\mathrm{val}_{\Theta^*}(c, r)$, and determinism
  of $\Mstar$ forces $r = 1$.
\end{enumerate}

\subsection{The all-zeros atom and the clash}

The body is $\Phi_Z = ((\lambda y.\,0)\cdot x = 0) \wedge a\cdot(x+1)$
(\isn{zPhi}), with $\wedge$ the abbreviation $\neg(\neg p \vee \neg q)$,
and $Z = \nu\Phi_Z$ (\isn{trmZ}).  $\Phi_Z$ is in the class
(\isn{nuGood\_zPhi}) and \isn{posP}-positive (\isn{posP\_zPhi}).

\begin{thm}[\isn{nuZ\_true}, \isn{nuZ\_true\_valueless}]
\label{thm:app-nuZ}
For every term $x$, $Z\cdot x$ is true in $\Mstar$.  Moreover $\bar 0
\in \Theta^*(\Phi_Z, 1)$, $Z\cdot 0$ is true in $\Mstar$, and for every
$n$ and $r$, $Z\cdot\bar n$ does not have the value $r$ in the least
model.
\end{thm}

\begin{thm}[\isn{inga\_Z}, \isn{inga\_omega\_unsound}]
\label{thm:app-ingaZ}
$\emptyset \ent Z\cdot 0$ is derivable in \inga.  Hence some $c$ has
$\emptyset \ent c$ derivable in \inga and is not true in the least
model.
\end{thm}

The valuelessness of the witness is the third clause of
\cref{thm:app-nuZ}; \isn{inga\_omega\_unsound} itself states ``not
true''.

\paragraph{Proof.}
\begin{enumerate}
\item \isn{FO\_zPhi\_UNIV}: at every oracle, $F^\Theta_{\Phi_Z}$ sends the
  set of all terms to itself.  The first conjunct is true by beta at
  every oracle (\isn{omvO\_s\_eq}), and every atom reads true because
  every argument lies in the set.  By \isn{gfp\_upperbound} every term
  lies in $\gfp\,F^{\Theta^*}_{\Phi_Z} = \Theta^*(\Phi_Z, 1)$, and
  \isn{omvO\_nu\_atomI} gives \isn{nuZ\_true}.  The membership of $\bar 0$
  is the inversion \isn{omvO\_nu\_atom\_inv}.
\item \emph{Valuelessness} (\isn{omv\_Z\_valueless\_gen}): a rule induction
  on the least model with a negative motive: no term with a value is
  one of the four shapes $Z\cdot\bar n$, $\neg(\neg e_n \vee \neg(Z\cdot
  \overline{n+1}))$, $\neg e_n \vee \neg(Z\cdot\overline{n+1})$ and
  $\neg(Z\cdot\overline{n+1})$, with $e_n$ the equation $(\lambda y.\,0)
  \cdot\bar n = 0$.  The four shapes are the unfolding of $Z\cdot\bar n$
  (\isn{gsub\_zPhi}) and its subterms.  At the beta case for $Z\cdot\bar
  n$ the head of $Z$ is unique (\isn{omh\_det}, \isn{omh\_trmNu}), so the
  body is the unfolding, which the hypothesis excludes.  The recursor
  and universal cases cannot match a $\lambda$-head.  The negation cases
  pass to the inner shape: the outer negation to the disjunction at the
  same $n$, the negated atom to $Z\cdot\overline{n+1}$.  A true disjunct $\neg e_n$ would
  make $e_n$ false, against \isn{omv\_s\_eq} and uniqueness; the other
  disjunction cases fall to the hypothesis.  In effect the induction is
  a well-founded descent: a value at $n$ needs a value at $n+1$.
\item \isn{inga\_Z}: one application of the rule (\isn{ir\_nu}) with
  certificate $c = (x = x)$.  The premises are \rga derivations
  (\isn{z\_prem1} to \isn{z\_prem3}): $\emptyset \ent 0 = 0$; in the context
  $i\,\mathbf N,\ i = i$, the body with $\lambda c$ for $a$, assembled
  from $(\lambda y.\,0)\cdot i = 0$ by conversion, $i+1 = i+1$ by
  reflexivity turned into $(\lambda c)\cdot(i+1)$ by conversion, and two
  double-negation introductions and a disjunction rule for the
  conjunction; and $\emptyset \ent 0\,\mathbf N$.
\item \isn{inga\_omega\_unsound}: step~2 at $n = 0$, with step~3.
\end{enumerate}

\begin{thm}[\isn{nu\_gfix\_clash}]
\label{thm:app-clash}
Let $S$ be a step relation on proof lists, contexts and conclusions
such that every \rga step is an $S$ step, $S$ is monotone in the proof
list (if the proof list $P$ is a subset of $P'$, as sets, and $S(P, \Gamma, c)$, then $S(P', \Gamma, c)$), every coinduction step is an $S$ step, and every
step of ground induction at $\nu$-terms (\isn{gfixNuS}) is an $S$ step.
Then $\emptyset \ent 0 = 1$ is derivable for $S$.
\end{thm}

\isn{gfixNuS} is the clause of \isn{pr\_gfix} with $\nu\Phi$ in place of
$\mu\Phi$: side conditions $\Phi$ \isn{posP}-positive, $\Phi$ and $c$
well-formed, and the three freshness conditions; premises $\Gamma \ent
\nu\Phi\cdot x$, $i\,\mathbf N,\ \Phi[\lambda c/a][i/x],\ \Gamma \ent
c[i]$ and $\Gamma \ent x\,\mathbf N$; conclusion $c[x]$.  The clause does
not ask $a$ to occur in $\Phi$.

\paragraph{Proof.}
\begin{enumerate}
\item \rga entailment embeds into $S$-entailment (\isn{entails\_step\_embed}).
\item $\emptyset \ent Z\cdot 0$ for $S$: merge the three embedded premises
  of step~3 of the preceding proof into one proof list
  (\isn{entails\_merge3}, using monotonicity), then one coinduction step.
\item \isn{clash\_prem2}, in \rga: the closure premise of ground induction
  at the point-free motive $0 = 1$.  Its hypothesis is the body with
  $\lambda(0 = 1)$ for $a$, $\neg(\neg e \vee \neg((\lambda(0=1))\cdot
  (i+1)))$.  A disjunction rule gives $\neg\neg((\lambda(0=1))\cdot(i+1))$,
  double-negation elimination gives the redex, and one conversion
  from the redex to its reduct (\isn{ir\_convE} at a beta step) gives
  $0 = 1$.
\item Merge step~2, step~3 and $\emptyset \ent 0\,\mathbf N$, and apply one
  \isn{gfixNuS} step; the conclusion $(0 = 1)[0]$ is $0 = 1$.
\end{enumerate}
The proof uses no $\omega$-model, neither the least model nor $\Mstar$: a build-time check confirms that no name
among its transitive dependencies contains one of the strings
\isn{omv}, \isn{omh}, \isn{Ostar}, \isn{nu\_model}, \isn{mjt} and
\isn{mtrue}, the names of the two models' relations and of their
soundness facts.

\subsection{A body read at both fixed points}

These three results concern the least model and the operators
$F_\Phi$ (\isn{gfixOp}) and $G_\Phi$ (\isn{gfixFOp}) of
\cref{thm:disjoint}, on sets of numbers; no oracle is involved.  The
readings consult an atom at a literal numeral: $\mathrm{rd}^+(S, \ldots)$
reads $a\cdot u$, with $a$ not occurring in $u$, true when
$\sigma\downarrow u = \bar m$ with $m \in S$ (\isn{os\_app}), and
$\mathrm{rd}^-(X, \ldots)$ reads such an atom false when
$\sigma\downarrow u = \bar m$ with $m \in X$ (\isn{fs\_app}).  Both
operators are monotone for every body (\isn{gfixOp\_mono},
\isn{gfixFOp\_mono}).

\begin{thm}[\isn{gfix\_omv\_false\_eq\_lfp}, \isn{gfp\_lfp\_disjoint}]
\label{thm:app-regions}
If $\Phi$ is \isn{pos}-positive at index $1$ and $\lambda a.\lambda
x.\Phi$ is closed, then $\{n \mid \mu\Phi\cdot\bar n$ is false in the
least model$\} = \lfp\,G_\Phi$.  For every term $\Phi$, $\gfp\,F_\Phi \cap
\lfp\,G_\Phi = \emptyset$.
\end{thm}

\begin{thm}[\isn{gfp\_eq\_compl\_iff}]
\label{thm:app-compl}
For every term $\Phi$, $\gfp\,F_\Phi = \mathbb N \setminus \lfp\,G_\Phi$
if and only if every $n$ lies in $F_\Phi(\mathbb N \setminus
\lfp\,G_\Phi)$ or in $\lfp\,G_\Phi$.
\end{thm}

\paragraph{Proof.}
\begin{enumerate}
\item \emph{Falsity set, one inclusion} (\isn{lfp\_gfix\_omv\_false}).  By
  \isn{lfp\_lowerbound} it suffices that the falsity set $\{n \mid
  \mu\Phi\cdot\bar n$ false$\}$ is closed under $G_\Phi$
  (\isn{fsem\_gfix\_omv}): \isn{fsem\_term}, at the predicate term
  $\mu\Phi$, gives the falsity of $\Phi[\mu\Phi/a][\bar n/x]$, and the
  value-level beta clause with \isn{omh\_trmGfix} folds that falsity back to
  $\mu\Phi\cdot\bar n$.
\item \emph{The other inclusion} (\isn{gfix\_omv\_false\_lfp}) applies
  \isn{gfix\_omv\_false\_ind\_set} with $X = \lfp\,G_\Phi$, whose closure
  hypothesis is the fixed-point equation (\isn{gfixFOp\_lfp}).
  \isn{gfix\_omv\_false\_ind\_set} is the falsity half of a transport
  lemma (\isn{gfix\_false\_transport}), a rule induction on the least
  model with a combined motive: false instances of $\mu\Phi$ lie in $X$,
  positive bodies read false, negative bodies read true, and the two
  quantifier-instance redexes likewise.  Positivity and closedness are
  used there.
\item \emph{Disjointness.}  Write $N = \gfp\,F_\Phi$.  By
  \isn{lfp\_lowerbound} it suffices that $\mathbb N \setminus N$ is
  closed under $G_\Phi$.  An $n$ in $G_\Phi(\mathbb N \setminus N) \cap N$
  would satisfy, by \isn{gfp\_fixpoint}, both the truth reading at $N$
  and the falsity reading at $\mathbb N \setminus N$.
  \isn{osem\_fsem\_excl}, an induction on the truth reading at disjoint
  sets, refutes the pair: at a leaf the falsity reading gives the
  opposite value (\isn{fsem\_nf\_only}), against uniqueness; at an atom
  the numeral would lie in both sets (\isn{fsem\_app\_F},
  \isn{num\_trm\_eq\_iff}); at a connective the inversions of the
  falsity reading apply.  No positivity is used.
\item \emph{Complement, left to right.}  With $L = \lfp\,G_\Phi$, the
  equation turns \isn{gfp\_fixpoint} into $F_\Phi(\mathbb N \setminus
  L) = \mathbb N \setminus L$, so every $n$ is in $F_\Phi(\mathbb N
  \setminus L)$ or in $L$.
\item \emph{Right to left.}  The hypothesis says $\mathbb N \setminus L
  \subseteq F_\Phi(\mathbb N \setminus L)$, a post-fixed point, so
  \isn{gfp\_upperbound} gives $\mathbb N \setminus L \subseteq
  \gfp\,F_\Phi$, and step~3 gives the converse inclusion.
\end{enumerate}
The oracle form of step~3, \isn{gfpO\_lfpO\_disjoint}, is step~3(b) of
the proof of \cref{thm:app-rga-nu-sound}; the oracle form needs a consistent oracle,
because its leaf case uses determinism of the oracle model.

\paragraph{Head steps never turn a $\nu$-atom into a $\mu$-atom.}
\isn{nu\_hsteps\_not\_gfix} (\isn{OIGA/GNu.thy}): for a body shaped as a formula (\isn{fshape}; every body of the $\nu$ class is, \isn{nuGood\_fshape}), no term reachable from $\nu\Phi\cdot u$ by head steps has the form $\mu\Psi\cdot v$.  The orbit of the atom (\isn{nuOrbit}) has five shapes, the atom, the three unfolding stages and the body instance, and is closed under head steps by inversion (\isn{hstep\_trmNu\_inv}, \isn{hstep\_app\_spine\_inv}, \isn{hstep\_beta\_inv}); in the instance case the formula shape leaves either a connective, which takes no step, or $\nu\Phi\cdot e'$, back in the orbit; induction on the reflexive transitive closure and a shape check finish it.  The shape condition is needed: the body that is itself a closed $\mu$-atom applied to the point reaches $\mu\Psi\cdot u$ in four steps (\isn{nu\_hsteps\_gfix\_unshaped}).

\paragraph{Lemmas a re-mechanization needs.}
\begin{itemize}
\item \isn{nu\_step1}--\isn{nu\_step4}, \isn{trmNu\_neq\_trmGfix}: the
  $\nu$-former unfolds in four head steps and differs from every $\mu$-atom.
\item \isn{omhO\_trmNu\_inv}, \isn{omvO\_nu\_atom\_inv}, \isn{omvO\_nu\_atomI}:
  a $\nu$-atom of the class takes values through the oracle alone.
\item \isn{omhO\_det\_excl}, \isn{osemO\_fsemO\_excl}, \isn{ocons\_Ostar}:
  determinism at a consistent oracle; $\Theta^*$ is consistent.
\item \isn{osemO\_term}, \isn{omvO\_pos\_neg\_transport}, \isn{nu\_coherent\_T},
  \isn{nu\_coherent\_F}, \isn{nu\_head\_val}: a $\nu$-atom has the value of
  its unfolding.
\item \isn{omO\_hstep\_sim}, \isn{omO\_hstep\_rev}: head steps preserve
  values in both directions, heads up to \isn{heqN}.
\item \isn{omO\_vcr\_transport}, \isn{Ocg\_gfp\_cong}, \isn{Ocg\_prefix},
  \isn{omvO\_vcr\_ex}: the value congruence of $\Mstar$, through the
  congruence interior.
\item \isn{nu\_model\_iff\_vc}, \isn{step\_prf\_mjt}, \isn{step\_M\_sound}:
  the interface reduced to two facts; soundness from step soundness.
\item \isn{gfix\_transportO\_par}, \isn{gfix\_step\_mjt}: ground induction
  is sound for $\Mstar$.
\item \isn{omvO\_pos\_neg\_transport\_P}, \isn{osemO\_restrict},
  \isn{osemO\_env\_off}, \isn{nu\_step\_mjt}: the coinduction step.
\item \isn{omv\_omvO\_Ostar}: the least model lies below $\Mstar$.
\item \isn{FO\_zPhi\_UNIV}, \isn{omv\_Z\_valueless\_gen}, \isn{clash\_prem2}:
  the all-zeros atom and the clash.
\item \isn{fsem\_gfix\_omv}, \isn{gfix\_omv\_false\_ind\_set},
  \isn{osem\_fsem\_excl}: the least model's falsity set; the two readings
  exclude each other at disjoint sets.
\end{itemize}

\section{The dichotomy}
\label{app:dich}

Files: \isn{IGA/IGA\_Dich.thy} (the dual, the premise, the swap identity,
the characterization of validity (\cref{thm:C}), the rule and its step
soundness, accessibility never false),
\isn{OIGA/GDich.thy} (the calculi with the dichotomy, the system
\isn{oinga\_dich}, whose entailment equals \dinga's, the price theorem,
the accessibility instances), \isn{OIGA/GDich2.thy} (exclusion,
the decided $\nu$-atom, the premise of the $\mathcal O$ body),
\isn{OIGA/GDcmp.thy} (totality of the decompiler, $\Dec(\mathcal O)$),
\isn{IGA\_Closing/CO\_Links.thy} (ATI from the dichotomy),
\isn{IGA\_Closing/CO\_Check.thy} (the proof checker of
\isn{oinga\_dich}) and \isn{IGA\_Closing/CO\_Dinga.thy} (\dinga, the
equality \isn{dinga\_eq\_oingad}, the re-exports).  The
results about $\Mstar$ are classical theorems in the locale \isn{rga};
the derivations are in the abstract calculi named with each theorem.

\paragraph{Definitions.}
The \emph{dual} (\isn{ndual}) is an inductive relation $\Psi \in
\partial_k(t)$ with five clauses: a term in which index $k$ does not
occur is its own dual; the dual of $v_k\cdot u$, with $u$ a numeral
expression in the variables below $k$ (\isn{pnexp}), is $\neg(v_k\cdot
u)$; and the dual commutes with disjunction, negation and the universal
(one index deeper).  The relation is functional (\isn{ndual\_fun}); a dual
exists for every \isn{posP}-positive or negative term
(\isn{ndual\_exists}); the dual of a positive term is negative and
conversely (\isn{ndual\_polar}), so if $\Psi$ is the dual of a positive
$\Phi$ then $\neg\Psi$ is positive (\isn{ndual\_posP}) and, when index
$1$ occurs in $\Phi$, in the $\nu$ class (\isn{ndual\_nuGood}); and the
dual commutes with substitution for parameters (\isn{ndual\_subf2}).
For a body $\Phi$ with dual $\Psi$, the $\nu$-body is $\neg\Psi$ and the
two atoms of the rule are $\mu\Phi$ and $\nu(\neg\Psi)$.

The \emph{decidedness premise} (\isn{nfdec}) is an inductive relation
$D \in \mathrm{prem}_k(t)$ with five clauses.  If $k$ does not occur in
$t$, then $D$ is the universal closure, over the $k$ indices below $k$,
of $t'\,\mathbf B$, where $t'$ is $t$ with the absent index removed
(\isn{lowr}; the parameters stay free).  For a recursive atom with a
numeral-expression argument, $D$ is $0 = 0$.  For a disjunction, $D$ is
the conjunction of a premise of each disjunct; for a negation, a premise
of its argument; for a universal, a premise of its body one index
deeper.  The relation is not functional: at a compound part in which
$k$ does not occur, the first clause and a structural clause both
apply, and the rule accepts any premise so obtained.  A premise exists
for every positive body (\isn{nfdec\_exists}).  For the accessibility
body one premise is $\Dec\,R \wedge (0 = 0)$ (\isn{nfdec\_accPhi}).

The rule (\isn{dichS}) has the side conditions $\Phi$
\isn{posP}-positive at $1$, index $1$ free in $\Phi$, $\Psi$ the dual of
$\Phi$, $D$ a premise of $\Phi$, $\Phi$ well-formed; its premises are
$\Gamma \ent D$ and $\Gamma \ent x\,\mathbf N$; its conclusion is $\Gamma
\ent \mu\Phi\cdot x \vee \nu(\neg\Psi)\cdot x$.  There is no eigenvariable,
and parameters are allowed.  The rule is \emph{valid} at $\Phi$ and
$\Psi$ (\isn{dichValid}) when, for every well-formed $\Gamma$ and
well-formed $x$ such that $\Gamma \ent x\,\mathbf N$ is valid in
$\Mstar$, $\Gamma \ent \mu\Phi\cdot x \vee \nu(\neg\Psi)\cdot x$ is valid
in $\Mstar$.  The premise $D$ does not occur in validity.

Three operators on sets of numbers, relative to $\Mstar$, carry the
proofs.  For a $\nu$-body $T$, $\hat F_T(S) = \{n \mid \bar n \in
F^{\Theta^*}_T(\bar S)\}$ and $\hat G_T(S) = \{n \mid \bar n \in
G^{\Theta^*}_T(\bar S)\}$ (\isn{FN}, \isn{GN}), where $\bar S$ is the set
of numerals of members of $S$.  For a $\mu$-body $\Phi$, $F^*_\Phi(S) =
\{n \mid \mathrm{rd}^*(S, 1, \sigma_n, \Phi, 1)\}$ (\isn{gfixOpM}), where
$\sigma_n$ sends $0$ to $\bar n$, $1$ to $\mu\Phi$ and $v+2$ to $v$, so
that it differs from $\mathrm{env}_{\neg\Psi}(\bar n)$ only at index $1$,
and the reading
$\mathrm{rd}^*$ takes non-recursive leaves from $\Mstar$ and reads an atom
true when its argument's instance is a numeral $\bar m$ with $m \in S$.
The three operators are those of \cref{thm:disjoint}, the truth and
falsity operators of a body, with the least model replaced by $\Mstar$.

\subsection{The swap identity}

\begin{thm}[\isn{swap\_identity}]
\label{thm:app-swap}
Let $\Psi$ be the dual of $\Phi$ at index $1$, with $\Phi$ and
$\neg\Psi$ \isn{posP}-positive at $1$ and $\neg\Psi$ in the $\nu$ class.
Then for every $n$: $\mu\Phi\cdot\bar n$ is true in $\Mstar$ if and only
if $\nu(\neg\Psi)\cdot\bar n$ is false in $\Mstar$.
\end{thm}

The swap identity is \cref{thm:swap}; the hypothesis that index $1$
occurs in $\Phi$ is not needed, since membership of $\neg\Psi$ in the
class already implies that index $1$ occurs (\isn{ndual\_free}).

\paragraph{Proof.}
\begin{enumerate}
\item \emph{The $\mu$-side is a least fixed point in $\Mstar$}
  (\isn{gfix\_omvO\_eq\_lfp}): for $\Phi$ \isn{posP}-positive, $\{n \mid
  \mu\Phi\cdot\bar n$ true in $\Mstar\} = \lfp\,F^*_\Phi$.  One inclusion
  (\isn{lfp\_gfix\_omvO}) is \isn{lfp\_lowerbound}: the truth set is closed
  under $F^*_\Phi$, because \isn{osemM\_term} turns a reading at the truth
  set into the truth of $\Phi[\mu\Phi/a][\bar n/x]$, and one backward
  beta step at the head of $\mu\Phi$ (\isn{omhO\_trmGfix}) gives
  $\mu\Phi\cdot\bar n$.  The other inclusion (\isn{gfix\_omvO\_lfp}) is the
  ground-induction property of $\Mstar$, \isn{gfix\_omvO\_ind\_set\_par},
  at $S = \lfp\,F^*_\Phi$, whose closure hypothesis is the fixed-point
  equation (\isn{gfixOpM\_lfp}).  The ground-induction property is the
  one place positivity of
  $\Phi$ enters.
\item \emph{The swap lemma} (\isn{ndual\_swap}): if $t'$ is the dual of
  $t$ at $k$, then for every $\sigma$, $b$ and $S$,
  $\mathrm{rd}^-_{\Theta^*}(\bar S, k, \sigma, t', b)$ iff
  $\mathrm{rd}^*(S, k, \sigma, t, b)$.  Induction on the dual, with
  $\sigma$ and $b$ arbitrary.  At a part without index $k$ both sides
  are the value of $\sigma\downarrow t$ in $\Mstar$.  At an atom with $b$
  true: the falsity reading of $\neg(v_k\cdot u)$ is true iff the atom
  reads false, iff $\sigma\downarrow u \in \bar S$, iff the $\mathrm{rd}^*$
  atom is true (\isn{fsemO\_app\_F}, \isn{osemM\_app\_T}).  At an atom with
  $b$ false neither side holds, since the falsity reading has no true
  atoms and $\mathrm{rd}^*$ no false ones (\isn{fsemO\_app\_T\_none},
  \isn{osemM\_app\_F\_none}).  Disjunction, negation and the universal
  follow from the inversion lemmas of both readings and the induction
  hypotheses.  The lemma needs no positivity and no condition on $\sigma$.
\item \emph{Numeral restriction} (\isn{nrestrict\_T}, \isn{nrestrict\_F}):
  for $T$ in the class and \isn{posP}-positive, the numerals in
  $\Theta^*(T, 1)$ are $\gfp\,\hat F_T$ and those in $\Theta^*(T, 0)$ are
  $\lfp\,\hat G_T$.  The engine is \isn{osemO\_restrict\_num} and
  \isn{fsemO\_restrict\_num}: at an environment that is numeral below the
  recursive index (\isn{envN}), a reading of a term of the positive or
  negative grammar at a set $S_1$ transfers to any $S_2$ containing the
  numerals of $S_1$, because every atom argument is then a numeral
  expression and evaluates to a numeral (\isn{pnexp\_envN}).  For falsity, forward:
  the set of terms $u$ such that $u = \bar n$ implies $n \in
  \lfp\,\hat G_T$ is closed under $G^{\Theta^*}_T$ (restrict, then use the
  fixed point), so $\lfp\,G^{\Theta^*}_T = \Theta^*(T, 0)$ lies inside
  that set; backward:
  $\{n \mid \bar n \in \lfp\,G^{\Theta^*}_T\}$ is closed under $\hat G_T$
  (\isn{GO\_mono2}).  Truth is the analog with \isn{gfp\_upperbound}.
\item \emph{The operators coincide} (\isn{GN\_dual}): $\hat G_{\neg\Psi} =
  F^*_\Phi$.  At a point $n$: $n \in \hat G_{\neg\Psi}(S)$ iff $\neg\Psi$
  reads false at $\mathrm{env}_{\neg\Psi}(\bar n)$, iff $\Psi$ reads true
  there (\isn{fsemO\_not\_F\_any}), iff $\Phi$ reads true in $\mathrm{rd}^*$
  there (step~2), iff $\Phi$ reads true at $\sigma_n$ (\isn{osemM\_env\_off}:
  the two environments differ only at index $1$, which the reading never
  consults), iff $n \in F^*_\Phi(S)$.
\item \emph{Assembly.}  By step~4 the two least fixed points agree.  By
  step~1 and step~3, $\mu\Phi\cdot\bar n$ is true iff $\bar n \in
  \Theta^*(\neg\Psi, 0)$, iff $\nu(\neg\Psi)\cdot\bar n$ is false
  (\isn{omvO\_nu\_atomI}, \isn{omvO\_nu\_atom\_inv}).
\end{enumerate}

\subsection{The characterization of validity and the step soundness of the rule}

\begin{thm}[\isn{thmC}]
\label{thm:app-C}
Let $\Phi$ be \isn{posP}-positive at $1$, with index $1$ free in $\Phi$,
and $\Psi$ its dual.  The rule is valid at $\Phi$ and $\Psi$ if and only
if, for every well-formed substitution $\sigma$ and every $n$, the atom
$\nu(\neg\Psi)_\sigma\cdot\bar n$ has a value in $\Mstar$, where
$(\neg\Psi)_\sigma$ is $\neg\Psi$ with $\sigma$ applied to its parameters
(the substitution pushed under two binders, \isn{subf\_lift}).
\end{thm}

The substitutions range over all well-formed ones, closing or not; the
identity substitution is one of them.

\begin{thm}[\isn{dich\_step\_mjt}]
\label{thm:app-dichstep}
Suppose the conclusion $H \ent c$ comes from one dichotomy step whose
premises are judgments of the proof list $P$; every judgment $H' \ent
c'$ in $P$ with $H'$ well-formed is valid in $\Mstar$; every judgment
$H' \ent c'$ in $P$ has $c'$ well-formed; and $H$ is well-formed.  Then $H \ent c$ is valid in
$\Mstar$.
\end{thm}

Together the characterization and the step soundness are \cref{thm:C}.  The step soundness assumes nothing
about the body: the premise $D$ does the work.

\paragraph{Shared lemmas.}
\begin{enumerate}
\item \emph{Instances} (\isn{dich\_inst}): at every $\sigma$, $\Phi_\sigma$
  is \isn{posP}-positive (\isn{posP\_subf2}), $\Psi_\sigma$ is its dual
  (\isn{ndual\_subf2}), $(\neg\Psi)_\sigma = \neg(\Psi_\sigma)$, and
  $\neg\Psi_\sigma$ is \isn{posP}-positive and in the class
  (\isn{ndual\_posP}, \isn{nuGood\_of\_posP}, \isn{free1\_subf2},
  \isn{ndual\_free}).
\item \emph{Swap at instances} (\isn{dich\_swap\_inst}): by
  \cref{thm:app-swap} at the instance pair, for every $\sigma$ and $n$,
  $\mu\Phi_\sigma\cdot\bar n$ is true iff $\nu(\neg\Psi)_\sigma\cdot\bar n$
  is false.
\item \emph{Substitution into the conclusion} (\isn{dich\_concl\_inst}):
  $\sigma$ applied to $\mu\Phi\cdot x \vee \nu(\neg\Psi)\cdot x$ is
  $\mu\Phi_\sigma\cdot(\sigma\downarrow x) \vee
  \nu(\neg\Psi)_\sigma\cdot(\sigma\downarrow x)$ (\isn{subf\_trmGfix},
  \isn{subf\_trmNu\_gen}).
\item \emph{The disjunction at a valued term} (\isn{disj\_at\_term}).
  Suppose that $\nu T\cdot\bar n$ has a value for every $n$, that
  $\mu\Phi\cdot\bar n$ is true iff $\nu T\cdot\bar n$ is false for every
  $n$, and that $y$ is well-formed with value $v$.  Then $\mu\Phi\cdot y
  \vee \nu T\cdot y$ is true in $\Mstar$.  The numeral $\bar v$ and $y$
  are congruent (\isn{ovcrO\_swap}).  If $\nu T\cdot\bar v$ is true, the
  value congruence (\isn{ovcrO\_appI}, \isn{omvO\_vcr\_ex}) moves the
  truth to $\nu T\cdot y$; if false, the swap makes $\mu\Phi\cdot\bar v$
  true and the congruence moves that truth to $\mu\Phi\cdot y$.  A
  disjunction clause closes either case.
\end{enumerate}

\paragraph{Proof of \cref{thm:app-C}, left to right.}
\begin{enumerate}
\item Assume validity; fix a well-formed $\sigma$ and $n$.  $\emptyset \ent
  \bar n\,\mathbf N$ is valid, since $\bar n$ has a value.
\item Validity at $\Gamma = \emptyset$ and $x = \bar n$ makes $\emptyset
  \ent \mu\Phi\cdot\bar n \vee \nu(\neg\Psi)\cdot\bar n$ valid, so its
  instance at $\sigma$ is true; by shared lemma~3, and since $\sigma$
  leaves the closed numeral alone, that instance is $\mu\Phi_\sigma\cdot
  \bar n \vee \nu(\neg\Psi)_\sigma\cdot\bar n$.
\item Invert the disjunction (\isn{omvO\_or\_inv}).  If the $\mu$-disjunct
  is true, shared lemma~2 makes the $\nu$-atom false; otherwise the
  $\nu$-atom is true.  Either way the $\nu$-atom has a value.  The
  swap at every instance needs \isn{ndual\_subf2}.
\end{enumerate}

\paragraph{Proof of \cref{thm:app-C}, right to left.}
Fix a well-formed $\Gamma$, a well-formed $x$ with $\Gamma \ent x\,\mathbf N$
valid, and a well-formed $\sigma$ satisfying $\Gamma$.  Then
$\sigma\downarrow x$ has a value (\isn{mtrue\_NE}) and is well-formed
(\isn{wsubf\_img\_wft}); shared lemma~4 with the assumption at $\sigma$,
shared lemma~2 and shared lemma~3 gives the truth of the conclusion's
instance.  The premise $D$ is not involved.

\paragraph{Proof of \cref{thm:app-dichstep}.}
The step replaces the assumption of the right-to-left direction by a
consequence of the premise, bivalence of the $\nu$-atom at every
instance that satisfies the context.
\begin{enumerate}
\item \emph{The leaf condition} (\isn{leafdec}): an inductive relation
  over the grammar, relative to a substitution $\tau$, whose leaf clause
  says that index $k$ does not occur in $t$ and that, for every $\nu :
  \mathbb N \to \mathbb N$, $t$ has a value in $\Mstar$ at the
  environment sending $v < k$ to $\overline{\nu(v)}$, $k$ to $v_0$ and $v >
  k$ to $\tau(v - k - 1)$ (\isn{penv}).
\item \emph{The premise gives the leaf condition} (\isn{nfdec\_leafdec}):
  if $D$ is a premise of $t$ at $k$ and $\tau\downarrow D$ is true, then
  $t$ satisfies the leaf condition at $\tau$.  Induction on the premise
  relation.  At the leaf clause the truth of the iterated universal
  (\isn{allqN\_true\_iff}, an induction on $k$) gives the truth of
  $t'\,\mathbf B$ at every numeral assignment of the quantified indices;
  that environment agrees with \isn{penv} on the free variables of $t$
  (\isn{envq\_lowr}); and $s\,\mathbf B$ is true iff $s$ has a value
  (\isn{omvO\_typeB\_true}).  At an atom the condition is immediate.  A
  true conjunction gives both conjuncts (\isn{omvO\_and\_true}); negation
  and the universal pass the hypothesis through.
\item \emph{The dual keeps the leaf condition} (\isn{ndual\_leafdec}), by
  induction on the dual.  For a compound part, split on whether $k$
  occurs: if not, the dual is the part itself (\isn{ndual\_nfree\_id});
  if so, only the structural clause of the leaf condition applies, and
  the induction hypotheses finish.  This split settles the overlap of
  the dual's first clause with its structural clauses.
\item \emph{From the leaf condition to determinacy}
  (\isn{leafdec\_dposE}).  The determinacy grammar at an environment
  $\rho$ (\isn{dposE}, \isn{dnegE}) admits a leaf when every environment
  that agrees with $\rho$ above $k$ and is numeral below $k$ gives the
  leaf a value (\isn{nfvalE}), admits recursive atoms with
  numeral-expression arguments as positive, and follows the grammar
  elsewhere, asking for a universal's body at every numeral extension
  of $\rho$.  If $t$ satisfies the leaf condition at $\sigma$ and $\rho$
  sends each index $v > k$ to $v_{v-k-1}$, the parameter that index
  names (\isn{idAbove}), a positive
  (negative) $t$ has its instance at $\sigma$ in the positive (negative)
  half.  Induction on the leaf condition with $\rho$ arbitrary; at a
  leaf, every admissible environment composed with the instance is a
  \isn{penv} instance of $t$ (\isn{inst\_leaf}).
\item \isn{dich\_premise\_dposE}: from the validity of $\Gamma \ent D$, a
  well-formed $\sigma$ satisfying $\Gamma$ makes $\sigma\downarrow D$
  true (\isn{mjtD}); steps~2 and~3 give the leaf condition for $\Phi$ and
  for $\neg\Psi$; step~4 at the environments $\mathrm{env}_T(u)$
  (\isn{idAbove\_nuenv}) puts the instance $T = (\neg\Psi)_\sigma$ in the
  positive half at every point.
\item \emph{Determinacy gives the union} (\isn{dposE\_det},
  \isn{dposE\_gdetM}, \isn{gdetM\_union}).  By mutual induction on the
  determinacy grammar, at an environment numeral below $k$: a positive
  part reads true at $\bar S$ or reads false at the numerals of the
  complement of $S$, and a negative part the reverse.  A leaf has a value,
  true or false; an atom's argument is a numeral $\bar m$
  (\isn{pnexp\_envN}), and classically $m \in S$ or not; the universal
  splits classically on whether every numeral instance reads true,
  otherwise the failing instance reads false.  At the environments of
  $T$ the determinacy says: for every $S$ and $n$, $n \in \hat F_T(S)$ or $n \in
  \hat G_T(\mathbb N \setminus S)$.  Let $L = \lfp\,\hat G_T$.  For $n
  \notin L$, the instance at $S = \mathbb N \setminus L$ gives $n \in
  \hat F_T(\mathbb N \setminus L)$, since $\hat G_T(L) \subseteq L$; so
  $\mathbb N \setminus L$ is a post-fixed point of $\hat F_T$, and
  \isn{gfp\_upperbound} gives $\gfp\,\hat F_T \cup \lfp\,\hat G_T =
  \mathbb N$.  Only this inclusion is used.
\item \emph{Bivalence} (\isn{thmC\_rhs\_forms}, \isn{dich\_premise\_bivalent}):
  by step~3 of the proof of \cref{thm:app-swap}, the union says that
  every numeral lies in $\Theta^*(T, 1) \cup \Theta^*(T, 0)$, so the
  $\nu$-atom of $T$ has a value at every numeral.
\item \emph{The step at a body} (\isn{dich\_valid\_at}): from the validity
  of $\Gamma \ent D$ and $\Gamma \ent x\,\mathbf N$ and a well-formed $x$,
  the right-to-left proof of \cref{thm:app-C} applies at each $\sigma$
  satisfying $\Gamma$, with step~7 in place of the assumption.
\item \emph{The theorem.}  Invert the step; well-formedness of the
  recorded judgment $\Gamma \ent x\,\mathbf N$ gives a well-formed $x$;
  the hypothesis on $P$ makes $\Gamma \ent D$ and $\Gamma \ent x\,\mathbf N$
  valid; step~8 concludes.
\end{enumerate}
The premise is a sufficient condition for bivalence: the rule is
valid at a body whenever the premise is valid in every well-formed
context (\isn{dich\_valid\_of\_premise}), and the step uses neither
\cref{thm:app-C} nor validity.  The controls are \isn{dichA\_valid}, the valid accessibility body at
the relation $R(y, x) \iff y = x+1$ (\isn{flipD}), and \isn{dichB\_invalid}, the invalid body $\bot
\vee a\cdot x$, whose premise is true at no substitution
(\isn{botD\_never\_true}).

\subsection{The system and its price}

\dinga is the locale \isn{dinga} (theory
\isn{IGA\_Closing/CO\_Dinga.thy}); its step relation
(\isn{dinga\_derives}) is \rga's step, a ground-induction step, a
coinduction step or a dichotomy step, with no ATI rule, and
\isn{dinga\_ent} is its entailment.  ATI is derived
(\isn{dinga\_ir\_allTI}), by \cref{thm:ati-dich} (\isn{g\_ati\_of\_dich},
in \isn{gcalc\_gfix\_dich}, with the sublocale edge
$\isn{gcalc\_gfix\_dich} \subseteq \isn{gcalc\_ati}$, in
\isn{IGA\_Closing/CO\_Links.thy}).  The system \isn{oinga\_dich}, whose
step relation (\isn{oingad\_derives}) is \oinga's or a dichotomy step,
so that ATI is primitive there, has the same entailment
(\isn{dinga\_eq\_oingad}; the inclusion \isn{oingad\_le\_dinga} runs
through \isn{g\_ati\_of\_dich}).  The system-level results below are
proved for \isn{oinga\_dich} and re-exported at \dinga through the
equality under \isn{dinga\_} names, and \dinga's own step decider is
primitive recursive without an ATI detector
(\isn{n2pr\_dinga\_derives}).  The results of the examples in
\isn{OIGA\_Coind} keep their \isn{oingad\_} names (that session is not
an ancestor of \isn{IGA\_Closing}); they are restated at \dinga by name in the junction theory
\isn{CD\_Junction} (session \isn{Coind\_Dinga}), which sees both
sessions: every result proved in an abstract calculus \isn{oingad\_X} as
\isn{dinga.dinga\_X}, with a build-time check that each restatement
agrees with its original, and the system-level facts the paper cites
(\isn{dinga\_stepmother\_not\_refuted}, \isn{dinga\_R1p\_not\_acc},
\isn{dinga\_aoTw\_not\_acc}, \isn{dinga\_trUp\_dich\_premise\_not\_derivable},
\isn{dinga\_trUp\_dich\_premise\_neg\_not\_derivable},
\isn{dinga\_tr\_dich\_prem\_open}, \isn{dinga\_tr\_dich\_unavailable}).
The same theory proves \oinga's derivability recursively enumerable
(\isn{oinga\_entails\_re}) and \isn{oinga\_not\_cert\_sound}.

The accessibility body
has the dual $\forall y.\,\neg R(y, x) \vee \neg a\cdot y$
(\isn{psiAcc}, \isn{ndual\_accPhi}), and $\Path_R = \nu(\neg\,\mathrm{psiAcc}\,R)$
(\isn{trmPath}, \isn{pathPhi}).

\begin{thm}[\isn{dinga\_dich\_price}, \isn{oinga\_dich\_price}]
\label{thm:app-price}
In \dinga:
(i) if $\Gamma \ent c$, $\Gamma$ is well-formed and every member of
$\Gamma$ is true in $\Mstar$, then $c$ is true in $\Mstar$;
(ii) $\emptyset \ent 0 = 1$ is not derivable;
(iii) for every $\Gamma$, $R$ and $x$, from $\Gamma \ent \Dec\,R$ and
$\Gamma \ent x\,\mathbf N$, $\Gamma \ent \Acc_R\cdot x \vee \Path_R\cdot x$;
(iv) for every well-formed $\Gamma$, and every $R$ and $x$, from $\Gamma \ent
\Dec\,R$ and $\Gamma \ent x\,\mathbf N$, $\Gamma \ent (\Path_R\cdot x)\,\mathbf B$;
(v) if $\Gamma \ent c$, $\Gamma$ is well-formed, every member of $\Gamma$
is true in the least model, and $c$ has some value $r$ in the least
model, then $c$ is true in the least model;
(vi) some $c$ has $\emptyset \ent c$ and is not true in the least model.
\end{thm}

Clause~(iii) has no well-formedness hypothesis on $\Gamma$.  The witness
of clause~(vi), $Z\cdot 0$, has no value in the least model by
\cref{thm:app-nuZ}.

\paragraph{Proof.}
\isn{dinga\_dich\_price} takes clauses (i), (ii), (v) and (vi) from the
corresponding \isn{oingad\_} facts through the two inclusions of
\isn{dinga\_eq\_oingad} (\isn{dinga\_nu\_sound},
\isn{dinga\_consistent}, \isn{dinga\_conservative},
\isn{dinga\_omega\_unsound}), and (iii) and (iv) directly
(\isn{dinga\_acc\_or\_path}, \isn{dinga\_dec\_path}).
The clauses for \isn{oinga\_dich}:
\begin{enumerate}
\item (i): \isn{oingad\_nu\_sound}, the skeleton \isn{step\_M\_sound} at
  $\Mstar$ with \isn{oingad\_wft} and the step soundness
  \isn{oingad\_step\_mjt}, which sends an \oinga step to
  \isn{oinga\_step\_mjt} and a dichotomy step to \cref{thm:app-dichstep}.
  (ii): \isn{oingad\_consistent}, by \isn{step\_M\_consistent}.
\item (iii): \isn{g\_acc\_or\_path}, proved in \isn{gcalc\_dich} and
  available through the sublocale $\isn{oinga\_dich} \subseteq
  \isn{gcalc\_ag\_nu\_dich}$.  Well-formedness of $R$ follows from $\Gamma
  \ent \Dec\,R$; the premise $\Dec\,R \wedge (0 = 0)$ is a conjunction
  introduction; one dichotomy step at the accessibility body concludes.
\item (iv): \isn{g\_dec\_path}, in \isn{gcalc\_gfix\_dich}, is the proof of
  \cref{thm:app-decnu} below at the accessibility body, with
  \isn{g\_acc\_not\_path} (in \isn{gcalc\_gfix}) for exclusion: transfinite
  induction at the proved accessibility (\isn{gi\_TI\_par}) with motive
  $\neg\Path_R\cdot y$, whose step is a double-negation introduction
  followed by folding the $\nu$-atom under negation, four backward head
  steps (\isn{g\_nu\_fold\_neg}).
\item (v): \isn{oingad\_conservative}, as for \cref{thm:app-inga-cons}: lift
  the hypotheses and the value $r$ to $\Mstar$ (\isn{omv\_le\_omvO\_Ostar}),
  apply (i), and use determinism.
\item (vi): \isn{oingad\_omega\_unsound}, with the derivation of $Z\cdot 0$
  embedded from \oinga (\isn{oingad\_Z}) and \cref{thm:app-nuZ}.
\end{enumerate}

\subsection{Exclusion and the decided $\nu$-atom}

\begin{thm}[\isn{g\_mu\_not\_nu}]
\label{thm:app-excl}
In \isn{gcalc\_gfix}: if $\Phi$ is \isn{posP}-positive at $1$, $\Psi$ is
its dual, $\Gamma$ is well-formed, $\Gamma \ent x\,\mathbf N$ and $\Gamma
\ent \mu\Phi\cdot x$, then $\Gamma \ent \neg(\nu(\neg\Psi)\cdot x)$.
\end{thm}

\begin{thm}[\isn{g\_dec\_nu}]
\label{thm:app-decnu}
In \isn{gcalc\_gfix\_dich}: if $\Phi$ is \isn{posP}-positive at $1$,
index $1$ is free in $\Phi$, $\Psi$ is its dual, $D$ is a premise of
$\Phi$, $\Gamma$ is well-formed, $\Gamma \ent D$, $\Gamma \ent x\,\mathbf
N$ and $\Phi$ is well-formed, then $\Gamma \ent (\nu(\neg\Psi)\cdot
x)\,\mathbf B$.
\end{thm}

Exclusion and the decided $\nu$-atom are \cref{thm:excl}.  Exclusion assumes neither that index $1$
occurs, nor any decidedness, nor any model.

\paragraph{Proof of \cref{thm:app-excl}.}
\begin{enumerate}
\item \emph{The dual induction} (\isn{g\_ndual\_beta\_env}, in \isn{gcalc}).
  Fix a well-formed $N$ and let $Q = \lambda y.\,\neg(N\cdot y)$.  Claim,
  by induction on the dual $t'$ of $t$ at $k$: for every well-formed
  $\Gamma$ and every substitution $\tau$ of well-formed terms with
  $\tau k = Q$, (a) $\Gamma \ent \tau\downarrow t$ implies $\Gamma \ent
  \tau[k \mapsto N]\downarrow t'$, and (b) $\Gamma \ent \neg(\tau\downarrow
  t)$ implies $\Gamma \ent \neg(\tau[k \mapsto N]\downarrow t')$.  The plain
  and negated statements travel together, so no polarity is needed.
  \begin{enumerate}
  \item Without index $k$: the two substitutions agree on $t$.
  \item Atom $v_k\cdot u$: with $w = \tau\downarrow u$, which the update
    does not change, $\tau\downarrow(v_k\cdot u) = Q\cdot w$, one beta step
    from $\neg(N\cdot w)$; (a) by \isn{ir\_convE}, (b) by
    \isn{ir\_conv\_negE}.
  \item Disjunction and negation: the elimination and introduction
    rules of each connective, with hypothesis (a) or (b) as the
    polarity requires; a negation's (a) uses its argument's (b).
  \item Universal: (a) by introduction at a fresh eigenvariable $j$
    (\isn{rp\_allI}): instantiate at $j$, apply hypothesis (a) at the
    substitution extended by $j$ (index $k+1$ still sent to $Q$), and fold
    back; (b) by \isn{rp\_not\_allE} at a refuted instance, hypothesis
    (b), one beta step under negation, and \isn{ir\_not\_allI}.
  \end{enumerate}
\item \emph{Body form} (\isn{g\_ndual\_beta}): from $\Gamma \ent
  \Phi[Q/a][x/x]$, $\Gamma \ent \Psi[N/a][x/x]$, by step~1(a) at the
  environment sending $0$ to $x$ and $1$ to $Q$.
\item \emph{The ground induction.}  Let $\nu = \nu(\neg\Psi)$ and the motive
  $c = \neg(\nu\cdot y)$ at the point $y$; then $\lambda c = Q$ at $N =
  \nu$.  Choose $i$ fresh for $\Gamma$, $\lambda c$ and $\lambda a.\lambda
  x.\Phi$, and let $\Delta$ be $i\,\mathbf N,\ \Phi[\lambda c/a][i/x],\
  \Gamma$.  Closure obligation: $\Delta \ent c[i]$.  In $\Delta$, step~2
  gives $\Psi[\nu/a][i/x]$; double-negation introduction gives
  $\neg((\neg\Psi)[\nu/a][i/x])$; folding under negation
  (\isn{g\_nu\_fold\_neg}) gives $\neg(\nu\cdot i)$, which is $c[i]$.
  Ground induction (\isn{ir\_gfix\_par}) with $\Gamma \ent \mu\Phi\cdot
  x$ and $\Gamma \ent x\,\mathbf N$ gives $c[x] = \neg(\nu\cdot x)$.  The
  $\nu$-atom is folded, never unfolded: the substance lies in the dual
  induction of step~1.
\end{enumerate}

\paragraph{Proof of \cref{thm:app-decnu}.}
\begin{enumerate}
\item The dichotomy (\isn{ir\_dich}) gives $\Gamma \ent M \vee P$ with $M =
  \mu\Phi\cdot x$ and $P = \nu(\neg\Psi)\cdot x$.
\item Eliminate the disjunction.  In the branch $M$, \cref{thm:app-excl} in
  the context $M, \Gamma$ gives $\neg P$; in the branch $P$ the
  hypothesis serves; a disjunction introduction gives $P \vee \neg P$
  in both branches, and \isn{ir\_boolI} gives $P\,\mathbf B$.
\end{enumerate}

\subsection{Decidedness of $\mathcal O$}

$\mathcal O = \mu\Phi_{\mathcal O}$ (\isn{trmO}, \isn{oPhi}) has seven
clauses, three for membership (zero, successor and limit codes) and four
for the order.  The limit clause reads the decoded sequence through the
graph atom $\mathrm{gr}(e, n, r) = \exists s.\,\mathrm{lock}\cdot(\mathrm{dc}
\cdot e)\cdot n\cdot s = r + 1$ (\isn{oGr}), where $\mathrm{dc}$ is the
decompiler of index trees (\isn{dcmp}, a general fixed point) and
$\mathrm{lock}$ the locked reader (\isn{lockT}): $\mathrm{lock}\cdot F\cdot
n\cdot s$ is the first nonzero value of $F\cdot n\cdot z$ for $z \le s$,
and $0$ if there is none.  $\Phi_{\mathcal O}$ is
\isn{posP}-positive, mentions $a$ (\isn{posP\_oPhi}, \isn{free1\_oPhi}), and
lies in the arithmetic class \isn{okP} (\isn{okP\_oPhi}): its
non-recursive parts are built by equations, disjunction, negation and
universals from terms built by variables, numerals, successor, pairing
and $\mathrm{lock}\cdot(\mathrm{dc}\cdot a)\cdot b\cdot c$ (\isn{okT},
\isn{okF}).

\begin{thm}[\isn{dcmp\_total}]
\label{thm:app-dcmp}
In \rga: from $\Gamma \ent e\,\mathbf N$, $\Gamma \ent u\,\mathbf N$ and
$\Gamma \ent v\,\mathbf N$, $\Gamma \ent (\mathrm{dc}\cdot e\cdot u\cdot
v)\,\mathbf N$, in every context.
\end{thm}

\begin{thm}[\isn{g\_dec\_O}, \isn{g\_O\_not\_coO}]
\label{thm:app-decO}
In \isn{gcalc\_ag\_dich}: if $\Psi$ is the dual of $\Phi_{\mathcal O}$,
$\Gamma$ is well-formed and $\Gamma \ent a\,\mathbf N$, then $\Gamma \ent
\mathcal O\cdot a \vee \nu(\neg\Psi)\cdot a$ and $\Gamma \ent
(\nu(\neg\Psi)\cdot a)\,\mathbf B$.  In \isn{gcalc\_gfix}: if $\Psi$ is the
dual of $\Phi_{\mathcal O}$, $\Gamma$ is well-formed, $\Gamma \ent
\mathcal O\cdot a$ and $\Gamma \ent a\,\mathbf N$, then $\Gamma \ent
\neg(\nu(\neg\Psi)\cdot a)$.
\end{thm}

The decidedness of $\mathcal O$ is \cref{thm:decO}, with $\mathrm{co}\mathcal O = \nu(\neg\Psi)$; a
dual exists (\isn{oPhi\_dual\_premise\_ex}).

\paragraph{Proof of \cref{thm:app-dcmp}.}
\begin{enumerate}
\item \isn{dct\_all}: $\emptyset \ent \forall x.\,\mathrm{dt}(x)$ with the
  closed motive $\mathrm{dt}(x) = \forall u\,v.\,(\mathrm{dc}\cdot x\cdot
  u\cdot v)\,\mathbf N$ (\isn{dctM}).  Course-of-values induction
  (\isn{ir\_cov\_ind\_acc}, derived from number induction over the
  accumulator ``$\mathrm{dt}$ holds below $x$'') reduces the goal to:
  from $x\,\mathbf N$ and the accumulator, derive $\mathrm{dt}(x)$; by
  \isn{dctM\_intro} it suffices to derive $(\mathrm{dc}\cdot x\cdot u\cdot
  v)\,\mathbf N$ at fresh numbers $u, v$ (\isn{dct\_core}).
\item \isn{dct\_core}.  The induction hypothesis gives totality at every
  $y < x$ (\isn{acc\_descent\_strict}).  Split $x$ into the four initial
  codes, at which the decompiler returns a total compiled initial
  function (\isn{dct\_init\_N}), and $x = K + 4$ (\isn{dct\_cases5}).  At
  $K + 4$ the children's codes are numbers at most $K$, hence below $x$
  (\isn{dct\_dG\_le}, \isn{dct\_dH1\_le}, \isn{dct\_dH2\_le}, \isn{dct\_lt4}),
  and the parity of $K$ chooses the node (\isn{dct\_dPar}).  At a composition node (\isn{dcmp\_disp\_sub})
  three applications of the hypothesis give totality.  At a recursion
  node (\isn{dcmp\_disp\_rec}) the result is a recursor term $R\cdot g
  \cdot h\cdot v$, total by number induction on the recursion argument
  (\isn{dct\_R\_N}), with base and step from the hypothesis.
\item The theorem instantiates \isn{dct\_all} at $e$ and the universals at
  $u, v$ (\isn{dctM\_use}).  \isn{g\_dcmp\_total} carries the theorem to every
  calculus.
\end{enumerate}

\paragraph{Proof of \cref{thm:app-decO}.}
\begin{enumerate}
\item \emph{The reader is total} (\isn{g\_lock\_total}, in \isn{gcalc}): from
  \cref{thm:app-dcmp} and numbers $e, n, s$, the locked reader's loop
  is a number by number induction on the step count (\isn{g\_lockK\_N}),
  and $\mathrm{lock}\cdot(\mathrm{dc}\cdot e)\cdot n\cdot s$ head-reduces
  to that loop (\isn{lockT\_unf}).
\item \emph{The premise is derivable} (\isn{g\_nfdec\_ok}, in
  \isn{gcalc\_ati}): given the reader's totality, if $D$ is a premise of
  $t$ at $k$, $t$ lies in \isn{okP} at $k$ and every free variable of $t$
  is at most $k$, then $\Gamma \ent D$ for every well-formed $\Gamma$.
  Induction on the premise relation.  At the leaf clause, $t$ is in
  \isn{okF}, and \isn{g\_allqN\_env} (an induction on $k$ by universal
  introduction) reduces the goal to $t'\,\mathbf B$ at every assignment of
  numbers to the quantified indices.  \isn{g\_okF\_B}, an induction on
  \isn{okF}, decides that formula: an equation between \isn{okT} terms is decided
  because both sides are numbers (\isn{g\_okT\_N}, an induction on
  \isn{okT} whose reader case is step~1, and \isn{ir\_eqTI}); disjunction
  and negation by \isn{ir\_orTI} and \isn{ir\_negTI}; a universal by
  deciding its body at a fresh number and then ATI (\isn{ir\_allTI}).  ATI
  enters only there; in $\Phi_{\mathcal O}$ the universals so decided
  include the existential over the step count in the graph atom.  At an
  atom the premise $0 = 0$ is reflexivity; a conjunction takes both
  hypotheses; negation and the universal pass through.
  \isn{g\_oPhi\_premise} instantiates at $\Phi_{\mathcal O}$
  (\isn{okP\_oPhi}, \isn{oPhi\_free\_le}).
\item \emph{Assembly} (\isn{g\_dec\_O\_total}, \isn{g\_dec\_O\_of\_dcmp\_total}):
  obtain a premise $D$ (\isn{nfdec\_exists}); step~2 derives it; one
  dichotomy step at $\Phi_{\mathcal O}$ (\isn{g\_dec\_O\_cond}) gives the
  disjunction, and \cref{thm:app-decnu} gives the decidedness.
\item \isn{g\_O\_not\_coO} is \cref{thm:app-excl} at $\Phi_{\mathcal O}$,
  with no premise, no ATI and no dichotomy.
\end{enumerate}

\paragraph{The characteristic term.}
Files: \isn{IGA\_Closing/GChiO.thy} (imports \isn{CO\_Dinga} and
\isn{IGA/IGA\_KleeneO\_Meas.thy}).  $\mathcal O_{\mathrm L}$
(\isn{oL}) is the least model's $\mathcal O$, the set of $n$ with
$\mathcal O\cdot\bar n$ true in the least model (\isn{oL\_iff}).  For a
dual $\Psi$ of $\Phi_{\mathcal O}$, $\chi_\Psi = \lambda
a.\,\neg(\nu(\neg\Psi)\cdot a)$ (\isn{chiO}).

\begin{thm}[\isn{O\_Mstar\_iff}, \isn{dinga\_chiO\_char}]
\label{thm:app-chiO}
(i) In \isn{rga}: $\mathcal O\cdot\bar n$ has the value $1$ in
$\Mstar$ if and only if $n \in \mathcal O_{\mathrm L}$.
(ii) In \isn{dinga}, for every dual $\Psi$ of $\Phi_{\mathcal O}$:
$\Gamma \ent \Tot\,\chi_\Psi$ at every well-formed $\Gamma$
(\isn{dinga\_Tot\_chiO}), and $\chi_\Psi\cdot\bar a$ has in $\Mstar$
the value $1$ if $a \in \mathcal O_{\mathrm L}$ and $0$ otherwise
(\isn{dinga\_chiO\_Mstar}).
(iii) Hence some $\chi$ has $\Gamma \ent \Tot\,\chi$ at every
well-formed $\Gamma$ and is, in $\Mstar$, the characteristic function
of $\mathcal O_{\mathrm L}$ (\isn{dinga\_chiO\_char}).
\end{thm}

\paragraph{Proof of \cref{thm:app-chiO}.}
\begin{enumerate}
\item \emph{From the least model to $\Mstar$}: a value in the least
  model persists (\isn{omvO\_bot\_le\_omv\_le\_Ostar}).
\item \emph{The leaves are valued} (\isn{oiga\_okF\_omv}): an \isn{okF}
  formula at an environment of numerals is decided in \oiga
  (\isn{g\_okF\_B}, with the reader's totality from
  \cref{thm:app-dcmp}), so by \oiga's soundness for the least model it
  has a value there.
\item \emph{The reading relative to $\Mstar$ is the least model's}
  (\isn{osemM\_osem\_okP}): at a body of \isn{okP} whose free variables
  lie at most at the point, under an environment of numerals, the
  nat-set reading relative to $\Mstar$ implies the reading relative to
  the least model, because every leaf has a least-model value by
  step~2, which persists to $\Mstar$ by step~1, where values are
  unique.
\item \emph{From $\Mstar$ to the least model} (\isn{O\_Mstar\_le\_oL}):
  ground induction in $\Mstar$ at a closed set of numbers
  (\isn{gfix\_omvO\_ind\_set\_par}) with the set $\mathcal O_{\mathrm L}$:
  a point whose body reads true over $\mathcal O_{\mathrm L}$ relative to
  $\Mstar$ reads true relative to the least model by step~3, so it lies
  in the operator's image of $\mathcal O_{\mathrm L}$, hence in
  $\mathcal O_{\mathrm L}$, a fixed point.  Steps 1 and 4 give (i).
\item \emph{Totality} (\isn{dinga\_Tot\_chiO}): at a number $j$, the
  co-set atom is decided (\isn{dinga\_dec\_O\_B}), so its negation is,
  and a beta step (\isn{chiO\_step}) carries the decidedness to
  $\chi_\Psi\cdot j$; a decided term equals $1$ or $0$
  (\isn{ir\_true\_eqI}, \isn{ir\_false\_eqI}), so it is a number.  No
  least-witness search is involved.
\item \emph{The value in $\Mstar$} (\isn{dinga\_chiO\_Mstar}).  If $a \in
  \mathcal O_{\mathrm L}$: from the hypothesis $\mathcal O\cdot\bar a$,
  \isn{g\_O\_not\_coO} and a conversion derive $\chi_\Psi\cdot\bar a$,
  and soundness for $\Mstar$ (\isn{dinga\_nu\_sound}) with (i) makes it
  true.  Otherwise the dichotomy's disjunction $\mathcal O\cdot\bar a
  \vee \nu(\neg\Psi)\cdot\bar a$ is true in $\Mstar$, its first
  disjunct is not by (i), so the co-set atom is true, and from it
  $\neg(\chi_\Psi\cdot\bar a)$ is derived under that hypothesis and is
  true, so $\chi_\Psi\cdot\bar a$ has the value $0$.
\item (iii) takes a dual from \isn{oPhi\_dual\_premise\_ex}.
\end{enumerate}
The reading of $\chi$ as a characteristic function is $\Mstar$'s: off
$\mathcal O_{\mathrm L}$ the co-set atom, and with it $\chi\cdot\bar a$,
has no value in the least model.
\subsection{Accessibility is never false}

\begin{thm}[\isn{accO\_false\_empty}, \isn{oingad\_acc\_not\_refuted}, \isn{dinga\_acc\_not\_refuted}]
\label{thm:app-accO}
If $R$ is closed, then for every term $x$, $\Acc_R\cdot x$ does not have
the value $0$ in $\Mstar$.  If $R$ is closed, then for every $x$, $\emptyset
\ent \neg(\Acc_R\cdot x)$ is not derivable in \dinga.
\end{thm}

Neither statement assumes decidedness of $R$ or well-formedness of $R$
or $x$; the underivability is for \dinga at the empty context.

\paragraph{Proof.}
\begin{enumerate}
\item \emph{Shapes.}  Write $B_x = \neg R(y, x) \vee \Acc_R\cdot y$ for the
  body at the point $x$ under the binder $y$ (\isn{accBody}).  The
  invariant \isn{accInvOF} says: if $t$ has the value $0$, then $t$ is
  none of (S1) $\cA\cdot\lambda y.B_x$, (S2) $(\lambda y.B_x)\cdot\bar k$,
  (S3) $\neg R(\bar k, x) \vee \Acc_R\cdot\bar k$, (S4) $\Acc_R\cdot x$.  For
  closed $R$, the head of $\Acc_R$ in every oracle model is the
  $\lambda$-head whose body, at the $\lambda$'s variable $x$, is
  $\cA\cdot\lambda y.B_x$ (\isn{omhO\_trmAcc}).
\item \isn{accO\_false\_inv}: every derivation of $\mathrm{val}_{\Theta^*}(t,
  r)$ satisfies the invariant.  Mutual rule induction on the oracle
  model, with the trivial motive for heads.  The cases with value $0$:
  \begin{enumerate}
  \item False disjunction: in shape (S3) the right disjunct
    $\Acc_R\cdot\bar k$ is shape (S4), excluded by the hypothesis.
  \item Beta, $t\cdot u$ with $t$ at a $\lambda$-head: (S1) would make
    $t = \cA$, which has no $\lambda$-head (\isn{omhO\_A\_inv}); in (S2)
    the head body is $B_x$ (\isn{omhO\_lam\_inv}) and the reduct
    $B_x[\bar k]$ is (S3) (\isn{bsub\_accBody}); in (S4) determinism of
    heads in $\Mstar$ with \isn{omhO\_trmAcc} makes the reduct (S1)
    (\isn{acc\_beta2}).  The hypothesis excludes each of these shapes.
  \item False universal: in (S1) the refuted instance is (S2).
  \item Recursor clauses and the false oracle reading: the head is a
    recursor head or a $\nu$-head, so no shape matches.  The false
    oracle clause (\isn{omvO\_nuF}), the only value-$0$ clause specific
    to $\Mstar$, is discharged the same way: a $\nu$-head is the head of
    no shape.
  \end{enumerate}
  \isn{accO\_false\_empty} is the invariant at shape (S4).
\item \isn{oingad\_acc\_not\_refuted}: a derivation of $\emptyset \ent
  \neg(\Acc_R\cdot x)$ would make $\neg(\Acc_R\cdot x)$ true in $\Mstar$
  by \cref{thm:app-price}(i), so $\Acc_R\cdot x$ would have the value $0$
  (\isn{omvO\_not\_inv}), against step~2.  \isn{dinga\_acc\_not\_refuted}
  is the same statement at \dinga, through \isn{dinga\_eq\_oingad}.
\end{enumerate}

\subsection{The dichotomy is not valid for the least model}

Files: \isn{OIGA/GDich.thy} (context \isn{rga}).  The least model reads the $\nu$-former by its unfolding, so a $\nu$-application is true there only through a finite derivation, which puts the point into the least fixed point of the body's truth operator.

\begin{thm}[\isn{nu\_omv\_le\_lfp}, \isn{dich\_not\_omega\_valid}; classical]
\label{thm:app-notvalid}
(i) For a positive body $\Phi$ with $\mu\Phi$ closed, if $\nu\Phi\cdot\bar n$ is true in the least model then $n \in \lfp\,F_\Phi$ (and $\mu\Phi\cdot\bar n$ is true as well).
(ii) For a closed well-formed order $R$ with $\mathrm{accRelN}\,R$ not well-founded, the dual $\Psi$ of the accessibility body and any premise $D$ of that body that is valid in the least model: it is not the case that, for every well-formed context $\Gamma$ and well-formed $x$ with $D$ and $x\,\mathbf N$ true under every substitution satisfying $\Gamma$, the conclusion $\Acc_R\cdot x \vee \nu(\neg\Psi)\cdot x$ is true as well.  The variant \isn{dich\_not\_omega\_valid\_decRel} takes the premise $\Dec\,R \wedge 0 = 0$ derived in \rga instead of valid.
\end{thm}

\begin{enumerate}
\item \emph{Transport at a head.}  \isn{hfix\_transport} replays the transport lemma of the IGA paper's soundness case with the hypothesis ``$T$ is the $\mu$-term'' replaced by ``$T$ reaches the head form with body $\Phi[T/a]$'' (\isn{omh}); the original proof used $T$ only through that head fact.  \isn{omh\_trmNu} supplies the head fact for $\nu\Phi$.
\item \emph{The closure condition.}  \isn{osem\_upd\_rec}: the oracle reading never consults the environment at the recursive index, so the closure condition at $\nu\Phi$ is the one for $F_\Phi$.  Induction on the derivation of the value then gives (i) (\isn{nu\_omv\_ind\_set}, \isn{nu\_omv\_le\_lfp}).
\item \emph{Path is nowhere true.}  \isn{pathPhi\_gfixOp\_empty}: the least fixed point of the path body's truth operator is empty, so $\Path_R\cdot\bar n$ is never true in the least model (\isn{path\_omv\_never\_true}); hence the disjunction is true exactly where $\Acc_R$ is, and $\forall n.\,\Acc_R\cdot\bar n \vee \Path_R\cdot\bar n$ true is equivalent to well-foundedness (\isn{acc\_or\_path\_omega\_iff\_wf}, from \isn{accbiv\_omega\_sound\_iff\_wf}).
\item \emph{The rule-validity form.}  Instantiate at the empty context and at a numeral where $\Acc_R$ fails (\isn{acc\_or\_path\_omega\_not\_wf}); \isn{ndual\_fun} identifies $\Psi$ with the dual the rule builds.  For the variant, \isn{rga\_sound\_oomjt} turns the \rga derivation of the premise into its validity.
\end{enumerate}

\paragraph{Lemmas a re-mechanization needs.}
\begin{itemize}
\item \isn{ndual\_fun}, \isn{ndual\_exists}, \isn{ndual\_posP},
  \isn{ndual\_subf2}: the dual is functional, exists in the grammar, is
  positive under negation and commutes with substitution.
\item \isn{gfix\_omvO\_eq\_lfp}: the $\mu$-atom's truth set in $\Mstar$ is
  $\lfp\,F^*_\Phi$.
\item \isn{ndual\_swap}, \isn{GN\_dual}: the falsity reading of the dual is
  the truth reading of the body, so $\hat G_{\neg\Psi} = F^*_\Phi$.
\item \isn{nrestrict\_T}, \isn{nrestrict\_F}: on numerals the oracle's
  sets are the fixed points of the numeral operators.
\item \isn{disj\_at\_term}: bivalence and the swap give the disjunction at
  any valued term.
\item \isn{nfdec\_leafdec}, \isn{ndual\_leafdec}, \isn{leafdec\_dposE}: a
  true premise puts the $\nu$-body's instances in the determinacy grammar.
\item \isn{dposE\_det}, \isn{gdetM\_union}, \isn{dich\_premise\_bivalent}:
  determinacy gives bivalence of the $\nu$-atom.
\item \isn{dich\_valid\_at}, \isn{oingad\_step\_mjt}: the step, and the
  system's step soundness.
\item \isn{g\_acc\_not\_path}, \isn{g\_nu\_fold\_neg}: accessibility excludes
  a path; a $\nu$-atom folds under negation.
\item \isn{g\_ndual\_beta\_env}: the dual induction under a negated
  predicate.
\item \isn{dct\_core}, \isn{ir\_cov\_ind\_acc}: totality of the decompiler
  by course-of-values induction.
\item \isn{g\_lock\_total}, \isn{g\_okF\_B}, \isn{g\_nfdec\_ok}: the premise of
  the $\mathcal O$ body, with ATI at the universals.
\item \isn{omhO\_trmAcc}, \isn{accO\_false\_inv}: the four-shape invariant
  for accessibility in $\Mstar$.
\end{itemize}

\section{\inga's certificate-grounded semantics}
\label{app:inga}

Files, all in the session \textsf{INGA\_PGSem}:
\isn{INGA\_PG0.thy} (the negative result and the reduction of the
interface to two properties), \isn{INGA\_PG1.thy} (the certified
evaluation, the oracle map and its least fixed point, determinism,
coherence, positive replacement and the semantic Park lemma),
\isn{NG\_Eval.thy} (completeness at a certified oracle),
\isn{NG\_Kit.thy} (substitution, value congruence and head steps for
the certified evaluation), \isn{NG\_Closure.thy} (the certified oracle
closed under substitution and value congruence), \isn{NG\_Step.thy}
(the rules of \rga), \isn{NG\_GfixStep.thy} (ground induction) and
\isn{NG\_Sound.thy} (soundness, the two properties and the instances).
All results are classical: HOL theorems in the locale \isn{inga}, or
in the locales below, which fix an arbitrary calculus.

The appendix proves \cref{thm:inga-pg,cor:inga-dpwp}.  The plan: the
interface of \cref{thm:nocert}'s hypotheses, read with truth as
derivability, reduces to the disjunction property and the witness
property (\cref{thm:app-inga-reduce}); a compositional evaluation with
an oracle for $\nu$-atoms, filtered by derivations, is complete at its
least fixed point (\cref{thm:app-ng-complete}); the coinduction step
(\cref{thm:app-ng-park}), the rules of \rga and ground induction
(\cref{thm:app-ng-step}) preserve validity for that evaluation, once
the evaluation is closed under substitution and value congruence
(\cref{thm:app-ng-kit,thm:app-ng-closure}); so \inga is sound for it,
and the two properties follow, with consistency
(\cref{thm:app-ng-sound}).

\paragraph{The calculus.}
The construction is stated over a calculus $S$, a step relation on
proof lists, contexts and conclusions in the sense of the \iga
paper's certificate semantics~\cite{iga-paper}: every \rga step is an
$S$ step, and $S$ is monotone in the proof list and preserves
well-formedness (the locale \isn{gsem}).  The locale \isn{ngsem} is
\isn{gsem} with the coinduction step (\isn{gstep\_nu}): $S$ contains
it.  ``$S$ derives $c$'' and
$\ent_S c$ mean that $\emptyset \ent c$ is derivable for $S$;
``$S$ refutes $c$'' means $\ent_S \neg c$.  \inga is the instance whose
step relation (\isn{inga\_derives}) is \iga's step or a coinduction
step (\isn{nuS}); \rga with the coinduction rule is the instance with
step relation \isn{rganu\_derives}, \rga's step or a coinduction step.

\paragraph{The interface.}
The hypotheses of \cref{thm:nocert} on a pair of predicates
$\mathrm{Tr}$ and $\mathrm{Fa}$ on terms form the interface
\isn{cert\_iface}, six conditions.
(allT) If $\mathrm{Tr}(\cA\cdot u)$, then $S$ derives the open instance
of the universal, $u$ applied to a fresh variable $v$ in the context
$v\,\mathbf N$.
(allF) If $\mathrm{Fa}(\cA\cdot u)$, then $\mathrm{Fa}$ holds at some
numeral instance $u\cdot\bar n$.
(orT) If $\mathrm{Tr}(a \vee b)$, then $\mathrm{Tr}(a)$ or
$\mathrm{Tr}(b)$.
(notT) If $\mathrm{Tr}(\neg a)$, then $\mathrm{Fa}(a)$.
(excl) Never both $\mathrm{Tr}(a)$ and $\mathrm{Fa}(a)$.
(rtrue) $\mathrm{Tr}(c)$ for every closed theorem $c$ of \rga.

\paragraph{\rga's certificate-grounded evaluation.}
The \rga and \iga papers' semantics for a calculus $S$ (\isn{gwhev},
\isn{gwosr}) is the least model's pair of relations, a head judgment
and a value judgment, with the two clauses for the universal
replaced.  $\cA\cdot u$ is true when $S$ derives the open instance of
the universal, by a proof list whose code is bounded by the
relation's natural-number index, and every numeral instance
$u\cdot\bar n$ is true; it is false when $S$ refutes an instance and
that instance is false.  A term is true or false when it has the
value $1$ or $0$ at some index.  The universal's truth needs a
certificate, its falsity a refuted false instance: the evaluation
meets the interface's conditions by construction.

\paragraph{The certified evaluation.}
The certified evaluation at an oracle $\Theta$ (\isn{gnhev},
\isn{gnosr}; written $\mathrm{hd}^{\mathrm c}_\Theta$ and
$\mathrm{val}^{\mathrm c}_\Theta$) is \rga's certificate-grounded
evaluation with the three changes that turned the least model into
the oracle model in Appendix~\ref{app:nu}: the head-level beta clause
does not unfold a $\nu$-term of the class; a $\nu$-term $\nu\Phi$ of
the class has the $\nu$-head $\nu_\Phi$ (\isn{nhNu}, \isn{gnhev\_nu});
and if $t$ has the $\nu$-head $\nu_\Phi$ and $u \in \Theta(\Phi, b)$,
then $t\cdot u$ has the value $b$ (\isn{gnosr\_nuT}, \isn{gnosr\_nuF}).
The oracle occurs only positively in the clauses, so the evaluation is
monotone in the oracle.  The \emph{readings} $\mathrm{rd}^\pm_{\Theta,
\mathrm c}$ (\isn{gosem}, \isn{gfsem}) are the oracle readings
$\mathrm{rd}^\pm_\Theta$ with non-recursive leaves read by the
certified evaluation and with the universal clauses of \rga's
evaluation, whose certificate is a derivation of the open instance of
the universal as substituted by the environment.  The two operators
are filtered (\isn{FOg}, \isn{GOg}):
\begin{align*}
F^\Theta_{\Phi,\mathrm c}(X) &= \{u \mid \ent_S \nu\Phi\cdot u,\
  \mathrm{rd}^+_{\Theta,\mathrm c}(X, 1, \mathrm{env}_\Phi(u), \Phi, 1)\},\\
G^\Theta_{\Phi,\mathrm c}(X) &= \{u \mid \ent_S \neg(\nu\Phi\cdot u),\
  \mathrm{rd}^-_{\Theta,\mathrm c}(X, 1, \mathrm{env}_\Phi(u), \Phi, 0)\}.
\end{align*}
The map $\mathrm{Om}_{\mathrm c}$ (\isn{OmG}) sends $\Theta$ to the
oracle with $\gfp\,F^\Theta_{\Phi,\mathrm c}$ as truth side and
$\lfp\,G^\Theta_{\Phi,\mathrm c}$ as falsity side at a body $\Phi$ of
the class, both empty elsewhere.  It is monotone, and the
\emph{certified oracle} $\Theta^*_{\mathrm c}$ (\isn{OstarG}) is its
least fixed point, so that at a body of the class
$\Theta^*_{\mathrm c}(\Phi, 1) = \gfp\,F^{\Theta^*_{\mathrm
c}}_{\Phi,\mathrm c}$ and $\Theta^*_{\mathrm c}(\Phi, 0) =
\lfp\,G^{\Theta^*_{\mathrm c}}_{\Phi,\mathrm c}$.  The filter is the
only difference from $\Theta^*$: a point enters a side only with a
derivation of the atom or of its negation.

A judgment $H \ent c$ is \emph{valid} at an oracle (\isn{gnowsjt})
when, for every well-formed substitution $\sigma$ whose instances of
the members of $H$ are true in the certified evaluation at the oracle,
$\sigma\downarrow c$ is true there.  An oracle is \emph{certified}
(\isn{ocert}) when every point of its truth side at $\Phi$ has $\ent_S
\nu\Phi\cdot u$ and every point of its falsity side has $\ent_S
\neg(\nu\Phi\cdot u)$.

\subsection{The negative result and the reduction}

\begin{thm}[\isn{inga\_not\_gsem\_sound}]
\label{thm:app-inga-notgsem}
There are no step relations $S$ and $X$ such that $S$ and $X$ form an
instance of the locale \isn{gsem\_sound} and every \inga step is an
$S$ step.
\end{thm}

The locale \isn{gsem\_sound} is the \iga paper's soundness of a
calculus for \rga's certificate-grounded evaluation: $S$ splits into
\rga's step and a further relation $X$ whose steps preserve validity.

\paragraph{Proof.}
\begin{enumerate}
\item Every instance of \isn{gsem\_sound} is sound for the least model
  (\isn{ng\_gs\_omega\_true}): \rga's certificate-grounded evaluation
  refines the least model, so a derivable sentence, true in the
  evaluation by soundness, is true in the least model.
\item If every \inga step were an $S$ step, the \inga derivation of the
  all-zeros atom $Z\cdot 0$ (\cref{thm:app-ingaZ}) would be an $S$
  derivation, and $Z\cdot 0$ would be true in the least model, where
  it has no value (\isn{inga\_Z\_valueless}, from the least model's
  fact alone; \cref{thm:app-nuZ}).
\end{enumerate}

\begin{thm}[\isn{inga\_cert\_iff\_props}, \isn{inga\_pg\_iff\_dp\_wp}]
\label{thm:app-inga-reduce}
Let $\mathrm{Tr}(c)$ say that \inga derives $c$ from no hypotheses and
$\mathrm{Fa}(c)$ that \inga derives $\neg c$ from no hypotheses.  The
interface holds for $\mathrm{Tr}$ and $\mathrm{Fa}$ if and only if
\inga has the disjunction property (from $\ent a \vee b$, $\ent a$ or
$\ent b$) and the witness property for refuted universals (from $\ent
\neg(\cA\cdot u)$, the refutation of some numeral instance of $u$).
Equivalently (\isn{inga\_pg\_iff\_dp\_wp}), the two properties hold
if and only if some pair $\mathrm{Tr}$, $\mathrm{Fa}$ meets the
interface with every derivable term in $\mathrm{Tr}$, every term in
$\mathrm{Tr}$ derivable, and every term in $\mathrm{Fa}$ refutable.
\end{thm}

\paragraph{Proof.}
\begin{enumerate}
\item (allT) holds outright (\isn{inga\_cert\_allT}): eliminate the
  universal at a fresh variable.
\item (excl) holds outright (\isn{inga\_cert\_excl}): from $\ent a$ and
  $\ent \neg a$, \inga derives $0 = 1$, against \inga's consistency.
  The theory \isn{INGA\_PG0} proves the reduction with consistency as
  a hypothesis (\isn{inga\_cert\_excl\_con},
  \isn{inga\_cert\_iff\_props\_con}); the hypothesis-free forms are
  taken in \isn{NG\_Sound} from the consistency that soundness for
  the certified evaluation gives (\isn{inga\_consistent\_mf},
  \cref{thm:app-inga-pg}), with no appeal to $\Mstar$.
\item (notT) is the definition: $\mathrm{Tr}(\neg a)$ is
  $\mathrm{Fa}(a)$.  (rtrue) is the embedding of \rga's derivations
  into \inga's (\isn{rga\_inga\_embed}).
\item (orT) is the disjunction property and (allF) the witness
  property, word for word.
\end{enumerate}

\subsection{The certified oracle}

\begin{thm}[\isn{ocons\_OstarG}, \isn{gnosr\_unique}]
\label{thm:app-ng-oracle}
For every calculus $S$: every point of $\Theta^*_{\mathrm c}(\Phi, 1)$
has $\ent_S \nu\Phi\cdot u$ (\isn{OstarG\_T\_derivable}), and every
point of $\Theta^*_{\mathrm c}(\Phi, 0)$ has $\ent_S \neg(\nu\Phi\cdot
u)$ (\isn{OstarG\_F\_refutable}); so $\Theta^*_{\mathrm c}$ is
certified (\isn{ocert\_OstarG}).  $\Theta^*_{\mathrm c}$ is consistent,
with no hypothesis on $S$ (\isn{ocons\_OstarG}; at \inga,
\isn{inga\_ocons\_OstarG}).  At a consistent oracle no term has two
values in the certified evaluation (\isn{gnosr\_unique}; at
$\Theta^*_{\mathrm c}$, \isn{gnosr\_OstarG\_unique\_mf}).
\end{thm}

\paragraph{Proof.}
\begin{enumerate}
\item The two inclusions are the fixed-point equations with the
  filter's first conjunct.  They make the certified oracle certified
  by its filter alone.
\item \emph{Determinism} is a simultaneous rule induction in the shape
  of \isn{omhO\_det\_excl} (step~3(a) of the proof of
  \cref{thm:app-rga-nu-sound}): a term with a head has no other head
  and no value, and a term with a value has no other value and no
  head.  The universal clauses are \rga's, and their cases are those
  of \rga's evaluation.  At the two oracle clauses, a true and a false
  reading of $t\cdot u$ at a $\nu$-head put $u$ on both sides of the
  oracle, which is consistent.  The beta guard separates a $\nu$-head
  from a $\lambda$-head.
\item \emph{The readings exclude each other}
  (\isn{gosem\_gfsem\_excl}).  At a consistent oracle $\Theta$, and
  sets $X$ and $Y$ with no common point, the truth reading over $X$
  and the falsity reading over $Y$ never give one body opposite
  values.  The proof is an induction on the truth reading.  At a leaf
  where the recursive variable is not free, both readings are values
  of one term in the certified evaluation (\isn{gfsem\_nf\_only}), and
  determinism (step~2) excludes opposite values.  At a recursive atom,
  the truth reading puts the point in $X$ and the falsity reading puts
  it in $Y$ (\isn{gfsem\_app\_F}).  The connectives and the universal
  invert the falsity reading (\isn{gfsem\_or\_F\_any},
  \isn{gfsem\_all\_F\_any} and their companions) and close by the
  induction hypothesis.
\item \emph{The two sides are disjoint} (\isn{gfpG\_lfpG\_disjoint}).
  At a consistent $\Theta$, write $N = \gfp\,F^\Theta_{\Phi,\mathrm c}$.
  The complement of $N$ is a pre-fixed point of $G^\Theta_{\Phi,\mathrm
  c}$: a point of $G^\Theta_{\Phi,\mathrm c}(\overline N)$ that lay in
  $N$ would lie in $F^\Theta_{\Phi,\mathrm c}(N)$, and its body would be
  read true over $N$ and false over $\overline N$, against step~3.  So
  $\lfp\,G^\Theta_{\Phi,\mathrm c} \subseteq \overline N$ by
  \isn{lfp\_lowerbound}.  The filters play no part.
\item \emph{The oracle's consistency} (\isn{ocons\_OstarG}).  By
  step~4 the map $\mathrm{Om}_{\mathrm c}$ preserves consistency
  (\isn{ocons\_OmG}); consistency holds at the empty oracle and passes
  to suprema of chains (\isn{ocons\_admissible}, the lemma of
  \isn{IGA\_Nu} used for $\Theta^*$).  Fixed-point induction along the
  iteration of $\mathrm{Om}_{\mathrm c}$ (\isn{fixp\_induct}) gives
  consistency at its least fixed point.  A second route,
  \isn{ocons\_OstarG\_filter}, needs the calculus's consistency: a
  point on both sides would let $S$ derive an atom and its negation.
\end{enumerate}
The oracle is read on both sides.  If the evaluation read only truth
from the oracle and still unfolded a $\nu$-atom to obtain its
falsity, the induction of step~2 would meet a $\nu$-atom true by the
oracle, a clause with no subderivation, and false through its
unfolding, and the two verdicts would be related only through the
fixed point (\S\ref{sec:price}; a reading, not a theorem).

\subsection{Completeness at a certified oracle}

\begin{thm}[\isn{gnhev\_gnosr\_complete}, \isn{ng\_complete}]
\label{thm:app-ng-complete}
At every certified oracle, a well-formed term that is true in the
certified evaluation is derivable for $S$, and a false one is
refutable (\isn{ng\_complete\_F}).  The forms used below are
\isn{ng\_complete\_val} (at a value: $\ent_S a = \bar r$),
\isn{ng\_complete\_open} (a well-formed term, open or closed, read
true at a certified oracle is derivable; applied below to the
$\sigma$-instances of a judgment's hypotheses) and
\isn{ng\_complete\_neg} (one read false is refutable).
\end{thm}

\paragraph{Proof.}
\begin{enumerate}
\item \isn{gnhev\_gnosr\_complete} replays \rga's completeness of its
  certificate-grounded evaluation (\isn{gwhev\_gwosr\_complete})
  clause for clause: a simultaneous rule induction whose motive
  relates a term with a value to a derivation that the term equals
  that value, and a term with a head to its canonical head term
  (\isn{rga.nhd\_canon}).  The head half of the motive is
  \isn{gnhev\_ok}: the canonical head term is well formed, and the
  term inherits every derivable spine equation of the canonical head
  term.
\item The universal clauses carry their certificates, as in \rga.
\item The two oracle clauses carry theirs through the oracle: the
  oracle is certified, so $u \in \Theta(\Phi, 1)$ gives $\ent_S
  \nu\Phi\cdot u$ and $u \in \Theta(\Phi, 0)$ gives $\ent_S
  \neg(\nu\Phi\cdot u)$.  \isn{ir\_true\_eqI} and \isn{ir\_false\_eqI}
  turn these derivations into the derivable equations of the atom
  with its value, which is the motive's conclusion.
\item At $\Theta^*_{\mathrm c}$, which is certified
  (\cref{thm:app-ng-oracle}), the theorem discharges the one
  assumption of the locale \isn{ngkit} (\isn{kit\_complete}, the
  completeness of the evaluation at $\Theta^*_{\mathrm c}$), so the
  locale \isn{ngsem} of the construction is a sublocale of \isn{ngkit},
  and \inga is an instance (the re-exports \isn{inga\_complete},
  \isn{inga\_nu\_true\_derivable}, \isn{inga\_nu\_false\_refutable},
  \isn{inga\_true\_not\_refutable}, \isn{inga\_false\_not\_derivable}).
\end{enumerate}

\paragraph{Coherence.}
\isn{ng\_coherent\_T} and \isn{ng\_coherent\_F}: for $\Phi$ in the
class, $u \in \Theta^*_{\mathrm c}(\Phi, 1)$ if and only if $\ent_S
\nu\Phi\cdot u$ and the unfolding $\Phi[\nu\Phi/a][u/x]$ has the value
$1$ in the certified evaluation at $\Theta^*_{\mathrm c}$; dually for
$0$, with $\ent_S \neg(\nu\Phi\cdot u)$.  No completeness is needed:
the proof is that of \isn{nu\_coherent\_T} (step~4 of the proof of
\cref{thm:app-rga-nu-sound}), through the fixed-point equation and the
transport of the reading into the evaluation, in both directions.  The
forms without the derivation conjunct, \isn{ng\_coherent\_T\_iff} and
\isn{ng\_coherent\_F\_iff}, hold for well-formed $\Phi$ and $u$, and
are where completeness enters: from a true unfolding, completeness
gives a derivation of the unfolding, and a fold
(\isn{g\_nu\_fold} on the truth side; on the falsity side four
conversions under negation through \isn{nu\_step1} to
\isn{nu\_step4}) turns it into $\ent_S \nu\Phi\cdot u$ or $\ent_S
\neg(\nu\Phi\cdot u)$, the derivation conjunct.  At \inga:
\isn{inga\_coherent\_T\_iff}, \isn{inga\_coherent\_F\_iff}.

\subsection{The semantic Park lemma}

\begin{thm}[\isn{nu\_step\_gnowsjt}; at \inga, \isn{inga\_nu\_step\_gnowsjt}]
\label{thm:app-ng-park}
In the locale \isn{ngkit}: if the three premises of a coinduction
step are valid at $\Theta^*_{\mathrm c}$, so is its conclusion.
\end{thm}

The proof follows the coinduction step for $\Mstar$
(\isn{nu\_step\_mjt}, step~3 of the proof of
\cref{thm:app-inga-sound}).  Two things are new: the region's points
must pass the filter, and the universals of the body must carry
certificates for the $\nu$-term rather than for the rule's
certificate.

\paragraph{Proof.}
Invert the step to obtain $\Phi$, $x$, $c$ and the eigenvariable $i$;
fix a well-formed $\sigma$ satisfying the hypotheses $H$, and write
$\Phi_\sigma$, $c_\sigma$, $P = \lambda c_\sigma$ and $x_\sigma$ as in
that proof.
\begin{enumerate}
\item \emph{The region.}
  \[
    E = \{u \mid u \text{ well-formed, with a value, and } P\cdot u
    \text{ true at } \Theta^*_{\mathrm c}\}.
  \]
  The set of valued points at which the certificate is true is the
  whole region; no union with the oracle's truth side is needed.
  $x_\sigma \in E$ by the validity of the first and third premises,
  with one beta step from $c_\sigma[x_\sigma]$ to $P\cdot x_\sigma$.
\item \emph{The filter} (\isn{ng\_rep}).  For $u \in E$, completeness
  (\cref{thm:app-ng-complete}) gives $\ent_S u\,\mathbf N$ and $\ent_S
  P\cdot u$.  The closure premise $i\,\mathbf N,\ c[i],\ \Gamma \ent
  \Phi[\lambda c/a][i/x]$ is instantiated by $\sigma$ with the
  eigenvariable renamed to one fresh for the target context, the body
  and the certificate.  The $\sigma$-instances of the members of
  $\Gamma$ are true, hence derivable by completeness, and cut removes
  them.  One coinduction step (\isn{s\_nu}, the coinduction step of
  the locale \isn{gstep\_nu}) with the three premises at $u$ gives
  $\ent_S \nu\Phi_\sigma\cdot u$.  The lemma is stated at every
  well-formed context $\Delta$ with $\Delta \ent u\,\mathbf N$ and
  $\Delta \ent P\cdot u$, which step~4 uses.
\item \emph{The body's truth.}  The closure premise's validity at
  $\sigma[i \mapsto u]$ makes $\Phi_\sigma[P/a][u/x]$ true at
  $\Theta^*_{\mathrm c}$, as in the proof for $\Mstar$.
\item \emph{Positive replacement} (\isn{ng\_posrep}).  An induction on
  \isn{posP}/\isn{negP} reads that truth directly into the region:
  $\mathrm{rd}^+_{\Theta^*_{\mathrm c},\mathrm c}(E, 1, \rho, \Phi_\sigma,
  1)$, with $\rho$ sending $0$ to $u$ and $1$ to $\nu\Phi_\sigma$.  At
  a recursive atom, the atom's argument is a numeral expression in the
  variables below $1$ and so has a value (\isn{pnexp\_ngval}), and the
  truth of $P\cdot v$ there is the remaining membership condition of
  $E$.  The new cases are the universals of the body.  The reading
  asks, at a true universal, for a derivation of the open instance
  with $\nu\Phi_\sigma$ in place of $a$, and the truth of
  $\Phi_\sigma[P/a][u/x]$ supplies one with $P$ in place of $a$.
  \isn{ng\_cert\_allT} moves the certificate: it is \iga's lemma
  \isn{gwgf\_cert\_allT} at environments of provable naturals
  (\isn{ngenv}), and it needs, at every context $\Delta$ that such an
  environment produces, the rule ``from $\Delta \ent v\,\mathbf N$ and
  $\Delta \ent P\cdot v$, $\Delta \ent \nu\Phi_\sigma\cdot v$''.  That
  rule is step~2 at $\Delta$.  \isn{ng\_cert\_allF} moves the
  refutation of an instance at a false universal in a negative
  position, dually.
\item \emph{Close.}  The environment $\rho$ is
  $\mathrm{env}_{\Phi_\sigma}(u)$, an equation of substitutions, so
  steps 2 and 4 give $u \in F^{\Theta^*_{\mathrm c}}_{\Phi_\sigma,\mathrm
  c}(E)$.  So $E$ is a post-fixed point, \isn{gfp\_upperbound} gives $E
  \subseteq \Theta^*_{\mathrm c}(\Phi_\sigma, 1)$, and \isn{gnosr\_nuT}
  makes $\nu\Phi_\sigma\cdot x_\sigma$ true.
\end{enumerate}
The theorem is proved in \isn{ngkit}, whose one assumption is
completeness; \cref{thm:app-ng-complete} discharges it, so the
theorem holds at \inga without hypothesis.

\subsection{Substitution, value congruence and head steps}

The rules of \rga need the facts that the interface of
Appendix~\ref{app:nu} supplied for $\Mstar$: invariance under
substitution, under value congruence and under head steps.  For the
certified evaluation they are proved over an oracle $\Theta_z$ in three
nested locales, the pattern of \isn{cohO}, \isn{cohO2} and
\isn{cohO3}: \isn{ngco} ($\Theta_z$ consistent, certified and coherent;
it holds at $\Theta^*_{\mathrm c}$ for every calculus,
\isn{ngco\_OstarG\_mf}, by \cref{thm:app-ng-oracle} and coherence),
\isn{ngsubst} (adding \isn{Oz\_osub}: $\Theta_z$ is closed under
substitution) and \isn{ngcong} (adding \isn{Oz\_ocg}: $\Theta_z$ is
closed under value congruence).

\begin{thm}[\isn{gnosr\_Oz\_subf}, \isn{gnosr\_Oz\_vcr}]
\label{thm:app-ng-kit}
For well-formed terms: in \isn{ngsubst}, a value of a term in the
certified evaluation at $\Theta_z$ is a value of its instance under a
well-formed substitution (heads: \isn{gnhev\_Oz\_subf}); in
\isn{ngcong}, a value is a value of every term congruent to it, the
congruence \isn{ngvcrd} being \rga's value congruence (\isn{gvcrd})
with values taken in the certified evaluation at $\Theta_z$ (heads up
to the relation \isn{nhrelG}: \isn{gnhev\_Oz\_vcr}; the oracle's
membership along the congruence: \isn{Oz\_vcr\_member}); and head
steps preserve values in both directions (\isn{ng\_hstep\_sim},
\isn{ng\_hstep\_rev}, \isn{ng\_nu\_unfold}).
\end{thm}

\paragraph{Proof.}
\begin{enumerate}
\item The proofs replay those for $\Mstar$ (step~5 of the proof of
  \cref{thm:app-rga-nu-sound}) and \rga's for its evaluation, with one
  addition at every oracle clause: the image of a point must carry a
  derivation.  The filter's derivations are moved by the system's own
  rules: along a substitution by the substitution rule
  (\isn{ng\_filter\_subf\_T}, \isn{ng\_filter\_subf\_F}), and along a
  value congruence by the proof-level transport
  (\isn{ir\_ngvcr\_transport}; \isn{ng\_filter\_vcr\_T},
  \isn{ng\_filter\_vcr\_F}), whose sites are provably equal because
  their common value is derivable at a certified oracle.
\item \emph{The hardest case} (\isn{ng\_tr\_case\_nuTF}): congruence
  at an oracle clause whose argument changes.  As in
  \isn{tr\_case\_nuTF}, the image is $t'\cdot u'$ with $t'$ at a
  $\nu$-head of a congruent body or at the $\lambda$-head of a mixed
  self-application.  The unfolding's value at the image comes from
  the congruence interior (\isn{OcgG}) at the source point.  The
  filter's derivation moves along the congruence
  (\isn{ng\_filter\_vcr\_T}, \isn{ng\_filter\_vcr\_F}, applied to the
  source atom's derivation or refutation), and \isn{Oz\_fold\_cert}
  takes the value and the derivation and puts the image point on the
  oracle's side at the congruent body, from which \isn{gnosr\_nuT} or
  \isn{gnosr\_nuF} gives the value.  At a mixed self-application there
  is no derivation step: the interior's value is read through beta
  (\isn{gnosr\_appLam}).
\item The interiors used in the next subsection, \isn{OsubG} for
  substitution and \isn{OcgG} for congruence, are those of $\Mstar$ in
  certified form; \isn{ng\_subf\_transport} and \isn{ng\_vcr\_transport}
  evaluate at them.
\end{enumerate}
The kit carries well-formedness hypotheses that the kit for $\Mstar$
does not.  The guard costs nothing at $\Theta^*_{\mathrm c}$, whose
points carry derivations and so are well formed; the theory's text
notes the shape a counterexample to the unguarded statement would
take, a beta step that drops an ill-formed argument.

\subsection{The certified oracle is closed}

\begin{thm}[\isn{OstarG\_osub}, \isn{OstarG\_ocg}]
\label{thm:app-ng-closure}
$\Theta^*_{\mathrm c}$ is closed under substitution and under value
congruence.  Hence the locale \isn{ngstar}, \isn{ngco} with $\Theta_z
= \Theta^*_{\mathrm c}$, is a sublocale of \isn{ngsubst} and of
\isn{ngcong}, and it holds at every calculus
(\isn{ngstar\_OstarG\_mf}).
\end{thm}

\paragraph{Proof.}
Both closures are pre-fixed-point arguments, as for $\Mstar$
(step~5(e)--(f) of the proof of \cref{thm:app-rga-nu-sound}): an
interior $\Theta^\circ \le \Theta^*_{\mathrm c}$ with
$\mathrm{Om}_{\mathrm c}(\Theta^\circ) \le \Theta^\circ$ contains
$\Theta^*_{\mathrm c}$ by \isn{lfp\_lowerbound}.
\begin{enumerate}
\item \emph{The readings move.}  \isn{gosem\_subf} and \isn{gfsem\_subf}
  (in \isn{ngco}) move the two readings along a well-formed
  substitution, the analogue of \isn{osemO\_subf} for $\Mstar$ with a
  well-formedness guard; \isn{gosem\_vcr} and \isn{gfsem\_vcr} (in
  \isn{ngsubst}) move them along the congruence, with environments
  related by \isn{envrelG} (\isn{envrel} with well-formed
  substitutions), the universals' certificates moved by the
  substitution rule and by \isn{ir\_ngvcr\_transport}.
\item \emph{Substitution.}  \isn{OsubG\_prefix}: $\mathrm{Om}_{\mathrm
  c}(\isn{OsubG}) \le \isn{OsubG}$, by step~1 and the filter's
  transport along substitution.  So $\Theta^*_{\mathrm c} \le
  \isn{OsubG}$, which is \isn{OstarG\_osub}, and \isn{ngstar} is a
  sublocale of \isn{ngsubst}.
\item \emph{Congruence, truth side} (\isn{OcgG\_gfp\_cong}): a
  coinduction whose region is the set of points congruent to a point
  of the interior's greatest fixed point, as \isn{Ocg\_gfp\_cong}.
\item \emph{Congruence, falsity side} (\isn{OcgG\_lfp\_cong}): a
  least-fixed-point induction, as for $\Mstar$, with one new
  obligation: the image point needs a refutation, which is the
  interior member's own refutation moved along the congruence
  (\isn{ng\_filter\_vcr\_F}).
\item \isn{OcgG\_prefix} combines steps 3 and 4, so $\Theta^*_{\mathrm
  c} \le \isn{OcgG}$, which is \isn{OstarG\_ocg}, and \isn{ngstar} is a
  sublocale of \isn{ngcong}.
\end{enumerate}
At $\Theta^*_{\mathrm c}$ the kit's facts hold without hypothesis on
the oracle (stated with the calculus's consistency as hypothesis in
\isn{gnosr\_OstarG\_subf} and \isn{gnosr\_OstarG\_vcr}, and at every
calculus through \isn{ngstar\_OstarG\_mf}), and at \inga without
hypothesis at all (\isn{inga\_gnosr\_subf},
\isn{inga\_gnosr\_vcr}, \isn{inga\_OstarG\_osub}, \isn{inga\_OstarG\_ocg}).

\subsection{The rules of \rga and ground induction}

\begin{thm}[\isn{ng\_rga\_step\_gnowsjt}, \isn{ng\_gfix\_step\_gnowsjt}]
\label{thm:app-ng-step}
In \isn{ngstar}: every \rga step, each of the forty-two rules,
number induction among them, preserves validity at
$\Theta^*_{\mathrm c}$.  In \isn{ngstar\_gfix}, \isn{ngstar} with the
calculus's ground-induction step (\isn{gsem\_gfix}): so does every
ground-induction step.  Since \isn{ngstar} holds at every calculus
(\cref{thm:app-ng-closure}), both hold at $\Theta^*_{\mathrm c}$ of
every calculus with the steps in question (stated with the calculus's
consistency as hypothesis in \isn{ng\_rga\_step\_OstarG} and
\isn{ng\_gfix\_step\_OstarG}), and at \inga without hypothesis
(\isn{inga\_rga\_step\_gnowsjt}, \isn{inga\_gfix\_step\_gnowsjt}).
\end{thm}

\paragraph{Proof.}
\begin{enumerate}
\item \emph{\rga's rules}, case by case, as in \rga's soundness for its
  certificate-grounded evaluation.  Five kinds of case use the new
  evaluation.
  \begin{enumerate}
  \item \emph{Universal introduction.}  The conclusion's truth needs a
    certificate for the $\sigma$-instance of the universal.  The
    premise is a derivation under the hypotheses $H$ at the
    eigenvariable; the proof surgery \isn{gS.ir\_subf1} substitutes
    $\sigma$ into it, the substituted hypotheses are true and so
    derivable by completeness at the certified oracle
    (\isn{ng\_complete\_open}), and cut discharges them.
  \item \emph{Conversions}, through the head steps in both directions
    (\isn{ng\_hstep\_gnosr}, \isn{ng\_hstep\_rev\_gnosr}), with the
    redex's well-formedness from \isn{hstep\_wft\_rev}.
  \item \emph{The equality rules}, by the value congruence
    (\isn{gnosr\_Oz\_vcr}).
  \item \emph{Explosion}: no substitution makes both $a$ and $\neg a$
    true, by determinism (\isn{gnosr\_Oz\_unique}).
  \item \emph{Refuting a universal}: the falsity clause of the
    universal asks for a refutation of the false instance, and
    completeness at the certified oracle supplies it
    (\isn{ng\_complete\_open}); \isn{gnosr\_allF} closes the case.
  \end{enumerate}
\item \emph{Ground induction} follows the proof for $\Mstar$ (step~2 of
  the proof of \cref{thm:app-inga-sound}).  \isn{ngf\_transport} is one
  rule induction on the certified evaluation with the five-part motive
  of the \iga paper's transport lemma (\isn{ngfA}, \isn{ngfPP},
  \isn{ngfNP}, \isn{ngfRPP}, \isn{ngfRNP}, together \isn{ngfInv}) over
  \iga's environments (\isn{wgfenv}), with the head facts
  \isn{gnhev\_trmGfix} and \isn{gnhev\_gfix\_head}.  The two new cases
  are the oracle clauses, and they are vacuous: a $\nu$-head is not
  the head of the $\mu$-atom, not $\cA$'s and not a $\lambda$-head
  (\isn{ngf\_spine}).  \isn{ngf\_ind\_par} is the
  resulting induction property, and \isn{ng\_gfix\_step\_gnowsjt} the
  step.
\end{enumerate}

\subsection{Soundness, the two properties and the instances}

The locale \isn{ngsound} is \isn{ngstar} with a further step relation
$X$, the split \isn{ng\_S\_split} (every $S$ step is an \rga step or an
$X$ step) and \isn{ng\_X\_sound} (every $X$ step preserves validity at
$\Theta^*_{\mathrm c}$).

\begin{thm}[\isn{gnsem\_sound}, \isn{ng\_open\_iff}]
\label{thm:app-ng-sound}
In \isn{ngsound}:
(i) every judgment of a valid proof list for $S$ from a well-formed
context is valid at $\Theta^*_{\mathrm c}$ (\isn{prf\_gnowsjt},
\isn{gnsem\_sound});
(ii) a derivable term is true and a refutable term false in the
certified evaluation at $\Theta^*_{\mathrm c}$
(\isn{ng\_provable\_true}, \isn{ng\_refutable\_false}), and with
completeness a well-formed term is derivable exactly when it is true,
and refutable exactly when it is false (\isn{ng\_open\_iff},
\isn{ng\_false\_iff});
(iii) at every body $\Phi$ of the class, the certified oracle consists
exactly of the derivable and the refutable $\nu$-atoms
(\isn{ng\_nu\_true\_iff}, \isn{ng\_nu\_false\_iff});
(iv) the disjunction property and the witness property for refuted
universals hold for $S$ (\isn{ng\_disj\_prop}, \isn{ng\_witness\_prop});
(v) a sentence whose decidedness $S$ derives is derivable or refutable
(\isn{ng\_no\_fict}; for $\nu$-atoms, \isn{ng\_nu\_no\_fict});
(vi) $S$ does not derive $0 = 1$ from no hypotheses
(\isn{ng\_consistent}).
\end{thm}

\paragraph{Proof.}
\begin{enumerate}
\item (i) is an induction on the proof list, the template of the \iga
  paper's \isn{prf\_gowsjt}: a new entry is an \rga step
  (\cref{thm:app-ng-step}) or an $X$ step (\isn{ng\_X\_sound}).
\item (ii): (i) at the empty context and the identity substitution;
  the converse is \cref{thm:app-ng-complete}.
\item (iii): a point of the oracle is a derivable or refutable atom by
  \cref{thm:app-ng-oracle}; conversely, a derivable atom is true by
  (ii), and an atom of the class is true only through the oracle.
\item (iv): a derivable $a \vee b$ is true, so a disjunct is true and,
  by completeness, derivable.  A refutable $\cA\cdot u$ is false, and
  the falsity clause of the universal carries a derivation of the
  refuted instance; this half needs no completeness.
\item (v) is the disjunction property at $c \vee \neg c$.
\item (vi): a derivable $0 = 1$ would be true by (ii).  An equation
  is true only when its two sides have equal values
  (\isn{gnosr\_eq\_inv}), and $0$ has the value $0$ while $1$ has a
  successor value.  Consistency follows from soundness, with no other
  input.
\end{enumerate}

\begin{thm}[\isn{inga\_pg}, \isn{inga\_pg\_not\_omega}]
\label{thm:app-inga-pg}
(i) \inga is an instance of \isn{ngsound}, with $X$ the
ground-induction step or the coinduction step, so every clause of
\cref{thm:app-ng-sound} holds at \inga (\isn{inga\_gnsem\_sound},
\isn{inga\_provable\_true}, \isn{inga\_gn\_open\_iff},
\isn{inga\_gn\_false\_iff}, \isn{inga\_nu\_true\_iff},
\isn{inga\_nu\_false\_iff}, \isn{inga\_disj\_prop},
\isn{inga\_witness\_prop}, \isn{inga\_no\_fict},
\isn{inga\_nu\_no\_fict}), and \inga is consistent
(\isn{inga\_consistent\_mf}).
(ii) With $\mathrm{Tr}$ and $\mathrm{Fa}$ as in
\cref{thm:app-inga-reduce}, the interface holds (\isn{inga\_pg};
\isn{inga\_cert\_grounded}).
(iii) \inga is proof-grounded in the sense of \S\ref{sec:price}, is not
sound for the least model, and is not an instance of
\isn{gsem\_sound} (\isn{inga\_pg\_not\_omega}).
(iv) \rga with the coinduction rule is an instance of \isn{ngsound}
with $X$ the coinduction step (\isn{rganu\_gnsem\_sound},
\isn{rganu\_nu\_true\_iff}, \isn{rganu\_disj\_prop},
\isn{rganu\_witness\_prop}), and it is consistent
(\isn{rganu\_consistent\_mf}).
\end{thm}

\paragraph{Proof.}
\begin{enumerate}
\item The sublocale $\isn{inga} \subseteq \isn{ngsound}$ (prefix
  \isn{ings}): \isn{ngstar} holds with no hypothesis
  (\isn{ngstar\_OstarG\_mf}); the split is the
  definition of \isn{inga\_derives}; $X$ is sound by
  \cref{thm:app-ng-step} for ground induction and
  \cref{thm:app-ng-park} for the coinduction step.
\item (ii) is (i)'s two properties with \cref{thm:app-inga-reduce}.
\item (iii) adds \isn{inga\_omega\_unsound\_mf}: \inga derives
  $Z\cdot 0$ (\cref{thm:app-ingaZ}), which has no value in the least
  model (\isn{inga\_Z\_valueless}); and \cref{thm:app-inga-notgsem}.
\item (iv): the sublocale $\isn{inga} \subseteq \isn{ngsound}$ at
  \isn{rganu\_derives} (prefix \isn{rgnu}), the split by definition,
  the coinduction step sound by \cref{thm:app-ng-park}.  \rga with the
  coinduction rule is consistent by its own soundness
  (\cref{thm:app-ng-sound}(vi)), with no appeal to \inga's
  consistency.
\end{enumerate}
No step of the construction assumes consistency.  The certified oracle
is consistent outright (\cref{thm:app-ng-oracle}), so the locale
\isn{ngstar} holds at every calculus, and consistency is an output of
soundness (\cref{thm:app-ng-sound}(vi)), from which exclusion in
\cref{thm:app-inga-reduce} takes it.  \inga's consistency thus has a
proof from its certified evaluation alone (\isn{inga\_consistent\_mf}),
beside the proof through $\Mstar$ (\isn{inga\_consistent},
\cref{thm:app-inga-sound}); the headlines of this appendix rest on no
fact about $\Mstar$.  The forms with the calculus's consistency
as hypothesis (\isn{ocons\_OstarG\_filter}, \isn{ngco\_OstarG},
\isn{ngstar\_OstarG} and the \isn{\_con} forms of \isn{INGA\_PG0}) are
kept in the development.

\paragraph{Lemmas a re-mechanization needs.}
\begin{itemize}
\item \isn{ng\_gs\_omega\_true}, \isn{inga\_Z\_valueless}: the locale
  of \rga's certificate semantics forces soundness for the least model,
  which \inga lacks.
\item \isn{inga\_cert\_allT}, \isn{inga\_cert\_excl},
  \isn{inga\_cert\_iff\_props}: the interface at derivability reduces
  to the two properties.
\item \isn{gnosr\_nuT}, \isn{gnosr\_nuF}, \isn{FOg}, \isn{GOg},
  \isn{OmG}, \isn{OstarG}: the evaluation reads both sides of a
  filtered oracle.
\item \isn{OstarG\_T\_derivable}, \isn{OstarG\_F\_refutable},
  \isn{ocert\_OstarG}, \isn{gnosr\_unique}: the certified oracle is
  certified by its filter; determinism at a consistent oracle.
\item \isn{gosem\_gfsem\_excl}, \isn{gfpG\_lfpG\_disjoint},
  \isn{ocons\_OmG}, \isn{ocons\_admissible}, \isn{ocons\_OstarG}: the
  certified oracle is consistent, with no hypothesis.
\item \isn{gnhev\_gnosr\_complete}, \isn{ir\_true\_eqI},
  \isn{ir\_false\_eqI}, \isn{ng\_complete}: completeness at a certified
  oracle.
\item \isn{ng\_coherent\_T\_iff}, \isn{ng\_coherent\_F\_iff}: coherence,
  backward through completeness and folding.
\item \isn{ng\_rep}, \isn{ng\_posrep}, \isn{ng\_cert\_allT},
  \isn{ng\_cert\_allF}, \isn{nu\_step\_gnowsjt}: the semantic Park lemma.
\item \isn{ng\_filter\_subf\_T}, \isn{ng\_filter\_vcr\_F},
  \isn{ir\_ngvcr\_transport}, \isn{ng\_tr\_case\_nuTF},
  \isn{Oz\_fold\_cert}: the filter moved along substitution and
  congruence; the hardest congruence case.
\item \isn{gosem\_subf}, \isn{gosem\_vcr}, \isn{OsubG\_prefix},
  \isn{OcgG\_lfp\_cong}, \isn{OcgG\_prefix}: the certified oracle is
  closed.
\item \isn{gS.ir\_subf1}, \isn{ng\_complete\_open},
  \isn{ng\_hstep\_gnosr}, \isn{gnosr\_Oz\_unique},
  \isn{ng\_rga\_step\_gnowsjt}: the rules of \rga.
\item \isn{ngf\_transport}, \isn{ngf\_spine}, \isn{ngf\_ind\_par},
  \isn{ng\_gfix\_step\_gnowsjt}: ground induction.
\item \isn{prf\_gnowsjt}, \isn{ng\_open\_iff}, \isn{ng\_disj\_prop},
  \isn{ng\_witness\_prop}, \isn{ng\_consistent}: soundness, the two
  properties and consistency.
\end{itemize}

\section{Streams, runs and branches}
\label{app:streams}

Files: \isn{OIGA\_Coind/GStream.thy} (prefix sets, the round trip,
the prefix tree, runs), \isn{OIGA\_Coind/GBranch.thy} (the run
from a path and the branch), \isn{OGA/OGA\_Coind\_Kit.thy} (runs, the
child relation, the replay, the decomposition \isn{fb\_surj}), \isn{OGA/OGA\_Koenig\_FB.thy} (sequence
codes), \isn{OGA/OGA\_Shell.thy} (membership, the search combinator,
\isn{srch\_from\_ex}), \isn{OIGA/GShell.thy} (the search lemma
\isn{mu\_choice}), \isn{OIGA/GDich.thy} (the dichotomy at the
accessibility body, the decided path atom, the negative fold),
\isn{OIGA/GNu.thy} (the unfolding) and \isn{OIGA/GPar.thy}
(ground induction with parameters).

\paragraph{Definitions: sequence codes and prefix sets.}
A finite sequence is coded as a number.  The empty sequence is $1$;
appending $a$ to $\sigma$ gives $\sigma{\frown}a = S(S\langle\sigma,
a\rangle)$ with Cantor pairing (\isn{fbSnoc}).  The parent, the last
element and the $j$-fold parent of a code are compiled primitive
recursive functions, written $\mathrm{par}$, $\mathrm{last}$ and
$\mathrm{drop}$ (\isn{fbPar}, \isn{fbLast}, \isn{fbDrop}), and a
code $\sigma$ \emph{has length} $n$ when $\mathrm{drop}(\sigma, n) =
1$.  The code $0$ is junk: its parent and every drop of it are $0$, so
it has no length (\isn{len\_zero\_absurd}).  A code of length $k+1$
is the parent with the last element appended (\isn{fb\_surj}).  For a
stream $s$, the prefix function $\mathrm{pfx}_s$ (\isn{pfx}) is
defined by primitive recursion from $1$: $\mathrm{pfx}_s\cdot 0 = 1$
and $\mathrm{pfx}_s\cdot(m+1) = (\mathrm{pfx}_s\cdot m){\frown}(s\cdot
m)$ (\isn{pfx\_0}, \isn{pfx\_S}; their premise is the totality of the
recursion's step function, which follows from $\Tot\,s$ at a
well-formed context, \isn{pfxStep\_Tot2}; \isn{pfx\_S} also assumes
$\Tot\,s$).  The prefix set
$\mathrm{pset}_s$ (\isn{pset}) is $\lambda\sigma.\,\exists n.\,\sigma =
\mathrm{pfx}_s\cdot n$.  For a set term $p$, $\sigma \in p$
abbreviates $p\cdot\sigma = 1$ (\isn{Mem}).  A set term $p$ is a
\emph{branch}, written $\mathrm{Branch}(p)$ (\isn{fbBranch}), when
five conditions hold, as one conjunction: $\Dec\,p$; \emph{weak
prefix closure} (\isn{fbPrefClosedS}), that every drop of a member is
$0$ or a member; \emph{a node at every length} (\isn{fbInfTree}), that
for every $n$ some member has length $n$; \emph{linearity}
(\isn{fbLin}), that two members of the same length are equal; and
\emph{codes only} (\isn{fbSeqs}), that every member has a length.  The
stream recovered from $p$ (\isn{spr}) is $\mathrm{spr}_p\cdot n =
\mathrm{last}(\mu_p\cdot n)$, where $\mu_p$ is the search combinator
at the matrix ``$v$ is a member of $p$ of length $n+1$''
(\isn{sprPhi}, \isn{muF}): it returns the least number satisfying the
matrix, searching from $0$.

\paragraph{Definitions: trees, runs and the replay.}
The child relation of a set $T$ (\isn{fbTreeRel}) is $C_T(y, x) \equiv
y \in T \wedge \mathrm{par}(y) = x$; only the child is required to be a
member.  $\Path_{C_T}$ is the $\nu$-atom of $\exists y.\,C_T(y, x)
\wedge a\cdot y$.  A \emph{run} of $R$, written $\mathrm{run}_R(s)$
(\isn{stRun}), is $\Tot\,s \wedge
\forall n.\,R(s\cdot(n+1), s\cdot n)$, in the orientation of the
accessibility body; the image of $s$ (\isn{stImg}) is $\lambda\sigma.\,
\exists k.\,\sigma = s\cdot k$; and $\mathrm{every}_Q(s)$
(\isn{stEvery}) is $\forall n.\,Q(s\cdot n)$.  The \emph{replay}
(\isn{gPath}) at a relation $C$ and a set $Q$ iterates one search
step from a start $a$ by primitive recursion: the step at $x$
(\isn{gStepBody}) is, if $Q(x)$, the least $y$ with $C(y,x) \wedge
Q(y)$, and otherwise $0$.  \emph{Progress} (\isn{gProg}) is $\forall
x.\,\neg Q(x) \vee \exists y.\,C(y,x) \wedge Q(y)$.  For a branch term
$b$: $\mathrm{in}_T(b)$ is $\forall n.\,b\cdot n \in T$
(\isn{fbPathIn}); $\mathrm{ext}(b)$ is $\forall n.\,\mathrm{par}(b\cdot
(n+1)) = b\cdot n$ (\isn{fbPathExt}); $\mathrm{len}(b)$ is $\forall
n.\,b\cdot n$ has length $n$ (\isn{fbPathLen}).

\begin{thm}[\isn{g\_RT\_s}, \isn{g\_RT\_p}; grounded]
\label{thm:app-roundtrip}
In the calculus \isn{gcalc\_ati}, with $\Gamma$ well-formed:
(i) if $\Gamma \ent \Tot\,s$, then $\Gamma \ent \forall n.\,
\mathrm{spr}_{\mathrm{pset}_s}\cdot n = s\cdot n$;
(ii) if $\Gamma \ent \mathrm{Branch}(p)$, then $\Gamma \ent \forall
\sigma.\,(\sigma \in \mathrm{pset}_{\mathrm{spr}_p} \to \sigma \in p)
\wedge (\sigma \in p \to \sigma \in \mathrm{pset}_{\mathrm{spr}_p})$.
The companion \isn{g\_pset\_Branch}, in \isn{gcalc\_ati}: if $\Gamma
\ent \Tot\,s$, then $\Gamma \ent \mathrm{Branch}(\mathrm{pset}_s)$.
\end{thm}

\paragraph{Proof of \cref{thm:app-roundtrip}.}
The prefix-function kit is proved in \rga and carried into every
calculus by cut: $\mathrm{pfx}_s$ is total; the parent of
$\mathrm{pfx}_s\cdot(m+1)$ is $\mathrm{pfx}_s\cdot m$
(\isn{pfx\_par\_S}); $\mathrm{pfx}_s\cdot n$ has length $n$
(\isn{pfx\_len}, number induction on $n$ with motive ``$\mathrm{pfx}_s
\cdot i$ has length $i$''); and that length is unique
(\isn{g\_pfx\_len\_uniq}: if $\mathrm{pfx}_s\cdot n$ has length $m$
then $m = n$, by number induction on $n$ with motive ``for every $m$,
$\mathrm{pfx}_s\cdot n$ does not have length $m$, or $m = n$'').
\begin{enumerate}
\item \emph{The prefix set is decided} (\isn{g\_pset\_Dec}).  It
  suffices, at every number $t$ of every extension, that $\exists
  n.\,t = \mathrm{pfx}_s\cdot n$ is decided (\isn{g\_Dec\_lamI}).  At
  a fresh number variable $v$, the equation $t = \mathrm{pfx}_s\cdot v$
  is decided, both sides being numbers, and so is its negation; ATI in
  its open form (\isn{ir\_allTI\_open}) decides $\forall v.\,\neg(t =
  \mathrm{pfx}_s\cdot v)$, and its negation is the existential.  This
  lemma is the only place where the theorem applies ATI; clause (ii)
  calls the lemma again in step 8.
\item \emph{The branch conditions} (\isn{g\_pset\_Branch}).  Decided:
  step 1.  Weak prefix closure (\isn{g\_pset\_PrefClosedS}): split on
  the membership of $\sigma$, decided by step 1; a member is
  $\mathrm{pfx}_s\cdot k$; for a drop count $j$, split $j \le k$ or
  $k < j$; the drop is $\mathrm{pfx}_s\cdot(k - j)$, a member
  (\isn{g\_pfx\_drop\_le}), or $0$ (\isn{g\_pfx\_drop\_gt}).  A node at
  every length: $\mathrm{pfx}_s\cdot n$.  Linearity (\isn{g\_pset\_Lin}):
  four decided case splits, on the membership of $\sigma$ and of $\tau$
  (step 1) and on the two length equations; in the remaining case
  $\sigma = \mathrm{pfx}_s\cdot a$ and $\tau = \mathrm{pfx}_s\cdot b$
  both have length $n$, so uniqueness gives $a = n = b$.  Codes only:
  split on membership (step 1); a member $\mathrm{pfx}_s\cdot k$ has
  length $k$.
\item \emph{The search premises} (\isn{g\_spr\_premises}).  From
  $\Dec\,p$ and a node at every length, at a fresh number variable
  $n$: the matrix ``member of length $n+1$'' is decided, from $\Dec\,p$
  and the decidedness of equality; and it has a witness, the instance
  at $n+1$ of ``a node at every length''.
\item \emph{The search} (\isn{mu\_choice}).  From the two premises of
  step 3: $\Tot\,\mu_p$, $\mu_p\cdot n$ equals the \emph{search value},
  the result of the least-number search from $0$ at the matrix for $n$
  (\isn{srchK}), and the search value satisfies the matrix; \isn{g\_mu\_choice\_all}
  packages the last clause as $\forall n.\,\mu_p\cdot n$ is a member
  of $p$ of length $n+1$.  The lemma is stated in \isn{gcalc}, the base calculus: its content is \rga's \isn{srch\_from\_ex} (a decided matrix with
  a witness makes the search from $0$ a number satisfying the matrix),
  whose two clauses the calculus cuts in directly by \isn{g\_via2}.  No ATI step occurs in that derivation, and
  no minimality of the search value is used anywhere in the theorem:
  uniqueness comes from \isn{g\_pfx\_len\_uniq} and from linearity.
\item \emph{Clause (i).}  Fix $n$ and let $c = \mu_{\mathrm{pset}_s}
  \cdot n$.  Step 4 at $p = \mathrm{pset}_s$, with step 1 for
  decidedness and \isn{g\_pset\_InfTree} (a node at every length) for
  the witness, makes $c$ a member of length $n+1$, so $c =
  \mathrm{pfx}_s\cdot j$ for some $j$ (\isn{g\_pset\_memE}), and
  uniqueness of length gives $j = n+1$.  Then $c = (\mathrm{pfx}_s\cdot
  n){\frown}(s\cdot n)$, whose last element is $s\cdot n$; one beta step
  (\isn{spr\_beta}) turns $\mathrm{spr}_{\mathrm{pset}_s}\cdot n$ into
  $\mathrm{last}(c)$.
\item \emph{The bridge} (\isn{g\_RTp\_bridge}).  Assume the branch
  conditions of $p$ and write $M = \mu_p$.  Then $M\cdot(k+1) =
  (M\cdot k){\frown}\mathrm{last}(M\cdot(k+1))$.  The code $x =
  M\cdot(k+1)$ has length $k+2$, so $x = \mathrm{par}(x){\frown}
  \mathrm{last}(x)$ (\isn{fb\_surj}); $\mathrm{par}(x)$ has length
  $k+1$ and is a member by weak prefix closure (the junk case is
  excluded by the length); $M\cdot k$ is a member of length $k+1$; and
  linearity identifies $\mathrm{par}(x)$ with $M\cdot k$ (\isn{g\_fbLinE}).
\item \emph{The prefixes of the recovered stream} (\isn{g\_RTp\_pfx}).
  $\mathrm{pfx}_{\mathrm{spr}_p}\cdot(k+1) = M\cdot k$, by number
  induction on $k$ with motive ``$\mathrm{pfx}_{\mathrm{spr}_p}\cdot
  (i+1) = M\cdot i$''.  Base: $M\cdot 0$ has length $1$, so it is
  $1{\frown}\mathrm{last}(M\cdot 0)$, which is
  $\mathrm{pfx}_{\mathrm{spr}_p}\cdot 1$.  Step (\isn{g\_RTp\_step}):
  unfold $\mathrm{pfx}$ once, rewrite with the hypothesis and
  \isn{spr\_beta}, and close with step 6.
\item \emph{First inclusion} (\isn{g\_RTp\_sub}).  At a number
  $\sigma$, split on $\sigma \in \mathrm{pset}_{\mathrm{spr}_p}$, which
  is decided by step 1 at the total stream $\mathrm{spr}_p$
  (\isn{g\_spr\_Tot}); this case split spends ATI.  A member is
  $\mathrm{pfx}_{\mathrm{spr}_p}\cdot k$; if $k = 0$ it is $1$, the
  drop of $M\cdot 0$ by one, a member by weak prefix closure; if $k =
  q+1$ it is $M\cdot q$ by step 7, a member by step 4.
\item \emph{Second inclusion} (\isn{g\_RTp\_sup}).  Split on $\sigma \in
  p$, decided by $\Dec\,p$ (no ATI).  A member has some length $n$
  (codes only).  If $n = 0$, $\sigma = 1 = \mathrm{pfx}_{\mathrm{spr}_p}
  \cdot 0$.  If $n = k+1$, $\sigma$ and $M\cdot k$ are members of the
  same length, so $\sigma = M\cdot k = \mathrm{pfx}_{\mathrm{spr}_p}
  \cdot(k+1)$.
\item \emph{Clause (ii).}  Universal introduction, the branch
  conditions unpacked (\isn{g\_fbBranchE}), and the conjunction of
  steps 8 and 9.
\end{enumerate}
The stronger prefix closure of the \oga development (\isn{fbPrefClosed},
every drop a member) fails at every prefix set: the drop of the root
$1$ is the junk code $0$ (\isn{g\_pset\_not\_fbPrefClosed}).  This is why
a branch is defined with the weak form.

\begin{thm}[\isn{g\_pset\_path}, \isn{g\_path\_pset\_iff}; grounded]
\label{thm:app-psetpath}
Let $\Gamma$ be well-formed, $\Gamma \ent \Tot\,s$, and write $P =
\mathrm{pset}_s$.
(i) In \isn{gcalc\_nu}: if $\Gamma \ent P\cdot x$ and $\Gamma \ent
x\,\mathbf N$, then $\Gamma \ent \Path_{C_P}\cdot x$.
(ii) In \isn{gcalc\_ati\_nu}: if $\Gamma \ent \sigma\,\mathbf N$,
$\Gamma \ent n\,\mathbf N$ and $\Gamma \ent \mathrm{drop}(\sigma, n) =
1$, then $\Gamma \ent (\Path_{C_P}\cdot\sigma \to P\cdot\sigma) \wedge
(P\cdot\sigma \to \Path_{C_P}\cdot\sigma)$.
\end{thm}

\paragraph{Proof of \cref{thm:app-psetpath}.}
\begin{enumerate}
\item \emph{Coinduction at a certificate set} (\isn{g\_path\_coind}, in
  \isn{gcalc\_nu}).  For a set term $Q$ and a relation $R$: from
  $\Gamma \ent Q\cdot x$, $\Gamma \ent x\,\mathbf N$, and the closure
  obligation ``at every $j$ not free in $\Gamma$, $j\,\mathbf N,\
  Q\cdot j,\ \Gamma \ent \exists y.\,R(y, j) \wedge Q\cdot y$'', conclude
  $\Gamma \ent \Path_R\cdot x$.  This lemma is one application of the
  coinduction rule at the body of $\Path_R$, with certificate $c(v) =
  Q\cdot v$ and an eigenvariable fresh for $\Gamma$, $\lambda c$ and the
  body at once.
\item \emph{Clause (i).}  Step 1 at $Q = P$ and $R = C_P$.  The region
  is the prefix set.  Closure: from $P\cdot j$, eliminate the
  existential, $j = \mathrm{pfx}_s\cdot k$; the child $w =
  \mathrm{pfx}_s\cdot(k+1)$ is in $P$ and its parent is $j$
  (\isn{g\_pfx\_par\_S}), so $C_P(w, j)$ (\isn{g\_TR\_I}), and
  \isn{g\_path\_child\_intro} introduces the existential of the body.
\item \emph{Forward half, no rule} (\isn{g\_path\_pset}, in \isn{gcalc}).
  Unfold $\Path_{C_P}\cdot\sigma$ once (\isn{g\_nu\_unfold}) and
  eliminate the existential: some $j$ with $C_P(j, \sigma)$.  Then
  \isn{g\_TR\_child\_pset}: a code with a length whose child is in $P$
  is itself in $P$.  The child is $\mathrm{pfx}_s\cdot k$; if $k = 0$
  the child is $1$, so $\sigma$, its parent, is the junk code $0$, which
  has no length, against the hypothesis on $\sigma$; if $k = m+1$, then
  $\sigma = \mathrm{pfx}_s\cdot m$ (\isn{g\_pfx\_par\_S}).
\item \emph{Non-members have no infinite branch}
  (\isn{g\_not\_pset\_not\_path}, in \isn{gcalc\_ati}).  If $\neg P\cdot
  \sigma$, then for every $y$ the atom $C_P(y,\sigma)$ is decided
  (\isn{g\_TR\_B}, from $\Dec\,P$ of step 1 of the proof of \cref{thm:app-roundtrip},
  which spends ATI); in the positive case \isn{g\_TR\_child\_pset} puts
  $\sigma$ in $P$, a contradiction.  So $\forall y.\,\neg C_P(y,\sigma)
  \vee \neg\Path_{C_P}\cdot y$, and the negative fold
  (\isn{g\_nu\_fold\_neg}) gives $\neg\Path_{C_P}\cdot\sigma$.
\item \emph{Clause (ii).}  $P\cdot\sigma$ is decided (step 1 of
  \cref{thm:app-roundtrip}).  For the left conjunct, split on $P\cdot\sigma$: the
  negative case gives $\neg\Path$ by step 4, the positive case the right
  disjunct.  For the right conjunct, implication introduction at the
  decided antecedent (\isn{ir\_impI\_bool}) with clause (i).  The same
  case split, producing $\Path \vee \neg\Path$, gives the decided form
  \isn{g\_path\_pset\_dec}.
\end{enumerate}
The interderivability form (\isn{g\_path\_pset\_hol\_iff}, in
\isn{gcalc\_nu}) is steps 2 and 3 as a statement about derivability:
$\Gamma \ent \Path_{C_P}\cdot\sigma$ if and only if $\Gamma \ent
P\cdot\sigma$, at a code with a length.  Its check forbids ATI and the
dichotomy.

\begin{thm}[\isn{g\_run\_acc}, \isn{g\_run\_path}; grounded]
\label{thm:app-runs}
Let $\Gamma$ be well-formed, $R$ a well-formed term and $\Gamma \ent
\mathrm{run}_R(s)$.
(i) In \isn{gcalc\_gfix}: if $\Gamma \ent \Acc_R\cdot x$ and $\Gamma
\ent x\,\mathbf N$, then $\Gamma \ent \forall k.\,\neg(s\cdot k = x)$.
(ii) In \isn{gcalc\_nu}: if $\Gamma \ent m\,\mathbf N$, then $\Gamma
\ent \Path_R\cdot(s\cdot m)$.
\end{thm}

\paragraph{Proof of \cref{thm:app-runs}.}
\begin{enumerate}
\item \emph{Clause (i): the motive.}  $c(z) \equiv \forall k.\,\neg(s
  \cdot k = z)$, with the stream as a parameter (its instances are
  rewritten by \isn{bsub\_runMotive}).
\item \emph{Ground induction} along $\Acc_R$ with parameters
  (\isn{gi\_TI\_par}, one application of \isn{ir\_gfix\_par} at the
  accessibility body, the eigenvariable fresh for the context, the
  motive and $R$).  The closure obligation: in every extension $\Delta$
  and at every number $z$, from $\forall y.\,\neg R(y, z) \vee c(y)$
  derive $c(z)$.
\item \emph{The closure.}  Fix $k$; $s\cdot k = z$ is decided, $s\cdot
  k$ being a number.  If $s\cdot k = z$, the run gives $R(s\cdot(k+1),
  s\cdot k)$, hence $R(s\cdot(k+1), z)$ by transport
  (\isn{g\_app\_transport0}); the hypothesis at $y = s\cdot(k+1)$ gives
  $c(s\cdot(k+1))$, whose instance at $k+1$ is $\neg(s\cdot(k+1) =
  s\cdot(k+1))$, against reflexivity.  Otherwise $\neg(s\cdot k = z)$ is the case hypothesis itself.
  The rule then yields $c(x)$.
\item \emph{Clause (ii)} (\isn{g\_run\_path}).  Step 1 of the
  proof of \cref{thm:app-psetpath} at $Q$ the image of $s$, which is the
  coinduction's region.  Start: $s\cdot m$ is in the image
  (\isn{g\_img\_memI}).  Closure: from $j = s\cdot k$
  (\isn{g\_img\_memE}), the state $s\cdot(k+1)$ is in the image and
  $R(s\cdot(k+1), j)$ by the run and transport; \isn{g\_path\_child\_intro}
  concludes.  One coinduction step.
\end{enumerate}

\begin{thm}[\isn{g\_run\_of\_path}; grounded]
\label{thm:app-runofpath}
In \isn{gcalc\_gfix\_dich}: if $\Gamma$ is well-formed, $\Gamma \ent
\Dec\,R$ (as a relation), $\Gamma \ent \Path_R\cdot x$ and $\Gamma
\ent x\,\mathbf N$, then, writing $P$ for the replay at $C = R$, $Q =
\Path_R$ and start $x$: $\Gamma \ent \Tot\,P$, $\Gamma \ent P\cdot 0
= x$, $\Gamma \ent \mathrm{run}_R(P)$ and $\Gamma \ent
\mathrm{every}_{\Path_R}(P)$.
\end{thm}

\paragraph{Proof of \cref{thm:app-runofpath}.}
\begin{enumerate}
\item \emph{The replay, in \rga} (\isn{gpb\_run\_of\_progress}, carried
  into every calculus as \isn{g\_run\_of\_progress}).  From $\Dec\,C$
  (as a relation), $\Dec\,Q$, progress, $a\,\mathbf N$ and $Q(a)$: the
  replay is total, starts at $a$, is a $C$-run and stays in $Q$.  The
  search matrix ``$C(y,x) \wedge Q(y)$'' is decided (\isn{gq\_Dec}); at
  $x \in Q$ the step equals the search value, which is a number, a
  $C$-child of $x$ and in $Q$ (\isn{gpb\_search}, by \isn{srch\_from\_ex}
  on the existential that progress supplies); the invariant $Q(P\cdot
  n)$ holds by number induction on $n$ with motive ``$Q(P\cdot v)$''
  (\isn{gpb\_inv}); and the edge $C(P\cdot(n+1), P\cdot n)$ follows from
  the recursion equation, the invariant and the search
  (\isn{gpb\_ext}).  The check forbids all four rules under the replay.
\item \emph{$\Path_R$ is decided} (\isn{g\_Dec\_path}), pointwise by
  \isn{g\_dec\_path}: the dichotomy at the accessibility body, and the
  exclusion $\Acc_R \ent \neg\Path_R$ (ground induction) for the
  accessible disjunct.
\item \emph{Progress} (\isn{g\_prog\_path}, in \isn{gcalc}, from $\Dec\,
  \Path_R$).  At $x$, split on $\Path_R\cdot x$.  In the positive case
  unfold once (\isn{g\_nu\_unfold}), eliminate the negated universal to
  a $j$ with $R(j,x)$ and $\Path_R\cdot j$, and conclude the
  existential.  No coinduction step: the path is a hypothesis.
\item Step 1 at $C = R$, $Q = \Path_R$, $a = x$, with steps 2 and 3.
\end{enumerate}

\begin{thm}[\isn{g\_tree\_branch}, \isn{g\_tree\_acc\_or\_branch}; grounded]
\label{thm:app-koenig}
In \isn{gcalc\_gfix\_dich}, with $\Gamma$ well-formed, $\Gamma \ent
\Dec\,T$, and $b$ the replay at $C = C_T$, $Q = \Path_{C_T}$ and start
$1$:
(i) if $\Gamma \ent 1 \in T$ and $\Gamma \ent \Path_{C_T}\cdot 1$, then
$\Gamma$ derives $\Tot\,b$, $b\cdot 0 = 1$, $\mathrm{in}_T(b)$,
$\mathrm{every}_{\Path_{C_T}}(b)$, $\mathrm{ext}(b)$ and
$\mathrm{len}(b)$;
(ii) if $\Gamma$ derives codes-only and weak prefix closure for $T$,
then $\Gamma \ent \Acc_{C_T}\cdot 1 \vee (\Tot\,b \wedge
\mathrm{in}_T(b) \wedge \mathrm{ext}(b) \wedge \mathrm{len}(b))$.
\end{thm}

\paragraph{Proof of \cref{thm:app-koenig}.}
\begin{enumerate}
\item $C_T$ is decided, from $\Dec\,T$ and the decidedness of equality
  (\isn{g\_TR\_decRel}, no rule beyond \rga's).
\item \emph{Clause (i).}  \cref{thm:app-runofpath} at $R = C_T$ and $x
  = 1$ gives totality, the start, the run and $\mathrm{every}$.  For
  $\mathrm{in}_T$, split $n = 0$ (the root, a member by hypothesis) or
  $n = m+1$ (the run's edge $C_T(b\cdot(m+1), b\cdot m)$, whose first
  conjunct is membership; \isn{g\_treeRun\_mem}).  $\mathrm{ext}$ is the
  second conjunct of each edge (\isn{g\_treeRun\_par}).  $\mathrm{len}$
  is number induction on $n$ with motive ``$b\cdot v$ has length $v$''
  (\isn{treeRun\_len}): the base is $b\cdot 0 = 1$, the step the parent
  equation and the drop law (\isn{dropF\_front}).
\item \emph{Removing the root hypothesis} (\isn{g\_tree\_root\_mem},
  in \isn{gcalc}).  With codes only and weak prefix closure,
  $\Path_{C_T}\cdot 1$ gives $1 \in T$: one unfolding gives a child
  $j$, a member with some length $n$; its $n$-fold drop is $1$, not
  junk, hence a member.  The corollary \isn{g\_tree\_branch\_seq}, in
  \isn{gcalc\_gfix\_dich}, is clause (i) without the root hypothesis.
\item \emph{Clause (ii).}  The dichotomy at the accessibility body of
  $C_T$, at $1$ (\isn{g\_acc\_or\_path}), gives $\Acc_{C_T}\cdot 1 \vee
  \Path_{C_T}\cdot 1$.  The accessible case is the left disjunct; in
  the path case, step 3 gives the conclusions of clause (i), and the
  conjunction of four of them (totality, $\mathrm{in}_T$,
  $\mathrm{ext}$, $\mathrm{len}$) is the right disjunct.
\end{enumerate}
The element stream (\isn{g\_tree\_elem}) sends $n$ to
$\mathrm{last}(b\cdot(n+1))$ and has the branch nodes as its prefixes,
by number induction with \isn{fb\_surj}; its prefix set lies in $T$
(\isn{g\_tree\_elem\_pset}, in \isn{gcalc\_ag\_dich}, where ATI decides
the prefix set; the ATI hypothesis is redundant there, since the
calculus has ground induction and the dichotomy and so derives ATI,
\cref{thm:ati-dich}).

\paragraph{Where ATI is spent.}
In this appendix ATI is applied as a rule only to decide a prefix set
(\isn{g\_pset\_Dec}, one \isn{ir\_allTI\_open}).  The round trip spends
ATI through \isn{g\_pset\_Dec} twice over: in clause (i) as the decidedness premise of the
search at $\mathrm{pset}_s$, in clause (ii) as the case split of the
first inclusion.  \cref{thm:app-psetpath}(ii) spends ATI through the
same lemma.  The search lemma \isn{mu\_choice} sits in the base calculus, since its derivation is \rga's search lemma, and the check for the round trip requires it by name; so a result that uses the search and decides nothing by ATI, such as player II's strategy in Appendix~\ref{app:games}, does not apply ATI.
\cref{thm:app-runs} is stated in calculi without ATI.
\cref{thm:app-runofpath,thm:app-koenig} are stated without ATI as a
primitive rule: their calculi have ground induction and the dichotomy
and so derive ATI (\cref{thm:ati-dich}), the proofs do not apply it,
and their theories are built where that derivation is not in
scope.

\paragraph{Lemmas a re-mechanization needs.}
\begin{itemize}
\item \isn{pfx\_0}, \isn{pfx\_S}, \isn{pfx\_par\_S}: the recursion
  equations of the prefix function and its parent law.
\item \isn{pfx\_len}, \isn{g\_pfx\_len\_uniq}: the prefix of index $n$
  has length exactly $n$.
\item \isn{g\_pfx\_drop\_le}, \isn{g\_pfx\_drop\_gt}: a drop of a prefix
  is a shorter prefix or junk.
\item \isn{g\_pset\_Dec}: a prefix set is decided (ATI).
\item \isn{g\_pset\_PrefClosedS}, \isn{g\_pset\_InfTree},
  \isn{g\_pset\_Lin}, \isn{g\_pset\_Seqs}: the four other branch
  conditions of a prefix set.
\item \isn{mu\_choice}, \isn{g\_mu\_choice\_all}: a decided matrix with
  a witness gives a total search function whose values satisfy it.
\item \isn{srch\_from\_ex}: the same in \rga, for one search.
\item \isn{fb\_surj}, \isn{g\_fbLinE},
  \isn{g\_fbPrefClosedS\_len}: decomposition of a nonempty code;
  elimination of linearity and of weak prefix closure at a code with a
  length.
\item \isn{g\_RTp\_bridge}, \isn{g\_RTp\_pfx}: consecutive search values
  are parent and child; the recovered stream's prefixes are the search
  values.
\item \isn{g\_path\_coind}: coinduction for $\Path_R$ at a certificate
  set.
\item \isn{g\_path\_child\_intro}: the body of $\Path_R$ from a child in
  the certificate.
\item \isn{g\_TR\_child\_pset}: a code with a length whose child is a
  prefix is a prefix.
\item \isn{g\_TR\_decRel}: the child relation of a decided set is
  decided.
\item \isn{gi\_TI\_par}: ground induction along $\Acc_R$ with a
  parametric motive.
\item \isn{gpb\_run\_of\_progress} (with \isn{gpb\_search},
  \isn{gpb\_inv}, \isn{gpb\_ext}): the replay is a run that stays in
  $Q$.
\item \isn{g\_prog\_path}: a decided $\Path_R$ has progress, by one
  unfolding.
\item \isn{g\_dec\_path}, \isn{g\_acc\_or\_path}: $\Path_R$ decided, and
  $\Acc_R \vee \Path_R$, at a decided relation.
\item \isn{g\_tree\_root\_mem}, \isn{treeRun\_len}: the root of a tree
  with codes only and weak prefix closure is a member when
  $\Path_{C_T}\cdot 1$ holds; a run of the child relation from
  the root has the right lengths.
\end{itemize}

\section{The two program examples}
\label{app:programs}

Files: \isn{OIGA\_Coind/GRand.thy} (random assignment and its
divergent twin), \isn{OIGA\_Coind/GTrace.thy} (the search body, its
two atoms, and the counting program's bodies),
\isn{OIGA\_Coind/GTraceModel.thy} (the model side of the counting
program), with the generic lemmas of \isn{OIGA/GDich.thy},
\isn{OIGA/GDich2.thy}, \isn{OIGA\_Coind/GStream.thy} and
\isn{OIGA\_Coind/GBranch.thy}.

\paragraph{Definitions: random assignment.}
A state is $0$ (before the assignment) or $m+1$ (the loop with $x =
m$).  The step relation $R$ (\isn{R1}) is $R(y, x) \equiv (x = 0
\wedge y \ne 0) \vee (y \ne 0 \wedge x = y+1)$: $0$ steps to every
$m+1$, $m+2$ steps to $m+1$, and $1$, the loop exit, has no successor.
The divergent twin $R'$ (\isn{R1p}) adds the step from $1$ to $0$:
$R'(y, x) \equiv R(y, x) \vee (x = 1 \wedge y = 0)$.  A term $V$ is a
\emph{variant} (\isn{r1Var}) when $\forall x\,\forall y.\,\neg R(y,x)
\vee V\cdot y < V\cdot x$, the order being the compiled test
\isn{ltOp}.  The certificate ``every state'' (\isn{allSt}) is $\lambda
v.\,0 = 0$.

\paragraph{Definitions: the search and the counting program.}
For a set term $q$, the search body (\isn{trSPhi}) is $q\cdot x \vee
a\cdot(x+1)$, and its dual (\isn{trSPsi}) is $q\cdot x \vee
\neg(a\cdot(x+1))$.  $\Term_q = \mu(\text{search body})$ (\isn{trTerm})
and $\Div_q = \nu(\neg\text{dual})$ (\isn{trDiv}), whose body reads
$\neg q\cdot x \wedge a\cdot(x+1)$.  The decidedness premise of the
search body is $\Dec\,q \wedge 0 = 0$ (\isn{trDecD},
\isn{nfdec\_trSPhi}).  $\mathrm{Ex}_q(x) \equiv \exists m.\,q\cdot(x+m)$
(\isn{trEx}) and $\mathrm{All}_q(x) \equiv \forall m.\,\neg q\cdot(x+m)$
(\isn{trAll}); the negation of $\mathrm{All}_q(x)$ is literally
$\mathrm{Ex}_q(x)$ (\isn{trAll\_neg}).  The counting program's state is
a pair $\langle x, y\rangle$, and its step (\isn{trStep}, source first)
is: if $y = 0$, go to $\langle x+1, x+1\rangle$; otherwise go to
$\langle x, y-1\rangle$.  The first component never decreases.
Reachability $\mathrm{Reach}_t$ (\isn{trReach}) is the $\mu$-atom of
$x = t \vee \exists s.\,\mathrm{step}(x, s) \wedge a\cdot s$: the run
from the point reaches $t$.  The body ``reaches $n$ and then goes up''
(the \emph{up body}, \isn{trUpPhi}) reads, at a point $p = \langle n, \sigma\rangle$: some
$\sigma'$ has $\mathrm{Reach}_{\sigma'}\cdot\sigma$, the first component
of $\sigma'$ is $n$, and some $\sigma''$ with $\mathrm{step}(\sigma',
\sigma'')$ has $a\cdot\langle n+1, \sigma''\rangle$ (through an
equality guard on a bound variable).  The reach atom is a leaf of that
body that does not mention the recursive variable.

\begin{thm}[\isn{g\_acc\_R1\_all}, \isn{g\_R1\_no\_variant}; grounded]
\label{thm:app-rand}
Let $\Gamma$ be well-formed.
(i) In \isn{gcalc}: if $\Gamma \ent x\,\mathbf N$, then $\Gamma \ent
\Acc_R\cdot x$.
(ii) In \isn{gcalc\_ati}: if $\Gamma \ent \Tot\,V$, then $\Gamma \ent
\neg\mathrm{Var}(V)$, where $\mathrm{Var}(V)$ is the variant property.
\end{thm}
Clause (ii) is a schema: for each term $V$ whose totality is
derivable, the negation is derivable.  No quantifier over terms is
involved.

\paragraph{Proof of \cref{thm:app-rand}.}
\begin{enumerate}
\item \emph{Accessibility introduction} (\isn{g\_acc\_intro0}, in
  \isn{gcalc}).  For a closed well-formed $R$, a well-formed $x$ and a
  well-formed $\Gamma$: if in every well-formed extension, at every
  number $j$, $\neg R(j, x) \vee \Acc_R\cdot j$ is derivable, then
  $\Gamma \ent \Acc_R\cdot x$.  The lemma is universal introduction
  followed by one fold of the $\mu$-atom (\isn{acc\_unfoldI},
  conversions), and a fold is not a fixed-point rule.
\item \emph{The loop states} (\isn{g\_acc\_R1\_loop}).  $\Acc_R\cdot
  (m+1)$, by number induction on $m$ with motive ``$\Acc_R\cdot(v+1)$''.
  Base: no $j$ has $R(j, 1)$, since $1 \ne 0$ refutes the left disjunct
  and $1 = j+1$ with $j \ne 0$ is refuted by a decided case on $j = 0$;
  step 1 concludes.  Step: at $m+2$ and a number $j$, split on the
  decided $j = m+1$.  If so, the hypothesis gives $\Acc_R\cdot j$ by
  transport; if not, $\neg R(j, m+2)$, because $m+2 \ne 0$ and
  successor is injective.  Step 1 concludes.
\item \emph{The initial state} (\isn{g\_acc\_R1\_init}).  One
  application of step 1 at $0$: split $j = 0$, where $R(0, 0)$ fails,
  or $j = k+1$, accessible by step 2.
\item \emph{Clause (i).}  Split $x = 0$ or $x = k+1$; steps 3 and 2.
\item \emph{A variant refutes itself, in \rga} (\isn{r1Var\_absurd}).
  From $\Tot\,V$ and $\mathrm{Var}(V)$, any well-formed $d$.  Number
  induction on $m$ with motive ``$m < V\cdot(m+1) + 1$'', that is $m \le
  V\cdot(m+1)$.  Base: $0 \le V\cdot 1$.  Step: the edge from $m+2$ to
  $m+1$ (\isn{R1\_edge\_dec}) gives $V\cdot(m+1) < V\cdot(m+2)$
  (\isn{r1Var\_E}); splitting $m \le V\cdot(m+1)$ into equality
  (transport) and $<$ (transitivity) gives $m < V\cdot(m+2)$
  (\isn{ir\_ltOp\_split}), so $m+1 \le V\cdot(m+2)$.  The
  instance at $m = V\cdot 0$ gives $V\cdot 0 \le V\cdot(V\cdot 0 + 1)$;
  the edge from $0$ to $V\cdot 0 + 1$ (\isn{R1\_edge\_init}) gives
  $V\cdot(V\cdot 0 + 1) < V\cdot 0$; both cases of $\le$ contradict
  irreflexivity of $<$.
\item \emph{The variant property is decided.}  At two fresh number
  variables the matrix $\neg R(v_0, v_1) \vee V\cdot v_0 < V\cdot v_1$ is
  decided; two applications of ATI (\isn{ir\_allTI\_open}) decide
  $\mathrm{Var}(V)$, hence its negation.
\item \emph{Clause (ii).}  Under the hypothesis $\neg\neg\mathrm{Var}
  (V)$, double-negation elimination and step 5 give $0 = 1$.
  \isn{neg\_decided\_explosive\_iff} (for a decided $a$, $\neg a \ent 0
  = 1$ holds exactly when $a$ is derivable), at the decided $a =
  \neg\mathrm{Var}(V)$, gives $\neg\mathrm{Var}(V)$.
\end{enumerate}
The contradiction from a variant is \rga's; ATI is spent only to
decide the variant property, which turns the hypothetical refutation
into a derivable negation.

\begin{thm}[\isn{g\_R1p\_div}, \isn{g\_R1p\_run}; grounded]
\label{thm:app-randdiv}
Let $\Gamma$ be well-formed.
(i) In \isn{gcalc\_nu}: $\Gamma \ent \Path_{R'}\cdot 0$.
(ii) In \isn{gcalc\_gfix\_nu\_dich}: writing $P$ for the replay at $C =
R'$, $Q = \Path_{R'}$ and start $0$ (Appendix~\ref{app:streams}),
$\Gamma$ derives $\Tot\,P$, $P\cdot 0 = 0$, $\mathrm{run}_{R'}(P)$ and
$\mathrm{every}_{\Path_{R'}}(P)$.
\end{thm}

\paragraph{Proof of \cref{thm:app-randdiv}.}
\begin{enumerate}
\item \emph{Clause (i).}  Coinduction at a certificate set
  (\isn{g\_path\_coind}) with $Q$ = ``every state''; the region is
  every number, and $Q\cdot 0$ holds by conversion.  Closure at a fresh
  $j$: split $j = 0$ or $j = k+1$, and $k = 0$ or $k = l+1$.  At $j =
  0$ the child is $1$, with $R'(1, 0)$ (\isn{R1\_edge\_init},
  \isn{R1p\_of\_R1}); at $j = 1$ the child is $0$, with the back edge
  (\isn{R1p\_back}); at $j = l+2$ the child is $l+1$
  (\isn{R1\_edge\_dec}).  Each case closes by \isn{g\_path\_child\_intro}.
  One coinduction step.
\item \emph{Clause (ii).}  $R'$ is decided in \rga
  (\isn{g\_R1p\_decRel}); \isn{g\_run\_of\_path}
  (Appendix~\ref{app:streams}) at $R'$, with the path of clause (i) at
  $x = 0$.
\end{enumerate}
Companions.  $\Acc_R\cdot x \vee \Path_R\cdot x$ and the same for $R'$
are \isn{g\_acc\_or\_path} at the two decided relations
(\isn{g\_R1\_acc\_or\_path}, \isn{g\_R1p\_acc\_or\_path}, in
\isn{gcalc\_dich}, the dichotomy alone).  The decided forms
(\isn{g\_R1\_dec\_path}, \isn{g\_R1p\_dec\_path}) are \isn{g\_dec\_path},
in \isn{gcalc\_gfix\_dich}, so they use ground induction as well, for
the exclusion.  Non-derivability (\isn{oinga\_R1p\_not\_acc}, classical):
an empty-context derivation of $\Acc_{R'}\cdot 0$ in \oinga gives
$\neg\Path_{R'}\cdot 0$ by the exclusion, which clause (i) refutes,
against the consistency of \oinga (\isn{oinga\_consistent\_notation});
\isn{oingad\_R1p\_not\_acc} is the same for \oinga with the
dichotomy, restated at \dinga as \isn{dinga\_R1p\_not\_acc} in the
junction theory \isn{CD\_Junction}.

\begin{thm}[\isn{tr\_Q1\_iff}, \isn{tr\_Q2\_iff}; grounded]
\label{thm:app-trace}
Let $\Gamma$ be well-formed and $\Gamma \ent x\,\mathbf N$.
(i) In \isn{gcalc\_gfix}, with $\Gamma \ent \Dec\,q$: $\Gamma \ent
\Term_q\cdot x$ if and only if $\Gamma \ent \mathrm{Ex}_q(x)$.
(ii) In \isn{gcalc\_nu}, for every $q$: $\Gamma \ent \Div_q\cdot x$ if
and only if $\Gamma \ent \mathrm{All}_q(x)$.
(iii) In \isn{gcalc\_gfix} (\isn{tr\_excl}): if $\Gamma \ent \Term_q
\cdot x$, then $\Gamma \ent \neg\Div_q\cdot x$.
\end{thm}
The two equivalences are between derivability facts, not internal
biconditionals.  Clause (ii) needs no decidedness of $q$.

\paragraph{Proof of \cref{thm:app-trace}.}
\begin{enumerate}
\item \emph{Backward induction} (\isn{tr\_back\_ind}, in \isn{gcalc}).
  For a predicate $P$: from $\Dec\,q$, $x\,\mathbf N$, $m\,\mathbf N$,
  $q\cdot(x+m)$ and, in every well-formed extension and at every
  number $y$, ``$q\cdot y$ gives $P(y)$'' and ``$P(y+1)$ gives
  $P(y)$'', conclude $P(x)$.  Number induction on $m$ with motive $K(k)
  \equiv \forall y.\,\neg q\cdot(y+k) \vee P(y)$ (\isn{trK}).  Base: at
  $y$, split on the decided $q\cdot(y+0)$; in the positive case
  \isn{tr\_pAdd\_0} transports it to $q\cdot y$, which gives $P(y)$.  Step: split on $q\cdot(y+k+1)$; in the positive case
  $q\cdot((y+1)+k)$ (\isn{tr\_pAdd\_shift}), so $K(k)$ at $y+1$ gives
  $P(y+1)$, hence $P(y)$.  The instance of $K(m)$ at $y = x$ finishes.
\item \emph{Clause (i), backward} (\isn{tr\_term\_of\_ex}, no rule
  beyond \rga's).  Eliminate the existential and apply step 1 with $P =
  \Term_q$: a witness gives the left disjunct of the body, and
  $\Term_q\cdot(y+1)$ the right one, each followed by one fold
  (\isn{gm\_gfix\_fold}).
\item \emph{Clause (i), forward} (\isn{tr\_ex\_of\_term}).  Ground
  induction (\isn{ir\_gfix\_par}) at the search body with motive $c(v)
  \equiv \mathrm{Ex}_q(v)$, the eigenvariable fresh for the context, the
  motive and the body.  Closure obligation: from $q\cdot v \vee
  \mathrm{Ex}_q(v+1)$, derive $\mathrm{Ex}_q(v)$.  In the left case the
  witness is $0$ (\isn{tr\_pAdd\_0}); in the right case a witness $m$
  for $v+1$ gives the witness $m+1$ for $v$ (\isn{tr\_pAdd\_shift}).
  The rule yields $c(x)$.  This direction does not use $\Dec\,q$.
\item \emph{Clause (ii), forward} (\isn{tr\_all\_of\_div}, no rule
  beyond \rga's).  First $\Div_q\cdot(x+m)$ for every $m$ (\isn{tr\_div\_add}),
  by number induction on $m$ with motive ``$\Div_q\cdot(x+v)$'', the step
  one unfolding and a transport along $x+(v+1) = (x+v)+1$.  Then, at
  each $m$, one unfolding at $x+m$ gives the conjunct $\neg q\cdot(x+m)$.
\item \emph{Clause (ii), backward} (\isn{tr\_div\_of\_all}).  The
  coinduction rule at the body $\neg q\cdot x \wedge a\cdot(x+1)$ with
  certificate $c(v) \equiv \mathrm{All}_q(v)$; the region is the set of
  points from which $q$ never holds.  The start $c(x)$ is the
  hypothesis.  Closure obligation: at a fresh $i$ with $\mathrm{All}_q
  (i)$, derive $\neg q\cdot i \wedge c(i+1)$.  The instance at $m = 0$
  gives $\neg q\cdot i$; the instance at $m+1$, transported, gives
  $\mathrm{All}_q(i+1)$; conversion reads that universal as the certificate at
  $i+1$.  One coinduction step.
\item \emph{Clause (iii).}  \isn{g\_mu\_not\_nu} at the search body and
  its dual: ground induction with motive ``$\neg\Div_q\cdot v$''.
\end{enumerate}

\begin{thm}[\isn{tr\_term\_of\_ndiv}, \isn{tr\_dich}; grounded]
\label{thm:app-tracegap}
Let $\Gamma \ent \Dec\,q$ and $\Gamma \ent x\,\mathbf N$.
(i) In \isn{gcalc\_ati\_nu}, with $\Gamma$ well-formed: if $\Gamma \ent
\neg\Div_q\cdot x$, then $\Gamma \ent \Term_q\cdot x$.
(ii) In \isn{gcalc\_dich}: $\Gamma \ent \Term_q\cdot x \vee \Div_q
\cdot x$.
The third route (\isn{tr\_ati\_disj}, in \isn{gcalc\_ati\_nu}, with
$\Gamma$ well-formed) derives the same disjunction.
\end{thm}

\paragraph{Proof of \cref{thm:app-tracegap}.}
\begin{enumerate}
\item \emph{The universal is decided} (\isn{tr\_all\_B}, in
  \isn{gcalc\_ati}).  At a fresh number variable $v$, $\neg q\cdot(x+v)$
  is decided by $\Dec\,q$; ATI (\isn{ir\_allTI\_open}) decides
  $\mathrm{All}_q(x)$.  Here ATI acts as the limited principle of
  omniscience.
\item \emph{Clause (i).}  Split on $\mathrm{All}_q(x)$.  If the universal holds,
  step 5 of the proof of \cref{thm:app-trace} gives $\Div_q\cdot x$, against the
  hypothesis.  If $\mathrm{All}_q(x)$ fails, its negation is literally $\mathrm{Ex}_q(x)$, and
  step 2 of the same proof gives $\Term_q\cdot x$.
\item \emph{Clause (ii).}  The premise $\Dec\,q \wedge 0 = 0$ by
  conjunction introduction; one application of the dichotomy at the
  search body, with its positivity, mention, dual and premise facts
  (\isn{posP\_trSPhi}, \isn{free1\_trSPhi}, \isn{ndual\_trSPhi},
  \isn{nfdec\_trSPhi}).
\item \emph{The third route.}  The case split of step 2, with
  $\Div_q\cdot x$ on the positive side and $\Term_q\cdot x$ on the
  negative.  The third route's decided form \isn{tr\_ati\_dec} uses the witness
  refutation \isn{tr\_ndiv\_of\_ex} on the negative side; the
  dichotomy's decided form \isn{tr\_dich\_dec} is \isn{g\_dec\_nu}.
\end{enumerate}

\begin{thm}[\isn{nfdec\_not\_false}, \isn{trUp\_premise\_valueless}; in $\Mstar$]
\label{thm:app-premise}
(i) For every $k$, body $t$ and premise $D$ with $\mathrm{nfdec}(k, t,
D)$, and every substitution $\tau$: $\tau(D)$ is not false in $\Mstar$.
(ii) For every $D$ with $\mathrm{nfdec}(1, \text{up body}, D)$ and every
$\tau$: $\tau(D)$ has no value in $\Mstar$.
(iii) In \dinga (\isn{oingad\_trUp\_dich\_premise\_not\_derivable},
\isn{oingad\_trUp\_dich\_premise\_neg\_not\_derivable}, restated at \dinga
as \isn{dinga\_trUp\_dich\_premise\_not\_derivable} and
\isn{dinga\_trUp\_dich\_premise\_neg\_not\_derivable} in \isn{CD\_Junction};
classical): for
every $D$ with $\mathrm{nfdec}(1, \text{up body}, D)$, neither $D$ nor
$\neg D$ is derivable from the empty context.
(iv) The combined forms (\isn{tr\_dich\_prem\_valueless},
\isn{oingad\_tr\_dich\_prem\_open}, \isn{oingad\_tr\_dich\_unavailable},
restated at \dinga as \isn{dinga\_tr\_dich\_prem\_open} and
\isn{dinga\_tr\_dich\_unavailable} in \isn{CD\_Junction};
classical): every such $D$ is neither true nor false in $\Mstar$ and
neither derivable nor refutable, and such a $D$ exists
(\isn{tr\_dich\_prem\_ex}), so the rule is unavailable at the up body.
\end{thm}
The premise relation is not functional: $D$ ranges over every premise
the grammar computes, and the theorem holds for each.

\paragraph{Proof of \cref{thm:app-premise}.}
\begin{enumerate}
\item \emph{The step in $\Mstar$.}  An equation $t = \bar v$ derivable
  in \rga from no hypotheses gives $t$ the value $v$ in $\Mstar$, by \isn{rga\_nu\_sound}
  (\isn{trmod\_val\_of\_eq}); so pairing, projections and predecessor
  have their arithmetic values at numerals.  A true step from $p$ to
  $q$ forces $q = \mathrm{succ}(p)$, where $\mathrm{succ}$ is the HOL
  function computing the program's next state (\isn{trmod\_step\_inv}), and the step to $\mathrm{succ}(p)$ is true
  (\isn{trmod\_step\_true}).
\item \emph{The reach atom is never false} (\isn{tr\_reach\_not\_false},
  at every target $T$ and numeral $p$).  Induction over the mutual
  inductive definition of the model's head and value relations, with
  the invariant \isn{trFInv}: no value judgment derived is one of
  ``$\mathrm{Reach}_T\cdot\bar p$ is false'', ``its body instance is
  false'', ``the negated universal of the body is false'', ``the
  universal $\forall s.\,\neg\mathrm{step}(\bar p, s) \vee
  \neg\mathrm{Reach}_T\cdot s$ is true'', for $q$ whose step from $p$ is
  not false, ``the $q$-instance of that universal, as a redex or
  contracted, is true'', and ``$\neg\mathrm{Reach}_T\cdot\bar q$ is
  true''.  The case that uses the program is the true universal: its truth includes the instance at $q
  = \mathrm{succ}(p)$, whose step is not false
  (\isn{trmod\_succ\_not\_false}), and that instance is excluded by the
  induction hypothesis.  Read operationally, a falsity derivation would
  have to follow the run forever.
\item \emph{The reach atom from $\langle 1,0\rangle$ to $\langle
  0,0\rangle$ is not true} (\isn{tr\_reach\_not\_true}).  The atom's truth
  set in $\Mstar$ is the least fixed point of its operator relative to
  $\Mstar$ (\isn{gfix\_omvO\_eq\_lfp}).  By least-fixed-point induction
  at the set ``first component $0$'' (\isn{tr\_reach\_lfp\_fst0}), every
  state from which $\langle 0,0\rangle$ is reachable has first
  component $0$: a state in the operator's image is $\langle 0,
  0\rangle$ itself, or steps truly to a member $q$ with first component
  $0$, and then $q = \mathrm{succ}(p)$ forces the first component of $p$
  to be $0$.  The state $\langle 1,0\rangle$ has first component $1$.
\item \emph{The premise is not true} (\isn{trUp\_premise\_not\_true}).  If
  $\tau(D)$ were true, its model-side reading (\isn{nfdec\_leafdec})
  would make every leaf valued at every environment of numerals;
  descending through the negation, the universal and the disjunction of
  the up body (\isn{trmod\_up\_leafdec}), the reach leaf would be valued
  at the environment sending $\sigma'$ to $\langle 0,0\rangle$ and the
  point to $\langle 0, \langle 1,0\rangle\rangle$.  Value congruence
  through the second projection turns that leaf into
  $\mathrm{Reach}_{\langle 0,0\rangle}\cdot\langle 1,0\rangle$, which
  step 3 makes not true and step 2 not false.
\item \emph{Clause (i).}  Induction on the premise grammar.  Leaf: the
  universal closure of ``$X$ is decided'' is never false
  (\isn{trmod\_allqN\_B\_not\_false}, by induction on the number of
  quantifiers: a false universal has a false instance; with none left,
  $X \vee \neg X$ false would make $X$ both true and false, against
  uniqueness of values).  Recursive atom: $0 = 0$ is not false.
  Disjunction: a false conjunction has a false conjunct.  Negation and
  the universal pass the premise through unchanged.
\item \emph{Clause (ii)}: a value of $\tau(D)$ would be true, against
  step 4, or false, against clause (i).
\item \emph{Clause (iii).}  An empty-context derivation of $D$ makes $D$
  true in $\Mstar$ (\isn{oingad\_nu\_sound}), against step 4 at the
  identity substitution.  A derivation of $\neg D$ makes $D$ false,
  against clause (i).
\end{enumerate}

\paragraph{Where ATI is spent.}
In this appendix ATI is applied as a rule only to decide a universal:
the variant property in \cref{thm:app-rand}(ii) and $\mathrm{All}_q(x)$
in \cref{thm:app-tracegap}(i) and the third route.  The termination of
random assignment, the divergence of its twin, the disjunctions $\Acc \vee \Path$ and
\cref{thm:app-trace} are stated in calculi without ATI; the decided
forms are stated without ATI as a primitive rule, in calculi with
ground induction and the dichotomy, which derive ATI
(\cref{thm:ati-dich}), but their proofs do not apply it and their
theories are built where that derivation is not in scope.  The checks
forbid it under \cref{thm:app-rand}(i), \isn{g\_R1p\_div}, the two
directions of \cref{thm:app-trace} that use a rule, and \isn{tr\_dich}.

\subsection{Nested random assignment at every fixed depth}

Files: \isn{OIGA\_Coind/GRandDepth.thy} (imports \isn{GRand}).  The relation of depth $k$ is a HOL function of $k$ producing a term (\isn{rdRel}): depth $0$ has no edges; depth $k+1$ runs on Cantor-pair codes $\langle s, c\rangle$ of a depth-$k$ state $s$ and a counter $c$, and one step either decrements the counter or steps $s$ along the depth-$k$ relation and resets the counter to any number.  Read lexicographically the relation has height $\omega^k$ (a reading).  Every number is a state, and the relation is decided at every $k$ (\isn{rdRel\_B}, HOL induction on $k$, no ATI).

\begin{thm}[\isn{rd\_ti}, \isn{acc\_rdRel}; grounded]
\label{thm:app-depth}
In the locale \isn{rga}, for every $k$: transfinite induction along the depth-$k$ relation holds for every closed motive that is progressive in the every-extension form of \isn{acc\_Rn\_I}; and for every well-formed $\Gamma$ and $\Gamma \ent s\,\mathbf N$, $\Gamma \ent \Acc_{R_k}\cdot s$.  In \isn{gcalc}: \isn{g\_acc\_rdRel}, \isn{g\_rdRel\_B}.
\end{thm}

\begin{enumerate}
\item \emph{HOL induction on $k$, motive arbitrary.}  Accessibility is not a strong enough motive for the induction on $k$; the schema is stated for every closed progressive $P$.
\item \emph{The step.}  Apply the depth-$k$ schema to the motive $\lambda s.\,\forall c.\,P\langle s, c\rangle$.  Its progressiveness is one internal number induction (\isn{ir\_nat\_ind}) on the counter, using \isn{allI\_cont}, the edge analysis \isn{rdB\_cases}, the pairing round trips with \isn{cmp\_argI} and \isn{app\_T\_cong}, and the outer hypothesis for a step over a depth-$k$ edge.
\item \emph{Accessibility.}  \isn{acc\_rdRel} is the schema at $P = \Acc$, made progressive by \isn{acc\_unfoldI} together with \isn{rdRel\_B}; \isn{acc\_rdRel\_3} is the instance $k = 3$.  The check at the end of the theory forbids ground induction, $\nu$, the dichotomy and ATI under the new theorems, with a blindness control on \isn{g\_R1\_not\_path}.
\end{enumerate}

\subsection{The Lattice rule and the WELL rule}

Files: \isn{OIGA\_Coind/GLattice.thy} (imports \isn{GTrace}).  For a closed step relation $T$, source first, and a target $G$, the \emph{eventually} atom $\mathrm{Ev}_G$ (\isn{trEv}) is the $\mu$-atom of the body $G\cdot x \vee \forall y.\,\neg T\cdot x\cdot y \vee a\cdot y$; a state without a successor satisfies it vacuously, so the step relation is taken total.  \emph{Leads-to} (\isn{ltTo}) is the sentence $\forall x.\,\neg P\cdot x \vee \mathrm{Ev}_Q\cdot x$.  The \emph{lattice hypothesis} (\isn{latHyp}) at an order $R$, a family $F$ and a target $G$ is the sentence $\forall c.\,\neg\Acc_R\cdot c \vee (F\cdot c \leadsto \lambda y.\,G\cdot y \vee \exists d.\,R\cdot d\cdot c \wedge F\cdot d\cdot y)$.

\begin{thm}[\isn{g\_lattice}, \isn{g\_well}; grounded]
\label{thm:app-lattice}
In \isn{gcalc\_gfix}, for closed well-formed $T$, $R$, $F$, $G$: from $\Gamma \ent \isn{latHyp}\,T\,R\,F\,G$, $\Gamma \ent \Acc_R\cdot c$ and $\Gamma \ent c\,\mathbf N$, $\Gamma \ent F\cdot c \leadsto G$; at a state, from $\Gamma \ent x\,\mathbf N$ and $\Gamma \ent F\cdot c\cdot x$, $\Gamma \ent \mathrm{Ev}_G\cdot x$ (\isn{g\_lattice\_at}).  The WELL rule: at the family ``$P\cdot x$ and $\delta\cdot x = c$'' (\isn{wellF}), from the hypothesis, $P\cdot x$, $x\,\mathbf N$, $(\delta\cdot x)\,\mathbf N$ and $\Acc_R\cdot(\delta\cdot x)$, $\Gamma \ent \mathrm{Ev}_G\cdot x$.
\end{thm}

\begin{enumerate}
\item \emph{Transitivity and monotonicity of eventually} (\isn{ev\_mono\_trans}, \isn{gcalc\_gfix}): from $\mathrm{Ev}_H\cdot x$ and, in every extension of the context, $H\cdot y$ giving $\mathrm{Ev}_G\cdot y$, conclude $\mathrm{Ev}_G\cdot x$; one ground induction on the outer atom, the unfolding and folding laws \isn{ev\_unfold}, \isn{ev\_fold}, \isn{ev\_now}, \isn{ev\_next} being \rga-level.
\item \emph{The motive.}  Ground induction on $\Acc_R\cdot c$ at the motive ``$\Acc_R\cdot c$ and $F\cdot c \leadsto G$'' (\isn{latMot}): the closure receives the body of the accessibility atom at the new point, not the atom, so the motive carries accessibility to meet the relativized hypothesis.  The schema form \isn{g\_lattice\_sch} takes the hypothesis at every extension; \isn{g\_lattice} instantiates the sentence.
\item \emph{The closure.}  Rebuild $\Acc_R\cdot z$ from the first conjunct of the induction hypothesis; apply the hypothesis at $z$; apply step~1 with $H$ the lattice target at $z$: a $G$-state folds directly (\isn{ev\_now}), and a witness $d \prec z$ gets $\mathrm{Ev}_G$ from the induction hypothesis' second conjunct.
\item \emph{Why the conclusion is a rule.}  The sentence $\forall c.\,\neg\Acc_R\cdot c \vee (F\cdot c \leadsto G)$ would need $\neg\Acc_R\cdot c$ at points the order leaves inaccessible, where accessibility is undecided; Lamport's existential conclusion is handled in the metatheory by \isn{g\_lattice\_at}.
\item \emph{The countdown} (\isn{cd\_term}): the step $x \mapsto x - 1$ stuttering at $0$, $F\cdot c\cdot x$ the equation $x = c$, $G$ the set $\{0\}$; every number is accessible along the predecessor order by one number induction (\isn{cd\_acc}) and the hypothesis sentence is derived without a fixed-point rule (\isn{cd\_hyp}); the Lattice rule at a state gives $\mathrm{Ev}_{\{0\}}\cdot c$ for every $c$.
\end{enumerate}

\subsection{AG\,AF\,$p$ with the inner leaf a negated path atom}

Files: \isn{OIGA\_Coind/GAgAf.thy} (imports \isn{OIGA/GDich2.thy}).
For closed relations $R$ and $Q$, child first as for $\Path$, the body
\isn{agafPhi}\,$R\,Q$ is $\Path_Q\cdot x \vee \neg\forall y.\,(\neg
R\cdot y\cdot x \vee \neg a\cdot y)$, with $a$ the recursive variable:
``$\Path_Q\cdot x$, or some $R$-successor of $x$ is in $a$''.  It is
positive in $a$ and mentions it (\isn{posP\_agafPhi},
\isn{free1\_agafPhi}).  Its decidedness premise (\isn{agafD}) is the
conjunction of $\forall z.\,(\Path_Q\cdot z)\,\mathbf B$, $\Dec\,R$ and
$0 = 0$ at the recursive atom (\isn{nfdec\_agafPhi}).  The intended
$Q$ is the restriction \isn{restrRel}\,$R\,p = \lambda y.\,\lambda
x.\,R\cdot y\cdot x \wedge \neg p\cdot x$, the steps of $R$ out of
the states where $p$ fails.  At that $Q$, $\Path_Q\cdot x$ is
EG\,$\neg p$, the body's $\mu$-atom is EF\,EG\,$\neg p$, and the
$\nu$-atom of the body's dual is AG\,AF\,$p$, a classical
identification.  $\mathrm{decSet}\,p$ (\isn{decSet}) says that $\neg
p\cdot a$ is decided at every number $a$.

\begin{thm}[\isn{g\_dec\_agaf}, \isn{g\_dec\_agaf\_restr}; grounded]
\label{thm:app-agaf}
In \isn{gcalc\_gfix\_dich}, at a well-formed $\Gamma$ and closed $R$:
(i) at closed $Q$, from $\Gamma \ent \Dec\,R$ and $\Gamma \ent
\Dec\,Q$, $\Gamma \ent \isn{agafD}\,R\,Q$ (\isn{g\_agaf\_premise});
and for every $\Psi$ with $\mathrm{ndual}(1, \isn{agafPhi}\,R\,Q,
\Psi)$ and $\Gamma \ent x\,\mathbf N$, $\Gamma \ent
(\nu(\neg\Psi)\cdot x)\,\mathbf B$ (\isn{g\_dec\_agaf}; with the dual
supplied, \isn{g\_dec\_agaf\_ex}).
(ii) In \isn{gcalc}, at closed $p$: from $\Gamma \ent \Dec\,R$ and
$\Gamma \ent \mathrm{decSet}\,p$, $\Gamma \ent
\Dec(\isn{restrRel}\,R\,p)$ (\isn{g\_restrRel\_decRel}).
(iii) Hence, at closed $p$, from $\Gamma \ent \Dec\,R$ and $\Gamma
\ent \mathrm{decSet}\,p$, the $\nu$-atom of the dual of
$\isn{agafPhi}\,R\,(\isn{restrRel}\,R\,p)$ is decided at every
number (\isn{g\_dec\_agaf\_restr}, \isn{g\_dec\_agaf\_restr\_ex}).
\end{thm}

\begin{enumerate}
\item \emph{The premise} (\isn{g\_agaf\_premise}).  The first conjunct
  is a universal introduction at every extension of the context
  (\isn{g\_allqN\_env}), each instance $(\Path_Q\cdot z)\,\mathbf B$
  given by \isn{g\_dec\_path}, the decided path atom of
  \cref{thm:excl}, at the decided $Q$; the second is the hypothesis,
  the third reflexivity; two conjunction introductions
  (\isn{g\_dich\_conjI}) assemble them.
\item \emph{The outer atom} (\isn{g\_dec\_agaf}): \isn{g\_dec\_nu}, the
  second clause of \cref{thm:excl}, at the body with its positivity,
  mention, dual and premise facts and step~1.
\item \emph{The restriction is decided} (\isn{g\_restrRel\_B},
  \isn{g\_restrRel\_decRel}, no rule beyond \rga's).  At numbers $y$
  and $x$, $R\cdot y\cdot x$ is decided by $\Dec\,R$ and $\neg p\cdot
  x$ by $\mathrm{decSet}\,p$ (\isn{g\_decSet\_inst}), so the
  conjunction is decided (\isn{g\_andTI\_neg}); two beta steps
  (\isn{restrRel\_step1}, \isn{restrRel\_step2}) carry the decidedness
  back to $\isn{restrRel}\,R\,p\cdot y\cdot x$ (\isn{ir\_bool\_conv}),
  and two universal introductions at every extension give
  $\Dec(\isn{restrRel}\,R\,p)$.
\item \emph{Clause (iii)} is clause (i) at $Q = \isn{restrRel}\,R\,p$,
  which is closed (\isn{closed\_restrRel}) and decided by step~3.
\end{enumerate}
The contrast is the trace example's up body
(\cref{thm:app-premise}): there the inner leaf is a $\mu$-atom whose
premise is valueless in $\Mstar$, here it is a path atom that the
dichotomy decides.

Fair termination, the fairness twin, the transformation as a term former, the alternating-bit protocol and the grounded fairness rule are treated in Appendix~\ref{app:fair}.

\paragraph{Lemmas a re-mechanization needs.}
\begin{itemize}
\item \isn{g\_acc\_intro0}: accessibility introduced by one fold from
  ``every successor is accessible''.
\item \isn{R1\_edge\_dec}, \isn{R1\_edge\_init}, \isn{R1p\_back},
  \isn{R1p\_of\_R1}: the edges of the two relations.
\item \isn{g\_R1\_decRel}, \isn{g\_R1p\_decRel}: both relations are
  decided in \rga.
\item \isn{r1Var\_E}, \isn{r1Var\_absurd}: the variant property's
  instance at an edge, and its refutation in \rga.
\item \isn{neg\_decided\_explosive\_iff}: at a decided $a$, $\neg a \ent
  0 = 1$ exactly when $a$ is derivable.
\item \isn{g\_path\_coind}, \isn{g\_path\_child\_intro}, \isn{g\_run\_of\_path}
  (Appendix~\ref{app:streams}).
\item \isn{tr\_back\_ind}: backward induction from a witness of $q$.
\item \isn{tr\_pAdd\_0}, \isn{tr\_pAdd\_S}, \isn{tr\_pAdd\_shift}: the
  addition equations used for transport.
\item \isn{tr\_div\_add}: divergence propagates along the search.
\item \isn{tr\_ndiv\_of\_ex}: a witness refutes divergence, by negative
  folds, with no rule.
\item \isn{tr\_all\_B}: $\mathrm{All}_q(x)$ decided (ATI).
\item \isn{g\_mu\_not\_nu}, \isn{g\_dec\_nu}: the generic exclusion and
  decided $\nu$-atom.
\item \isn{trmod\_step\_inv}, \isn{trmod\_step\_true}: the step relation
  in $\Mstar$ is the successor function.
\item \isn{tr\_false\_inv}: the falsity invariant of the reach atom.
\item \isn{tr\_reach\_lfp\_fst0}: every state that reaches $\langle
  0,0\rangle$ has first component $0$.
\item \isn{nfdec\_leafdec}, \isn{trmod\_up\_leafdec}: a true premise
  values every leaf, here the reach leaf.
\item \isn{trmod\_allqN\_B\_not\_false}: a closure of a decidedness
  atom is never false.
\item \isn{agafPhi}, \isn{agafD}, \isn{nfdec\_agafPhi},
  \isn{restrRel}, \isn{g\_restrRel\_decRel}, \isn{g\_agaf\_premise}:
  AG\,AF\,$p$ with a path-atom leaf, decided at a decided relation.
\end{itemize}

\section{The two games}
\label{app:games}

Files: \isn{OIGA\_Coind/GDet.thy} (the game body, determinacy, the
Baire instance and player II's strategy),
\isn{OIGA\_Coind/GDetModel.thy} (the game atom in $\Mstar$),
\isn{OIGA\_Coind/GBucket.thy} (the bucket game),
\isn{OIGA\_Coind/GBucketModel.thy} (the Stepmother's legal pour and
the model side), \isn{OGA/OGA\_Coind\_Kit.thy} (tuples, folds over five
buckets, multiplication and the dispatch on a phase).

\paragraph{Definitions: the game body.}
A game has three parameters: a target set $U$, a move relation $M$
for player I and a reply relation $\mathit{Rp}$ for player II, both
child-first: $M(q, p)$ says $q$ is a legal move of I at $p$, and
$\mathit{Rp}(r, q)$ that $r$ is a legal reply of II at $q$.  Player I's
body (\isn{gmPhi}) is
\[
  U(p) \vee \exists q.\,\bigl(M(q, p) \wedge \forall r.\,(\neg
  \mathit{Rp}(r, q) \vee a\cdot r)\bigr),
\]
with the existential written $\neg\forall\neg$ and the parameters
lifted past the two binders.  Its dual (\isn{gmPsi}) replaces $a\cdot r$
by $\neg(a\cdot r)$.  $\WinI = \mu(\text{body})$ (\isn{gmWinI}) and
$\WinII = \nu(\neg\text{dual})$ (\isn{gmWinII}).  Unfolded, the body of
$\WinII$ at $p$ reads $\neg U(p) \wedge \forall q.\,(M(q,p) \to \exists
r.\,\mathit{Rp}(r, q) \wedge a\cdot r)$: II is outside the target and
answers every move with a reply back in the region.  The body's
decidedness premise (\isn{gmDec}) is $\Dec\,U \wedge (\Dec\,M \wedge
(\mathit{Rp}_3 \wedge 0 = 0))$, where $\mathit{Rp}_3$ (\isn{gmRp3}) is
the reply leaf decided at every triple of numbers, the outermost
variable unused.  The structural facts \isn{posP\_gmPhi},
\isn{free1\_gmPhi}, \isn{ndual\_gmPhi} and \isn{nfdec\_gmPhi} (body
positive, mentions $a$, has that dual and that premise) are proved in
\rga's metatheory by one grammar rule per connective.

\paragraph{Definitions: the Baire instance and the strategy.}
Positions are sequence codes, the root is the empty code $1$, and a
move or reply appends one number: $\mathrm{Ext}(q, p) \equiv q =
p{\frown}\mathrm{last}(q)$ (\isn{gmExt}), decided in \rga
(\isn{g\_gmExt\_decRel}).  In this instance $M = \mathit{Rp} =
\mathrm{Ext}$.  Write $W$ for $\WinII$ at the target $U$.  A code $q$
is \emph{good} for a set $V$ (\isn{gmGood}) when $V(\mathrm{par}(q))
\wedge \mathrm{Ext}(q, \mathrm{par}(q))$: $q$ is a move of I from a
position of $V$.  The strategy's matrix (\isn{gmSigB}) is
$\phi(q, n) \equiv W(q{\frown}n) \vee \neg\mathrm{good}_W(q)$, and the
strategy $\sigma$ (\isn{gmSig}) is the search combinator at $\phi$: the
least $n$ satisfying $\phi(q, n)$.  The second disjunct makes the
matrix satisfiable at every $q$.

\paragraph{Definitions: the bucket game.}
Five buckets stand in a circle, each of capacity $c$ units.  At each
round the Stepmother (player I) pours one unit, split as she likes, and
Cinderella (player II) empties two adjacent buckets.  Amounts are
numerators over a common denominator $d$.  A state is an eight-tuple
$\langle \mathit{ph}, d, k, n_0, \ldots, n_4\rangle$ of nested Cantor
pairs (\isn{csTup}): the phase of the round-robin clock, the
denominator, the last refinement factor and the five numerators, read
by \isn{csPh}, \isn{csDen}, \isn{csFac} and \isn{csBk}.  The initial
state $\mathit{init}$ (\isn{csInit}) has phase $0$, denominator and factor $1$, and
empty buckets.  The target is empty: $U(x) \equiv (1 = 0)$ (\isn{csU}).
The Stepmother's move $M(q, p)$ (\isn{csM}) is the conjunction of: $k_q
\ne 0$; $d_q = k_q\cdot d_p$; $\mathit{ph}_q = \mathit{ph}_p$; $k_q\cdot
n_i(p) \le n_i(q)$ for every $i < 5$ (each bucket scales and may grow);
and the five increments $n_i(q) - k_q\cdot n_i(p)$ sum to $d_q$.  The
order is truncated subtraction ($a \le b$ is $a \mathbin{\dot{-}} b = 0$); the
refinement factor is stored in the successor, so no quantifier occurs.
Cinderella's reply $\mathit{Rp}_c(r, q)$ (\isn{csRp}) is the
conjunction of: no overflow at $q$, that every bucket of $q$ is at most
$c\cdot d_q$ (\isn{csNoOvf}, with $c\cdot d$ the $c$-fold sum
\isn{cdSum}); $d_r = d_q$; $k_r = k_q$; the phase advances, $4$ to $0$
and $a$ to $a+1$ otherwise (\isn{csAdv}); and for some $j < 5$, $r$
empties buckets $j$ and $j+1 \bmod 5$ of $q$ and copies the rest
(\isn{csEmpB}).  So overflow is tested after the pour and before the
emptying; at an overflowing state II has no legal reply, and the move
clause of I's body holds there vacuously.  The schedule empties at phase
$r$ the pair starting at $\pi(r) = 2r \bmod 5$ (\isn{csPi}).  The table
$G(i, r)$ (\isn{csG}) has rows $[2,0,1,0,1]$, $[1,0,1,2,0]$,
$[1,2,0,1,0]$, $[0,1,0,1,2]$, $[0,1,2,0,1]$ indexed by the bucket $i$;
it satisfies $G(i, r+1 \bmod 5) = 0$ if $i$ is in the pair of phase
$r$ and $G(i, r) + 1$ otherwise (\isn{cs\_T1}), and $G(i, r) + 1 \le 3$
(\isn{cs\_T2}).  The certificate $\mathit{Cert}$ (\isn{csW}) holds at
$p$ when, for some $r_0 < 5$, the phase of $p$ is $r_0$ and every
bucket satisfies $n_i(p) \le G(i, r_0)\cdot d_p$ (\isn{csWr}); it need
not be decided.  The strategy term $\mathrm{sig}$ (\isn{csSig}) dispatches on the phase
$r_0$ of $q$ (\isn{cdSel}) to the tuple with phase $r_0 + 1 \bmod 5$,
the denominator and factor of $q$, and the buckets of $q$ with the pair
of phase $r_0$ emptied (\isn{csSigma}).  The two regions at capacity $c$
are $\WinI$ and $\WinII$ at these three parameters, written
$\WinI^{(c)}$ and $\WinII^{(c)}$ (\isn{csWinI}, \isn{csWinII}).

\begin{thm}[\isn{gm\_det}, \isn{gm\_decII}; grounded]
\label{thm:app-det}
Let $\Gamma$ be well-formed, with $\Gamma \ent \Dec\,U$, $\Gamma \ent
\Dec\,M$, $\Gamma \ent \Dec\,\mathit{Rp}$ (the relations decided at
every pair) and $\Gamma \ent x\,\mathbf N$.
(i) In \isn{gcalc\_dich}: $\Gamma \ent \WinI\cdot x \vee \WinII\cdot x$.
(ii) In \isn{gcalc\_gfix\_dich}: $\Gamma \ent (\WinII\cdot x)\,\mathbf B$.
Companions: in \isn{gcalc\_gfix}, from $\Gamma$ well-formed, $x\,\mathbf
N$ and $\WinI\cdot x$, $\neg\WinII\cdot x$ (\isn{gm\_excl}, no
decidedness hypotheses); in \isn{gcalc}, the body's decidedness premise
(\isn{gmDec}) from the three decidedness hypotheses (\isn{g\_gmDec}).
\end{thm}

\paragraph{Proof of \cref{thm:app-det}.}
\begin{enumerate}
\item \emph{The premise} (\isn{g\_gmDec}, no rule beyond \rga's).  Three nested
  universal introductions at every extension (\isn{g\_allq\_env})
  introduce $j_1, j_2, j_3$; at the innermost level $\Dec\,\mathit{Rp}$
  instantiated at $j_3, j_2$ (\isn{g\_decRel\_inst}) is the reply leaf
  ($j_1$ unused).  Then $0 = 0$, and conjunction introduction with
  $\Dec\,U$ and $\Dec\,M$.  The check forbids all four rules here.
\item \emph{Clause (i).}  One application of the dichotomy at the game
  body, with dual \isn{gmPsi}, premise \isn{gmDec} from step 1, the four
  structural facts and $x\,\mathbf N$.  That application is the whole
  proof; the check forbids $\nu$, ATI and ground induction.
\item \emph{Exclusion} (\isn{gm\_excl}): \isn{g\_mu\_not\_nu} at the
  body.  The lemma's proof is ground induction at $\WinI\cdot x$ with motive
  $c(v) \equiv \neg\WinII\cdot v$.  Closure obligation: at a fresh $i$
  with $i\,\mathbf N$ and the body at $i$ with $c$ for $a$, derive
  $c(i)$.  \isn{g\_ndual\_beta} turns that body instance into the dual's
  instance at $\WinII$; double negation and the negative fold
  \isn{g\_nu\_fold\_neg} give $\neg\WinII\cdot i$.
\item \emph{Clause (ii)} (\isn{g\_dec\_nu}).  Step 2 gives $\WinI\cdot x
  \vee \WinII\cdot x$.  In the left case step 3 gives $\neg\WinII\cdot
  x$; in the right case $\WinII\cdot x$; either way $\WinII\cdot x \vee
  \neg\WinII\cdot x$, which is $(\WinII\cdot x)\,\mathbf B$
  (\isn{ir\_boolI}).  The check requires the dichotomy and ground
  induction and forbids $\nu$ and ATI.  The calculus has no $\nu$ and
  no ATI as a primitive rule; it derives ATI (\cref{thm:ati-dich}),
  the proof does not apply it, and \isn{GDet} is built where that
  derivation is not in scope.
\end{enumerate}
The Baire instances (\isn{gm\_baire\_det}, \isn{gm\_baire\_excl},
\isn{gm\_baire\_decII}) are the three results at $M = \mathit{Rp} =
\mathrm{Ext}$.  In the two concrete Baire games, I wins at the root in
one round when the target is ``player I's last move was $0$, at a
position other than the root''
(\isn{gm\_I\_one\_round}, two folds), and II wins at the root when the
target is ``the last two moves agree, at a position other than the
root'' by answering every move $m$ with
$m+1$, the certificate being the complement of the target
(\isn{gm\_II\_copies}, one coinduction step).

\begin{thm}[\isn{gm\_sigII}; grounded]
\label{thm:app-sigII}
In \isn{gcalc\_gfix\_dich}, at the Baire instance, with $\Gamma$ well-formed and $\Gamma \ent \Dec\,U$:
(i) $\Gamma \ent \Tot\,\sigma$;
(ii) if $\Gamma \ent W\cdot p$, $\Gamma \ent p\,\mathbf N$ and $\Gamma
\ent m\,\mathbf N$, then $\Gamma \ent W\cdot\bigl((p{\frown}m){\frown}
(\sigma\cdot(p{\frown}m))\bigr)$.
\end{thm}

\paragraph{Proof of \cref{thm:app-sigII}.}
Work in the context extended by a number variable $v$, the argument of
the strategy, with $U$ and $W$ lifted (\isn{lift\_gmWinII}).
\begin{enumerate}
\item \emph{The region is decided.}  At every extension and number $y$,
  $(W\cdot y)\,\mathbf B$: \cref{thm:app-det}(ii) at $\Dec\,U$ and the
  decidedness of $\mathrm{Ext}$.  This step is where the dichotomy and
  ground induction enter.
\item \emph{Goodness is decided}, from step 1 at $\mathrm{par}(v)$, the
  decidedness of $\mathrm{Ext}$, and conjunction.
\item \emph{The matrix is decided} (\isn{g\_Dec\_lamI}): at every number
  $t$, $W\cdot(v{\frown}t)$ is decided by step 1 and $\neg\mathrm{good}_W
  (v)$ by step 2.
\item \emph{A witness exists.}  Split on the decided
  $\mathrm{good}_W(v)$.  If not good, $n = 0$ satisfies the matrix by
  its second disjunct.  If good, $W\cdot\mathrm{par}(v)$ unfolds once
  (\isn{g\_nu\_unfold}, no rule) to II's body at $\mathrm{par}(v)$, and
  $v$ is a move from $\mathrm{par}(v)$; the body's move clause
  (\isn{gm\_psiE\_M}) supplies $r$ with $\mathrm{Ext}(r, v)$ and $W\cdot
  r$; then $r = v{\frown}\mathrm{last}(r)$ (\isn{gmExt\_eq}), so $W\cdot
  (v{\frown}\mathrm{last}(r))$ by transport, and $\mathrm{last}(r)$
  witnesses the matrix.
\item \emph{Clause (i).}  The search lemma \isn{mu\_choice} on steps 3
  and 4 gives $\Tot\,\sigma$.
\item \emph{The search satisfies the matrix}: \isn{g\_mu\_choice\_all}
  gives $\forall q.\,\phi(q, \sigma\cdot q)$.
\item \emph{Clause (ii).}  Let $q = p{\frown}m$.  Instantiate step 6 at
  $q$.  Since $\mathrm{par}(q) = p$ (\isn{fbPar\_snoc}), the hypothesis
  $W\cdot p$ and $\mathrm{Ext}(q, p)$ (\isn{gmExt\_snoc}) make $q$
  good; disjunction elimination against goodness (\isn{gmSigB\_E}) gives
  $W\cdot(q{\frown}\sigma\cdot q)$.
\end{enumerate}
No coinduction step occurs: the region is a hypothesis and is only
unfolded, and the check forbids $\nu$.  No ATI step occurs either:
decidedness comes from the dichotomy and ground induction (step 1),
and the search lemma \isn{mu\_choice} is \rga's, stated in the base
calculus (Appendix~\ref{app:streams}).

\begin{thm}[\isn{gm\_winI\_not\_false}; in $\Mstar$]
\label{thm:app-winInotfalse}
Let $S$ be a set of naturals with $\mathrm{Surv}_{M,\mathit{Rp}}(S)$
(\isn{gmSurvI}): for every $p \in S$ some $q$ has $M(\bar q, \bar p)$
not false in $\Mstar$, and every $r$ with $\mathit{Rp}(\bar r, \bar q)$
true in $\Mstar$ lies in $S$.  Then for every $p \in S$ and every $U$,
$\WinI\cdot\bar p$ is not false in $\Mstar$.
\end{thm}

\paragraph{Proof of \cref{thm:app-winInotfalse}.}
The proof is an induction over the mutual inductive definition of the
model's head and value relations (\isn{gm\_false\_inv}), with the
trivial motive on head judgments and, on value judgments $(t, v)$, the
invariant \isn{gmFInv}: $(t, v)$ is none of ten shapes.  Write $L_1(p)$
for the instance at $p$ of the outer universal's matrix and $L_2(q, p)$
for the inner universal's matrix (\isn{gmL1}, \isn{gmL2}).  For $p \in S$: (1)
$\WinI\cdot\bar p$ false; (2) the body instance at $p$ false; (3) the body instance's
negated universal $\neg\forall\lambda L_1(p)$ false; (4) $\forall\lambda
L_1(p)$ true.  For $p, q$ with $M(\bar q, \bar p)$ not false and every
true reply of $q$ in $S$: (5) $(\lambda L_1(p))\cdot\bar q$ true; (6) its
beta reduct $\neg M(\bar q, \bar p) \vee \neg\forall\lambda L_2(q,p)$
true.  For $p, q$ with every true reply of $q$ in $S$: (7)
$\neg\forall\lambda L_2(q,p)$ true; (8) $\forall\lambda L_2(q,p)$ false;
(9) $(\lambda L_2(q,p))\cdot\bar r$ false; (10) its reduct $\neg
\mathit{Rp}(\bar r, \bar q) \vee \WinI\cdot\bar r$ false.  The cases:
\begin{enumerate}
\item Negation true from its argument false: only (7) has that form,
  and its argument would be (8).  Negation false from its argument
  true: (3) would need (4).
\item Disjunction true from its left disjunct: (6) would need $\neg M(
  \bar q, \bar p)$ true, so $M(\bar q, \bar p)$ false
  (\isn{omvO\_not\_inv}), against the shape's hypothesis.  From its right
  disjunct: (6) would need (7).
\item Disjunction false from both disjuncts false: (2) would need (3);
  (10) would need $\mathit{Rp}(\bar r, \bar q)$ true, so $r \in S$, and
  $\WinI\cdot\bar r$ false, which is (1) at $r$.
\item Application of an abstraction: the head of $\WinI$ is the
  abstraction of the body with $\WinI$ for $a$ (\isn{gm\_head\_W}, by
  determinism of heads), so (1) reduces to (2) (\isn{gm\_L0\_eq}), (5) to
  (6) (\isn{gm\_L1\_eq}) and (9) to (10) (\isn{gm\_L2\_eq}).
\item Universal true from every instance true: (4) would need
  $(\lambda L_1(p))\cdot\bar q$ true at the move $q$ that survival
  supplies, which is (5).  This case is the only use of the existential
  in the survival hypothesis.  Universal false from one instance false:
  (8) would need (9).
\item Recursor heads, $\nu$-heads and the remaining clauses produce no
  shape.
\end{enumerate}
The theorem is shape (1).  Read operationally: a falsity derivation of
$\WinI$ is well-founded, and at a surviving position it would have to
pass through shapes (1) to (10) and return to (1) at a new position of
$S$, forever.

\begin{thm}[\isn{cs\_cind}, \isn{cs\_cind\_init}; grounded]
\label{thm:app-cind}
In \isn{gcalc\_nu}, for a capacity $c \ge 3$ (a HOL natural), with
$\Gamma$ well-formed: if $\Gamma \ent x\,\mathbf N$ and $\Gamma \ent
\mathit{Cert}\cdot x$, then $\Gamma \ent \WinII^{(c)}\cdot x$.  In
particular $\Gamma \ent \WinII^{(c)}\cdot\mathit{init}$.
\end{thm}

\paragraph{Proof of \cref{thm:app-cind}.}
\begin{enumerate}
\item \emph{One coinduction step} (\isn{gm\_winII\_coind}): the
  coinduction rule at II's body with certificate $\mathit{Cert}\cdot v$,
  so the region is the set of certificate states, and an eigenvariable
  fresh for the context, the certificate's abstraction and the body at
  once.  Closure
  obligation: at every $j$ not free in $\Gamma$, from $j\,\mathbf N$ and
  $\mathit{Cert}\cdot j$, derive II's body at $j$ with the certificate
  for $a$.
\item \emph{Introducing II's body} (\isn{gm\_psiI}, in \isn{gcalc}).
  From $\neg U(j)$, $\Dec\,M$, $j\,\mathbf N$, and, at every extension
  and number $q$ with $M(q, j)$, a number $r$ with $\mathit{Rp}(r, q)$ and
  $\mathit{Cert}\cdot r$, the body follows: universal introduction over
  $q$, a case split on the decided $M(q, j)$, and in the move case the
  existential introduced at the supplied $r$.  Here $\neg U(j)$ is
  \isn{csU\_never} and $\Dec\,M$ is \isn{g\_csM\_decRel}; the
  decidedness of the move relation serves this case split, not a
  dichotomy.
\item \emph{The reply.}  Take $r = \mathrm{sig}\cdot q$, a number
  (\isn{csSig\_N}).  The step theorem \isn{cs\_sig\_step} (in \isn{gcalc},
  at $c \ge 3$), from $\mathit{Cert}\cdot j$ and $M(q, j)$, gives
  $\mathit{Rp}_c(\mathrm{sig}\cdot q, q)$ and $\mathit{Cert}\cdot(
  \mathrm{sig}\cdot q)$.  The step theorem is the \rga lemma \isn{cs\_sig\_step0}
  carried by cut:
  \begin{enumerate}
  \item The certificate at $j$ is a five-way disjunction over $r_0$
    (\isn{csAny\_E}); fix the case $r_0$: the phase of $j$ is $r_0$ and
    $n_i(j) \le G(i, r_0)\cdot d_j$.  The move keeps the phase, so $q$
    has phase $r_0$ (\isn{cs\_M\_parts}).
  \item The strategy evaluates at phase $r_0$ to the tuple for phase $r_0$
    (\isn{cs\_sig\_eval}: one beta step and the dispatch at a numeral
    phase), whose fields are read off (\isn{cs\_sigma\_acc}).
  \item \emph{One pour adds at most one unit} (\isn{cs\_pour\_q}): from
    $n_i(j) \le g\cdot d_j$ and the move, $n_i(q) \le (g+1)\cdot d_q$.
    With $k = k_q$: $k\cdot n_i(j) \le n_i(q)$; the increment is at most
    the sum of the five increments (\isn{cs\_summand\_le}), which is $d_q
    = k\cdot d_j$; the arithmetic core \isn{cs\_pour\_bound} then gives
    $n_i(q) = k\cdot n_i(j) + (n_i(q) - k\cdot n_i(j)) \le k\cdot(g\cdot
    d_j) + k\cdot d_j = (g+1)\cdot(k\cdot d_j)$, using $k\cdot(g\cdot d)
    = g\cdot(k\cdot d)$ (\isn{cs\_pMul\_cdSum}, by HOL induction on $g$).
  \item \emph{The reply is legal} (\isn{cs\_step\_Rp}).  No overflow:
    $n_i(q) \le (G(i, r_0) + 1)\cdot d_q \le 3\cdot d_q \le c\cdot d_q$,
    by step (c), \isn{cs\_T2} and $c \ge 3$; this is the only use of the
    capacity hypothesis.  Denominator and factor are copied.  The phase
    advances, by cases on $r_0 = 4$.  The emptied pair is $\pi(r_0)$.
  \item \emph{The certificate holds at the reply} (\isn{cs\_step\_W}) at
    phase $r_0 + 1 \bmod 5$: an emptied bucket is $0$; any other bucket
    equals $n_i(q) \le (G(i, r_0) + 1)\cdot d_q$, which is $G(i, r_0 + 1
    \bmod 5)\cdot d_r$ by \isn{cs\_T1} and $d_r = d_q$.
  \end{enumerate}
  One beta step turns $\mathit{Cert}\cdot r$ into the certificate's
  abstraction applied to $r$, as step 2 requires.
\item \emph{The initial state}: $\mathit{init}$ is a number and
  satisfies the certificate at phase $0$ with empty buckets
  (\isn{cs\_W\_init}).
\end{enumerate}
The coinduction rule is used once, in step 1, and every other step is
in \isn{gcalc}; the check forbids the dichotomy, ATI and ground
induction under the theorem, all four rules under the strategy, the
step and the kit, and ATI and \isn{mu\_choice} under every result of the
theory.  The strategy is total in \rga (\isn{cs\_sig\_total}).

\begin{thm}[\isn{cs\_rr\_tight}; grounded]
\label{thm:app-tight}
In \rga, at every context $\Gamma$.  Let $\mathrm{play}(0) =
\mathit{init}$ and $\mathrm{play}(n+1) = \mathrm{sig}\cdot
\mathrm{pour}(\mathrm{play}(n))$ (\isn{csPlay}), the play in which the
Stepmother pours the whole unit into bucket $0$ every round
(\isn{csPour}, \cref{thm:app-pour}) and Cinderella answers by the
round-robin strategy.
(i) For every $n$, $\Gamma \ent M(\mathrm{pour}(\mathrm{play}(n)),
\mathrm{play}(n))$: every pour of the play is legal.
(ii) For every $n < 5$, $\Gamma \ent \mathit{Rp}_2(\mathrm{play}(n+1),
\mathrm{pour}(\mathrm{play}(n)))$: the first five replies are legal at
capacity $2$.
(iii) $\Gamma \ent n_0(\mathrm{pour}(\mathrm{play}(5))) = 3$.
(iv) $\Gamma \ent d_{\mathrm{pour}(\mathrm{play}(5))} = 1$.
(v) For every capacity $c < 3$ and every well-formed term $r$:
$\mathit{Rp}_c(r, \mathrm{pour}(\mathrm{play}(5))),\ \Gamma \ent 0 = 1$.
\end{thm}
So after the sixth pour bucket $0$ holds three units over the
denominator one, and at capacity two or less every candidate reply of
Cinderella's is refuted: the Stepmother has won that play.  The bound
$c \ge 3$ of \cref{thm:app-cind} is therefore attained against the
strategy it certifies.  The theorem says nothing about other
strategies or about the game's threshold.  The proof is by
computation: the fields of each state of the play are read off by
\rga's rules (\isn{cs\_play\_acc}), and for clause (v) the proof
unfolds the reply relation and refutes its no-overflow conjunct, which
concerns the poured position (bucket $0$ holds three units, above the
capacity), so the refutation holds whatever the reply.

\begin{thm}[\isn{cs\_stepmother\_refuted}; grounded]
\label{thm:app-refuted}
In \isn{gcalc\_gfix\_nu}, for $c \ge 3$ and $\Gamma$ well-formed:
$\WinI^{(c)}\cdot\mathit{init},\ \Gamma \ent 0 = 1$.
\end{thm}

\paragraph{Proof of \cref{thm:app-refuted}.}
\begin{enumerate}
\item The extended context is well-formed and proves $\mathit{init}\,
  \mathbf N$.
\item Exclusion (\isn{cs\_excl}, step 3 of the proof of \cref{thm:app-det} at the
  bucket parameters) under the hypothesis gives $\neg\WinII^{(c)}\cdot
  \mathit{init}$.  Ground induction.
\item \cref{thm:app-cind} in the extended context gives $\WinII^{(c)}\cdot
  \mathit{init}$.  One coinduction step.
\item Negation elimination gives $0 = 1$.
\end{enumerate}
The check requires ground induction and $\nu$ and forbids the dichotomy
and ATI.  \oinga reads the theorem through the edge from
\isn{gcalc\_ag\_nu} to \isn{gcalc\_gfix\_nu}.  At every capacity, with no
capacity hypothesis, every state is won by one of the players
(\isn{cs\_det}, in \isn{gcalc\_dich}) and Cinderella's region is decided
(\isn{cs\_decII}, in \isn{gcalc\_gfix\_dich}, so with ground induction as
well): \cref{thm:app-det} at the decidedness facts \isn{g\_csU\_Dec},
\isn{g\_csM\_decRel} and \isn{g\_csRp\_decRel}.

\begin{thm}[\isn{cs\_pour\_legal}, \isn{cs\_init\_surv}; grounded and in $\Mstar$]
\label{thm:app-pour}
(i) In \rga: if $\Gamma \ent p\,\mathbf N$, then $\Gamma \ent
M(\mathrm{pour}(p), p)$, where $\mathrm{pour}(p)$ (\isn{csPour}) keeps the
phase and denominator of $p$, has factor $1$, adds $d_p$ to bucket $0$
and copies the other buckets.
(ii) In $\Mstar$, at every capacity $c$:
$\mathrm{Surv}_{M,\mathit{Rp}_c}(\mathbb N)$.
\end{thm}

\paragraph{Proof of \cref{thm:app-pour}.}
\begin{enumerate}
\item \emph{Clause (i).}  Read the fields of $\mathrm{pour}(p)$
  (\isn{cs\_acc\_tup}) and use $1\cdot y = y$ (\isn{cs\_pMul\_oneL}, number
  induction).  The five conjuncts: $1 \ne 0$; $d_p = 1\cdot d_p$; the
  phase is copied; $n_0(p) \le n_0(p) + d_p$ and $n_i(p) \le n_i(p)$; the
  increments are $d_p$ and four zeros, whose five-fold sum collapses to
  $d_p$ (\isn{cs\_sum5\_unfold}).  The pour has no overflow precondition
  and is legal at every number $p$.
\item \emph{Clause (ii).}  Fix $p$; with $S = \mathbb N$ the reply
  condition is trivial.  Clause (i) at $\bar p$ is an empty-context
  derivation in \rga, so $M(\mathrm{pour}(\bar p), \bar p)$ is true in
  $\Mstar$ (\isn{rga\_nu\_sound}).  $\mathrm{pour}(\bar p)$ is derivably a
  number, so it has a numeral value $\bar q$ (\isn{cs\_N\_num\_rel}); value
  congruence under application (\isn{ovcrO\_appI}) and transport
  (\isn{omvO\_vcr\_ex}) make $M(\bar q, \bar p)$ true; by uniqueness of
  values it is not false.
\end{enumerate}

\begin{thm}[\isn{cs\_winI\_not\_false}, \isn{dinga\_stepmother\_not\_refuted}; in $\Mstar$, classical]
\label{thm:app-notrefuted}
(i) At every capacity $c$, $\WinI^{(c)}\cdot\mathit{init}$ is not false in
$\Mstar$.
(ii) At every capacity $c$, \dinga does not derive $\neg\WinI^{(c)}\cdot
\mathit{init}$ from the empty context (\isn{oingad\_stepmother\_not\_refuted}
for \oinga with the dichotomy, restated at \dinga in the junction
theory \isn{CD\_Junction}).
\end{thm}

\paragraph{Proof of \cref{thm:app-notrefuted}.}
\begin{enumerate}
\item \emph{Clause (i).}  $\mathit{init}$ is derivably a number, so it
  has a numeral value $\bar q$ (\isn{cs\_N\_num\_rel}); if
  $\WinI^{(c)}\cdot\mathit{init}$ were false, value congruence would
  make $\WinI^{(c)}\cdot\bar q$ false; but
  \cref{thm:app-winInotfalse} with \cref{thm:app-pour}(ii) and $q \in
  \mathbb N$ says that atom is not false.
\item \emph{Clause (ii).}  An empty-context derivation of the negation
  makes $\neg\WinI^{(c)}\cdot\mathit{init}$ true in $\Mstar$
  (\isn{oingad\_nu\_sound}); by \isn{omvO\_not\_inv} the atom then has a
  value $b$ with $\neg b$ true, so the atom is false, against clause (i).
\end{enumerate}
Overflow plays no part here: the Stepmother survives everywhere
because the legal pour exists at every position.

\paragraph{Where ATI is spent.}
No result of this appendix applies ATI as a rule.  The results stated
in a calculus with ground induction and the dichotomy (\isn{gm\_decII},
\isn{gm\_sigII}, \isn{cs\_decII}) have ATI as a derived rule
(\cref{thm:ati-dich}) but not as a primitive one, and their theories
are built where that derivation is not in scope.  \cref{thm:app-sigII} uses the search lemma \isn{mu\_choice}, which is \rga's search carried into the base calculus by cut; the check for the strategy requires \isn{mu\_choice} and the dichotomy by name and forbids ATI and $\nu$.
The bucket game's check forbids ATI and \isn{mu\_choice} under every
result of its theory.

\paragraph{A stuck position is lost.}
\isn{gm\_winI\_stuck\_false} (\isn{OIGA\_Coind/GDetModel.thy}): if the target is false at $p$ in $\Mstar$ and every move from $p$ is false in $\Mstar$, then $\WinI\cdot\bar p$ is false in $\Mstar$.  One unfolding: at every $q$ the first disjunct of the move clause, ``the move is not legal'', is true (\isn{omvO\_notT}, \isn{omvO\_orT1}, \isn{gm\_L1\_eq}); so the universal over moves is true (\isn{omvO\_A\_lamI}) and its negation, the existential, false; so the body is false (\isn{omvO\_orF}, \isn{gm\_L0\_eq}), and \isn{omvO\_appLam} with \isn{omhO\_trmGfix} carries the value to the atom.

\paragraph{Lemmas a re-mechanization needs.}
\begin{itemize}
\item \isn{posP\_gmPhi}, \isn{free1\_gmPhi}, \isn{ndual\_gmPhi},
  \isn{nfdec\_gmPhi}: the side conditions of the dichotomy at the game
  body.
\item \isn{g\_gmDec}: the premise from the decidedness of the three
  parameters, in \rga.
\item \isn{g\_allq\_env}, \isn{g\_decRel\_inst}: universal introduction
  at every extension; a decided relation's instance.
\item \isn{g\_mu\_not\_nu}, \isn{g\_ndual\_beta}, \isn{g\_dec\_nu}: the
  generic exclusion and decided $\nu$-atom.
\item \isn{gm\_psiE\_M}, \isn{gm\_psiI}: eliminating and introducing II's
  body at a move.
\item \isn{gm\_winII\_coind}: coinduction for II's region at a
  certificate set.
\item \isn{gmExt\_eq}, \isn{gmExt\_snoc}, \isn{g\_gmExt\_decRel}: the Baire
  extension relation.
\item \isn{mu\_choice}, \isn{g\_mu\_choice\_all}: the search
  (Appendix~\ref{app:streams}).
\item \isn{gmSigB\_E}: the matrix at a good code gives the region.
\item \isn{gm\_head\_W}, \isn{gm\_L0\_eq}, \isn{gm\_L1\_eq}, \isn{gm\_L2\_eq},
  \isn{gm\_false\_inv}: the head of $\WinI$ in $\Mstar$, the three
  instance steps, and the falsity invariant.
\item \isn{csU\_never}, \isn{g\_csU\_Dec}, \isn{g\_csM\_decRel},
  \isn{g\_csRp\_decRel}: the target never holds, and the three bucket
  parameters are decided in \rga.
\item \isn{cs\_sig\_eval}, \isn{cs\_sigma\_acc}, \isn{cs\_M\_parts}: the
  strategy at a numeral phase, its fields, and the move's parts.
\item \isn{cs\_pour\_q}, \isn{cs\_pour\_bound}, \isn{cs\_pMul\_cdSum}: one
  pour adds at most one unit of the new denominator.
\item \isn{cs\_step\_Rp}, \isn{cs\_step\_W}, \isn{cs\_T1}, \isn{cs\_T2},
  \isn{cs\_W\_init}: the reply is legal and re-establishes the
  certificate; the table laws; the certificate at the start.
\item \isn{cs\_pour\_legal}, \isn{cs\_N\_num\_rel}, \isn{ovcrO\_appI}: the
  legal pour; numeral values of derivable numbers; value congruence
  under application.
\end{itemize}

\section{Fair termination and eventual delivery through random assignment}
\label{app:fair}

Files, all in \isn{OIGA\_Coind}: \isn{GFairCorr.thy} (imports
\isn{Main}; runs, fairness, the exact and the weak edge and their
correspondences), \isn{GFairX.thy} (imports \isn{GFairCorr}; the relaxed
masked edge, its correspondence, the eventuality transfer),
\isn{GFair.thy} (imports \isn{GLattice}; the two-component loop),
\isn{GFairModel.thy} (imports \isn{GFair}, \isn{GFairCorr} and
\isn{GTraceModel}; the loop in the least model), \isn{GFairTwin.thy}
(imports \isn{GFair}; the fairness twin), \isn{GFairTwinModel.thy}
(imports \isn{GFairTwin} and \isn{GFairModel}), \isn{GFairGen.thy}
(imports \isn{GFairTwin}; the term former, its kit, the motive pack,
the helpful direction, the regression), \isn{GFairGenModel.thy}
(imports \isn{GFairGen}, \isn{GFairX} and \isn{GFairModel}; the generic
least-model side), \isn{GFairAbp.thy} (imports \isn{GFairGen}; the
alternating-bit protocol), \isn{GFairAbpModel.thy} (imports
\isn{GFairAbp} and \isn{GFairGenModel}), \isn{GFairLift.thy} (imports
\isn{GFairGen}; the grounded fairness rule and the eventuality rule).
The generic lemmas come from \isn{GLattice.thy}, \isn{GStream.thy},
\isn{GBranch.thy}, \isn{GRand.thy}, \isn{OIGA/GDich.thy},
\isn{OIGA/GShell.thy} and the tuple kit \isn{OGA/OGA\_Coind\_Kit.thy};
the least-model facts from \isn{RGA/RGA\_Omega\_Sound.thy},
\isn{IGA/IGA\_Sound.thy} and \isn{IGA/IGA\_Lfp.thy}.

\paragraph{Conventions.}
For a closed step relation $T$, $T\cdot y\cdot x$ says that $y$ is a
successor of $x$, the orientation of accessibility.  The eventually atom
$\mathrm{Ev}_G$ (\isn{trEv}) takes its relation source first, so the
term former below produces one relation in two orientations.  The order
$a \le b$ is the equation $a \mathbin{\dot-} b = 0$ for truncated
subtraction (\isn{gapSub}), and $a < b$ is $a+1 \le b$.  Codes are
nested Cantor pairs built by the tuple kit of \isn{OGA\_Coind\_Kit}:
$\langle t_0, \ldots, t_{n-1}\rangle$ (\isn{csTup}), read by field
(\isn{csGet}), with the kit's finite conjunction and disjunction over a
HOL-indexed family of terms (\isn{csAll}, \isn{csAny}).  A premise holds
\emph{at every extension} when it is required in every well-formed
context that derives everything $\Gamma$ derives; universal introduction
and the case lemmas deliver their cases in such extensions.  Results of
\rga (locale \isn{rga}) hold in every calculus by \isn{g\_viaLW}; the
rows below are those transfers.

\subsection{The classical correspondence}

\paragraph{Definitions.}
A program has $n$ components, with guards $g_i$ and one-step relations
$a_i$ (source first) on a state type.  A run is a state sequence $\rho$
with a schedule $\sigma$ such that at every $k$ the component $\sigma_k$
is enabled at $\rho_k$ and steps from $\rho_k$ to $\rho_{k+1}$
(\isn{aoRun}); only infinite runs exist in this vocabulary.  A run is
\emph{strongly fair} (\isn{aoSF}) when every component enabled
infinitely often is scheduled infinitely often, and \emph{weakly fair}
(\isn{aoWF}) when every component enabled from some point on is
scheduled infinitely often (\isn{aoSF\_aoWF}: strong implies weak).  A
\emph{mask} marks the components owed fairness, and \isn{aoSFM}
restricts strong fairness to them.  A transformed state is a pair $(s,
z)$ of a program state and a counter vector.  In each of the three edges
below, the scheduled component $i$ is enabled at $s$ and steps $s$ to
$s'$, and its own counter is unconstrained, Apt and Olderog's random
assignment $z_i \mathrel{:=} {?}$.
\begin{itemize}
\item The \emph{exact edge} (\isn{aoEdgeAt}, \isn{aoEdge}): every other
  enabled component has a positive counter (the strict guard) and its
  counter decreases by exactly one; the counter of every other disabled
  component is kept.
\item The \emph{weak edge} (\isn{aoEdgeW}): the counter of a disabled
  component other than $i$ is reassigned at random instead.
\item The \emph{relaxed masked edge} (\isn{aoEdgeX}): for every masked
  $j \ne i$, the counter strictly decreases if $j$ is enabled, which
  forces it to have been positive, and does not increase if $j$ is
  disabled; unmasked counters are free.  The exact edge is a relaxed edge
  for the full mask (\isn{aoEdge\_aoEdgeX}).
\end{itemize}
The \emph{pending set} (\isn{aoPend}) of $j$ at time $k$ is the set of
steps $m \ge k$ at which $j$ is enabled and such that $j$ is not
scheduled at any step from $k$ through $m$; the counter (\isn{aoCnt}) is its size, and the
lift (\isn{aoLift}) of a run at $k$ is $\rho_k$ with the counters at $k$.
The weak counter (\isn{aoWCnt}) counts the steps before $j$ is next
scheduled or next disabled; the masked lift (\isn{aoXLift}) sets
unmasked counters to $0$.

\begin{thm}[\isn{ao\_sound}, \isn{ao\_complete}; classical]
\label{thm:app-fair-corr}
(i) Every infinite chain of the exact edge projects onto an infinite
strongly fair run (\isn{ao\_sound}); of the weak edge, onto a weakly fair
one (\isn{ao\_wsound}); of the relaxed masked edge, onto one strongly
fair for the mask (\isn{ao\_xsound}).
(ii) Every infinite strongly fair run lifts to an infinite chain of the
exact edge with the same projection (\isn{ao\_complete}); likewise for
the weak edge (\isn{ao\_wcomplete}) and the relaxed edge
(\isn{ao\_xcomplete}).
(iii) If no infinite chain of the exact edge starts at $(s, z)$ for any
$z$, no infinite strongly fair run starts at $s$ (\isn{ao\_fair\_term});
likewise for the weak and the relaxed edge (\isn{ao\_fair\_wterm},
\isn{ao\_fair\_xterm}), and conversely a chain gives a fair run
(\isn{ao\_fair\_div}, \isn{ao\_fair\_wdiv}, \isn{ao\_fair\_xdiv}).
(iv) If every infinite relaxed chain from every $(s, z)$ reaches a state
whose program part satisfies $P$, every infinite run from $s$ strongly
fair for the mask reaches $P$ (\isn{ao\_fair\_xev}).
\end{thm}
No \rga term occurs; the checks of \isn{GFairCorr} and \isn{GFairX}
forbid every dependency whose name mentions \isn{rga}.  The
correspondence holds for infinite runs only: the strict guard blocks
where two enabled components both have counter $0$, and the program
does not.

\paragraph{Proof of \cref{thm:app-fair-corr}.}
\begin{enumerate}
\item \emph{Projection} (\isn{ao\_sound}).  Choose $\sigma_k$ from each
  edge (choice in HOL).  Along a chain on which $j$ is not scheduled from
  $k_0$ on, $z_j(k_0+d)$ plus the number of steps in $[k_0, k_0+d)$ at
  which $j$ is enabled equals $z_j(k_0)$ (\isn{ao\_cnt\_invariant},
  induction on $d$).  If $j$ is enabled infinitely often, that count
  exceeds $z_j(k_0)$ for some $d$ (\isn{ao\_card\_unbounded}).  The weak
  and masked forms use $z_j(k_0+d) + d = z_j(k_0)$ on a stretch where $j$
  is always enabled (\isn{ao\_wcnt\_invariant}) and an inequality
  (\isn{ao\_xcnt\_invariant}).
\item \emph{Finiteness of the pending set} (\isn{aoPend\_finite}), the
  one place a classically undecidable question is split.  If $j$ is
  scheduled at some $N \ge k$, the set lies below $N$; if not, strong
  fairness makes $j$ enabled only finitely often, and the set lies below
  the last such step.  The definition of the counter makes no case split.
\item \emph{One step} (\isn{aoPend\_step}, \isn{aoCnt\_step},
  \isn{aoCnt\_pos}).  When $\sigma_k \ne j$, the pending set at $k$ is
  the set at $k+1$ with $k$ added if $j$ is enabled at $k$, and the same
  set otherwise: the counter decreases by one or is kept, and an enabled
  $j$ has a positive counter.
\item \emph{Lifting} (\isn{ao\_complete}).  By step~3 each lifted step is
  an exact edge.  The weak set is finite by weak fairness
  (\isn{aoWPend\_finite}); the masked lift needs finiteness only from the
  fairness clause of $j$ itself (\isn{aoPend\_finite\_at},
  \isn{aoCnt\_xstep}, \isn{aoCnt\_xpos}).
\item \emph{Transfers.}  Lift a fair run from $s$ at the counter vector
  it induces at time $0$, or project a chain.  For (iv), the lift keeps
  the program part (\isn{ao\_fair\_xev}).
\end{enumerate}

\subsection{The two-component loop}

\paragraph{Definitions.}
The program is $b \mathrel{:=} \mathit{true};\ x \mathrel{:=} 0;\
\mathbf{do}\ b \to x \mathrel{:=} x+1 \mathrel{\square} b \to b
\mathrel{:=} \mathit{false}\ \mathbf{od}$, after Apt and Olderog and
Grumberg, Francez, Makowsky and de Roever: always scheduling component~$1$
diverges, and every strongly fair run terminates.  The transformed
relation \isn{aoT} is written by hand at the \rga level, on codes
$\langle \mathit{pc}, x, z_1, z_2\rangle$ ($\mathit{pc}$ $0$ initial,
$1$ loop, $2$ exit; fields \isn{aoPc}, \isn{aoXv}, \isn{aoZ1},
\isn{aoZ2}).  Its body \isn{aoB} has three disjuncts, written by
projection equations on the successor $y$ and the source $x$: the
initial edge ($\mathit{pc}$ from $0$ to $1$, $x_y = 0$, both counters
free); component~$1$ ($\mathit{pc}$ stays $1$, $x_y = x_x+1$, $z_2(y) +
1 = z_2(x)$, $z_1$ free); component~$2$ ($\mathit{pc}$ from $1$ to $2$,
$x$ kept, $z_1(y) + 1 = z_1(x)$, $z_2$ free).  The strict guard and the
decrement are one equation.  Both components are enabled at every loop
code, so the update is exact.

\begin{thm}[\isn{g\_acc\_ao\_all}, \isn{aomod\_fair\_term}; grounded and classical]
\label{thm:app-fair-loop}
(i) In \isn{gcalc}: for every well-formed $\Gamma$ and $\Gamma \ent
x\,\mathbf N$, $\Gamma \ent \Acc_{\isn{aoT}}\cdot x$; and $\Gamma \ent
\mathrm{decRel}(\isn{aoT})$ (\isn{g\_ao\_decRel}).
(ii) Classically, through \rga's soundness for the least model: the HOL
loop has no infinite strongly fair run from $b = \mathit{true}$, $x =
0$.
\end{thm}

\paragraph{Proof of \cref{thm:app-fair-loop}.}
\begin{enumerate}
\item \emph{Decidedness and the case lemma.}  The edge is decided at
  numbers (\isn{aoT\_B}), with no ATI.  From $\Gamma \ent
  \isn{aoB}\,y\,x$ and three case premises, each at every extension and
  receiving the equations of one disjunct, any well-formed $d$ follows
  (\isn{aoB\_cases}); no tuple rebuild lemma is needed, because the
  motive below quantifies over every code with the given fields.
  Impossible cases close by $a = 0$ against $a = b+1$ (\isn{ao\_ne01})
  and $a = 1$ against $a = 2$ (\isn{ao\_ne12}).
\item \emph{Introduction and the exit} (\isn{acc\_ao\_I},
  \isn{acc\_ao\_exit}).  If at every extension every number $j$ with
  $\isn{aoB}\,j\,x$ is accessible, so is $x$: one fold of the
  $\mu$-atom by conversion, universal introduction, the decided edge
  split.
  At an exit code every case is refuted.
\item \emph{The loop step} (\isn{acc\_ao\_loop\_step}).  A loop code $u$
  is accessible if every loop code $j$ with $z_2(j)+1 = z_2(u)$ is: the
  initial case is refuted, component~$1$ is the hypothesis, and
  component~$2$ lands at an exit code.
\item \emph{The loop induction} (\isn{acc\_ao\_loop}).  The closed motive
  $\isn{aoMot}\cdot m \equiv \forall u.\,\neg(\mathit{pc}\,u = 1) \vee
  \neg(z_2(u) = m) \vee \Acc_{\isn{aoT}}\cdot u$ (\isn{aoMot\_I},
  \isn{aoMot\_E}), and one internal number induction (\isn{ir\_nat\_ind})
  on $m$.  At $0$ every component-$1$ successor would have $z_2 + 1 = 0$;
  at $m+1$ it has $z_2 = m$, covered by the hypothesis.  The universal
  over the whole code absorbs the free $z_1$ and the change of $x$.
\item \emph{Every code} (\isn{acc\_ao\_all}).  One introduction with no
  case split on the source: every successor is a loop code, accessible at
  its own $z_2$ (\isn{acc\_ao\_loop\_at}), or an exit code.
\item \emph{The rows.}  In \isn{gcalc}: \isn{g\_acc\_ao\_all},
  \isn{g\_acc\_ao\_init}, \isn{g\_ao\_decRel}.  In \isn{gcalc\_gfix}: no
  path (\isn{g\_ao\_not\_path}), and no run of \isn{aoT}, the transformed
  relation and not the loop, visits the initial code
  (\isn{g\_ao\_no\_run}).  In \isn{gcalc\_dich}: $\Acc \vee \Path$
  (\isn{g\_ao\_acc\_or\_path}); in \isn{gcalc\_gfix\_dich}: the path atom
  decided (\isn{g\_ao\_dec\_path}).  \oinga and \dinga readings.
\item \emph{The bridge.}  HOL states $(b, x)$, HOL component $0$ the
  increment and $1$ the exit, and the encoding \isn{aoEnc} at
  $\mathit{pc} = 1$ or $2$ by $b$ (\isn{aomod\_tup\_eq},
  \isn{aomod\_fields}).  A HOL exact edge is derivable in \rga at the
  encoded numerals (\isn{aomod\_edge\_der}) and a non-edge is refuted
  (\isn{aomod\_edge\_neg}); through soundness for the least model
  (\isn{aomod\_true}, by \isn{rga\_omega\_sound}), \isn{aoT} is true at
  two encoded numerals exactly when the HOL edge holds
  (\isn{aomod\_edge\_iff}).
\item \emph{Clause (ii)} (\isn{aomod\_fair\_term}).  By
  \isn{ao\_fair\_term}, exclude for each $z$ an infinite chain of exact
  edges from $((\mathit{true}, 0), z)$.  \isn{acc\_ao\_init} at the empty
  context puts the numeral of the initial code in HOL's accessible part
  of the relation true in the least model (\isn{acc\_omv\_sound}); the
  initial edge of \isn{aoT} reaches the encoded start code at every $z$,
  so every start code is accessible; a chain of HOL edges descends along
  that relation (step~7), and an accessible point has no such chain
  (\isn{aomod\_acc\_chain}).  Only the one-directional soundness of
  accessibility is used; \isn{aomod\_acc\_iff} is an extra.  For this
  program the HOL statement also has a direct proof; the point is the
  pipeline.
\end{enumerate}
The rank at the initial code is $\omega$, finite at every loop code (a
reading of the relation and of the proof).  The measure
is the postponed component's counter $z_2$, never the counter $z_1$
that component~$1$ resets.

\subsection{The fairness twin}

\paragraph{Definitions.}
The program is $b \mathrel{:=} \mathit{true};\ c \mathrel{:=}
\mathit{false};\ x \mathrel{:=} 0;\ \mathbf{do}\ b \to x \mathrel{:=}
x+1;\ c \mathrel{:=} \neg c \mathrel{\square} b \wedge c \to b
\mathrel{:=} \mathit{false}\ \mathbf{od}$.  A run that never schedules
component~$2$ enables it at every other step: infinitely often, never
continuously.  Codes are $\langle \mathit{pc}, x, z_1, z_2, c\rangle$.
Four edge kinds: the initial edge (\isn{aoTwI}: $\mathit{pc}$ to $1$, $x
= c = 0$, counters free); component~$1$ at $c = 0$, where component~$2$
is disabled (\isn{aoTwA}: $x$ incremented, $c$ set to $1$, counters
free); component~$1$ at $c = 1$ (\isn{aoTwC}: $x$ incremented, $c$ set
to $0$, $z_2(y)+1 = z_2(x)$); component~$2$ at $c = 1$ (\isn{aoTwD}:
$\mathit{pc}$ to $2$, $z_1(y)+1 = z_1(x)$).  The strong body
\isn{aoTsB} adds $z_2(y) = z_2(x)$ to the second kind, keeping the
disabled component's counter as the exact edge does; the weak body
\isn{aoTwB} omits that conjunct, as the weak edge does.  That conjunct is
the whole difference between \isn{aoTs} and \isn{aoTw}, both relations decided
(\isn{aoTs\_B}, \isn{aoTw\_B}).  The certificate \isn{aoCert} holds at a
code that is initial, a loop code with $c = 0$, or a loop code with $c
= 1$ and $z_2$ a successor.

\begin{thm}[\isn{g\_acc\_aoTs\_all}, \isn{g\_aoTw\_div}; grounded and classical]
\label{thm:app-fair-twin}
Let $\Gamma$ be well-formed.
(i) In \isn{gcalc}: every number is accessible for \isn{aoTs}
(\isn{g\_acc\_aoTs\_all}; the initial code, \isn{g\_acc\_aoTs\_init}).
(ii) In \isn{gcalc\_nu}: $\Gamma \ent \Path_{\isn{aoTw}}\cdot
\isn{aoTwInit}$.
(iii) In \isn{gcalc\_gfix\_nu\_dich}: a run of \isn{aoTw}, a term total in $\Mstar$, from
the initial code, every value on the path (\isn{g\_aoTw\_run}).
(iv) Classically: neither \oinga nor \dinga derives
$\Acc_{\isn{aoTw}}\cdot\isn{aoTwInit}$ from the empty context
(\isn{oinga\_aoTw\_not\_acc}, \isn{oingad\_aoTw\_not\_acc}, restated at
\dinga as \isn{dinga\_aoTw\_not\_acc} in \isn{CD\_Junction}), nor does
\rga (\isn{aomod\_tw\_weak\_not\_der}).
(v) Classically: the HOL twin has no infinite strongly fair run from its
initial state, and has an infinite weakly fair one
(\isn{aomod\_tw\_separation}).
\end{thm}

\paragraph{Proof of \cref{thm:app-fair-twin}.}
\begin{enumerate}
\item \emph{Clause (i).}  A four-way case lemma (\isn{aoTsB\_cases}).  A
  loop code with $c = 1$ is accessible if every loop code with $z_2$ one
  less is (\isn{acc\_aoTs\_c1\_step}: the second kind is refuted by $c$,
  the third is the hypothesis, the fourth exits).  Any loop code is, under
  the same hypothesis (\isn{acc\_aoTs\_loop\_step}): a successor of the
  second kind keeps $z_2$ and has $c = 1$, and the hypothesis transports
  along $z_2(y) = z_2(u)$.  Then one number induction on $z_2$ at a motive
  of the shape of \isn{aoMot} (\isn{acc\_aoTs\_loop}); the case on $c$
  comes from the case lemma, not from a second induction.
\item \emph{Clause (ii)} (\isn{g\_aoTw\_div}).  One coinduction step
  (\isn{g\_path\_coind}) at the certificate, from \isn{aoCert\_init}.
  Closure at a fresh variable: split the certificate's three disjuncts
  and name a successor in each disjunct by \isn{g\_path\_child\_intro}.  The
  successors are $\langle 1,0,0,0,0\rangle$ from the initial code; $\langle
  1, x+1, 0, 1, 1\rangle$ from $c = 0$, the counter of the disabled
  component~$2$ reassigned to $1$; $\langle 1, x+1, 0, p, 0\rangle$ from
  $c = 1$ and $z_2 = p+1$ (\isn{aoTw\_succ0}, \isn{aoTw\_succ1},
  \isn{aoTw\_succ2}, \isn{aoCert\_succ}).  The pattern is that of
  \isn{g\_R1p\_div}.
\item \emph{Clause (iii).}  \isn{g\_run\_of\_path} at the decided relation
  (\isn{g\_aoTw\_decRel}) and the path of clause (ii).
\item \emph{Clause (iv).}  An empty-context derivation of the
  accessibility, with \isn{g\_acc\_not\_path} and clause (ii), would derive
  $0 = 1$, against \isn{oinga\_consistent\_notation} (itself proved
  through $\Mstar$), or \isn{oingad\_consistent\_notation}.  For \rga:
  the weakly fair run of clause (v) lifts by \isn{ao\_wcomplete} to a
  chain of weak edges, every edge an \rga edge of \isn{aoTw} at the encoded
  numerals (\isn{aomod\_tw\_wedge\_der}) and so true in the least model,
  so the start is not accessible there (\isn{aomod\_tw\_weak\_not\_acc}),
  and soundness gives the non-derivability.
\item \emph{Clause (v).}  No strongly fair run: the chain of
  \cref{thm:app-fair-loop}(ii), with \isn{acc\_aoTs\_init} and the bridge for the exact
  edge (\isn{aomod\_tw\_edge\_der}, \isn{aomod\_tw\_fair\_term}).  The
  weakly fair run always schedules component~$1$, $\rho_k =
  (\mathit{true}, k, \mathit{odd}\,k)$; the run is weakly and not strongly
  fair (\isn{aomod\_tw\_weak\_run}), a fact about the HOL program alone,
  with no least model on that half.
\end{enumerate}
\isn{ao\_wcomplete} is used only for the non-derivability in \rga;
\isn{ao\_wsound} has no downstream use.  The twin's bridge is
one-directional; no refutation of a non-edge is proved or needed.  One
reading is withheld: the path atom of clause (ii) read in $\Mstar$ as an
infinite HOL chain, which would need the truth set of a path atom
turned into a sequence by dependent choice, with the edge decided at
numerals in $\Mstar$; the development has no such lemma, and the
classical weak divergence is read off the HOL run instead.  The rank at
the initial code under \isn{aoTs} is $\omega$ (a reading).

\subsection{The transformation as a term former}

\paragraph{Definitions.}
A program is a HOL list of components $(G_i, A_i, f_i)$: a closed guard
($G_i\cdot s$ a sentence about a program-state code), a closed step
($A_i\cdot s'\cdot s$, successor first) and a HOL fairness flag, the
mask (\isn{aoGd}, \isn{aoStp}, \isn{aoFr}; \isn{aoOk} says every guard
and step is closed and well-formed).  A transformed code is $\langle s,
z_0, \ldots, z_{n-1}\rangle$, read by \isn{aoSt} (written $\mathrm{st}$)
and \isn{aoZc} (written $z_j$).  The counter law (\isn{aoCtr}) of a
masked component, from $x$ to $y$, is
\[
  (G\cdot\mathrm{st}(x) \wedge z(y)+1 \le z(x)) \vee (\neg G\cdot
  \mathrm{st}(x) \wedge z(y) \le z(x)),
\]
and $0 = 0$ for an unmasked one.  The edge of component $i$
(\isn{aoCompB}) is $G_i\cdot\mathrm{st}(x) \wedge A_i\cdot\mathrm{st}(y)
\cdot\mathrm{st}(x)$ and the counter law of every other component, a
kit conjunction; the body \isn{aoRelB} is the kit disjunction over $i$.
The \emph{term former} is the HOL function \isn{aoRel} from the list to
the closed term $\lambda y\,x.\,\isn{aoRelB}\,y\,x$; \isn{aoRelF} is
the same relation in the forward orientation, for the eventually atom.
A program has \emph{decided components} (\isn{aoDecCs}) when, in every
context, every guard is decided at every number and every step at every
pair of numbers.  The \emph{motive pack} is
\[
  \isn{aoMotL}\,P\,F\,Q\cdot m = \forall u.\,\neg P\cdot u \vee \neg(F
  \cdot u \le m) \vee Q\cdot u,
\]
with $P$, $F$, $Q$ lifted as wholes, so they may have free variables (a
parametric target, or a predicate that mentions an eigenvariable); the
pack is not a closed motive.  \isn{aoPrd} reads a program-state
predicate on the program part, and \isn{aoZf} turns a counter into a
function.

\begin{thm}[\isn{aoRel\_cases}, \isn{acc\_aoRel\_helpful}; grounded]
\label{thm:app-fair-gen}
In \rga, for every list with \isn{aoOk} and every well-formed $\Gamma$:
(i) From $\Gamma \ent \isn{aoRelB}\,y\,x$ and, at every extension and
every HOL index $i$, a case premise receiving $G_i\cdot\mathrm{st}(x)$,
$A_i\cdot\mathrm{st}(y)\cdot\mathrm{st}(x)$ and the counter law of every
masked $j \ne i$, the cases' conclusion follows.
(ii) With decided components, the relation is decided (\isn{aoRel\_B},
\isn{g\_aoRel\_decRel}).
(iii) For well-formed $P$, $F$, $Q$: if $P\cdot v$ is decided and $F\cdot v$ a number at every
extension, and $P$ is progressive for $Q$ along $F$ (at every extension,
a $u$ with $P\cdot u$ has $Q\cdot u$ once every $v$ with $P\cdot v$ and
$F\cdot v < F\cdot u$ has $Q\cdot v$ at every further extension), then
every number $u$ with $P\cdot u$ has $Q\cdot u$ (\isn{aoMotL\_wf}).
(iv) Let the components be decided, $h$ masked, and $I$, $E$
well-formed, with, for numbers $s$ and $s'$:
\begin{itemize}
\item[IB] $I\cdot s$ is decided, in every context;
\item[P1] $I\cdot s$ gives $G_h\cdot s$, in every context;
\item[P2] for every $i \ne h$, $I\cdot s$, $G_i\cdot s$ and $A_i\cdot
  s'\cdot s$ give $I\cdot s' \vee E\cdot s'$, in every context;
\item[P3] $I\cdot s$ and $A_h\cdot s'\cdot s$ give $E\cdot s'$, in every
  context;
\item[P4] at every extension of $\Gamma$, every code whose program part
  satisfies $E$ is accessible.
\end{itemize}
Then every number $u$ with $\Gamma \ent I\cdot\mathrm{st}(u)$ has
$\Gamma \ent \Acc_{\isn{aoRel}}\cdot u$.
\end{thm}
Clause (iv) is the one-level case of the rule of Grumberg, Francez,
Makowsky and de Roever; its rank is $\omega$ (a reading).

\paragraph{Proof of \cref{thm:app-fair-gen}.}
\begin{enumerate}
\item \emph{The kit.}  Closedness and head steps of both orientations,
  edge introduction and elimination (\isn{aoRel\_I}, \isn{aoRel\_E},
  \isn{aoRelF\_I}, \isn{aoRelF\_E}), and \isn{aoRelB\_I}, which builds an
  edge from a HOL index, its guard, its step and the other counter laws.
  The counter law is introduced when enabled with a strict decrease
  (\isn{aoCtr\_I\_en}) or disabled without increase (\isn{aoCtr\_I\_dis});
  eliminated, an enabled component strictly decreases (\isn{aoCtr\_E\_en})
  and no counter increases (\isn{aoCtr\_E\_le}).
\item \emph{Clauses (i) and (ii).}  Kit disjunction elimination, then
  \isn{aoCompB\_D} and kit conjunction elimination; cases indexed by HOL
  numbers need one step lemma per component.  Decidedness by HOL
  induction over the list, with no ATI.  Accessibility is introduced from
  every successor by the decided edge split (\isn{acc\_aoRel\_I}), and a
  code with no enabled component is accessible (\isn{acc\_aoRel\_idle}).
\item \emph{Clause (iii).}  One internal number induction on $m$ at the
  pack, with a fresh variable.  At $0$ every hypothesis case is refuted,
  since $F\cdot v + 1 \le 0$ fails (\isn{ao\_suc\_le0}); at $m+1$ the
  hypothesis case reduces to $F\cdot v \le m$, closed by \isn{aoMotL\_E}.
  The instance at $m = F\cdot u$ finishes: well-founded induction on a
  counter at a decided predicate costs one number induction.
\item \emph{Clause (iv).}  Clause (iii) at $P = \isn{aoPrd}\,I$, $F =
  \isn{aoZf}\,h$, $Q = \Acc$.  At a code $x$ with $I$, introduce and split
  the edge to $y$ by clause (i).  If $i = h$, P3 and P4 close.  If $i \ne
  h$, P2 gives $I$ or $E$ at $y$; $E$ closes by P4; in the $I$ branch $h$
  is enabled at the source (P1) and masked, so \isn{aoCtr\_E\_en} gives
  $z_h(y) + 1 \le z_h(x)$ (\isn{aoZf\_eq}), and the induction hypothesis
  closes.
\item \emph{The rows}, at every program with decided components: in
  \isn{gcalc}, \isn{g\_aoRel\_decRel}; in \isn{gcalc\_gfix}, from a
  derived accessibility, no path (\isn{g\_aoRel\_not\_path}) and no run
  through the code (\isn{g\_aoRel\_no\_run}); in \isn{gcalc\_dich}, $\Acc
  \vee \Path$ (\isn{g\_aoRel\_acc\_or\_path}); in \isn{gcalc\_gfix\_dich},
  the path atom decided (\isn{g\_aoRel\_dec\_path}).
\item \emph{The regression} (\isn{aoLp}).  The loop as a list over
  two-field program codes $\langle \mathit{pc}, x\rangle$: both guards
  $\mathit{pc} = 1$ (\isn{aoLpG}), steps \isn{aoLpA1} and \isn{aoLpA2},
  both masked, exit set $\mathit{pc} = 2$ (\isn{aoLpE}).  Exit codes are
  idle; loop codes are clause (iv) at $h$ the exit component
  (\isn{acc\_aoLp\_loop}); every code by a decided split on the guard
  (\isn{g\_acc\_aoLp\_all}).  The regression has no initial edge, and the
  relation \isn{aoT} is not re-derived literally: its code is flat, with
  an initial edge and the exact update.  The twin is not re-derived, and
  its weak transformation lies outside the term former, which has only
  the strong, masked law.  The regression has no least-model reading.
\end{enumerate}

\subsection{The least-model side of the term former}

For HOL guards $g$, steps $a$ and an encoding $\mathit{enc}$ of program
states, \isn{aoGnEnc} encodes $(s, z)$ as the nested Cantor pair of
$\mathit{enc}(s), z_0, \ldots, z_{n-1}$ (\isn{aoGnCode}).  The
\emph{bridge premises}, for every component $i$, with derivations in \rga
from no hypotheses at encoded numerals, are: BG, if $g_i(s)$ then $G_i
\cdot\overline{\mathit{enc}(s)}$ is derivable; BN, if $i$ is masked and
$\neg g_i(s)$ then $\neg G_i\cdot\overline{\mathit{enc}(s)}$ is
derivable; BA, if $g_i(s)$ and $a_i(s, s')$ then the step is derivable.

\begin{thm}[\isn{aomod\_gen\_fair\_term}, \isn{aomod\_gen\_fair\_ev}; classical]
\label{thm:app-fair-genmodel}
For every list with \isn{aoOk}, under BG, BN and BA:
(i) if \rga derives, from no hypotheses, $\Acc_{\isn{aoRel}}$ at the
encoded code of $(s_0, z)$ for every $z$, the HOL program has no infinite
run from $s_0$ strongly fair for the mask;
(ii) for a closed target $G$ whose truth in the least model at an
encoded code implies $P$ of the program part, if \rga derives, from no
hypotheses, $\mathrm{Ev}_G$ along \isn{aoRelF} at the encoded code of
$(s_0, z)$ for every $z$, every infinite run from $s_0$ strongly fair
for the mask reaches a state satisfying $P$.
\end{thm}
The pipeline reads \rga derivations from the empty context only,
through \isn{rga\_omega\_sound}; a derivation that uses ground induction
in \oinga is not read by the pipeline.  No decidedness premise appears.

\paragraph{Proof of \cref{thm:app-fair-genmodel}.}
\begin{enumerate}
\item \emph{The code and the edge at numerals.}  The kit tuple of
  numerals is the numeral of the code, read back field by field
  (\isn{aomod\_gn\_tup\_eq}, \isn{aomod\_gn\_fields}).  A relaxed masked
  HOL edge gives an \rga derivation of \isn{aoRelB} at the encoded
  numerals, in both orientations (\isn{aomod\_gen\_edge\_der}): the guard
  by BG, the step by BA, and per masked $j \ne i$ a HOL case on $g_j$,
  \isn{aoCtr\_I\_en} with BG or \isn{aoCtr\_I\_dis} with BN.  The
  refutation of a non-edge is not needed.
\item \emph{Clause (i).}  By \isn{ao\_fair\_xterm}, exclude an infinite
  relaxed chain from each $(s_0, z)$.  The derived accessibility is true
  (\isn{aomod\_true}), so the start lies in HOL's accessible part of the
  true relation (\isn{acc\_omv\_sound}, which needs closedness and
  well-formedness); each
  chain edge is true by step~1, against \isn{aomod\_acc\_chain}.
\item \emph{The eventually atom in the least model}
  (\isn{aomod\_ev\_chain}).  For closed $T$ and $G$, if $\mathrm{Ev}_G\cdot
  \bar x$ is true, every sequence from $x$ with true edges meets $G$.  The
  atom's truth set is the least fixed point of its operator
  (\isn{gfix\_omv\_eq\_lfp}, positivity \isn{posP\_evPhi}).  The set of
  points from which every true chain meets $G$ is closed under the
  operator: the $G$ branch is index $0$; in the universal branch, the
  instance at the chain's second point has a false first disjunct, since
  that edge is true, so the recursive atom puts the second point in the
  set (\isn{aomod\_os\_varT}, \isn{aomod\_num\_inj}).
\item \emph{Clause (ii)} (\isn{aomod\_gen\_fair\_ev}).  \isn{ao\_fair\_xev}
  reduces the claim to relaxed chains; the derived atom is true, each
  edge is true in the forward orientation (step~1), step~3 gives an index
  where $G$ is true, and the hypothesis on $G$ gives $P$.
\end{enumerate}

\subsection{The alternating-bit protocol}

\paragraph{Definitions.}
The program-state code is $\langle i, a, b, p, q\rangle$: the message
index, the sender's bit, the receiver's expected bit, and the data and
acknowledgement channels, each $0$ (empty) or $c+1$ (holding bit $c$).
The invariant (\isn{aoAbInvH}) on the control part $(a, b, p, q)$ is $(a
= b \wedge q \ne a+1) \vee (a \ne b \wedge p \ne b+1)$: a stale bit is
never taken for a fresh one; its $24$ tuples of $36$ are listed
(\isn{aoAbInvL}, \isn{aoAbInvL\_char}).  Six components: send ($p
\mathrel{:=} a+1$), lose on the data channel, deliver on the data
channel ($p \mathrel{:=} 0$, and $b \mathrel{:=} \neg b$ at a match $p =
b+1$), acknowledge ($q \mathrel{:=} \neg b + 1$), lose on the
acknowledgement channel, deliver on the acknowledgement channel ($q \mathrel{:=} 0$, and $a
\mathrel{:=} \neg a$, $i \mathrel{:=} i+1$ at a match $q = a+1$); the
mask \isn{aoAbFrH} leaves the two losses unconstrained.  The grounded
guard of $j$ (\isn{aoAbGd}) is membership of the control in the table
of invariant tuples enabling $j$ (\isn{aoAbL}), so every guard contains
the invariant; steps follow the HOL tables \isn{aoAbStH} and
\isn{aoAbIncH}.  Phase~$2$ ($a \ne b$) splits into the \emph{full-good}
tuples with $q = a+1$ (\isn{aoAbFG2}), the \emph{empty} ones with $q =
0$ (\isn{aoAbEM2}) and the \emph{stale} ones with $q = b+1$
(\isn{aoAbST2}); phase~$1$ ($a = b$) likewise on $p$.  The program
\isn{aoAbPg} has decided components (\isn{aoAbPg\_dec}).
\isn{aoAbPc}\,$X\,k$ holds at a code with index $k$, an arbitrary term,
and control in the table $X$.  The target $\isn{aoAbG}\cdot K\cdot y$ is
$K \le i_y$.  The sentence \isn{aoAbSent} is
\[
  \forall K.\,\forall u.\,\neg(\isn{aoAbInv}\cdot\mathrm{st}(u) \wedge
  i(\mathrm{st}(u)) = 0) \vee \mathrm{Ev}_{\isn{aoAbG}\cdot K}\cdot u,
\]
the atom along $\isn{aoRelF}\,\isn{aoAbPg}$.

\begin{thm}[\isn{aoAb\_every}, \isn{aomod\_ab\_deliver}; grounded and classical]
\label{thm:app-abp}
(i) In \rga: from $\Gamma$ well-formed, $\Gamma \ent K\,\mathbf N$,
$\Gamma \ent u\,\mathbf N$ and $\Gamma \ent \isn{aoAbInv}\cdot
\mathrm{st}(u)$, $\Gamma \ent \mathrm{Ev}_{\isn{aoAbG}\cdot K}\cdot u$,
at any index and any counters.
(ii) In \rga: $\Gamma \ent \isn{aoAbSent}$ (\isn{aoAb\_sent}; in every
calculus, \isn{g\_aoAb\_sent}).
(iii) Classically, through \rga's soundness for the least model: every
infinite run of the HOL protocol from index $0$, both bits $0$ and both
channels empty, strongly fair for send, acknowledge and the two
deliveries, with the losses unconstrained, reaches an index at least
$K$, for every $K$.
\end{thm}
The sentence covers the codes with the invariant and index $0$; clause
(i) covers every invariant code, as a derivation at each $K$ and $u$,
and in every calculus at a numeral $K$ (\isn{g\_aoAb\_sent\_at}).  Read
through the atom, clause (i) says that every run of the transformed
protocol from $u$ reaches an index at least $K$; the atom also holds
vacuously at a code without successor.

\paragraph{Proof of \cref{thm:app-abp}.}
\begin{enumerate}
\item \emph{The tables} (HOL, one \isn{simp} each).  Every invariant
  tuple enabling $j$ steps into the invariant (\isn{aoAbInv\_presH}), and
  the invariant is the union of the two phases (\isn{aoAbInv\_split}).  At
  full-good phase-$2$ tuples the acknowledgement delivery increments the
  index and lands in the invariant, and every other step keeps the index
  and lands at full-good or empty tuples; at empty tuples acknowledge
  lands at full-good ones and the acknowledgement delivery is disabled; a
  stale delivery lands at empty tuples (\isn{aoAbFG2\_tab},
  \isn{aoAbEM2\_tab}, \isn{aoAbST2\_tab}).  Phase~$1$ mirrors these with
  the data delivery, which at full-good tuples lands in phase~$2$ with
  the index kept, the progress edge.
\item \emph{The successor lemma} (\isn{aoAb\_succ}).  An edge from a code
  whose control is the HOL tuple $t$ hands its case, at every extension,
  a component $j$ with $t$ in its table, the successor control, the index
  law and the counter laws (\isn{aoRel\_cases}, \isn{aoAbStp\_cases}; a
  different tuple at the same code is refuted by \isn{aoAbIs\_uniq}).
  Every later case analysis is the successor lemma plus a table fact.
\item \emph{The nested induction} (\isn{aoAb\_nest}).  For a decided
  closed $T$, an arbitrary target $G$, decided sets FG and EM, a delivery
  counter $F_d$, a refill counter $F_r$ and a progress set $I_k$ giving
  the atom at every extension: if every edge from FG leads to $I_k$, or to
  FG or EM with smaller $F_d$, and every edge from EM leads to FG with
  $F_d$ not larger, or to EM with $F_d$ not larger and $F_r$ smaller, then
  every number in FG or EM satisfies $\mathrm{Ev}_G$.  An outer
  \isn{aoMotL\_wf} on $F_d$; inside the outer induction, the full codes by \isn{ev\_next};
  then the empty codes by an inner \isn{aoMotL\_wf} on $F_r$ at the
  predicate ``EM and $F_d$ at most $F_d\cdot u$'' (\isn{aoAbLe}), which
  mentions the outer code $u$: a parametric motive.
\item \emph{Phase~$2$} (\isn{aoAb\_phase2}).  If at every extension every
  invariant code with index $k+1$ satisfies the atom, so does every
  phase-$2$ code with index $k$.  The full-good and empty codes by step~3,
  with $F_d$ the acknowledgement delivery's counter, $F_r$ acknowledge's,
  and the case lemmas \isn{aoAb\_cFG2} (by \isn{aoCtr\_E\_en}) and
  \isn{aoAb\_cEM2} (by \isn{aoCtr\_E\_le}).  The stale codes by one more
  \isn{aoMotL\_wf} on the delivery counter: a stale delivery lands at an
  empty code with its counter reset at random, covered for every counter.
\item \emph{Phase~$1$ and the round} (\isn{aoAb\_phase1},
  \isn{aoAb\_round}).  Phase~$1$ is the mirror, with the data delivery's
  and send's counters, \isn{aoAb\_cFG1}, \isn{aoAb\_cEM1}, and step~4 as
  the progress continuation.  The round splits the invariant codes at
  index $k$ by phase.
\item \emph{Clause (i)} (\isn{aoAb\_every}).  \isn{aoMotL\_wf} at the
  invariant and the distance $K \mathbin{\dot-} i$ (\isn{aoAbDk}), $K$ a
  parameter of the distance and of the target.  At $x$, split the
  distance: at $0$, $K \le i_x$ (\isn{aoAbG\_I}) and the atom folds
  (\isn{ev\_now}); at $e+1$, the round at $k = i_x$, whose hypothesis is
  the induction hypothesis, the distance at index $i_x+1$ being $e$.
\item \emph{Clause (ii)} (\isn{aoAb\_sent}).  Two universal
  introductions, a decided split on the invariant, clause (i).  Rows: in
  \isn{gcalc}, \isn{g\_aoAb\_sent}, \isn{g\_aoAb\_sent\_at},
  \isn{g\_aoAb\_decRel}; in \isn{gcalc\_dich}, \isn{g\_aoAb\_acc\_or\_path};
  in \isn{gcalc\_gfix\_dich}, \isn{g\_aoAb\_dec\_path}.
\item \emph{The bridge.}  HOL states $(i, (a, b, p, q))$; the HOL guard
  \isn{aoAbGrdH} is membership in the table, so every state of a run
  satisfies the invariant; \isn{aoAbActH} applies the step and index
  tables.  BG holds (\isn{aomod\_ab\_BG}); BN holds for every component,
  masked or not (\isn{aomod\_ab\_BN}, from \isn{aoAbIn\_neg}); BA holds
  (\isn{aomod\_ab\_BA}).  Truth of the target at an encoded code implies
  $K \le i$ (\isn{aomod\_ab\_GP}, with \isn{aomod\_gap\_pos}).
\item \emph{Clause (iii)} (\isn{aomod\_ab\_deliver}).  Clause (i) at the
  empty context, the numeral $\bar K$ and every encoded start code, not
  the sentence; then \isn{aomod\_gen\_fair\_ev} with $P$ the order $K \le
  i$.
\end{enumerate}
The inductions number seven: two in each nested induction, one per phase
for the stale codes, and the distance induction (a count read off the
proof script; the check records presence).  The rank is $\omega^2$ per
phase, outer on the delivery counter and inner on the refill counter (a
reading).  The index increments exactly at a matching
acknowledgement delivery, which the matching data delivery precedes in
the phase order; no HOL statement about data deliveries is made, and no
theorem shows the protocol failing under weak fairness.

\subsection{The grounded fairness rule}

\paragraph{Definitions.}
A stream parameter $r$, with $\Tot\,r$, carries a program state and the
scheduled component: $\rho_k = \mathrm{fst}(r\cdot k)$ and $\sigma_k =
\mathrm{snd}(r\cdot k)$.  The \emph{run sentence} \isn{aoRunS} is $\Tot
\,r \wedge \forall k.\,\bigvee_i (\sigma_k = \bar\imath \wedge G_i\cdot
\rho_k \wedge A_i\cdot\rho_{k+1}\cdot\rho_k)$, a kit disjunction.  For a
component $j$, Fa (\isn{aoFaS}) is $\exists k_0.\,\forall m.\,\neg(k_0
\le m) \vee \neg G_j\cdot\rho_m$, ``never enabled after some point'', and
Fb (\isn{aoFbS}) is $\forall k.\,\exists m.\,k \le m \wedge \sigma_m =
\bar\jmath$, ``scheduled infinitely often''.  The \emph{fairness
sentence} \isn{aoFairS} is the kit conjunction of Fa $\vee$ Fb over the
masked components: strong fairness in its implication form written as a
disjunction.  At a model-true stream the two sentences say that the
stream is a run strongly fair for the mask (a reading; no least-model
reading of the sentences is formalized).  The \emph{lifted stream}
\isn{aoLiftT} at bound streams $B_j$ has the value $\langle \rho_k,
B_0\cdot k \mathbin{\dot-} k, \ldots, B_{n-1}\cdot k \mathbin{\dot-} k
\rangle$ at $k$ (\isn{aoLiftV}): the counter of $j$ is $c_j(k) = B_j\cdot
k \mathbin{\dot-} k$.

\begin{thm}[\isn{g\_ao\_fair\_rule}, \isn{g\_ao\_run\_ev}; grounded]
\label{thm:app-fair-grounded}
In \isn{gcalc\_gfix}, for every list with \isn{aoOk} and every
well-formed $\Gamma$:
(i) with decided components, if at every extension every code $x$ with
$\mathrm{st}(x) = s_0$ is accessible for \isn{aoRel}, then from $\Gamma
\ent s_0\,\mathbf N$, $\Gamma \ent \Tot\,r$, $\Gamma \ent \isn{aoRunS}
\,r$ and $\Gamma \ent \isn{aoFairS}\,r$, $\Gamma \ent \neg(\rho_0 =
s_0)$;
(ii) for a well-formed $G$, from $\Gamma \ent \mathrm{stRun}_{\isn{aoRel}}(s)$ and $\Gamma \ent
\mathrm{Ev}_G\cdot(s\cdot 0)$, the atom along \isn{aoRelF}, $\Gamma \ent
\exists k.\,G\cdot(s\cdot k)$ (\isn{aoEvTgt}).
\end{thm}
Clause (i) is a rule, a schema in the stream term $r$, not a sentence;
clause (ii) needs no decidedness of the components.  The regression at
\isn{aoLp} and $s_0 = \langle 1, 0\rangle$ is \isn{g\_aoLp\_fair\_rule},
its premise from \isn{g\_acc\_aoLp\_all}.

\paragraph{Proof of \cref{thm:app-fair-grounded}.}
\begin{enumerate}
\item \emph{The counter laws} (in \rga).  From (L1) $B\cdot(k+1) \le
  B\cdot k$, $c(k+1) \le c(k)$ (\isn{ao\_ctr\_le}); with also (L2) $k <
  B\cdot k$, $c(k+1) < c(k)$ (\isn{ao\_ctr\_lt}).  The gap arithmetic is
  \isn{ao\_gap\_mono1}, by the decided split \isn{gapSub\_le\_or\_lt}, and
  \isn{ao\_gap\_anti2}.
\item \emph{The bound in the Fb branch.}  $B\cdot k$ is the least $m$
  with $k \le m$ and $\sigma_m = \bar\jmath$, the search (\isn{muF}) at a
  decided matrix (\isn{ao\_fb\_dec}) whose existence is Fb
  (\isn{ao\_fb\_open}); $B$ is total and meets its matrix
  (\isn{ao\_fb\_Tot}, \isn{ao\_fb\_all}, by \isn{gcalc.mu\_choice} and
  \isn{g\_mu\_choice\_all}).  When $\sigma_k \ne \bar\jmath$: (L2), since
  $B\cdot k \le k$ would force $B\cdot k = k$ and $\sigma_k = \bar\jmath$
  (\isn{ao\_fb\_L2}, no guard premise); and (L1), since $B\cdot k$ is a
  witness at $k+1$ and the search is minimal (\isn{ao\_fb\_L1}, by
  \isn{gcalc.mu\_choice\_min}).
\item \emph{The bound in the Fa branch.}  The constant stream at the
  existential's witness $e$, an eigenvariable (\isn{ao\_Tot\_const}).
  (L1) is trivial (\isn{ao\_fa\_L1}); (L2) uses the guard: from Fa at $e$
  and $G_j\cdot\rho_k$, $k < e$ by the decided split of $e$ and $k$
  (\isn{ao\_fa\_L2}).
\item \emph{The branch continuation} (\isn{ao\_fair\_branches}).
  \isn{aoBnd} says that every bound stream is total and, for masked $j$,
  at every extension and $k$ with $\sigma_k \ne \bar\jmath$, satisfies
  (L1), and (L2) when $G_j\cdot\rho_k$.  HOL induction on the remaining
  components (\isn{ao\_fair\_branches\_aux}): an unmasked component gets
  the constant $0$; a masked one has Fa $\vee$ Fb, eliminated by cases;
  the Fa existential, written $\neg\forall\neg$, yields a fresh
  eigenvariable (\isn{rp\_not\_allE}, \isn{ir\_dnegE}).  The derivation
  has $2^m$ leaves for $m$ masked components; the HOL text is linear.
  The result: any goal that follows at every extension and every family
  with \isn{aoBnd} follows from \isn{aoFairS}.
\item \emph{The lifted run} (\isn{ao\_lift\_step}, \isn{ao\_lift\_run}).
  At $k$, the run disjunction (\isn{ao\_gAny\_E}) gives the scheduled $i$
  with its guard and step.  For every masked $j \ne i$, $\sigma_k \ne
  \bar\jmath$ by numeral distinctness, so (L1) holds; the decided guard
  $G_j\cdot\rho_k$ (\isn{aoDecCsD}) is split: enabled gives (L2),
  \isn{ao\_ctr\_lt} and \isn{aoCtr\_I\_en}, disabled \isn{ao\_ctr\_le} and
  \isn{aoCtr\_I\_dis}.  Then \isn{g\_aoRelB\_I} and \isn{aoRel\_I}; with
  totality (\isn{ao\_lift\_Tot}), the lifted stream is a run of
  \isn{aoRel}.
\item \emph{Clause (i).}  Split on the decided $\rho_0 = s_0$.  In the
  positive branch, step~4 with the continuation: the lifted stream $u$ is
  a run (step~5) whose first value has program part $s_0$, so the premise
  makes $u\cdot 0$ accessible; \isn{g\_run\_acc} (ground induction at the
  parametric motive ``no value of the stream is the point'') gives
  $\forall k.\,\neg(u\cdot k = u\cdot 0)$, refuted at $k = 0$.  The
  negative branch is the conclusion.
\item \emph{Clause (ii)} (\isn{g\_ao\_run\_ev}).  One ground induction
  (\isn{ir\_gfix\_par}) on the eventually atom, at the motive
  \isn{aoEvMot}, $\lambda z.\,\forall k.\,\neg(s\cdot k = z) \vee \exists
  k.\,G\cdot(s\cdot k)$, the stream a parameter, the eigenvariable $v$
  fresh for the context, the motive and the body.  Closure at $v$ under
  the body, at every extension and $k$, by the decided split $s\cdot k =
  v$: in the $G$ branch, the witness $k$ (\isn{ao\_exqI}); in the
  universal branch, the run's edge from $s\cdot k$ to $s\cdot(k+1)$
  (\isn{g\_runE}) turned forward (\isn{aoRel\_E}, \isn{aoRelF\_I}) gives
  the motive at $s\cdot(k+1)$, whose instance at $k+1$ is the target.  The
  final instance is at $s\cdot 0$ and $k = 0$.
\end{enumerate}

\paragraph{Where a decision would hide.}
The classical lifting uses one uniform counter (\isn{aoCnt}), whose
finiteness proof (\isn{aoPend\_finite}) splits cases on ``$j$ is
scheduled again from $k$'', a $\Sigma^0_1$ question, and in the
negative case obtains from strong fairness that $j$ is never enabled
from some point, the negation of a $\Pi^0_2$ question.  A grounded
transcription of that counter would have to decide these questions (a reading of
the classical proof).  The grounded proof builds the counter inside the
branches of each fairness disjunction, by search or as a constant.  Every
case split of the grounded proof is on a decided atom: the start equation, a guard at
a number, the order of two numbers, and in clause (ii) the equation
$s\cdot k = v$; the check forbids ATI, the dichotomy and $\nu$ under both
clauses and the regression.  A grounded ``every fair run meets the
target'', combining steps~4 and 5 with clause (ii), is not proved; the
fairness eventuality is obtained only classically
(\cref{thm:app-fair-genmodel}, \cref{thm:app-abp}).  The negated
sentence ``no stream is a fair run'' is not stated, since the calculus has no
quantifier over stream terms.  The no-variant theorem of random assignment reaches its
negation by spending ATI to decide its hypothesis; for ``no stream is
a fair run'' the corresponding hypothesis, the fairness sentence, has
$\Pi^0_2$ leaves (a reading; no non-derivability is proved), and the
grounded reading of fairness is the disjunctive one.

\paragraph{Where the rung rules are spent.}
For each result: its calculus, and what the check at the end of its
theory requires and forbids among ground induction ($\mu$), coinduction
($\nu$), the dichotomy (Dich) and ATI.  Every result of the four model
theories is classical, and the checks of those theories forbid all four rules
under the fair-termination and delivery theorems.
\begin{center}
\footnotesize
\begin{tabular}{@{}p{0.48\textwidth}lll@{}}
\toprule
result & calculus & check needs & check forbids \\
\midrule
loop, twin, term former, protocol: decidedness, case lemmas, accessibility, \isn{aoMotL\_wf}, \isn{acc\_aoRel\_helpful}, \isn{aoAb\_every}, \isn{aoAb\_sent} & \isn{gcalc} & & all four \\
no path, no run: \isn{g\_ao\_not\_path}, \isn{g\_aoTs\_no\_run}, \isn{g\_aoRel\_no\_run} & $\mu$ & $\mu$ & $\nu$, Dich, ATI \\
$\Acc \vee \Path$: \isn{g\_ao\_acc\_or\_path}, \isn{g\_aoTw\_acc\_or\_path}, \isn{g\_aoAb\_acc\_or\_path} & Dich & Dich & $\mu$, $\nu$, ATI \\
the path atom decided: \isn{g\_aoRel\_dec\_path}, \isn{g\_aoAb\_dec\_path} & $\mu$, Dich & $\mu$, Dich & $\nu$, ATI \\
weak divergence of the twin: \isn{g\_aoTw\_div} & $\nu$ & $\nu$ & $\mu$, Dich, ATI \\
the run of the weak twin: \isn{g\_aoTw\_run} & $\mu$, $\nu$, Dich & $\nu$, Dich & \\
bound streams, counter laws, branches, lift: \isn{ao\_fair\_branches}, \isn{ao\_lift\_run} & \isn{gcalc} & & all four \\
\isn{g\_ao\_fair\_rule}, \isn{g\_aoLp\_fair\_rule} & $\mu$ & $\mu$, the search & $\nu$, Dich, ATI \\
\isn{g\_ao\_run\_ev} & $\mu$ & $\mu$ & $\nu$, Dich, ATI \\
\bottomrule
\end{tabular}
\end{center}
ATI is used nowhere in this appendix.  The search enters
\isn{g\_ao\_fair\_rule} through \isn{mu\_choice\_min} and
\isn{g\_mu\_choice\_all}, both derived rules of \rga, carried into the calculus by cut.
The non-derivability results rest on \isn{oinga\_consistent} and
\isn{oingad\_consistent}.

\paragraph{Lemmas a re-mechanization needs.}
\begin{itemize}
\item \isn{aoPend\_finite}: the pending set is finite under strong
  fairness, the one classically split question.
\item \isn{aomod\_acc\_chain}: no infinite descending chain from an
  accessible point (HOL).
\item \isn{aoB\_cases}, \isn{aoTsB\_cases}, \isn{aoRel\_cases}: case
  lemmas at an eigenvariable code from projection equations alone, each
  case at every extension; no tuple rebuild lemma.
\item \isn{aoCtr\_E\_en}, \isn{aoCtr\_E\_le}: the counter law eliminated,
  a strict decrease when enabled and no increase in either case.
\item \isn{aoMotL\_E}, \isn{aoZf\_eq}: the pack at $m$ eliminated at a
  code with $P$ and $F \le m$, and a counter read as a function.
\item \isn{aoAbIs\_uniq}: a code carries one control tuple, so every
  protocol case analysis is the successor lemma plus a HOL table fact.
\item \isn{aoAbIn\_neg}: a control tuple outside a table refutes
  membership, which gives BN for every component.
\item \isn{ao\_gap\_mono1}, \isn{ao\_gap\_anti2}: monotonicity and
  antitonicity of truncated subtraction, the arithmetic of the
  grounded proof's counters.
\item \isn{g\_run\_acc}, \isn{g\_runE}, \isn{g\_path\_coind},
  \isn{g\_path\_child\_intro}, \isn{g\_run\_of\_path}, \isn{mu\_choice},
  \isn{mu\_choice\_min}, \isn{g\_mu\_choice\_all}, \isn{ev\_now},
  \isn{ev\_next}, \isn{gfix\_omv\_eq\_lfp}, \isn{acc\_omv\_sound}: the
  generic machinery (Appendix~\ref{app:streams}).
\end{itemize}

\clearpage
\section{The Theorems and Their Isabelle Names}
\label{sec:table}

Sessions: \textsf{RGA} (the base and the model interface),
\textsf{IGA} (the $\nu$-former, the model, the rules, the
dichotomy), \textsf{OIGA} (the abstract calculi and the systems),
\textsf{OIGA\_Coind} (the stream layer and the examples, written
Coind below), \textsf{OGA} (the replay and K\H{o}nig),
\textsf{IGA\_PGSem} (certificate-grounded semantics, PGSem below),
\textsf{IGA\_Closing} (ATI from the dichotomy, \dinga, the closed-off
oracle; Closing below), \textsf{Coind\_Dinga} (the junction theory
\isn{CD\_Junction}, which restates every example result
\isn{oingad\_X} at \dinga as \isn{dinga.dinga\_X}; Junction below)
and \textsf{INGA\_PGSem} (\inga's certificate-grounded semantics,
INGA-PG below).
The last group of rows cites three sessions of the companion papers
outside the coinduction chain: \textsf{IGA\_ConPA\_Checker} (the
guarded proof checker for Peano Arithmetic and its completeness for
the derivations of Peano Arithmetic), \textsf{IGA\_ConPA\_Main} (the
checker's consistency sentence in \iga) and \textsf{OIGA\_Fold}
(\rga's consistency sentence in \oiga).  A name
qualified by a calculus, such as \isn{gcalc\_nu.cs\_cind}, is the
theorem in that abstract calculus; the systems carry instances under
their own prefixes.  View: grounded, a derivation in a calculus;
classical, a HOL theorem about a calculus or the least model;
$\Mstar$, a HOL theorem about the $\nu$-model; def., a definition.
A HOL biconditional between two derivability judgments is listed as
grounded and classical: an interderivability.

{\footnotesize
\setlength{\tabcolsep}{4pt}
\begin{longtable}{@{}p{0.34\textwidth}p{0.40\textwidth}p{0.09\textwidth}p{0.11\textwidth}@{}}
\toprule
Statement & Isabelle name & Session & View \\
\midrule
\endhead
\multicolumn{4}{@{}l}{Both fixed points (\S\ref{sec:regions})}\\
falsity set is the lfp; disjointness (Thm.~\ref{thm:disjoint}) & \isn{gfix\_omv\_false\_eq\_lfp}, \isn{gfp\_lfp\_disjoint} & IGA & classical \\
never false & \isn{gfp\_never\_false} & IGA & classical \\
complement iff decided (Thm.~\ref{thm:compl}) & \isn{gfp\_eq\_compl\_iff}, \isn{gdet\_gfp\_eq}, \isn{dpos\_gdet} & IGA & classical \\
determinacy not necessary & \isn{liarOr\_identity}, \isn{liarOr\_not\_gdet} & IGA & classical \\
truth-teller; Liar conjunct & \isn{tt\_gfp}, \isn{tt\_valueless}, \isn{liar\_gap} & IGA & classical \\
accessibility's $\nu$-reading is everything & \isn{acc\_gfp\_UNIV}, \isn{acc\_lfp\_false\_empty} & IGA & classical \\
survival game; the up-by-one game & \isn{win\_gfp\_eq}, \isn{win\_surv\_iff}, \isn{stepU\_valueless} & IGA & classical \\
$\mathcal O$ at gfp & \isn{theorem\_D}, \isn{o\_const\_gfp\_not\_lfp} & IGA & classical \\
\midrule
\multicolumn{4}{@{}l}{The former, the rule, the model (\S\ref{sec:nu})}\\
the former; distinct from $\mu$; unfolding & \isn{trmNu}, \isn{trmNu\_neq\_trmGfix}, \isn{nu\_step1}, \isn{nu\_step2}, \isn{nu\_step3}, \isn{nu\_step4} & IGA & def., classical \\
the clash (Thm.~\ref{thm:clash}) & \isn{nu\_gfix\_clash}, \isn{oinga\_gfixNu\_inconsistent} & IGA, OIGA & classical \\
head steps never turn a $\nu$-atom into a $\mu$-atom & \isn{nu\_hsteps\_not\_gfix}, \isn{nu\_hsteps\_not\_gfix\_good} & OIGA & classical \\
the rule; the calculus with it & \isn{nuS}; \isn{gcalc\_nu} (assumption \isn{ir\_nu}) & IGA, OIGA & def. \\
positive replacement; monotonicity; up to union & \isn{g\_pos\_replace}, \isn{g\_nu\_mono}, \isn{g\_nu\_upto} & OIGA & grounded \\
the $\nu$ class; invariance under value congruence & \isn{nuGood}, \isn{ovcrdO\_INV} & IGA & def., classical \\
oracle model; oracle map; $\Theta^*$ & \isn{omvO}, \isn{Om}, \isn{Ostar}, \isn{ocons\_Ostar} & IGA & def., classical \\
coherence & \isn{nu\_coherent\_T}, \isn{nu\_coherent\_F} & IGA & $\Mstar$ \\
least model persists & \isn{omv\_le\_omvO\_Ostar} & IGA & classical \\
the interface; its soundness theorem; \rga sound for $\Mstar$ & \isn{omega\_model}, \isn{omega\_model.step\_M\_sound}, \isn{rga\_nu\_sound} & RGA, IGA & classical \\
coherent oracles; the interface at every one closed under substitution and congruence; the closed-off oracle & \isn{cohO}, \isn{cohO3.nu\_model\_omega\_model}, \isn{OstarC}, \isn{closed\_off\_omega\_model} & IGA, Closing & def., classical \\
ground induction at $\Mstar$ & \isn{omega\_model\_gfix}, \isn{gfix\_omvO\_ind\_par} & IGA & $\Mstar$ \\
the rule's step is sound & \isn{nu\_step\_mjt} & IGA & $\Mstar$ \\
\inga sound, consistent, conservative (Thm.~\ref{thm:inga}) & \isn{inga\_nu\_sound}, \isn{inga\_consistent}, \isn{inga\_conservative} & IGA & classical \\
\oinga sound, consistent, conservative (Thm.~\ref{thm:inga}); \oiga sound for $\Mstar$ & \isn{oinga\_nu\_sound}, \isn{oinga\_consistent}, \isn{oinga\_conservative}, \isn{oiga\_nu\_sound} & OIGA & classical \\
unsound for the least model by design (Thm.~\ref{thm:unsound}) & \isn{inga\_omega\_unsound}, \isn{nuZ\_true}, \isn{nuZ\_true\_valueless}; \isn{inga\_Z} & IGA & classical; grounded \\
\midrule
\multicolumn{4}{@{}l}{The dichotomy (\S\ref{sec:dich})}\\
the dual; its laws & \isn{ndual}, \isn{ndual\_posP}, \isn{ndual\_exists}, \isn{ndual\_fun}, \isn{ndual\_subf2} & IGA & def., classical \\
the premise & \isn{nfdec}, \isn{nfdec\_accPhi} & IGA & def. \\
swap identity (Thm.~\ref{thm:swap}) & \isn{swap\_identity} & IGA & $\Mstar$ \\
validity characterized; the premise gives bivalence (Thm.~\ref{thm:C}) & \isn{dichValid}, \isn{thmC}, \isn{dich\_premise\_bivalent} & IGA & $\Mstar$ \\
the rule; its step is sound (Thm.~\ref{thm:C}) & \isn{dichS}, \isn{dich\_step\_mjt} & IGA & def., $\Mstar$ \\
controls & \isn{dichB\_invalid}, \isn{botD\_never\_true}, \isn{dichA\_valid} & IGA & $\Mstar$ \\
ATI from the dichotomy (Thm.~\ref{thm:ati-dich}) & \isn{g\_ati\_of\_dich}, \isn{gcalc\_gfix\_dich} $\subseteq$ \isn{gcalc\_ati} & Closing & grounded \\
\dinga; its derivability equals \oinga's with the dichotomy & \isn{dinga}, \isn{dinga\_derives}, \isn{dinga\_eq\_oingad}, \isn{n2pr\_dinga\_derives} & Closing & def., classical \\
the price (Thm.~\ref{thm:price}) & \isn{dinga\_dich\_price}, \isn{oinga\_dich\_price} & Closing, OIGA & classical \\
no certificate-grounded semantics with ATI (Thm.~\ref{thm:nocert}) & \isn{cert\_sem\_gap}, \isn{ati\_cert\_unsound}, \isn{dich\_cert\_unsound}, \isn{dich\_not\_cert\_sound} & PGSem, Closing & classical \\
\dinga and \oiga have none (Cor.~\ref{cor:nocert}) & \isn{dinga\_entails\_re}, \isn{dinga\_pi1\_ouniv}, \isn{dinga\_not\_cert\_sound}, \isn{oiga\_not\_cert\_sound} & Closing, PGSem & classical \\
\oinga sound for $\Mstar$ and has none (Cor.~\ref{cor:nocert}) & \isn{oinga\_entails\_re}, \isn{oinga\_not\_cert\_sound} & Junction & classical \\
same price as bivalence in the least model & \isn{acc\_or\_mupath\_omega\_iff\_wf}, \isn{accbiv\_iff\_acc\_or\_mupath} & OIGA & classical \\
the dichotomy not valid for the least model, with the $\nu$-atom itself & \isn{nu\_omv\_le\_lfp}, \isn{acc\_or\_path\_omega\_iff\_wf}, \isn{dich\_not\_omega\_valid} & OIGA & classical \\
the same, with the premise derived in \rga & \isn{dich\_not\_omega\_valid\_decRel} & OIGA & classical \\
true in $\Mstar$ at every decided relation; the contrast & \isn{oingad\_dich\_decided\_true}, \isn{flipD\_contrast} & OIGA & $\Mstar$, classical \\
exclusion; decided $\nu$-atom (Thm.~\ref{thm:excl}) & \isn{g\_mu\_not\_nu}, \isn{g\_dec\_nu}, \isn{g\_acc\_not\_path}, \isn{g\_dec\_path} & OIGA & grounded \\
$\Acc \vee \Path$ & \isn{g\_acc\_or\_path} & OIGA & grounded \\
accessibility never false in $\Mstar$; $\neg\Acc$ underivable & \isn{accO\_false\_empty}, \isn{oingad\_acc\_not\_refuted}, \isn{dinga\_acc\_not\_refuted} & IGA, OIGA, Closing & $\Mstar$, classical \\
\midrule
\multicolumn{4}{@{}l}{\inga's certificate-grounded semantics (\S\ref{sec:price}, Appendix~\ref{app:inga})}\\
\inga is proof-grounded and not sound for the least model (Thm.~\ref{thm:inga-pg}); not an instance of the locale (Thms.~\ref{thm:app-inga-notgsem}, \ref{thm:app-inga-pg}) & \isn{inga\_cert\_grounded}, \isn{inga\_pg}, \isn{inga\_pg\_not\_omega} & INGA-PG & classical \\
\inga and \rga with $\nu$ consistent, from the certified evaluation (Cor.~\ref{cor:inga-con-mf}) & \isn{ng\_consistent}, \isn{inga\_consistent\_mf}, \isn{rganu\_consistent\_mf} & INGA-PG & classical \\
the disjunction and witness properties; no fictions (Cor.~\ref{cor:inga-dpwp}) & \isn{inga\_disj\_prop}, \isn{inga\_witness\_prop}, \isn{inga\_no\_fict}, \isn{inga\_nu\_no\_fict} & INGA-PG & classical \\
the locale entails the least model; \inga not an instance (Thm.~\ref{thm:app-inga-notgsem}); the all-zeros atom valueless, from the least model alone & \isn{ng\_gs\_omega\_true}, \isn{inga\_not\_gsem\_sound}, \isn{inga\_Z\_valueless}, \isn{inga\_omega\_unsound\_mf} & INGA-PG & classical \\
the interface; its reduction to the two properties (Thm.~\ref{thm:app-inga-reduce}) & \isn{cert\_iface}, \isn{inga\_cert\_allT}, \isn{inga\_cert\_excl}, \isn{inga\_cert\_iff\_props}, \isn{inga\_pg\_iff\_dp\_wp} & INGA-PG & def., classical \\
the certified evaluation; the filtered map; the certified oracle; validity & \isn{gnhev}, \isn{gnosr}, \isn{gosem}, \isn{gfsem}, \isn{FOg}, \isn{GOg}, \isn{OmG}, \isn{OstarG}, \isn{ocert}, \isn{gnowsjt} & INGA-PG & def. \\
the oracle certified and consistent with no hypothesis; determinism (Thm.~\ref{thm:app-ng-oracle}) & \isn{OstarG\_T\_derivable}, \isn{OstarG\_F\_refutable}, \isn{ocert\_OstarG}, \isn{gfpG\_lfpG\_disjoint}, \isn{ocons\_OmG}, \isn{ocons\_OstarG}, \isn{gnosr\_unique}, \isn{gnosr\_OstarG\_unique\_mf} & INGA-PG & classical \\
completeness at a certified oracle (Thm.~\ref{thm:app-ng-complete}); coherence & \isn{gnhev\_gnosr\_complete}, \isn{ng\_complete}, \isn{ng\_complete\_F}, \isn{ng\_coherent\_T\_iff}, \isn{ng\_coherent\_F\_iff} & INGA-PG & classical \\
the semantic Park lemma (Thm.~\ref{thm:app-ng-park}) & \isn{nu\_step\_gnowsjt}, \isn{inga\_nu\_step\_gnowsjt}, \isn{ng\_posrep}, \isn{ng\_cert\_allT}, \isn{ng\_rep} & INGA-PG & classical \\
substitution, congruence, head steps (Thm.~\ref{thm:app-ng-kit}) & \isn{ngco}, \isn{ngsubst}, \isn{ngcong}, \isn{gnosr\_Oz\_subf}, \isn{gnosr\_Oz\_vcr}, \isn{ng\_hstep\_sim}, \isn{ng\_hstep\_rev}, \isn{ng\_tr\_case\_nuTF} & INGA-PG & def., classical \\
the certified oracle closed (Thm.~\ref{thm:app-ng-closure}) & \isn{ngstar}, \isn{OsubG\_prefix}, \isn{OstarG\_osub}, \isn{OcgG\_prefix}, \isn{OstarG\_ocg}, \isn{ngstar\_OstarG\_mf} & INGA-PG & def., classical \\
the rules of \rga and ground induction (Thm.~\ref{thm:app-ng-step}) & \isn{ng\_rga\_step\_gnowsjt}, \isn{ng\_gfix\_step\_gnowsjt}, \isn{ngf\_transport}, \isn{ngf\_spine} & INGA-PG & classical \\
soundness; derivable iff true; the oracle; consistency (Thm.~\ref{thm:app-ng-sound}) & \isn{ngsound}, \isn{gnsem\_sound}, \isn{ng\_open\_iff}, \isn{ng\_false\_iff}, \isn{ng\_nu\_true\_iff}, \isn{ng\_nu\_false\_iff}, \isn{ng\_consistent} & INGA-PG & def., classical \\
at \inga; at \rga with $\nu$ (Thm.~\ref{thm:app-inga-pg}) & \isn{inga\_gnsem\_sound}, \isn{inga\_gn\_open\_iff}, \isn{inga\_nu\_true\_iff}; \isn{rganu\_gnsem\_sound}, \isn{rganu\_disj\_prop}, \isn{rganu\_witness\_prop} & INGA-PG & classical \\
\midrule
\multicolumn{4}{@{}l}{Streams, runs, branches (\S\ref{sec:streams})}\\
prefix function and set; elementwise view & \isn{pfx}, \isn{pset}, \isn{g\_pset\_elem\_iff} & Coind & def., grounded \\
branches; the prefix set a branch (Thm.~\ref{thm:branch}); round trip (Thm.~\ref{thm:roundtrip}) & \isn{fbBranch}, \isn{g\_pset\_Dec}, \isn{g\_pset\_Branch}, \isn{g\_RT\_s}, \isn{g\_RT\_p} & Coind & grounded \\
prefix set and branches (Thm.~\ref{thm:pset-path}) & \isn{g\_pset\_path}, \isn{g\_path\_pset\_hol\_iff}, \isn{g\_path\_pset\_iff}, \isn{g\_path\_pset\_dec} & Coind & grounded, classical \\
generic certificate coinduction & \isn{g\_path\_coind} & Coind & grounded \\
runs (Thms.~\ref{thm:runs}, \ref{thm:runofpath}) & \isn{stRun}, \isn{g\_run\_acc}, \isn{g\_run\_path}, \isn{g\_run\_of\_path} & Coind, OGA & def., grounded \\
the replay & \isn{gPath}, \isn{gProg}, \isn{gpb\_run\_of\_progress} & OGA & def., grounded \\
well-founded or explicit branch (Thm.~\ref{thm:koenig}) & \isn{g\_tree\_acc\_or\_branch}, \isn{g\_tree\_branch}, \isn{g\_tree\_elem} & Coind & grounded \\
K\H{o}nig's lemma as a corollary & \isn{koenig\_fb\_inst}, \isn{koenig\_fb\_cor} & OGA & grounded \\
the Kleene tree: decided, infinite, no recursive path & \isn{kTree\_Dec}, \isn{kTree\_InfTree}, \isn{kleene\_no\_recursive\_path} & OGA & grounded, classical \\
\midrule
\multicolumn{4}{@{}l}{Examples (\S\ref{sec:showpieces})}\\
random assignment (Thm.~\ref{thm:rand}) & \isn{g\_acc\_R1\_all}, \isn{g\_R1\_no\_variant} & Coind & grounded \\
nested random assignment at every fixed depth (\S\ref{sec:rungmap}) & \isn{rdRel}, \isn{rd\_ti}, \isn{acc\_rdRel}, \isn{g\_acc\_rdRel} & Coind & def., grounded \\
fair termination through random assignment (Thm.~\ref{thm:fair}) & \isn{aoT}, \isn{g\_acc\_ao\_all}, \isn{g\_acc\_ao\_init},\newline \isn{g\_ao\_not\_path}, \isn{g\_ao\_acc\_or\_path} & Coind & def., grounded \\
the correspondence of infinite runs; the least-model reading (Thms.~\ref{thm:app-fair-corr}, \ref{thm:app-fair-loop}) & \isn{ao\_sound}, \isn{ao\_complete}, \isn{ao\_fair\_term},\newline \isn{aomod\_fair\_term} & Coind & classical \\
the relaxed masked edge: correspondence and transfers (Thm.~\ref{thm:app-fair-corr}) & \isn{aoEdgeX}, \isn{ao\_xsound}, \isn{ao\_xcomplete},\newline \isn{ao\_fair\_xterm}, \isn{ao\_fair\_xev} & Coind & def., classical \\
the fairness twin: strong termination, weak divergence (Thm.~\ref{thm:fair-twin}) & \isn{aoTs}, \isn{aoTw}, \isn{g\_acc\_aoTs\_init}, \isn{g\_aoTw\_div},\newline \isn{oinga\_aoTw\_not\_acc} & Coind & def., grounded, classical \\
the fairness twin: every strong state accessible; the weak run; non-derivability (Thm.~\ref{thm:app-fair-twin}) & \isn{g\_acc\_aoTs\_all}, \isn{g\_aoTw\_run}, \isn{oingad\_aoTw\_not\_acc},\newline \isn{aomod\_tw\_weak\_not\_der}; at \dinga \isn{dinga\_aoTw\_not\_acc} & Coind, Junction & grounded, classical \\
the weak correspondence; the fairness twin's separation & \isn{ao\_wsound}, \isn{ao\_wcomplete}, \isn{aomod\_tw\_fair\_term},\newline \isn{aomod\_tw\_weak\_run}, \isn{aomod\_tw\_separation} & Coind & classical \\
the transformation as a term former; its case lemma; decidedness; the motive pack (Thm.~\ref{thm:app-fair-gen}) & \isn{aoRel}, \isn{aoRel\_cases}, \isn{g\_aoRel\_decRel},\newline \isn{aoMotL\_wf} & Coind & def., grounded \\
the loop as a regression of the term former & \isn{aoLp}, \isn{g\_acc\_aoLp\_all}, \isn{g\_aoLp\_fair\_rule} & Coind & def., grounded \\
one helpful direction; the generic least-model pipeline (Thm.~\ref{thm:app-fair-genmodel}) & \isn{acc\_aoRel\_helpful}, \isn{aomod\_gen\_fair\_term},\newline \isn{aomod\_ev\_chain}, \isn{aomod\_gen\_fair\_ev} & Coind & grounded, classical \\
the alternating-bit protocol (Thm.~\ref{thm:abp}) & \isn{aoAbPg}, \isn{aoAbInv}, \isn{aoAb\_phase2}, \isn{aoAb\_round},\newline \isn{aoAb\_every}, \isn{aoAb\_sent}, \isn{g\_aoAb\_sent}, \isn{g\_aoAb\_sent\_at},\newline \isn{aomod\_ab\_deliver} & Coind & def., grounded, classical \\
the grounded fairness rule (Thm.~\ref{thm:fair-grounded}) & \isn{aoRunS}, \isn{aoFairS}, \isn{ao\_fair\_branches}, \isn{ao\_lift\_run},\newline \isn{g\_ao\_fair\_rule}, \isn{g\_ao\_run\_ev} & Coind & def., grounded \\
divergent twin (Thm.~\ref{thm:rand-div}) & \isn{g\_R1p\_div}, \isn{g\_R1p\_run}, \isn{g\_R1\_acc\_or\_path}, \isn{g\_R1\_dec\_path} & Coind & grounded \\
non-derivability & \isn{oinga\_R1p\_not\_acc}, \isn{oingad\_R1p\_not\_acc}; at \dinga \isn{dinga\_R1p\_not\_acc} & Coind, Junction & classical \\
rank $\omega^2$ & \isn{acc\_Rn\_init}, \isn{g\_acc\_Rn\_init} & Coind & grounded \\
termination and divergence (Thm.~\ref{thm:trace}) & \isn{tr\_Q1\_iff}, \isn{tr\_Q2\_iff}, \isn{tr\_excl} & Coind & grounded, classical \\
the gap closed (Thm.~\ref{thm:trace-gap}) & \isn{tr\_term\_of\_ndiv}, \isn{tr\_dich}, \isn{tr\_dich\_dec}, \isn{tr\_ati\_disj} & Coind & grounded \\
AG\,AF\,$p$ with a path-atom leaf: the body, its premise, the restriction; decided (Thms.~\ref{thm:agaf}, \ref{thm:app-agaf}) & \isn{agafPhi}, \isn{agafD}, \isn{nfdec\_agafPhi}, \isn{restrRel}, \isn{g\_restrRel\_decRel}, \isn{g\_agaf\_premise}, \isn{g\_dec\_agaf}, \isn{g\_dec\_agaf\_restr} & Coind & def., classical, grounded \\
the eventually atom, leads-to; the Lattice and WELL rules (Thm.~\ref{thm:lattice}); the countdown & \isn{trEv}, \isn{ltTo}, \isn{ev\_mono\_trans}, \isn{g\_lattice}, \isn{g\_lattice\_at}, \isn{g\_well}, \isn{cd\_term} & Coind & def., grounded \\
reaches $n$ then goes up & \isn{tr\_up0} & Coind & grounded \\
the reach atom valueless at $\langle 1,0\rangle$; the premise valueless, underivable, irrefutable & \isn{tr\_reach\_not\_true}, \isn{tr\_reach\_not\_false}, \isn{trUp\_premise\_valueless}, \isn{oingad\_trUp\_dich\_premise\_not\_derivable}, \isn{tr\_dich\_prem\_valueless}, \isn{oingad\_tr\_dich\_unavailable}; at \dinga \isn{dinga\_trUp\_dich\_premise\_not\_derivable}, \isn{dinga\_tr\_dich\_prem\_open}, \isn{dinga\_tr\_dich\_unavailable} & Coind, Junction & $\Mstar$, classical \\
a decidedness atom is never false in $\Mstar$ & \isn{nfdec\_not\_false} & Coind & $\Mstar$ \\
determinacy (Thm.~\ref{thm:det}) & \isn{gm\_det}, \isn{gm\_excl}, \isn{gm\_decII}, \isn{g\_gmDec} & Coind & grounded \\
concrete games & \isn{gm\_I\_one\_round}, \isn{gm\_II\_copies} & Coind & grounded \\
player II's strategy (Thm.~\ref{thm:sigII}) & \isn{gm\_sigII} & Coind & grounded \\
bucket game (Thm.~\ref{thm:bucket}) & \isn{cs\_cind}, \isn{cs\_cind\_init}, \isn{cs\_det}, \isn{cs\_decII}, \isn{cs\_stepmother\_refuted}, \isn{cs\_sig\_total}, \isn{cs\_sig\_step} & Coind & grounded \\
the round-robin play lost below capacity three (Thm.~\ref{thm:app-tight}) & \isn{csPlay}, \isn{cs\_rr\_tight} & Coind & def., grounded \\
the game atom not false in $\Mstar$ where player I survives & \isn{gmSurvI}, \isn{gm\_winI\_not\_false} & Coind & def., $\Mstar$ \\
the game atom false at a stuck position & \isn{gm\_winI\_stuck\_false} & Coind & $\Mstar$ \\
the Stepmother's legal pour; her loss underivable & \isn{cs\_pour\_legal}, \isn{cs\_init\_surv}, \isn{cs\_winI\_not\_false}, \isn{oingad\_stepmother\_not\_refuted}; at \dinga \isn{dinga\_stepmother\_not\_refuted} & Coind, Junction & grounded, $\Mstar$, classical \\
$\Dec(\mathcal O)$ (Thm.~\ref{thm:decO}) & \isn{g\_dec\_O}, \isn{dinga\_dec\_O}, \isn{g\_O\_not\_coO}, \isn{dcmp\_total} & OIGA, Closing & grounded \\
the $\mathcal O$-atom in $\Mstar$ is Kleene's $\mathcal O$; the characteristic term, total in \dinga, characteristic in $\Mstar$ (Thm.~\ref{thm:chiO}) & \isn{O\_Mstar\_iff}, \isn{oiga\_okF\_omv}, \isn{chiO}, \isn{dinga\_Tot\_chiO}, \isn{dinga\_chiO\_Mstar}, \isn{dinga\_chiO\_char} & Closing & $\Mstar$, def., grounded \\
Cantor--Bendixson shape & \isn{g\_cb\_dich}, \isn{g\_kernel\_perfect}, \isn{g\_cb\_rank\_induct}, \isn{g\_cb\_dec\_kernel}; at \dinga \isn{dinga\_cb\_dich}, \isn{dinga\_kernel\_perfect} & OIGA, Closing & grounded \\
\midrule
\multicolumn{4}{@{}l}{Facts of the companion papers cited here (\S\S\ref{sec:bg}, \ref{sec:dich}, \ref{sec:concl})}\\
\oga decides a universal it neither derives nor refutes & \isn{oga\_godel\_fict} & OGA & classical \\
the class rung: consistent, sound for the closed-off model; $\mu$-bivalence at every positive body derives $0=1$ & \isn{mubivD\_consistent}, \isn{imb\_closed\_off\_sound}, \isn{g\_mubiv\_contra} & Closing & classical \\
\iga derives the consistency sentence of the guarded checker for Peano Arithmetic; for every derivation of Peano Arithmetic, \rga derives that the checker accepts a code of the derivation; \iga's derivations embed in \inga, \oinga, \dinga & \isn{gi\_ConPA}; \isn{pa\_ChkPAr\_complete}; \isn{gi\_inga\_embed}, \isn{inga\_oinga\_embed}, \isn{inga\_dinga\_embed} & IGA\_\allowbreak ConPA\_\allowbreak Main, IGA\_\allowbreak ConPA\_\allowbreak Checker, IGA, OIGA, Closing & grounded; classical \\
\oiga derives \rga's consistency sentence; its derivations embed in \oinga and, through the equality, \dinga & \isn{oiga\_ConRGA}; \isn{oiga\_oinga\_embed}, \isn{oiga\_oingad\_embed}, \isn{dinga\_eq\_oingad} & OIGA\_\allowbreak Fold, OIGA, Closing & grounded; classical \\
\bottomrule
\end{longtable}}

\end{document}